\documentclass[12pt, a4paper]{article}
\usepackage[english]{babel}
\usepackage[top=1in, bottom=1in, left=1in, right=1in]{geometry}

\usepackage{amssymb}
\usepackage{amsmath}
\usepackage{amsthm}
\usepackage{amsfonts}
\usepackage{hyperref}
\usepackage{enumitem}
\usepackage{titlesec}
\usepackage{secdot}
\sectiondot{subsection}
\sectiondot{subsubsection}
\setlist[enumerate]{topsep=0pt, itemsep=.5ex, parsep=0ex}
\usepackage{tkz-tab}
\usepackage{fancyhdr}
\usepackage{color}
\usepackage[colorinlistoftodos]{todonotes}

\hypersetup{colorlinks=true, allcolors = blue}
\usepackage{tcolorbox}
\usepackage{subcaption}
\usepackage{titlesec}
\usepackage{relsize}
\usepackage{scalerel}
\usepackage{tabularx}
\usepackage[bitstream-charter]{mathdesign}
\usepackage{xcolor}
\usepackage{mathtools}
\usepackage{multicol}
\theoremstyle{remark}
\theoremstyle{definition}

\newtheorem{remark}{Remark}

\numberwithin{hq}{section}

\usepackage{pgf,tikz}
\usepackage{setspace}
\usepackage{mathrsfs}

\usepackage{makecell}
\usepackage{algorithm}
\usepackage{algpseudocode}
\usepackage{array, booktabs}
\newcommand{\bu}{\pmb{u}}

\newcommand{\br}{\pmb{r}}

\newcommand{\bv}{\pmb{v}}
\newcommand{\bn}{\pmb{n}}
\newcommand{\bD}{\pmb{D}}
\newcommand{\bM}{\pmb{M}}

\newcommand{\bSig}{\pmb{\Sigma}}

\newcommand{\bK}{\pmb{K}}
\newcommand{\Ph}[1]{\textcolor{black}{#1}}

\begin{document}
\begin{center}
\Large
\textbf{Joint State-Parameter Estimation for the Reduced Fracture Model via the United Filter} 
\end{center}
\begin{center}
Toan Huynh, Thi-Thao-Phuong Hoang, Guannan Zhang, Feng Bao
\end{center}
\begin{abstract}
In this paper, we introduce an effective United Filter method for jointly estimating the solution state and physical parameters in flow and transport problems within fractured porous media. Fluid flow and transport in fractured porous media are critical in subsurface hydrology, geophysics, and reservoir geomechanics. Reduced fracture models, which represent fractures as lower-dimensional interfaces, enable efficient multi-scale simulations. However, reduced fracture models also face accuracy challenges due to modeling errors and uncertainties in physical parameters such as permeability and fracture geometry. To address these challenges, we propose a United Filter method, which integrates the Ensemble Score Filter (EnSF) for state estimation with the Direct Filter for parameter estimation. EnSF, based on a score-based diffusion model framework, produces ensemble representations of the state distribution without deep learning. Meanwhile, the Direct Filter, a recursive Bayesian inference method, estimates parameters directly from state observations. The United Filter combines these methods iteratively: EnSF estimates are used to refine parameter values, which are then fed back to improve state estimation. Numerical experiments demonstrate that the United Filter method surpasses the state-of-the-art Augmented Ensemble Kalman Filter, delivering more accurate state and parameter estimation for reduced fracture models. This framework also provides a robust and efficient solution for PDE-constrained inverse problems with uncertainties and sparse observations.
\end{abstract}
%
%
%

\section{Introduction}
Fluid flow and transport problems in fractured porous media play a crucial role in various fields, such as subsurface hydrology, geophysics, and reservoir geomechanics. Therefore, accurate and robust numerical simulations are necessary for addressing such problems. However, the presence of the fractures often causes significant challenges in designing the required algorithms. In particular, depending on its permeability, a fracture can act as either a fast conduit or a geological barrier relative to the surrounding rock matrix. Additionally, the widths of the fractures are significantly smaller than the size of the computational domain and any practical spatial mesh size. Consequently, after introducing a global mesh over the entire domain, the local mesh around the fractures requires refinement. This process is known to be computationally inefficient. An effective way to tackle this issue is to take advantage of the small widths of the fractures and treat them as interfaces, which avoids local refinement around the fractures. This approach reformulates the original problem into a new one that accounts for the interaction between the fractures and the surrounding rock matrix (see~\cite{Alboin1999, Alboin2002, Jaffre2005, Angot2009, Fuma2011, Morales2012, SHLee2020, Amir2021,  Gander2021} and the references therein). Models which consider fractures as low-dimensional objects are referred to as reduced fracture models or mixed-dimensional models. Moreover, these models enable the application of different time scales in the fracture and in the rock matrix through a multi-scale approach, such as global-in-time domain decomposition methods~\cite{Hoang2013, Hoang2016, Hoang2017, Toan2023a, Toan2023b, Toan2024}.

Despite the success of reduced fracture models in capturing fracture-matrix interactions and enabling multi-scale simulations, their application to practical problems often leads to unexpected model errors, such as fracture representation errors from treating fractures as lower-dimensional interfaces, fracture-matrix connectivity errors caused by incorrect assumptions about the interaction between fractures and the surrounding porous medium, and discretization errors in fracture geometry due to oversimplifications in fracture shape and orientation. These errors primarily arise from incomplete knowledge or neglected physics, which introduce intrinsic modeling inaccuracies and uncertainties that can significantly impact simulation reliability. 

Even within well-studied domains where the reduced fracture model has been thoroughly validated, deriving the model requires key characteristics of the rock matrix and fractures, such as permeability, porosity, and fracture widths. These properties not only influence fluid flow patterns but also play a critical role in determining suitable numerical algorithms. However, in practical settings, such as modeling underground fluid flow or contaminant transport, these properties are often not directly observable. Instead, available data are typically sparse and noisy measurements of subsurface flow and transport processes, making the development of an accurate numerical solver challenging.

To mitigate model errors and uncertainties, an inverse problem must often be solved. The study of inversion schemes for state and parameter estimation in flow and transport within (fractured) porous media has attracted great attention of researchers \cite{Neuman1979, Mauldon1993, Renshaw1996, Gwo2001, Ameur2002, Dai2004,LeGoc2010,  Krause2013,  Ameur2018, Wei2023,  Noii2021, Noii2022, Toan2025}. However, most of these studies focus exclusively on either estimating the state of the solution or the characteristics of the fractures, rather than addressing both simultaneously. Furthermore, only a limited number of studies specifically consider state and/or parameter estimation for the reduced fracture model \cite{Ameur2018, Toan2025}.

\vspace{0.5em}

In this work, we introduce a data assimilation framework for joint state-parameter estimation to simultaneously improve the accuracy of numerical solutions and refine the model. The mathematical foundation of data assimilation is the ``optimal filtering problem'', which seeks to determine the ``optimal'' estimate for the target state variable and construct an approximation of its conditional probability density function (PDF) given observational data -- known as the ``filtering density (distribution)". The optimal state estimation is then formulated as the corresponding conditional expectation. Key approaches for solving the optimal filtering problem include Bayesian filters such as Kalman-type filters \cite{Kalman1961, Evensen1994,  Houtekamer1998, Julier2004} and particle filters \cite{Gordon1993, Pitt1999, Snyder2008, Chorin2009,   Leeuwen2010, Andrieu2010,  Kang2018}, along with other optimal filtering methods that rely on solving stochastic partial differential equations (SPDEs) \cite{Zakai1969, Bao2014b, Bao2016, Bao2017, Bao2020,  Bao2021b}. A standard method for joint state-parameter estimation is the Augmented Ensemble Kalman Filter (AugEnKF) \cite{Kitagawa1998}, which extends the state variable by incorporating unknown model parameters and updates them simultaneously using the Ensemble Kalman Filter (EnKF). However, the effectiveness of the AugEnKF is heavily dependent on EnKF's performance, which tends to deteriorate when dealing with nonlinear optimal filtering problems. Additionally, since the AugEnKF indirectly estimates parameters by relying on state updates to guide parameter adjustments, it often yields less accurate parameter estimates, as state estimation takes priority in the joint estimation process \cite{Kantas2015}. To overcome these limitations, a more robust and effective approach, i.e., the United Filter method, was introduced in \cite{Bao2024b}. The United Filter method integrates the novel Ensemble Score Filter \cite{Bao2024a} for state estimation with the Direct Filter method \cite{Bao2019b} for online parameter estimation, forming an iterative algorithm that sequentially refines both the state and model parameters.

The Ensemble Score Filter (EnSF) is a generative artificial intelligence enabled method for optimal filtering, which adopts the score-based diffusion model framework to produce a large ensemble of samples representing the filtering density of the target state variable. As a crucial class of generative machine learning methods, diffusion models utilize noise injection to gradually perturb data and then learn to reverse this process, which enables the generation of new samples. Diffusion models have been widely applied in image synthesis \cite{Song2019, Ho2020,  Cai2020, Dhariwal2021, Ho2022, Meng2022}, image denoising \cite{Sohl2015,  Ho2020, Kawar2021, Luo2021}, image enhancement \cite{Kim2021, Li2022, Whang2022, Saharia2023} and natural language processing \cite{Austin2021, Hoogeboom2021, LLi2022, Savinov2022, PYu2022}. These models rely on the ``score'' to store information about the data. When applied to optimal filtering, score-based diffusion models encode filtering densities within the score models \cite{Bao2024c}. In the EnSF framework, we employ an ensemble approximation for the score model without training it through deep neural networks \cite{Bao2024a}. With the powerful generative capabilities of diffusion models, it has been demonstrated that EnSF can produce samples from complex target distributions \cite{YSong2021, YLiu2024}, and it can track stochastic nonlinear dynamical systems in very high dimensional spaces with superior accuracy and efficiency \cite{Bao2024a, Bao2025a}.

The Direct Filter method, on the other hand, estimates parameters through recursive Bayesian inference, directly mapping observational data into the parameter space~\cite{Bao2019a, Bao2019b}. Unlike the AugEnKF, which merges state and parameter estimation into a single augmented problem, the Direct Filter separately exploits state dynamics and observations to derive the most probable parameter values. This method has demonstrated accurate parameter estimation for high-dimensional problems and practical applications~\cite{Bao2021,  Dyck2021, Bao2022, Bao2023, Bao_DF_2025}. However, a notable drawback is that the Direct Filter requires complete observations of the state variable. When only partial state observations are available, an auxiliary state estimation step is necessary. To resolve this issue, the United Filter method iteratively couples the EnSF and Direct Filter methods. The EnSF estimates the state from available observations, and the Direct Filter uses this estimated state to update model parameters. This iterative process progressively refines both state and parameter estimates, which can significantly improve the accuracy of the reduced fracture model’s numerical solutions.

In this paper, we apply the United Filter method to jointly estimate states and parameters for reduced fracture models of flow and transport problems. Our forward solver is a discretized reduced fracture model based on the backward Euler scheme and mixed finite element methods~\cite{Brezzi1991, Roberts1991, Boffi2013}, which are mass conservative and effective in handling heterogeneous and anisotropic diffusion tensors. We present several numerical experiments to evaluate the performance of the United Filter method and compare its results with those obtained from the AugEnKF. The numerical results demonstrate the superior accuracy and reliability of the proposed United Filter approach for joint state-parameter estimation, particularly in scenarios with nonlinear and sparse observations.

The rest of this paper is organized as follows: Section~\ref{sec2} introduces the reduced fracture model, derived from the pure diffusion equation, and describes its spatial and temporal discretization. We shall discuss the data assimilation framework in Section~\ref{sec3}, where the Ensemble Score Filter and Direct Filter methods are introduced sequentially, along with their implementation. These two methods are then combined to derive the United Filter algorithm. In Section~\ref{sec4}, several numerical experiments are carried out and compared with the AugEnKF to illustrate the performance of the proposed method. Finally, the paper is closed with a conclusion section.

\section{Formulation of the Reduced Fracture Model}\label{sec2}
In this section, we briefly review the derivation of the reduced fracture model for porous media with a single fracture and use the pure diffusion equation as an example for convenience of presentation \footnote{The formula with multiple (possibly intersecting) fractures can be derived in a similar manner (see. e.g., \cite{Alboin2002}).}. We then introduce the corresponding numerical PDE solver and present its compact formulation for the data assimilation task, enabling the joint estimation of state and parameters. 

\begin{figure}[h!]
\vspace{-0.4cm}
\centering
\includegraphics[scale=0.55]{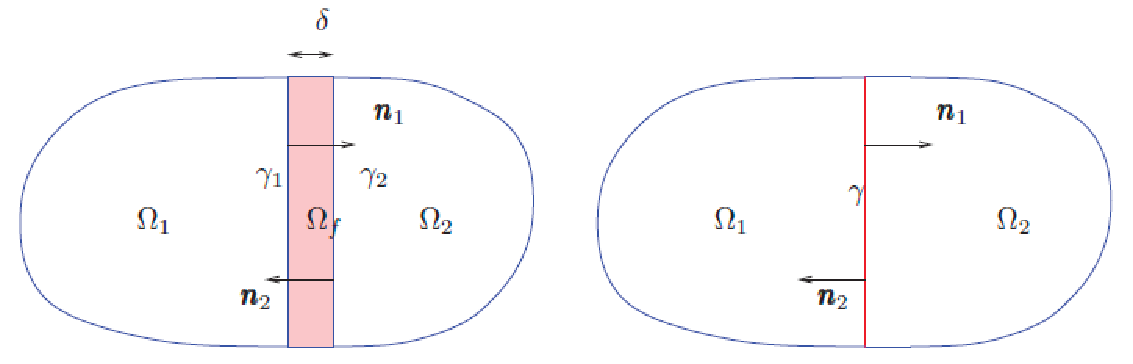}
\caption{The domain $\Omega$ with the fracture $\Omega_{f}$ (left) and the fracture-interface~$\gamma$~(right).}
\label{old_new_domain}  
\end{figure}

We let $\Omega$ be a bounded domain in $\mathbb{R}^{2}$ with Lipschitz boundary $\partial\Omega$, and $T>0$ be some fixed time.  Consider the flow problem of a single phase, compressible fluid written in mixed form as follows:
\begin{equation}
\label{original_problem}
\begin{array}{clll}
\partial_t{p} + \text{div } \bu & = & q &\text{ in } \Omega \times (0, T), \\
\bu & = & -\bK \nabla{p} & \text{ in } \Omega \times (0, T), \\
p & = & 0 & \text{ on } \partial\Omega \times (0, T), \\
p(\cdot, 0) & = & p_0 & \text{ in } \Omega,
\end{array}
\end{equation}
where $p$ is the pressure, $\bu$ the velocity, $q$ the source term, $\phi$ the storage coefficient, and $\bK$ a symmetric, time-independent, hydraulic, conductivity tensor. Suppose that the fracture $\Omega_f$ is a subdomain of $\Omega$, whose thickness is $\delta$, that separates $\Omega$ into two connected subdomains: $\Omega \backslash \overline{\Omega}_f = \Omega_1 \cup \Omega_2, $ and $\Omega_1 \cap \Omega_2 = \emptyset.$ For simplicity, we assume further that $\Omega_f$ can be expressed as \vspace{-0.2cm}
$$
\Omega_f  = \left\{ \textbf{\textit{x}} \in \Omega: \textbf{\textit{x}} = \textbf{\textit{x}}_{\gamma} + s\bn \text{ where } \textbf{\textit{x}}_{\gamma} \in \gamma \text{ and } s \in \left({-\dfrac{\delta}{2}, \dfrac{\delta}{2}}\right)\right\}, \vspace{-0.1cm}
$$
where $\gamma$ is the intersection between a line with $\Omega$ (see Figure~\ref{old_new_domain}). 

We denote by $\gamma_i$ the part of the boundary of $\Omega_i$ shared with the boundary of the fracture $\Omega_f$: $\, \gamma_i = (\partial\Omega_i \cap \partial\Omega_f) \cap \Omega$,  for $i =1,2$.  Let $\bn_i$ be the unit, outward pointing, normal vector field on $\partial\Omega_i$, where $\bn= \bn_1 = - \bn_2$. For $i=1,\; 2,\; f$, and for any scalar, vector, or tensor valued function $\varphi$ defined on $\Omega$, we denote by $\varphi_i$ the restriction of $\varphi$ to $\Omega_i$. The original problem \eqref{original_problem} can be rewritten as the following transmission problem: \vspace{-0.2cm}
\begin{equation}
\label{multidomain_problem}
\begin{array}{cllll}
\partial_t{p_i} + \text{div }\bu_i & = & q_{i} & \text{ in } \Omega_i \times (0, T), & i = 1, 2, f, \\
\bu_i & = & - \bK_i \nabla{p}_i & \text{ in } \Omega_i \times (0, T), & i = 1, 2, f, \\
p_i & = & 0 & \text{ on } \left(\partial\Omega_i \cap \partial\Omega\right) \times (0, T), & i=1, 2, f, \\
p_i & = & p_f & \text{ on } \gamma_i \times (0, T), & i=1, 2, \\
\bu_i\cdot \bn_i & = & \bu_f\cdot\bn_i & \text{ on } \gamma_i \times (0, T), & i=1, 2, \\
p_i(\cdot, 0) & = & p_{0, i} & \text{ in } \Omega_i, & i=1, 2, f.
\end{array} \vspace{-0.2cm}
\end{equation}
The reduced fracture model that we consider in this section was first proposed in~\cite{Alboin1999, Alboin2002} under the assumption that the fracture has larger permeability than that in the rock matrix. The model is obtained by averaging across the transversal cross sections of the two-dimensional fracture $\Omega_f$. We use the notation $\nabla_{\tau}$ and $\text{div}_{\tau}$ for the tangential gradient and tangential divergence, respectively. The reduced model consists of equations in the subdomains, \vspace{-0.2cm}
\begin{equation}
\label{reduced_subdomain}
\left.\begin{array}{rcll}
\partial_t{p_i}+\text{div }\bu_i&=&q_{i} &\text{ in } \Omega_i\times (0, T), \\
\bu_i&=&-\bK_i\nabla{p_i} &\text{ in } \Omega_i\times (0, T), \\
p_i&=&0 &\text{ on } \left(\partial\Omega_i \cap \partial\Omega\right) \times (0, T), \\
p_i&=&p_{\gamma} &\text{ on } \gamma \times (0, T), \\
p_i(\cdot, 0)&=&p_{0, i} &\text{ in } \Omega_i, 
\end{array}\right.  \vspace{-0.2cm}
\end{equation}
for $i=1,2,$ and equations in the fracture-interface $\gamma$, \vspace{-0.2cm}
\begin{equation}
\label{reduced_fracture}
\begin{array}{rcll}
\partial_t{p_{\gamma}}+\text{div}_{\tau }\bu_{\gamma}&=& q_{\gamma} +\sum\limits^{2}_{i=1}\left( \bu_i \cdot \bn_i\right)_{\vert \gamma} & \text{ in } \gamma \times (0, T), \\
\bu_{\gamma} &=&{-\bK_{f, \tau}\delta\nabla_{\tau}p_{\gamma}} & \text{ in } \gamma \times (0, T), \\
p_{\gamma}&=&0 &\text{ on } \partial\gamma \times (0, T), \\
p_{\gamma}(\cdot, 0)&=&p_{0, \gamma} & \text{ in } \gamma,
\end{array} \vspace{-0.2cm}
\end{equation}
where $\pmb{K}_{f, \tau}$ is the tangential component of $K_f$, \Ph{$p_{\gamma}$, $\bu_{\gamma}$, and $q_{\gamma}$ are the reduced pressure, flux, and the source term, respectively, which are given by} 
\begin{align*}
\fontsize{10pt}{10pt}\selectfont
\begin{array}{l}
{p_{\gamma} = \dfrac{1}{\delta}\mathlarger{\int}^{\delta/2}_{-\delta/2} p_{f}(x_{\gamma}+s\bn)ds}, \; \; \; \; \Ph{\bu_{\gamma} = \mathlarger{\int}^{\delta / 2}_{-\delta / 2}\bu_{f, \tau}(x_{\gamma}+s\bn)ds}, \; \; \; \; \Ph{q_{\gamma} = \mathlarger{\int}^{\delta / 2}_{-\delta / 2} q_{f}(x_{\gamma}+s\bn)ds},
\end{array}
\end{align*}
{where $\bu_{f, \tau}$ is the tangential component of $\bu_f$}. 

\vspace{0.5em}

Next, we derive the weak formulation of \eqref{reduced_subdomain}-\eqref{reduced_fracture} and its fully discrete version. The discretized equations are then reformulated into a matrix system to serve as the forward model for the data assimilation framework. To this end, we introduce the following notation: we use the convention that if $V$ is a space of functions, then $\pmb{V}$ is a space of vector functions having each component in $V$.  For arbitrary domain $\mathcal{O}$, we denote by $\left(\cdot, \cdot\right)_{\mathcal{O}}$ the inner product in $L^2\left(\mathcal{O}\right)$ or $\mathbf{\textbf{\textit{L}}^2\left(\mathcal{O}\right)}$. 
We next define the following Hilbert spaces: \vspace{-0.2cm}
\begin{align*}
\begin{array}{rl}
M &= \left\{v  = \left( v_1, v_2, v_{\gamma}\right) \in L^2\left(\Omega_1\right) \times L^2\left(\Omega_2\right) \times L^2\left(\gamma\right) \right\},  \vspace{0.2cm}\\
\Sigma &= \left\{\bv = \left(\bv_1, \bv_2, \bv_{\gamma}\right) \in \textbf{\textit{L}}^2\left(\Omega_1\right) \times \textbf{\textit{L}}^2\left(\Omega_2\right) \times \textbf{\textit{L}}^2\left(\gamma\right): \; \text{div} \; \bv_i \in L^2\left(\Omega_i\right), \; i=1,2, \right. \\
& \qquad \left.\text{and } \text{div}_{\tau} \; \bv_{\gamma} - \sum\limits^2_{i=1} \bv_i\cdot\bn_{i \vert \gamma} \in L^2(\gamma)\right\}.
\end{array}
\end{align*}
We define the bilinear forms $a(\cdot, \cdot)$, $b(\cdot, \cdot)$ and $c(\cdot, \cdot)$ on $\Sigma \times \Sigma$,  $\Sigma \times M$, and $M \times M$, respectively, and the linear form $L_q$ on $M$ by \vspace{-0.2cm}
\begin{equation}
\label{BilinearForm}
\begin{array}{l}
{a\left(\bu, \bv \right)} = \sum\limits^2_{i=1} \left(\bK^{-1}_i\bu_i, \bv_i\right)_{\Omega_i} + \left(\left({\bK_{f, \tau}\delta}\right)^{-1}\bu_{\gamma}, \bv_{\gamma}\right)_{\gamma}, \vspace{0.1cm} \\
b\left(\bu, \mu\right) = \sum\limits^2_{i=1}\left(\text{div}\; \bu_i, \mu_i\right)_{\Omega_i} + \left(\text{div}_{\tau} \; \bu_{\gamma} - \sum\limits^2_{i=1} \bu_i\cdot\bn_{i \vert \gamma}, \mu_{\gamma}\right)_{\gamma},\\
c(\eta, \mu) = \sum\limits^2_{i=1}\left(\eta_i, \mu_i\right)_{\Omega_i} + \left(\eta_{\gamma}, \mu_{\gamma}\right)_{\gamma}, \quad L_{q}(\mu) = \sum\limits^2_{i=1}(q_i, \mu_i)_{\Omega_i}+(q_\gamma, \mu_\gamma)_{\gamma}. \vspace{-0.1cm}
\end{array}
\end{equation}
The weak form of \eqref{reduced_subdomain}-\eqref{reduced_fracture} can be written as follows: 

Find $p \in H^1(0, T; M)$ and $\bu \in L^2(0, T;\Sigma)$ such that \vspace{-0.1cm}
\begin{equation}
\label{weak_reduced}
\begin{array}{rcll}
{a\left(\bu, \bv \right)} - b\left(\bv, p \right) & = & 0 & \forall \bv \in \Sigma, \\
c\left(\partial_{t}p, \mu\right) + b\left(\bu, \mu\right) &=& L_q(\mu) & \forall \mu \in M,
\end{array} \vspace{-0.1cm}
\end{equation}
together with the initial conditions: \vspace{-0.1cm}
\begin{equation}
\label{initial_weak_reduced}
p_i(\cdot, 0)  =  p_{0, i}, \; \text{in} \; \Omega_i, \; i=1,2,  \quad \text{and} \quad p_{\gamma}(\cdot, 0)  =  p_{0, \gamma}, \; \text{in} \; \gamma.\vspace{-0.1cm}
\end{equation}
The well-posedness of problem \eqref{weak_reduced}-\eqref{initial_weak_reduced} was proved in~\cite{Hoang2016}.
To find the numerical solutions to \eqref{weak_reduced}-\eqref{initial_weak_reduced}, we discretize the problem in space using mixed finite element method (see, e.g., \cite{Boffi2013, Brezzi1991, Roberts1991}) and in time using backward Euler method. To this end, let $\mathcal{K}_{h, i}$ be a finite element partition of $\Omega_i \;(i=1, 2)$ into triangles. We denote by $\mathcal{G}_{h, i}$ the set of the edges of elements $\mathcal{K}_{h, i}$ lying on the interface~$\gamma$. Since $\mathcal{K}_1$ and $\mathcal{K}_2$ coincide on $\gamma$, the spaces $\mathcal{G}_{h, 1}$ and  $\mathcal{G}_{h, 2}$ are identical, thus we set $\mathcal{G}_h := \mathcal{G}_{h, 1} = \mathcal{G}_{h, 2}$. For $i=1, \, 2$, we consider the lowest order Raviart-Thomas mixed finite element spaces $M_{h, i} \times \Sigma_{h, i} \subset L^2(\Omega_i) \times H(\text{div}, \Omega_i)$:
\begin{align*}
\fontsize{8pt}{8pt}\selectfont
\vspace{-0.2cm}
&M_{h, i}  =\left\{\mu_{h, i} \in L^2(\Omega_i): \; \mu_{h, i\vert K_i} = \text{const}, \; \; \forall K_i \in \mathcal{K}_{h, i} \right\}, \\
&\begin{array}{l}
\Sigma_{h, i} = \left\{\bv_{h, i} \in H(\text{div}, \Omega_i): \; \bv_{h, i\vert K_i} = (b_{K, i} + a_{K, i}x,  \; c_{K, i} + a_{K, i}y), \right. \vspace{0.1cm} \\
\hspace{4cm}\left.\left(a_{K, i}, \; b_{K, i}, \; c_{K, i}\right) \in \mathbb{R}^3, \; \forall K_i \in \mathcal{K}_{h, i}\right\}.
\end{array}
\end{align*}
Similarly for the fracture, {let }$M_{h, \gamma} \times \Sigma_{h, \gamma} \subset L^2(\gamma) \times H(\text{div}_{\tau}, \gamma)$ {be the lowest order Raviart-Thomas spaces in one dimension}:
\begin{align*}
& M_{h, \gamma} = \left\{ \mu_{h, \gamma} \in L^2(\gamma): \;  \mu_{h, \gamma\vert E} = \text{const}, \; \forall E \in \mathcal{G}_h \right\}, \\
& \Sigma_{h, \gamma} = \left\{ \bv_{h, \gamma} \in H(\text{div}_{\tau}, \gamma): \; \bv_{h, \gamma\vert E} = az + b, \; (a, \;b) \in \mathbb{R}^2, \; \forall E \in \mathcal{G}_h \right\}.
\end{align*}
Finally, we introduce the products of these finite element spaces:
\begin{align*}
\begin{array}{l}
\bM_h  = M_{h, 1} \times M_{h, 2} \times M_{h, \gamma}, \; \; \; \bSig_h = \Sigma_{h, 1} \times \Sigma_{h, 2} \times \Sigma_{h, \gamma}.
\end{array}
\end{align*}

For the discretization in time, we consider a uniform partition of $\Ph{(0, \; T)}$ into $N$ sub-intervals $\Ph{(t^{n}, \; t^{n+1})}$ of length $\Delta{t}=t^{n+1}-t^{n}$, for $n=0, \hdots, N-1$. The time derivatives are approximated by the backward difference quotients
$\bar{\partial}c^n = \dfrac{c^n-c^{n-1}}{\Delta{t}}$, $n= 1, \hdots, N$,
where the superscript $n$ indicates the evaluation of a function at the discrete time $t=t^n$. Altogether, the fully-discrete version of \eqref{weak_reduced}-\eqref{initial_weak_reduced} reads as follows: 

For $n=1, \hdots, N$, find $(p^n_h, \bu^n_h) \in \bM_h \times \bSig_h$ satisfying
\begin{equation}
\label{weak_discrete_reduced}
\begin{array}{rcll}
{a\left(\bu^n_h, \bv_h\right)} - b\left(\bv_h, p^n_h \right) & = & 0 & \forall \bv_h \in \bSig_h, \vspace{0.2cm} \\
c\left(\bar{\partial}p^n_h, {\mu}_h\right) + b\left(\bu^n_h, {\mu}_h\right) &=& L_q^n({\mu}_h) & \forall \mu_h \in \bM_h,
\end{array} \vspace{-0.2cm}
\end{equation}
together with the initial conditions: \vspace{-0.1cm}
\begin{equation}
\label{initial_weak_discrete_reduced}
p^0_{h, i\vert K_i} := \dfrac{1}{\vert{K_i}\vert }\int_{K_i}p_{0, i}, \; \forall K_i \in \mathcal{K}_{h, i} \; i=1,2, \; \text{and} \; p^0_{h, \gamma\vert E} = \dfrac{1}{\vert{E} \vert} \int_{E}p_{0, \gamma}, \; \forall E \in \mathcal{G}_h.\vspace{-0.1cm}
\end{equation}

\begin{remark}
For the case when the fracture has lower permeability than the rock matrix and acts as a barrier between the two subdomains, the assumption that the pressure is continuous across the fracture (fourth equation of \eqref{reduced_subdomain}) is no longer valid. Moreover, the normal component of $\pmb{K}_{f}$ needs to be taken into account. We refer the reader to \cite{Jaffre2005} for the derivation of the reduced fracture model in this case. This version of the reduced model is more general compared to the one derived in \cite{Alboin2002} and can be applied to both types of fracture. We will present the model in detail in Section~\ref{sec4}, where we conduct numerical experiments with different types of fracture.
\vspace{-0.1cm}
\end{remark}

\begin{remark}
A similar process can be applied to obtain the forward solver used to solve the advection-diffusion equations. We refer the readers to \cite{Alboin2002, Toan2023b} for the derivation of the corresponding reduced model and its space-time discretization into fully discrete weak formulation. 
\end{remark}

\vspace{0.5em}

In practice, PDEs used to model physical systems often contain intrinsic model errors. For the reduced fracture model studied in this paper, several sources of errors can arise. For example, fracture representation errors may result from treating fractures as lower-dimensional interfaces, which can fail to fully capture the complexity of real fracture networks, including roughness and aperture variations. Fracture-matrix connectivity errors occur when incorrect assumptions about how fractures interact with the surrounding porous medium lead to inaccuracies in flow simulations. Also, discretization errors in fracture geometry stem from oversimplifications in fracture shape and orientation, potentially introducing inaccuracies in transport predictions. Additionally, errors due to missing or incorrect physical assumptions can arise, such as simplified flow models or the neglect of nonlinear effects from coupled thermo-hydro-mechanical-chemical (THMC) processes. Since the precise mathematical representation of these errors is generally unknown due to limited physical knowledge and computational capability, they are often modeled as random perturbations to the PDE system. As a result, computational methods capable of filtering out uncertainties associated with these model errors are necessary.

Beyond model errors, the PDE system for the reduced fracture model often involves unknown parameters, such as diffusion coefficients and fracture transmissivities. To mitigate the impact of these uncertainties and obtain reliable estimates of key quantities -- such as pressure and velocity fields, along with model parameters -- a robust data assimilation framework is essential.
\vspace{0.5em}

To facilitate the incorporation of data and the execution of the data assimilation task, which will be detailed in the next section, we rewrite the system \eqref{weak_discrete_reduced}-\eqref{initial_weak_discrete_reduced} in a more compact form. For clarity and consistency with the preceding discussion, we consider the case where the permeability is isotropic, that is, $\bK_i = k_i \pmb{I}, \; i=1, 2$, and $\bK_{f, \tau} = k_f \pmb{I}$ and assume that the model depends on the diffusion coefficients $k_i$ on the subdomains $\Omega_i, \; i=1, 2$, and the fracture transmissivity given by $\alpha_{\gamma} := k_{f}\delta$. Since the model parameters are typically unobservable, we treat them as unknown and collectively represent them by a generic parameter vector $\pmb{\theta}$, where $\pmb{\theta} = \left(k_1, \; k_2, \alpha_{\gamma}\right)$. To simplify the formulation, we also assume that all parameters are space-independent and denote 
\begin{equation}
\begin{array}{l}
a\left(\cdot, \cdot, \pmb{\theta}\right) := a\left(\cdot, \cdot\right).
\end{array}
\vspace{-0.1cm}
\end{equation}
Next, let $\bar{\pmb{M}}_h$ and $\bar{\bSig}_h$ be the finite collections of basis functions in $\pmb{M}_h$ and $\bSig_h$, respectively. We introduce the following matrices resulted in the bi-linear forms given in \eqref{BilinearForm}:
\begin{align}
\begin{array}{l}
\pmb{A}_h(\pmb{\theta}) = \left(a(\br_h, \bv_h; \pmb{\theta})\right)_{\br_h, \bv_h \in \bar{\bSig}_h}, \; \pmb{B}_h = \left(-b(\bv_h, {\eta}_h)\right)_{\bv_h\in \bar{\bSig}_h, \eta_h \in \bar{\pmb{M}}_h}, \; \pmb{C}_{h, \phi} = \left(-c_{\phi}({\eta}_h, {\lambda}_h)\right)_{{\eta}_h, {\lambda}_h \in \bar{\pmb{M}}_h}. 
\end{array}
\end{align}
By replacing $v_h$ and $\mu_h$ in \eqref{weak_discrete_reduced} with the basis functions and expressing $\bu^n_h$ and $p^n_h$ in terms of those basis functions, we obtain the following matrix form for \eqref{weak_discrete_reduced}
\begin{equation}
\label{MatForm_v1}
\begin{array}{l}
\pmb{\Lambda}_h(\pmb{\theta})X_n = F(X_{n-1}, \; \pmb{G}_h),
\end{array}
\end{equation}
where
$\begin{array}{l}
\pmb{\Lambda}_h(\pmb{\theta}) =\begin{bmatrix}
\pmb{A}_h(\pmb{\theta}) & \pmb{B}_h \vspace{0.1cm} \\
\Delta{t}\left(\pmb{B}_h\right)^T  &\pmb{C}_{h, \phi}
\end{bmatrix}
\end{array}$
with $\Delta{t}$ the chosen timestep size, $\pmb{G}_h$ a vector containing the boundary conditions, and $X_n = \left(\bu^{n}_h, p^n_h\right)$ . Equivalently, we have the following recursive relation
\begin{equation}
\label{StateEqs_Discrete}
\begin{array}{l}
X_n = \pmb{\Phi}(X_{n-1}, \; \pmb{\theta}), \; n=1, 2, \hdots,
\end{array}
\end{equation} 
where $\pmb{\Phi}(X_{n-1}, \pmb{\theta}) = \left(\pmb{\Lambda}(\pmb{\theta})\right)^{-1}F(X_{n-1}, \;\pmb{G}_h).$ \textit{The operator $\pmb{\Phi}$ will be treated as the forward solver in the setting of the data assimilation problem.}

\vspace{0.5em}

In the next section, we formulate a data assimilation problem for the reduced fracture model and develop the necessary tools to address the joint state - parameter estimation problem within this framework.

\section{A Data Assimilation Framework for Joint State - Parameter Estimation}\label{sec3}
Our goal is to address the intrinsic errors and uncertainties in the reduced fracture model and accurately estimate its state variables, such as the fluid's pressure $p$ and velocity $\bu$ in \eqref{weak_reduced}, along with the parameter vector $\pmb{\theta}$. We assume that the observational data $Y_{1:n} := \left\{Y_i\right\}^{n}_{i=1}$ are available to provide partial noisy measurements up to time instant $n$ to assist in these estimations. To achieve this, we begin with a brief overview of data assimilation for joint state-parameter estimation, and we adopt the recursive Bayesian filter framework as our computational approach. 

\subsection{General framework: The joint state-parameter estimation problem}
We consider the following state-space model for a parameterized stochastic dynamical system:
\begin{equation}
\label{StateModel}
\begin{array}{l}
X_{n+1} = \pmb{\Phi}(X_n, \pmb{\theta})+\omega_n,
\end{array}
\end{equation}
where $\pmb{\Phi}: \mathbb{R}^l \times \mathbb{R}^k \rightarrow \mathbb{R}^l$ is a dynamical model parameterized by $\pmb{\theta} \in \mathbb{R}^k$, which is defined in \eqref{StateEqs_Discrete}; $X_n \in \mathbb{R}^l$ is the state of the dynamics at a time instant $n$ produced by the forward solver $\pmb{\Phi}$. The random variable $\omega_n \in \mathbb{R}^l$ accounts for model uncertainties and errors -- in addition to those arising from unknown parameters in the model's domain decomposition. For this work, $\omega$ is chosen as a standard $l-$dimensional Gaussian random variable with covariance matrix $\Theta$. More general uncertainties, including those with unknown distributions, can be handled similarly (e.g., see \cite{Bao2024d}).

In many practical scenarios, the true state $X_n$ is often not fully observed, and the parameter $\pmb{\theta}$ is unknown. To estimate both state $X_n$ and parameter $\pmb{\theta}$, we obtain partial noisy observations of $X_n$ through the following observation process
\begin{equation}
\label{ObserProcess}
\begin{array}{l}
Y_{n+1} = g(X_{n+1})+\epsilon_{n+1},
\end{array}
\end{equation}
where $g$ is the observation function, and $\epsilon_{n+1}$ is a Gaussian noise that representing observational errors with a given covariance $\Upsilon$. The goal of the joint state-parameter estimation is to determine the best estimates for both $X_n$ and $\pmb{\theta}$ based on the observations $\left\{Y_{i}\right\}^{n}_{i=1}$.

To enable an online mechanism for dynamically estimating the model state and parameters, we represent the target parameter as a sequential process, which reads as follows
\begin{equation}
\label{ParaProcess}
\begin{array}{l}
\pmb{\theta}_{n+1} = \pmb{\theta}_{n}+\xi_{n},
\end{array}
\end{equation}
where $\xi_n$ is a pre-defined noise term that controls the exploration rate of the sequence $\{\pmb{\theta}_n\}$ in the parameter space over time, and $\pmb{\theta}_n$ represents the estimate of $\pmb{\theta}$ at time instant $n$.

For state estimation, the best estimation $\bar{X}_{n+1}$ of $X_{n+1}$ is computed as the conditional expectation 
\begin{align*}
\begin{array}{l}
\bar{X}_{n+1} := \mathbb{E}\left[X_{n+1} \vert Y_{1:n+1}\right],
\end{array}
\end{align*}
conditioned on the observational data $Y_{1:n+1}$. The estimated state $\bar{X}_{n+1}$ is then used to compute the updated estimate $\pmb{\theta}_{n+1}$ for the parameter $\pmb{\theta}$. The mathematical technique that computes $\bar{X}_{n+1}$ is called ``optimal filtering'', which aims to
find the conditional distribution of the state variable $X_{n+1}$ conditioning on the observations $\left\{Y_{1:n+1}\right\}$, i.e., $P(X_{n+1} \vert Y_{1:n+1})$. In the context of the optimal filtering problem, this conditional distribution is commonly referred to as the ``filtering density (distribution)''. Even when the state dynamics parameter $\pmb{\theta}$ in \eqref{StateModel} is known, estimating the state remains a highly challenging task. One of the main difficulties lies in the complexities of high-dimensional approximation and Bayesian inference for nonlinear systems. The problem becomes even more prohibitive when there is a mismatch between the predicted state distribution (the prior) and the observational data (the likelihood), making the computation of the posterior distribution in high-dimensional spaces particularly demanding. Furthermore, when the model parameter is unknown, such discrepancies become more substantial, which makes the application of Bayesian inference more challenging.

In this study, we employ the recursive Bayesian filter framework, which is the standard approach for solving the optimal filtering problem, to dynamically estimate the filtering densities for $X_{n+1}$ and $\pmb{\theta}$. In the next section, two important methods will be discussed: the first one is the Ensemble Score Filter method (EnSF) for state estimation, while the second method is the Direct Filter for parameter estimation. Then, a United Filter framework will be introduced to effectively combine the above two procedures for the joint state-parameter estimation problem.
\subsection{The Ensemble Score Filter for state estimation} \label{EnSF}

We begin with the discussion of the Ensemble Score Filter method. In the derivation of this method, the parameter $\pmb{\theta}$ is assumed to be known, therefore, our focus is on the state estimation. 

In the Bayesian filter, given the filtering density, denoted by $p(X_n \vert Y_{1:n})$, for $X$ at time instant $n$, we generate predicted filtering density, i.e. $p(X_{n+1}\vert Y_{1:n})$, through the following Chapman-Kolmogorov formula
\begin{equation}
\label{PredictionEnSF}
\begin{array}{l}
p(X_{n+1}\vert Y_{1:n}) = \mathlarger{\int} p(X_{n}\vert Y_{1:n}) p(X_{n+1}\vert X_n) \; dX_n,
\end{array}
\end{equation}
where $p(X_{n+1}\vert X_n)$ is the transition probability governed by the dynamical model $\pmb{\Phi}$ and the random variable $\omega_n$ in \eqref{StateModel}.

Upon receiving the observational data $Y_{n+1}$, we apply the following Bayesian inference procedure to integrate the data information to the prior in the form of likelihood and obtain the updated filtering density 
\begin{equation}
\label{UpdateEnSF}
\begin{array}{l}
p(X_{n+1} \vert Y_{1:n+1}) \propto p(X_{n+1} \vert X_{n}) p(Y_{n+1} \vert X_{n+1}),
\end{array}
\end{equation}
where the likelihood $p(Y_{n+1} \vert X_{n+1})$ is defined by
\begin{equation}
\label{likelihoodEnSF}
\begin{array}{l}
p(Y_{n+1} \vert X_{n+1}) \propto \text{exp} \left[-\dfrac{1}{2}\left(g(X_{n+1})-Y_{n+1}\right)^T \Upsilon^{-1} \left(g(X_{n+1})-Y_{n+1}\right) \right],
\end{array}
\end{equation}
with $\Upsilon$ denoting the covariance of the observational noise $\epsilon_{n}$ introduced in \eqref{ObserProcess}.

The filtering density can be propagated over time by recursively performing the prediction procedure \eqref{PredictionEnSF} and the update procedure \eqref{UpdateEnSF}. However, when the problem is in high dimensions, it is very challenging to accurately characterize the filtering density and effectively incorporate new data into the filtering density.

To tackle this difficulty, the author in \cite{Bao2024c} utilized the idea of the score-based diffusion model framework to establish a connection between a filtering density and the standard normal distribution. The information of this filtering density is then stored and represented by a score function. More specifically, in \cite{Bao2024c}, two score functions, denoted by $S_{n+1\vert \; n}(\cdot, \cdot)$ and $S_{n+1 \vert n+1}(\cdot, \cdot)$ are constructed, which store the information regarding the filtering densities $p\left(X_{n+1} \vert \; Y_{1:n}\right)$ and $p\left(X_{n+1} \vert \; Y_{1:n+1}\right)$. The score-based model consists of a forward stochastic differential equation (SDE)
\begin{equation}
\label{ForwardSDE}
\begin{array}{l}
dZ_t = b(t) Z_t dt +\sigma(t) dW_t,
\end{array}
\end{equation}
a reverse-time SDE defined in a pseudo-time domain $t \in [0, T]$
\begin{equation}
\label{BackwardSDE}
dZ_t = \left[b(t)Z_t - \sigma^2(t) S(Z_t, t)\right]\; dt +\sigma(t) d\overleftarrow{W}_t,
\end{equation}
where $W_t$ is  a standard $p$- dimensional Brownian motion, $\int \cdot dW_t$ denotes a standard It\^{o} type stochastic integral, $\int \cdot d\overleftarrow{W}_t$ is a backward It\^{o} integral, and $b$ and $\sigma$ are the drift coefficient and the diffusion coefficient, respectively. 

In this work, we let the terminal time $T$ in the diffusion model to be $T=1$, and the score function $S(\cdot, \cdot)$ in \eqref{BackwardSDE} is defined by
\begin{align*}
S(z, t) := \nabla\log (Q_t(z)), \; z \in \mathbb{R}^p,
\end{align*}
where $Q_t$ is the probability distribution of solution $Z_t$ of the forward SDE \eqref{ForwardSDE}. The drift and the diffusion in \eqref{ForwardSDE} are given by 
\begin{align*}
b(t) = \dfrac{d\log \alpha_t}{dt}, \; \; \; \sigma^2(t) = \dfrac{d\beta^2_t}{dt} - 2\dfrac{d\log \alpha_t}{dt} \beta^2_t,
\end{align*}
with $\alpha_t = 1-t$ and $\beta_t = t$ for $t \in [0, 1]$. With these two functions, the forward SDE \eqref{ForwardSDE} transforms any target distribution $Q_0$ of $Z_0$ to the standard Gaussian distribution $\mathcal{N}(0, I_p)$. Moreover, we have the opposite direction, that is, from any sample drawn from $\mathcal{N}(0, I_p)$, we can generate as many sample from $Z_0$ as we want by solving the reverse SDE \eqref{BackwardSDE}. For a detailed explanation on this behavior, we refer the readers to \cite{Bao2024c}. The score function also fully characterizes the information of the initial distribution $Q_0(Z_0)$. More precisely, we have
\begin{equation}
\label{AnalyticScore}
\begin{array}{rl}
S(Z_t, t) &= \nabla \log(Q_t(Z_t)) = \nabla \log\left(\mathlarger{\int}_{\mathbb{R}^p}Q_t(Z_t \vert Z_0)Q_0(Z_0)\; dZ_0\right) \vspace{0.1cm} \\
&= \mathlarger{\int}_{\mathbb{R}^p} -\dfrac{Z_t-\alpha_tZ_0}{\beta^2_t}w_t(Z_t, Z_0)Q_0(Z_0)\; dZ_0,
\end{array}
\end{equation}
where 
\begin{align*}
w_t(Z_t, Z_0):= \dfrac{Q_t(Z_t\vert \; Z_0)}{\int_{\mathbb{R}^p}Q_t(Z_t\vert\; Z'_0)Q_0(Z'_0)\; dZ'_0}
\end{align*}
is a weight factor.

Based on \cite{Bao2024c}, we define the following score functions which are corresponding to the density $p\left(X_{n+1}\vert Y_{1:n}\right)$ and $p\left(X_{n+1}\vert Y_{1:n+1}\right)$, respectively, as follows:
\vspace{0.25em}

\begin{enumerate}[label = \roman*)]
    \item Prior filtering score $S_{n+1\vert \; n}$: set $Z_0 = X_{n+1} \vert Y_{1:n}$ and hence $Q_0(Z_0) = p(X_{n+1}\vert \; Y_{1:n})$,
    \item Posterior filtering score: $S_{n+1 \vert \; n+1}$: set $Z_0 = X_{n+1} \vert \; Y_{1:n+1}$ and hence $Q_0(Z_0) = p(X_{n+1}\vert \; Y_{1:n+1})$.
\end{enumerate}
\vspace{0.25em}

To obtain the score $S_{n+1 \vert \; n}$, data samples for $Z_0 = X_{n+1} \vert \; Y_{1:n}$ are needed. To this end, we first generate samples $\{x^{j}_{n \vert \; n}\}^J_{j=1}$ from the posterior filtering scores $\bar{S}_{n \vert \; n}$ at time instant $n$, where $\bar{S}_{n \vert \; n}$ is an approximation of $S_{n \vert \; n}$. The samples $\{x^{j}_{n \vert \; n}\}^J_{j=1}$ are then input into the state model \eqref{StateModel} to generate the samples for $Z_0 = X_{n+1} \vert Y_{1: n}$:
\begin{equation}
\label{ParaPredict}
\begin{array}{l}
x^{j}_{n+1 \vert \; n} = \pmb{\Phi}\left(x^{j}_{n \vert \; n}, \pmb{\theta}\right) + \omega^j_n, \; j=1, \hdots, J.
\end{array}
\end{equation}
The remaining task is to estimate the score function $S_{n+1\vert \; n}$. In \cite{Bao2024c}, the score $S_{n+1\vert \; n}$ was obtained by using deep neural networks. However, there is one major drawback, that is, the neural network
used to learn the score function needs to be re-trained at every filtering step after assimilating new observational data. To tackle this difficulty, a training-free ensemble approach, namely the Ensemble Score Filter, was derived in \cite{Bao2024a} for approximating the score $S_{n+1\vert \; n}$. This approach uses the analytic expression \eqref{AnalyticScore} for the score model to obtain an approximation of $S_{n+1\vert \; n}$ as follows:
\begin{equation}
\label{ApproPrior}
\begin{array}{l}
S_{n+1\vert \; n}(z, t) \approx \bar{S}_{n+1\vert \; n}(z, t) := \sum\limits^{O}_{o=1} - \dfrac{z-\alpha_t \cdot x^{j_o}_{n+1\vert \; n}}{\beta^2_t} \bar{w}_t\left(z, x^{j_o}_{n+1\vert \; n}\right),
\end{array}
\end{equation}
where $\{x^{j_o}_{n+1\vert \; n}\}^O_{o=1}$ is a mini-batch of samples in the sample set $\{x^{j}_{n+1\vert \; n}\}^J_{j=1}$ introduced by \eqref{ParaPredict}, and $\bar{w}_t$ is an approximation to weight $w_t$, which is defined as follows:
\small{
\begin{align*}
\bar{w}_t\left(z, x^{j'_o}_{n+1\vert \; n}\right):= \dfrac{Q_t\left(z\vert x^{j'_o}_{n+1\vert n}\right)}{\sum\limits^O_{o=1}Q_t\left(z\vert x^{j_o}_{n+1\vert \; n}\right)}.
\end{align*}
}

Finally, we provide an approximation of the score $S_{n+1\vert \; n+1}$. The main challenge of this process lies in the absence of the samples from the posterior filtering density. To address this issue, we analytically incorporate the likelihood information into the currect estimated score $\bar{S}_{n+1 \vert \; n}$. As a result, an approximation for $S_{n+1 \vert \; n+1}$ is given as follows:
\begin{equation}
\label{ApproPoster}
\begin{array}{l}
\bar{S}_{n+1\vert \; n+1}(z, t):= \bar{S}_{n+1\vert \; n}(z, t) + h(t) \nabla \log p\left(Y_{n+1} \vert \; z\right),
\end{array}
\end{equation}
where $\bar{S}_{n+1\vert \; n}$ is the training-free approximation introduced in \eqref{ApproPrior}, $p\left(Y_{n+1}\vert \; \cdot\right)$ is the likelihood given by \eqref{likelihoodEnSF} and $h$ is the damping function satisfying:
$h(\cdot) \text{ monotonically decreases in } [0, 1] \text{ with } h(0)=1 \text{ and } h(1) =0$.

%

With the approximated score $\bar{S}_{n+1 \vert \; n+1}$ derived in \eqref{ApproPoster}, we can generate a collection of samples $\{x^{j}_{n+1 \vert \; n+1}\}^{J}_{j=1}$ that follows the posterior filtering density $p\left(X_{n+1} \vert \; Y_{n+1}\right)$ from the Gaussian distribution $\mathcal{N}\left(0, I_p\right)$ through the reverse-time SDE \eqref{BackwardSDE}, and the estimated state can be obtained by
\begin{align*}
\begin{array}{l}
\bar{X}_{n+1} \approx \dfrac{1}{J}\sum\limits^J_{j=1} x^{j}_{n+1 \vert \; n+1}.
\end{array}
\end{align*}

In the above discussion, we assume that the parameter $\pmb{\theta}$ is known. To dynamically estimate the unknown parameter in the model $f$ upon receiving the observational data, in the next part, we explain the framework of the Direct Filter algorithm developed in \cite{Bao2019a}.
\subsection{The Direct Filter for dynamical parameter estimation}
In the setting of Direct Filter, we assume that the estimates $\{\bar{X}_n\}_n$ of $\{X_n\}_n$ are available, and these estimates shall be used to derive an estimate for the parameter $\pmb{\theta}$. More specifically, the following conditional expectation
\small{
\begin{align*}
\bar{\pmb{\theta}}_{n+1} = \mathbb{E}\left[ \pmb{\theta}_{n+1} \vert \bar{X}_{1:n+1}\right]
\end{align*}
}
is considered as an estimate for $\pmb{\theta}_{n+1}$. To this end, we aim to provide an estimate of the parameter distribution $p\left(\pmb{\theta}_{n+1} \vert \bar{X}_{1:n+1}\right)$ for $\pmb{\theta}_{n+1}$. It is done by using the following Bayesian formula 
\begin{align*}
\begin{array}{l}
p\left(\pmb{\theta}_{n+1} \vert \bar{X}_{1:n+1}\right) \propto p\left(\pmb{\theta}_{n+1}\vert \bar{X}_{1:n}\right)p\left(\bar{X}_{n+1} \vert \pmb{\theta}_{n+1}\right),
\end{array}
\end{align*}
where $p\left(\pmb{\theta}_{n+1}\vert \bar{X}_{1:n}\right)$ is the prior parameter distribution obtained through the pseudo process \eqref{ParaProcess} and the previous distribution $p\left(\pmb{\theta}_n \vert \bar{X}_{1:n}\right)$, and $p\left(\bar{X}_{n+1}\vert \pmb{\theta}_{n+1}\right)$ is the likelihood that measures the discrepancy between the predicted parameter and the estimated state $\bar{X}_{n+1}$.

In this work, we utilize the Direct Filter approach to compute the distribution $p\left(\pmb{\theta}_{n+1}\vert \bar{X}_{1:n+1}\right)$. Assume that we have a collection of parameter particles, denoted by $\left\{\pmb{\theta}^m_n\right\}^M_{m=1}$, from the estimated distribution $p\left(\pmb{\theta}_n\vert \bar{X}_{1:n}\right)$, i.e.,
\small{
\begin{align*}
\begin{array}{l}
p\left(\pmb{\theta}_n\vert \bar{X}_{1:n}\right) \approx \dfrac{1}{M}\sum\limits^M_{m=1}\delta_{\pmb{\theta}^{m}_n}(\pmb{\theta}_n),
\end{array}
\end{align*}}
where $M$ is a pre-chosen positive integer, and $\delta$ is the Dirac delta function. The particles $\left\{\pmb{\theta}^m_n\right\}^M_{m=1}$ is then perturbed using \eqref{ParaProcess} to create a set of predicted particles $\{\tilde{\pmb{\theta}}^m_{n+1}\}^M_{n=1}$, which is defined as,
\begin{align*}
\begin{array}{l}
\tilde{\pmb{\theta}}^m_{n+1} = \pmb{\theta}^m_{n}+\xi^m_{n},
\end{array}
\end{align*}
where $\xi^m_{n} \sim \mathcal{N}(0, \Gamma)$ is a Gaussian type random variable with a pre-determined covariance matrix $\Gamma$. As the model parameters are not directly observed, we input each particle $\tilde{\pmb{\theta}}^m_{n+1}$ into the state model \eqref{StateModel} to obtained a parameter-dependent state 
\begin{align*}
\begin{array}{l}
\bar{X}^{\tilde{\pmb{\theta}}^m_{n+1}}_{n+1} = \pmb{\Phi}\left(\bar{X}_n, \tilde{\pmb{\theta}}^m_{n+1}\right)+\omega^m_{n},
\end{array}
\end{align*}
and use these states to compute the likelihood function
\begin{equation}
\label{LikelihoodDF}
p\left(\bar{X}_{n+1}\vert \tilde{\pmb{\theta}}^m_{n+1}\right) \propto \text{exp}\left[ -\dfrac{1}{2}\left(\bar{X}^{\tilde{\pmb{\theta}}^m_{n+1}}_{n+1}-\bar{X}_{n+1}\right)^T \Theta^{-1}\left(\bar{X}^{\tilde{\pmb{\theta}}^m_{n+1}}_{n+1}-\bar{X}_{n+1}\right)\right],
\end{equation}
where $\Theta$ is the covariance of the model uncertainty $\{\omega_n\}$ in \eqref{StateModel}. This function compares $\bar{X}^{\tilde{\pmb{\theta}}^m_{n+1}}_{n+1}$ with the estimated state $\bar{X}_{n+1}$. The likelihood $p\left(\bar{X}_{n+1}\vert \tilde{\pmb{\theta}}^m_{n+1}\right)$ then assigns a weight to the parameter particle $\tilde{\pmb{\theta}}^m_{n+1}$, which can be used to give an approximation to $p\left(\pmb{\theta}_{n+1}\vert \bar{X}_{n+1}\right)$ through the following weighted sample scheme
\begin{align*}
\vspace{-0.5em}
p\left(\pmb{\theta}_{n+1}\vert \bar{X}_{n+1}\right) \approx \sum\limits^M_{m=1} w_m\delta_{\tilde{\pmb{\theta}}^m_{n+1}}(\pmb{\theta}_{n+1}),
\vspace{-0.5em}
\end{align*}
where the weight $w_m \propto p\left(\bar{X}_{n+1}\vert \tilde{\pmb{\theta}}^m_{n+1}\right)$.

To avoid degeneracy, which means that a lot of the weights of the particles will be ignored, one needs a re-sampling step for the updated measure. The purpose of the resampling step is to generate a set of equally weighted samples to avoid the degeneracy issue (see, e.g., ~\cite{Bao2019a, Bao2019b, Bao2014, Bao2020, Bao2018, Bao2017}). In this work, we simply use the importance sampling method~(\cite{Doucet2001b, Morzfeld2018}) to generate additional duplicates of particles in $\{\tilde{\pmb{\theta}}^m_{n+1}\}^M_{m=1}$, denoted as $\{\pmb{\theta}^m_{n+1}\}^M_{m=1}$, where particles with smaller weights are replaced by particles with larger weights. This can be done by sampling with replacement $M$ particles $\{\tilde{\theta}^{m}_{n+1}\}^M_{m=1}$ according to the weights $\left\{w_{m}\right\}^M_{m=1}$ using the distribution $p\left(\pmb{\theta}_{n+1}\vert \bar{X}_{n+1}\right)$. Then, we obtain the following estimation for $\pmb{\theta}_{n+1}$:
$\bar{\pmb{\theta}}_{n+1} = \dfrac{1}{M}\sum\limits^M_{n=1}\pmb{\theta}^m_{n+1}$,
which is an updated parameter for the next filtering step.
\begin{remark}
Throughout Section~\ref{sec4}, the values of the perturbation noise ${\epsilon_n}$ in \eqref{ParaProcess} will be selected empirically. However, there exists a systematic approach for determining these values. Specifically, the method proposed in \cite{Liu2001} utilizes the kernel smoothing approach with location shrinkage introduced in \cite{West1993} to establish a relationship between the variance matrix of the perturbation noise $\epsilon_{n+1}$ at time step $n+1$ and the variance matrix of the particle parameters $\theta_n$ at the previous time step $n$. This relationship provides a basis for systematically selecting the noise values.
\end{remark}

\subsection{The United Filter framework for joint state-parameter estimation}
In the following, we introduce an iterative algorithm, namely the United Filter \cite{Bao2024b}, which unifies the EnSF and the direct filter method to solve the joint state-parameter estimation problem \eqref{StateModel}-\eqref{ObserProcess} numerically.

At time instant $n$, we assume that an estimated state $\bar{X}_n$ is available, accompanied by a posterior score $\bar{S}_{n\vert n}$ corresponding to the posterior filtering density $p(X_n \vert Y_{1:n})$. We assume further that we have a set of parameter particles $\{\pmb{\theta}^m_n\}^M_{m=1}$ describing the parameter $\pmb{\theta}_n$, and we denote the mean of this set as $\bar{\pmb{\theta}}_n$. The United Filter algorithm for joint state-parameter estimation at time instant $n+1$ consists of two stages: $(\text{Stage I})$ state estimation with the current approximated parameter, and $(\text{Stage II})$ parameter estimation based on the updated state estimate. To further refine the results from the filtering process at time instant $n+1$, these two stages are executed within a nested iteration procedure, indexed by $r = 0, 1, \hdots, R$. This nested loop is initialized by setting the initial estimate for $\pmb{\theta}_{n+1}$ as $\bar{\pmb{\theta}}^{(0)}_{n+1} = \bar{\pmb{\theta}}_{n}$ with parameter particles chosen as $\{\pmb{\theta}^{m, (0)}_{n+1}\}^M_{m=1} : = \{\pmb{\theta}^m_{n=1}\}^M_{m=1}$. Similar to the perturbation noise $\epsilon_n$ in \eqref{ParaProcess}, the number $R$ will be chosen empirically in our numerical experiments presented in Section~\ref{sec4}. We observe numerically that increasing  $R$ leads to better results. However, a rigorous theoretical justification for this observation remains an open question and will be explored in future research.

\textbf{Stage I.} \textbf{State estimation with the current estimated parameter} $\pmb{\bar{\pmb{\theta}}^{(r)}_{n+1}}$ \textbf{to obtain} $\pmb{\bar{X}^{(r)}_{n+1}}$.

At iteration $r$ for the time period $n$ to $n+1$, we have a set of $J$ states $\{x^j_{n \vert n}\}^J_{j=1}$ generated from the reverse-time SDE with score model $\bar{S}_{n \vert n}$. Each sample from this collection is then propagated through that state dynamics \eqref{StateModel} based on the current estimated parameter $\bar{\pmb{\theta}}^{(r)}_{n+1}$ to obtain a set of predicted state samples $\{x^{(r), j}_{n \vert n}\}^J_{j=1}$ for the prior filtering density, i.e, 
\begin{align*}
\begin{array}{l}
x^{(r), j}_{n+1\vert \; n} = \pmb{\Phi}\left(x^{j}_{n \vert n}, \bar{\pmb{\theta}}^{(r)}_{n+1}\right) +\omega^{(r)}_{n, j}.
\end{array}
\end{align*}

The equation \eqref{ApproPrior} is then applied to estimate the prior score $\bar{S}^{(r)}_{n+1\vert n}$. To integrate the observational data $Y_{n+1}$, the formula \eqref{ApproPoster} is utilized to approximate the posterior score $\bar{S}^{(r)}_{n+1\vert n+1}$. As a consequence, the posterior state samples $\{x^{(r), j}_{n+1 \vert n+1}\}^J_{j=1}$ can be generated from the score $\bar{S}^{(r)}_{n+1 \vert \; n+1}$. These samples are then used to give an updated estimate for $X_{n+1}$, denoted by $\bar{X}^{(r)}_{n+1}$. 

\textbf{Stage II.} \textbf{Parameter estimation to obtain} $\pmb{\bar{\pmb{\theta}}^{(r+1)}_{n+1}}$ \textbf{with the estimated state} $\pmb{\bar{X}^{(r)}_{n+1}}$.

The direct filter based parameter estimation relies on the current estimate $\bar{X}^{(r)}_{n+1}$ for the state variable $X_{n+1}$, and the scheme is slightly modified by replacing $\bar{X}_{n+1}$ in \eqref{LikelihoodDF} by $\bar{X}^{(r)}_{n+1}$. The predicted parameter particles for the next $(r+1)-$th iteration state are generated by using the pseudo process \eqref{ParaProcess} 
\begin{align*}
\begin{array}{l}
\tilde{\pmb{\theta}}^{m, (r+1)}_{n+1} = \pmb{\theta}^{m, (r)}_{n}+\xi^{m, (r+1)}_{n}, \; m=1, \hdots, M,
\end{array}
\end{align*}
and the likelihood in \eqref{LikelihoodDF} is rewritten as
\begin{equation}
p\left(\bar{X}^{(r)}_{n+1}\vert \tilde{\pmb{\theta}}^{m, (r+1)}_{n+1}\right) \propto \text{exp}\left[ -\dfrac{1}{2}\left(\bar{X}^{\tilde{\pmb{\theta}}^{m, (r+1)}_{n+1}}_{n+1}-\bar{X}^{(r)}_{n+1}\right)^T \Theta^{-1}\left(\bar{X}^{\tilde{\pmb{\theta}}^{m, (r+1)}_{n+1}}_{n+1}-\bar{X}^{(r)}_{n+1}\right)\right],
\end{equation}
to compare $\bar{X}^{\tilde{\pmb{\theta}}^{m, (r+1)}_{n+1}}_{n+1}$ with $\bar{X}^{(r)}_{n+1}$, and to assign a weight to the particle $\tilde{\pmb{\theta}}^{m, (r+1)}_{n+1}$. After the resampling procedure, we obtain a set of equally weighted parameter particles $\{\pmb{\theta}^{m, (r+1)}_{n+1}\}^M_{m=1}$, and the estimated parameter at the $(r+1)-$th iteration stage is 
\begin{align*}
\bar{\pmb{\theta}}^{(r+1)}_{n+1} = \dfrac{1}{M}\sum\limits^M_{m=1}\pmb{\theta}^{m, (r+1)}_{n+1}.
\end{align*}

The two stages $(\text{Stage I})$ and $(\text{Stage II})$ are repeated sequentially throughout $R$ iterations, and the observational data $Y_{n+1}$ is appropriately assimilated into both the estimated state and the estimated parameter. Subsequently, the estimated state and the estimated parameter at time instant $n+1$ are defined as
\begin{align*}
\bar{X}_{n+1}:= \bar{X}^{(R)}_{n+1},
\end{align*}
and  
\begin{align*}
\bar{\pmb{\theta}}_{n+1}:= \bar{\pmb{\theta}}^{(R)}_{n+1},
\end{align*}
respectively.
\subsection{Summary of the United Filter algorithm for the reduced fracture model}
We conclude Section~\ref{sec3} by presenting the pseudo algorithm of the United Filter method for the joint state-parameter estimation for the reduced fracture model. The entire algorithm is summarized in Algorithm~\ref{UnitedFilter} below.
\begin{algorithm}[H]
    \caption{The United Filter Algorithm}\label{UnitedFilter}    
\begin{algorithmic}[1]
\State \textbf{Input:} the numerical solver $\pmb{\Phi}$ for the reduced fracture PDE system as the state dynamics, the observation function $g(X)$, distribution $P(X_0)$ as the guess for the initial state of the PDE, and the initial parameter particles $\left\{\pmb{\theta}^m_0\right\}^{M}_{m=1}$;
\State Generate $J$ samples $\left\{x^{j}_{0\vert \; 0}\right\}^J_{j=1}$ from the prior $P(X_0)$;
\For{$n=0, \hdots,$}
\State Set $\bar{\pmb{\theta}}^{(0)}_{n+1} = \bar{\pmb{\theta}}_n$ and $\left\{\pmb{\theta}^{m, (0)}_{n+1}\right\}^M_{m=1} = \left\{\pmb{\theta}^m_{n}\right\}^M_{m=1}$;
\For{$r =0, \hdots, R-1$}
\State Perform the state estimation using the EnSF algorithm with the current estimated parameter $\bar{\pmb{\theta}}^{(r)}_{n+1}$ to obtain the estimated state $\bar{X}^{(r)}_{n+1}$ as well as $\bar{S}^{(r)}_{n+1 \vert \; n+1}$;
\State Carry out the Direct Filter algorithm for the parameter estimation using the updated estimated state $\bar{X}^{(r)}_{n+1}$ to obtain the updated estimated parameter $\bar{\pmb{\theta}}^{(r+1)}_{n+1}$ as well as the new particles $\left\{\pmb{\theta}^{m, (r+1)}_{n+1}\right\}^M_{m=1}$;
\EndFor
\State Execute the EnSF based state estimation using the estimated parameter $\bar{\pmb{\theta}}^R_{n+1}$ to obtain $\bar{X}^{(R)}_{n+1}$ as well as $\bar{S}^{(R)}_{n+1 \vert \; n+1}$;
\State Set $\bar{\pmb{\theta}}_{n+1} = \bar{\pmb{\theta}}^{(R)}_{n+1}$, \; $\left\{\pmb{\theta}^{m}_{n+1}\right\}^M_{m=1} = \left\{\pmb{\theta}^{m, (R)}_{n+1}\right\}^M_{m=1}, \; \bar{X}_{n+1} = \bar{X}^{(R)}_{n+1}$, and $\bar{S}_{n+1\vert \; n+1} = \bar{S}^{(R)}_{n+1\vert \; n+1}$.
\EndFor
\end{algorithmic}
\end{algorithm}
\section{Numerical Experiments}\label{sec4}
In this section, we conduct numerical experiments to demonstrate the effectiveness of the proposed United Filter approach in producing accurate state estimates and robust parameter estimates. Furthermore, we compare the United Filter with a state-of-the-art online joint state-parameter estimation method, namely the augmented ensemble Kalman filter (AugEnKF), to showcase the superior performance of the United Filter in solving the state-parameter estimation problem for reduced fracture models.

Specifically, we perform three numerical experiments: Test Cases 1 and 2 focus on the reduced fracture model for the pure diffusion equation, while Test Case 3 addresses the advection-diffusion equation. In Test Case 1, the fracture permeability is isotropic, whereas in Test Case 2, the fracture is anisotropic and heterogeneous. Across all experiments, we use 200 ensembles for the EnSF with 200 time steps in the reverse-time SDE. The number of particles for the Direct Filter algorithm and the iterations for the United Filter process are determined based on the specific requirements of each test case. 
\begin{remark}
The proposed method is implemented in Python. All the numerical results presented in this section can be reproduced using the code on GitHub. The source code is publicly available at \url{https://github.com/Toanhuynh997/UF_ReducedFrtModel}. To plot the estimated solutions in Test case 1 and Test case 2, we adapt the codes from the package DarcyLite developed in \cite{JLiu2016}.
\end{remark}
\subsection{Test Case 1: Pure diffusion equation with a single fracture}
\begin{figure}[h!]
\centering
\begin{minipage}{.5\textwidth}
  \centering  \includegraphics[width=0.85\linewidth]{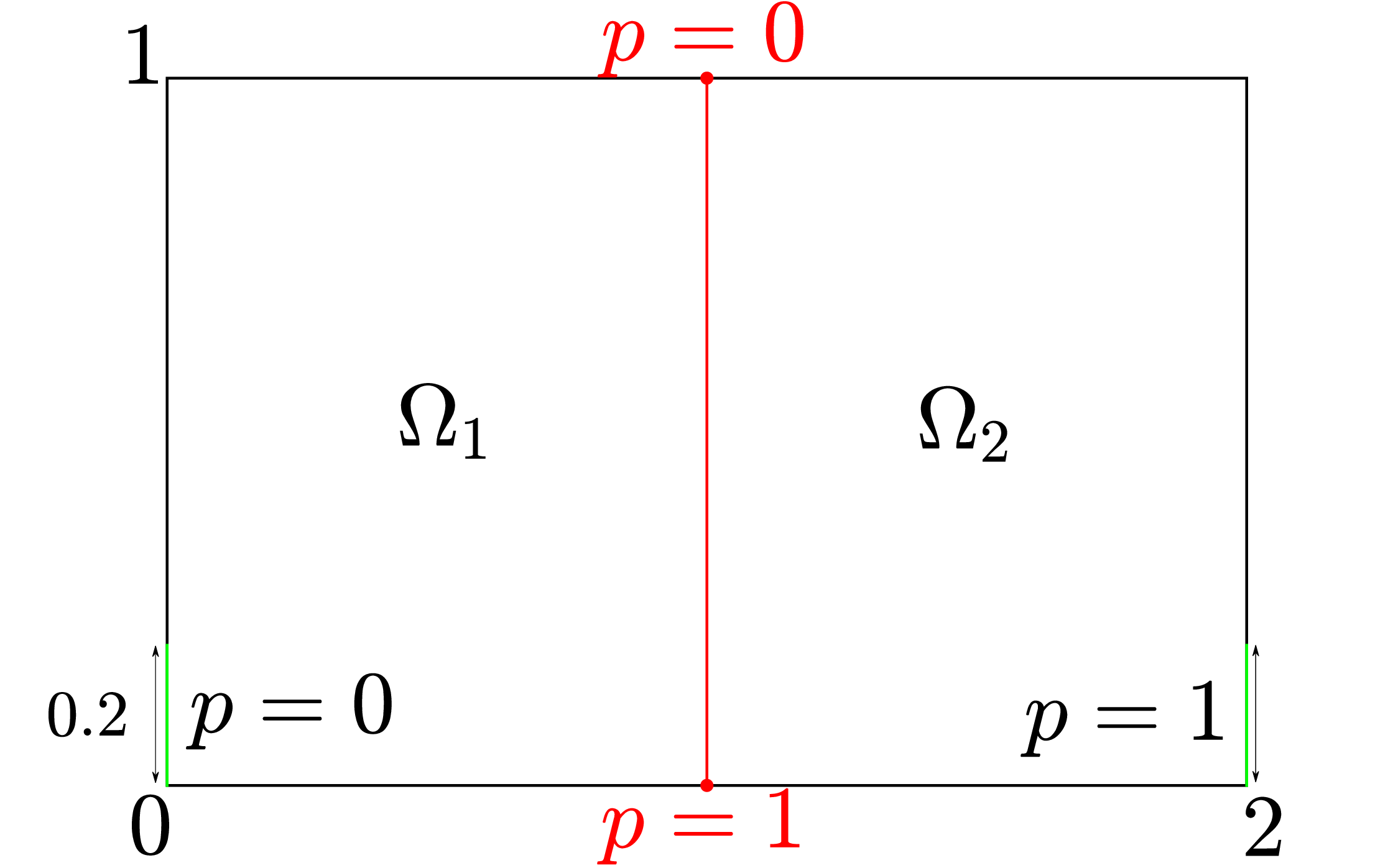}
  \label{fig:test1}
\end{minipage}%
\begin{minipage}{.5\textwidth}
  \centering
  \includegraphics[width=0.85\linewidth]{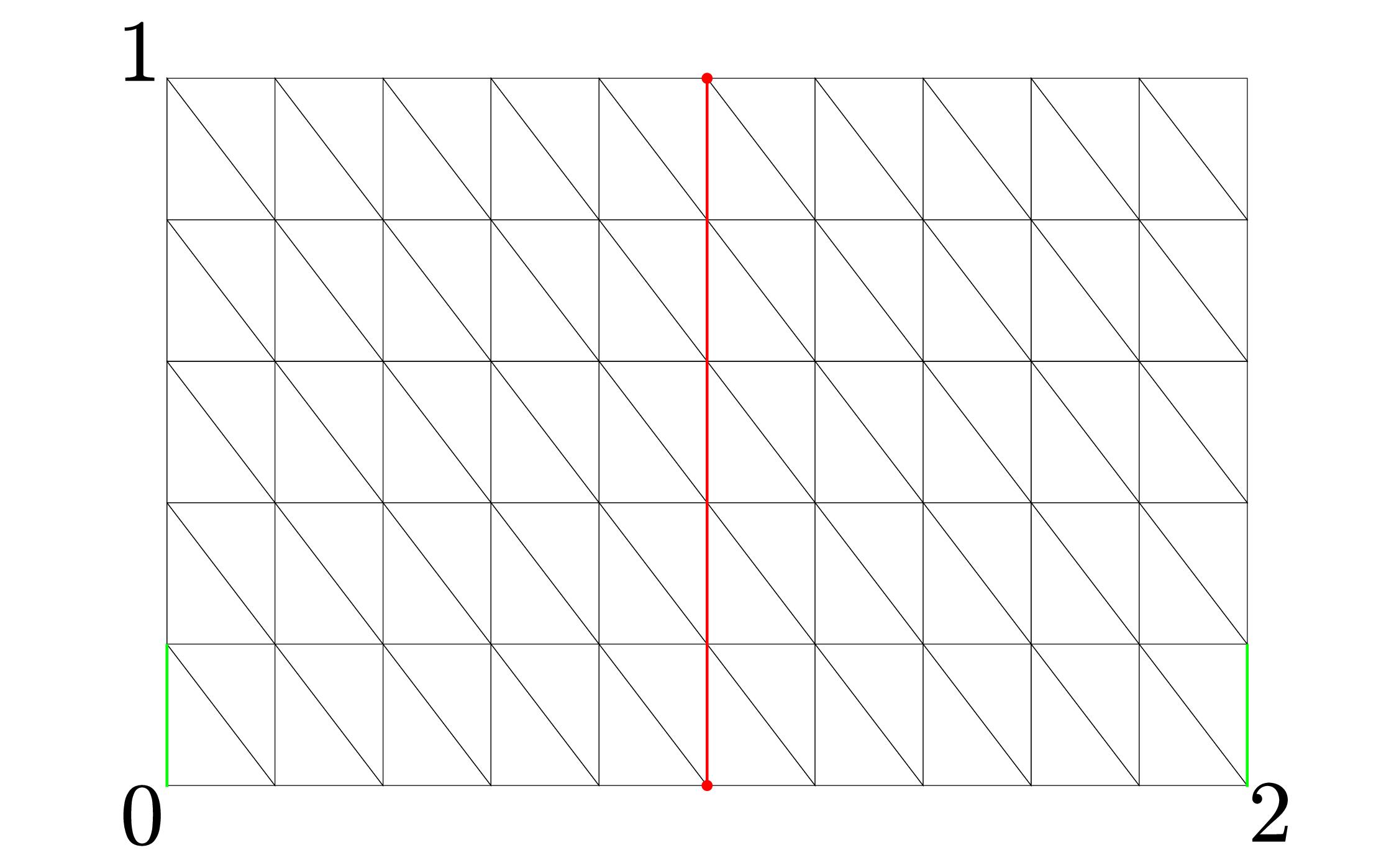}
  \label{fig:test2}
\end{minipage} 
\caption{\small [Test Case 1] (Left) Geometry and boundary conditions of the test case.  (Right) Example of an uniform triangular mesh for spatial discretization.}
\label{Non_Immersed} \vspace{-0.4cm}
\end{figure}

For the first test case, the domain of calculation $\Omega=(0,2) \times (0,1)$ is divided into two equally sized subdomains by a fracture of width ${\delta = 0.001}$ parallel to the $y$-axis (see Figure~\ref{Non_Immersed}).  For the boundary conditions, we impose $p =1$ at the bottom and $p=0$ at the top of the fracture.  On the external boundaries of the subdomains, a no flow boundary condition is imposed except on the lower fifth (length 0.2) of both lateral sides where a Dirichlet condition is imposed: $p = 1$ on the right and $p = 0$ on the left. The partition of the domain of calculation into triangular mesh is shown in the second figures of Figure~\ref{Non_Immersed}.
\begin{figure}[h!]
\centering
\begin{minipage}{0.35\textwidth}
\includegraphics[scale=0.25]{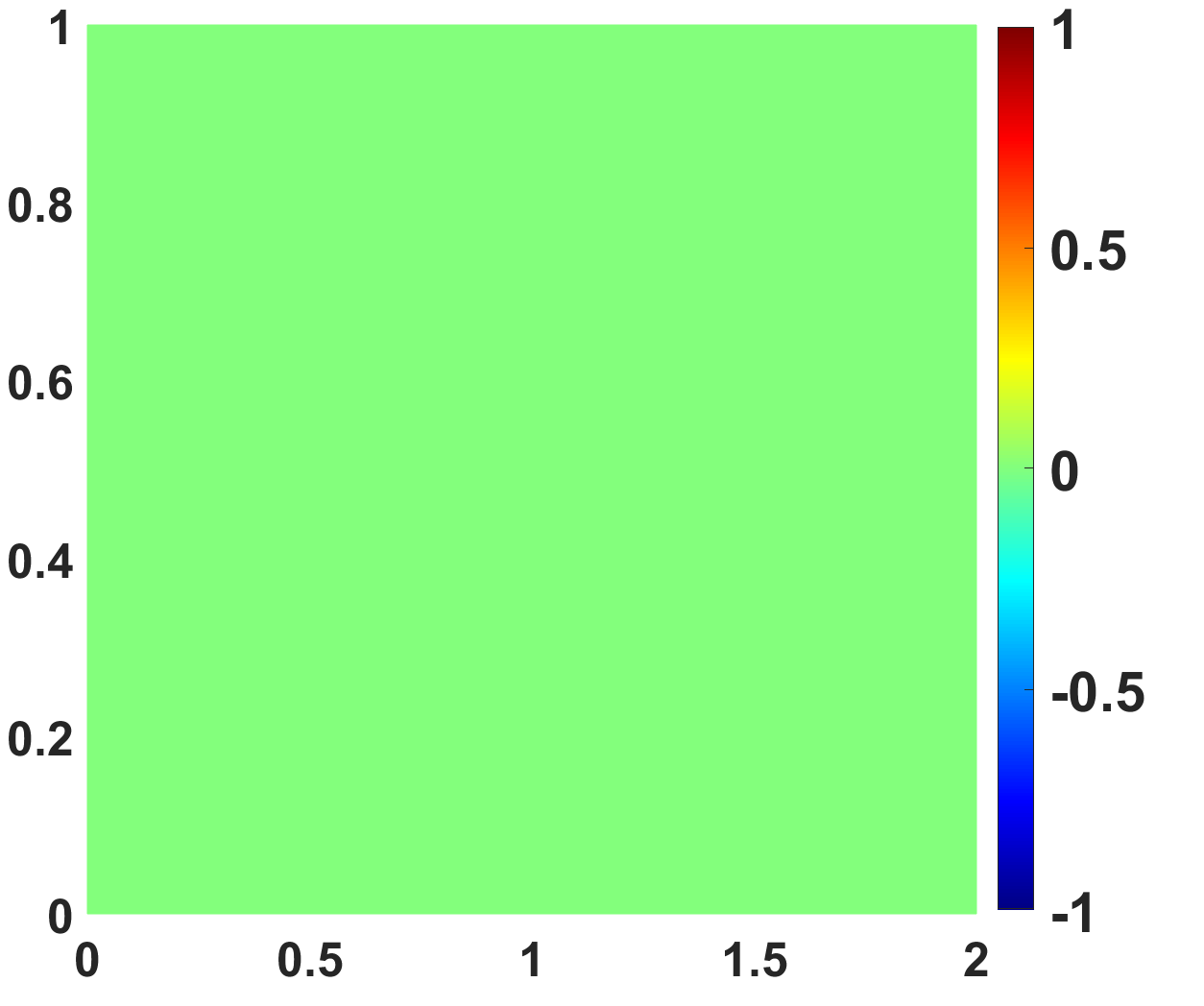}
\end{minipage}%
\begin{minipage}{0.35\textwidth}
\hspace{0.9cm}\includegraphics[scale=0.25]{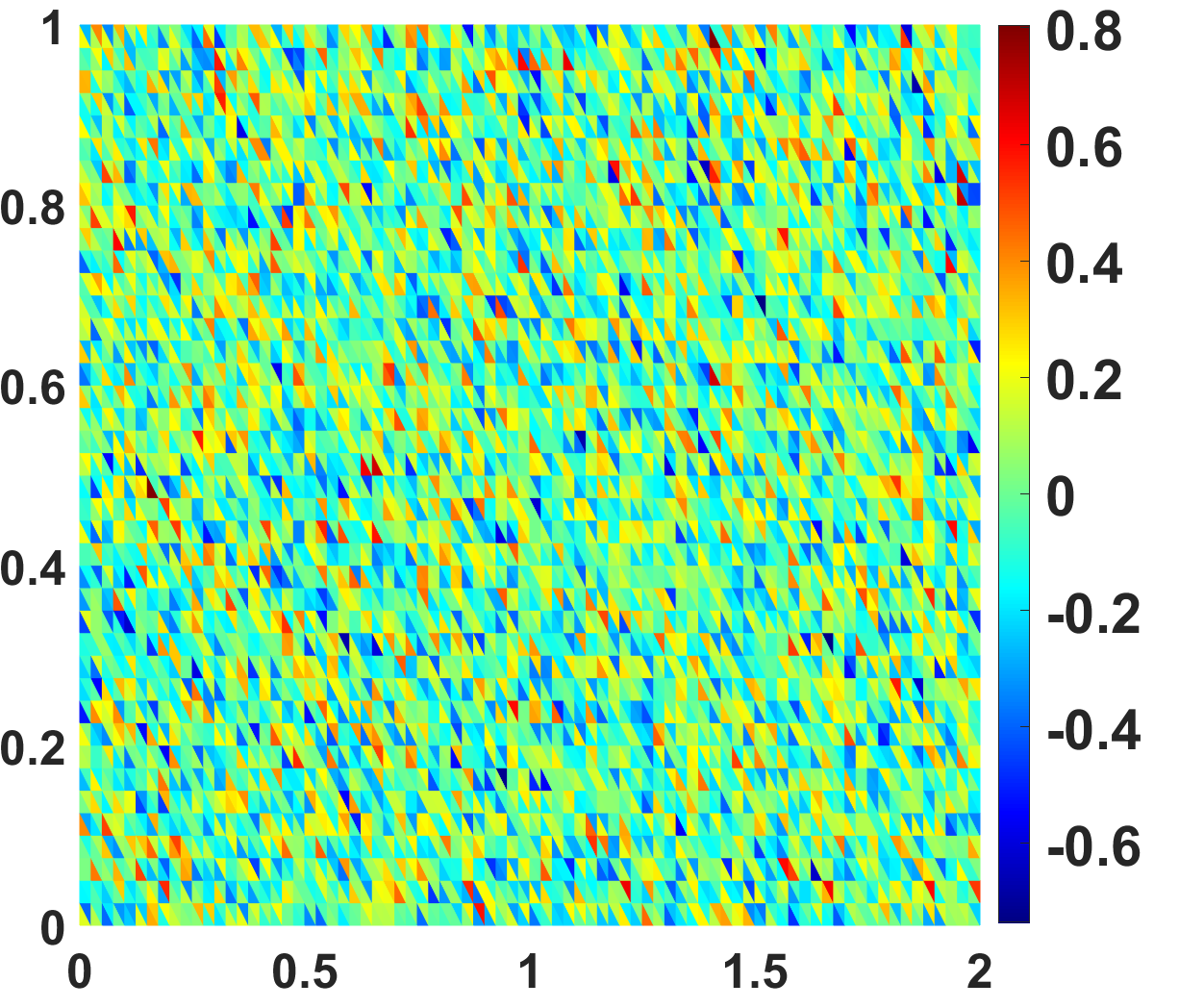}
\end{minipage}
\caption{\small [Test Case 1] The initial condition for the 2D equation. (Left) The true initial condition. (Right) The initial condition for the state dynamics in data assimilation.}
\label{InitialCond_Pure}\vspace{-0.4cm}
\end{figure}

Since data assimilation relies heavily on observational data, the accuracy of the assimilation algorithm may depend on the observational function provided. We aim to illustrate the robustness of the United Filter by considering three situations with varying levels of data availability. For the first case, we assume that we have full observation of the data, meaning that the function $g$ in \eqref{ObserProcess} is the identity operator. In the second case, only $75\%$ of the spatial points of the solution are observed, which are selected randomly. Finally, in the third case, we increase the problem's complexity and introduce additional nonlinearity by observing just $50\%$ of the state, with half of these observations transformed using the arctangent operation. We further assume that the actangent operator is applied for the velocity part of the solution as the magnitude of the velocity for this case is much larger than that of the pressure. The noise $\omega$ in \eqref{StateModel} is chosen to be $0.001\sqrt{\Delta{t}_{\text{Filter}}}\varepsilon$ where $\varepsilon \sim \mathcal{N}(0, \pmb{I}_l)$ with $\pmb{I}_l$ being the identity operator, and $\Delta{t}_{\text{Filter}}$ is the filtering timestep size. For all three cases, we choose the number of particles for the Direct Filter algorithm to be $M=30$, and the number of iterations in the United Filter setting to be $R=3$.
\begin{figure}[h!]
\centering
\includegraphics[scale=0.25]{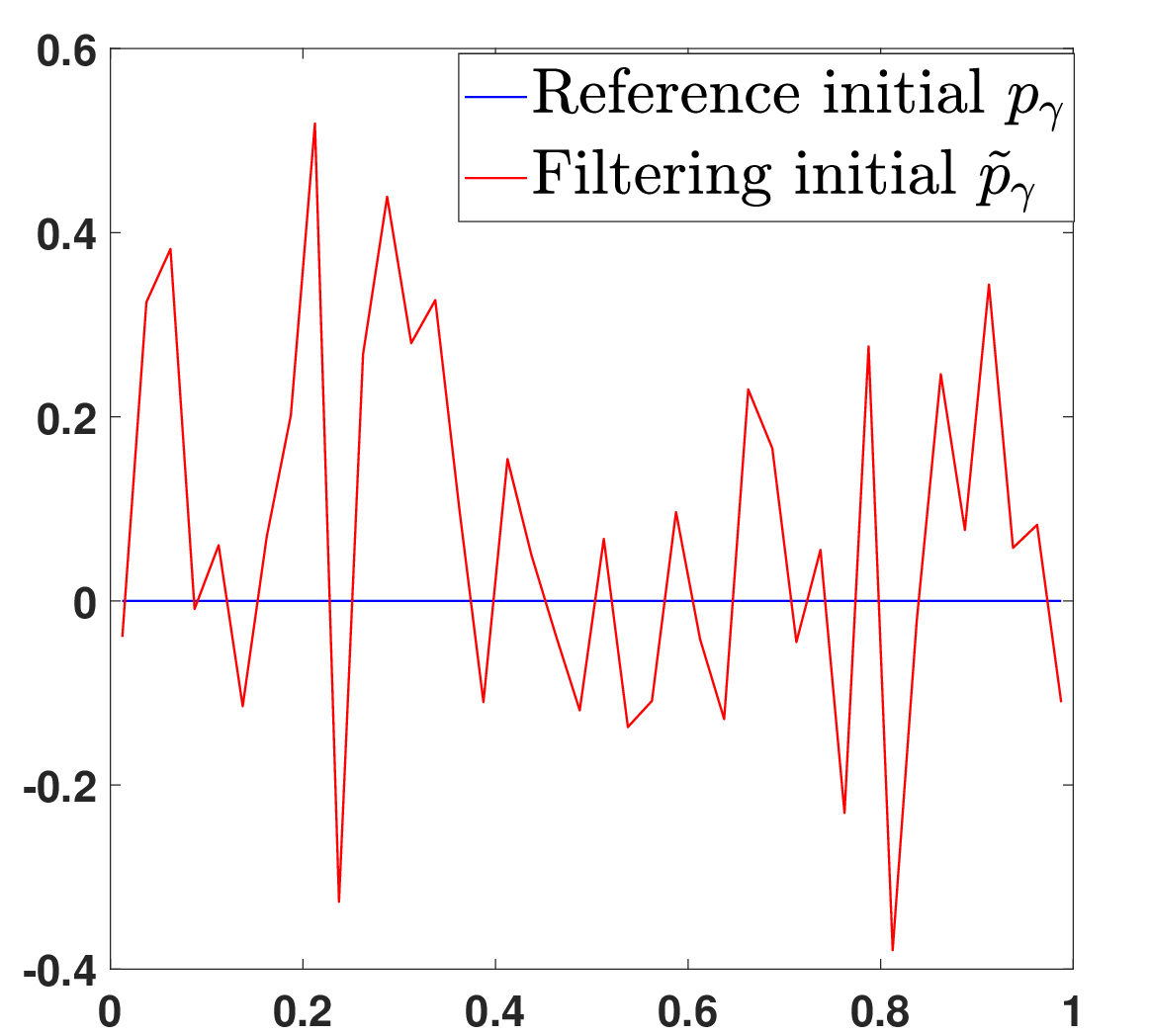}
\caption{[Test Case 1] The initial condition for the 1D equations. The straight blue line represents the true initial condition. The red curve depicts the initial condition used in the United Filter method for data assimilation.}
\label{InitialCond_1DPure}\vspace{-0.4cm}
\end{figure}

\begin{figure}[h!]
\centering
\begin{minipage}{0.25\textwidth}
\includegraphics[scale=0.22]{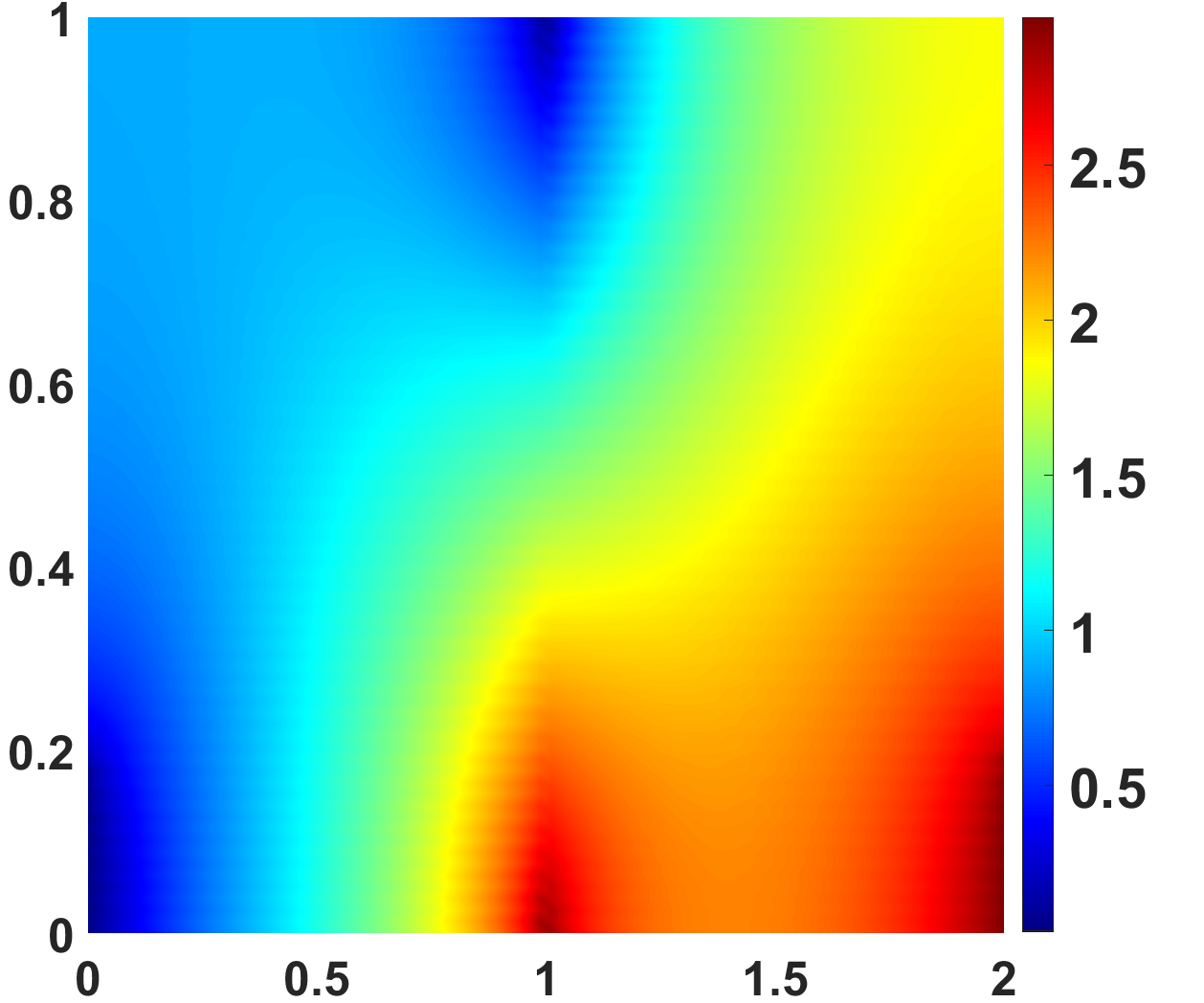}
\end{minipage}%
\begin{minipage}{0.226\textwidth}
\hspace{0.1cm}
\includegraphics[scale=0.22]{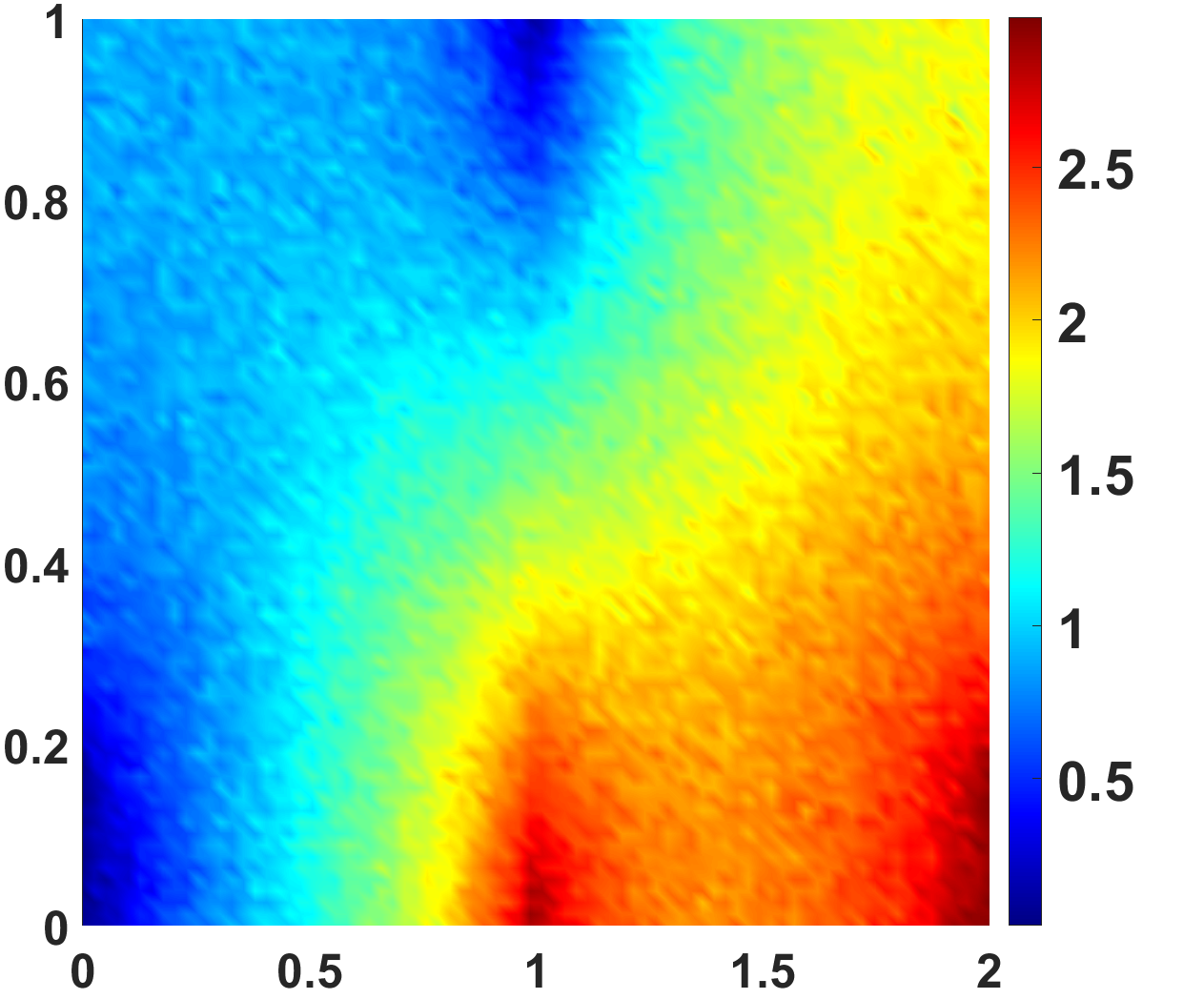}
\end{minipage} 
\hspace{0.1cm}
\begin{minipage}{0.235\textwidth}
\includegraphics[scale=0.22]{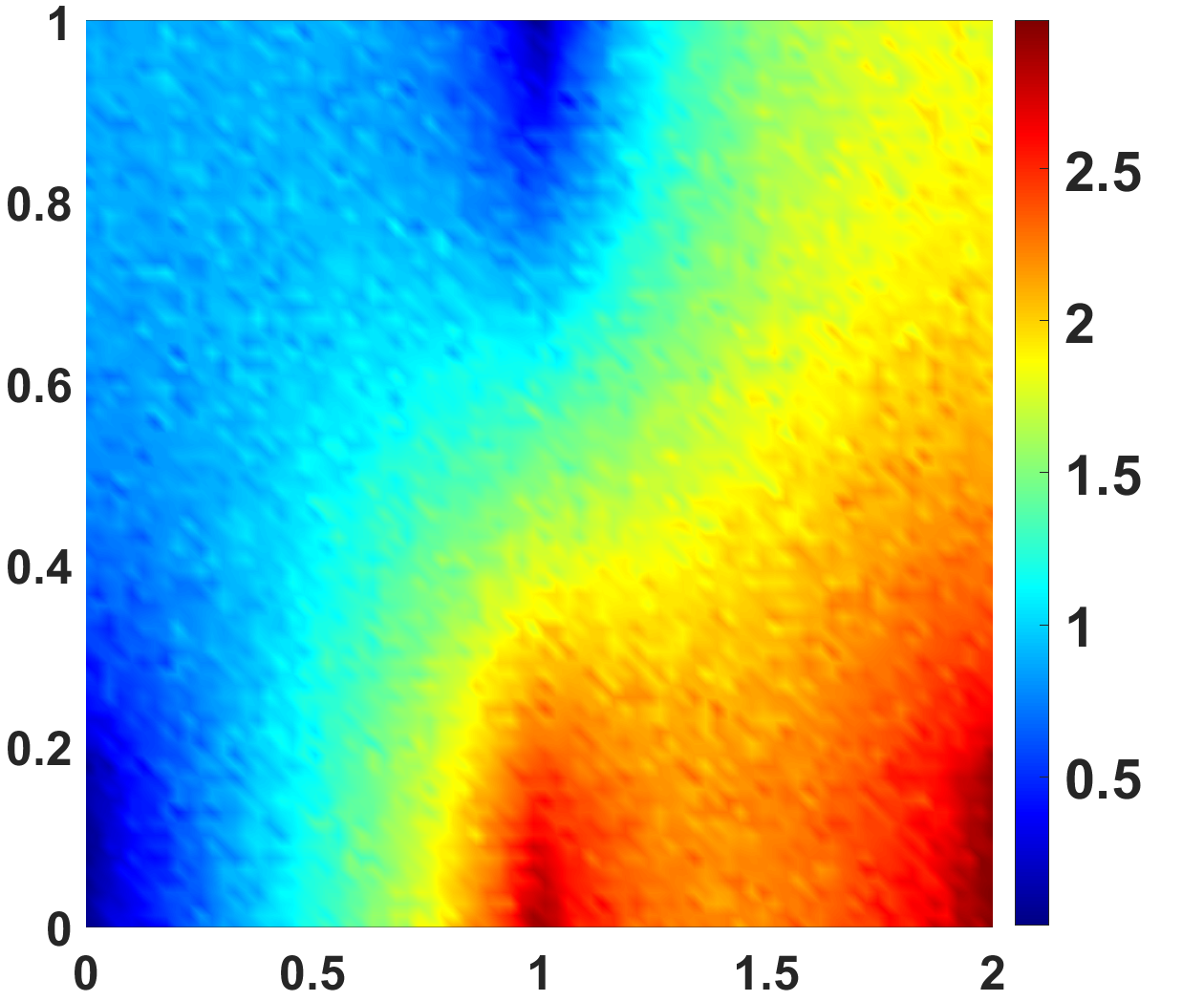}
\end{minipage}%
\hspace{0.1cm}
\begin{minipage}{0.23\textwidth}
\includegraphics[scale=0.22]{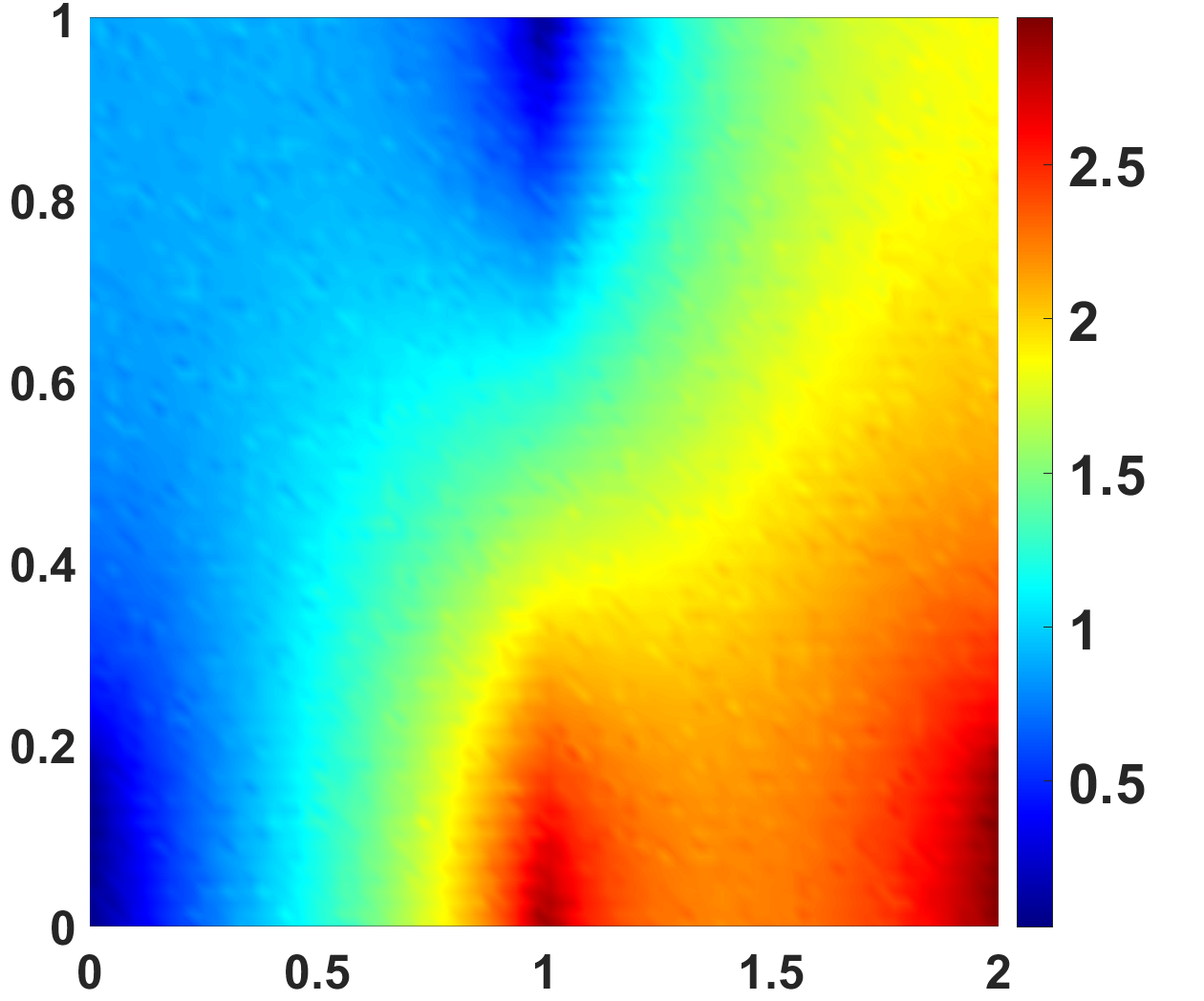}
\end{minipage}
\caption{[Test Case 1] Heat map illustrating the accuracy of the United Filter's pressure state estimation. (First) Reference pressure field at final time $T$. (Second) Estimated pressure field state with $100\%$ direct observation.
(Third) With $75\%$ direct observation. (Fourth) With $50\%$ mixed observation.}
\label{Test1PresField_Heat}
\vspace{-0.3cm}
\end{figure}

To generate the synthetic observational data, we first solve the problem \eqref{weak_discrete_reduced} with a spatial mesh size $h = 1/40$ and a very fine temporal mesh $\Delta{t} = T/800$ to obtain the reference solution. We then extract the parts of the reference solution at the following temporal location $\bar{t}_{i} = t_0 + i\frac{T-t_0}{Nt_{\text{Filter}}}, i=1, \cdots, Nt_{\text{Filter}}, $ where $Nt_{\text{Filter}} = 50$,
and apply the process \eqref{ObserProcess} to these parts to obtain the observation $Y_{1:Nt_{\text{Filter}}}$. The true parameters corresponding to the reference solution are $k_1 = 1, \; k_2=1$ and $\alpha_{\gamma} = k_{f}\delta = 2$. 

For all experiments, the United Filter is carried out on a coarser temporal mesh size $\Delta{t}_{\text{Filter}} = T/Nt_{\text{Filter}}$ with the same spatial mesh size. The reference initial conditions are set to zero across the domain, while those for the data assimilation algorithm are chosen randomly and are shown in Figure \ref{InitialCond_Pure} and Figure \ref{InitialCond_1DPure}. The initial guesses for the parameter estimation are $[8, 8, 8]^T$.
\begin{figure}[h!]
\centering
\begin{minipage}{0.248\textwidth}
\includegraphics[scale=0.22]{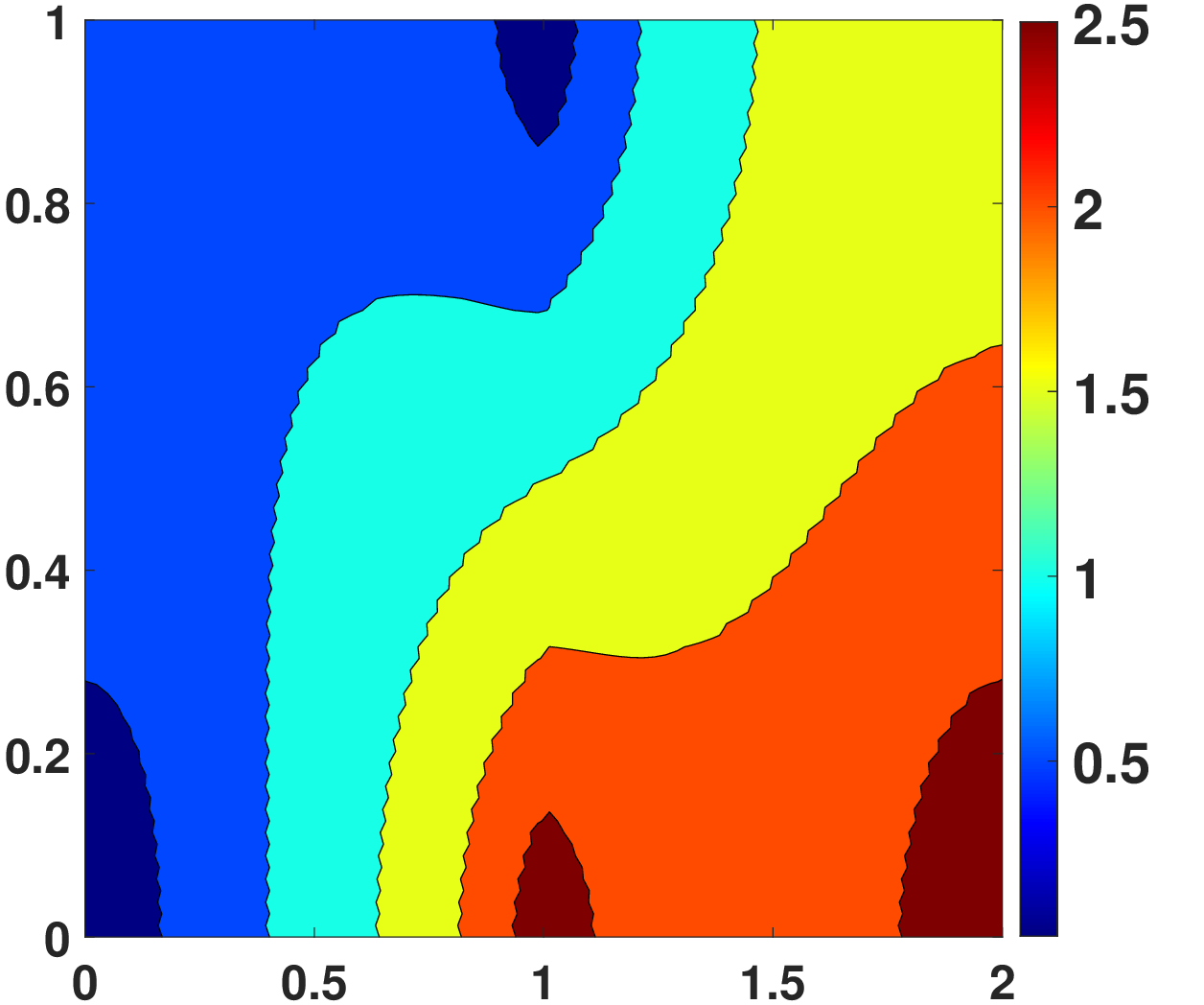}
\end{minipage}%
\begin{minipage}{0.228\textwidth}
\hspace{0.1cm}\includegraphics[scale=0.22]{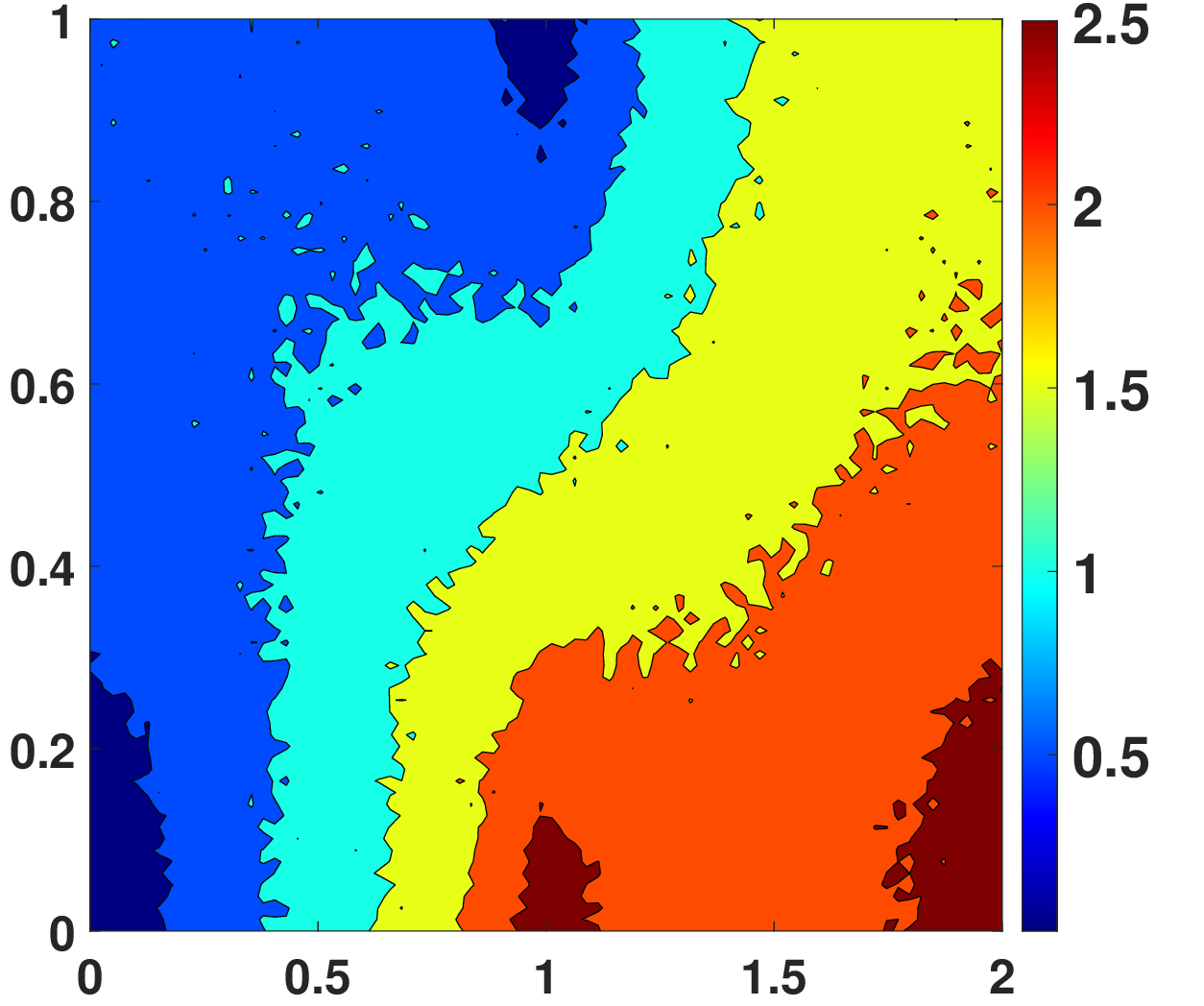}
\end{minipage}
\hspace{0.1cm}
\begin{minipage}{0.246\textwidth}
\includegraphics[scale=0.22]{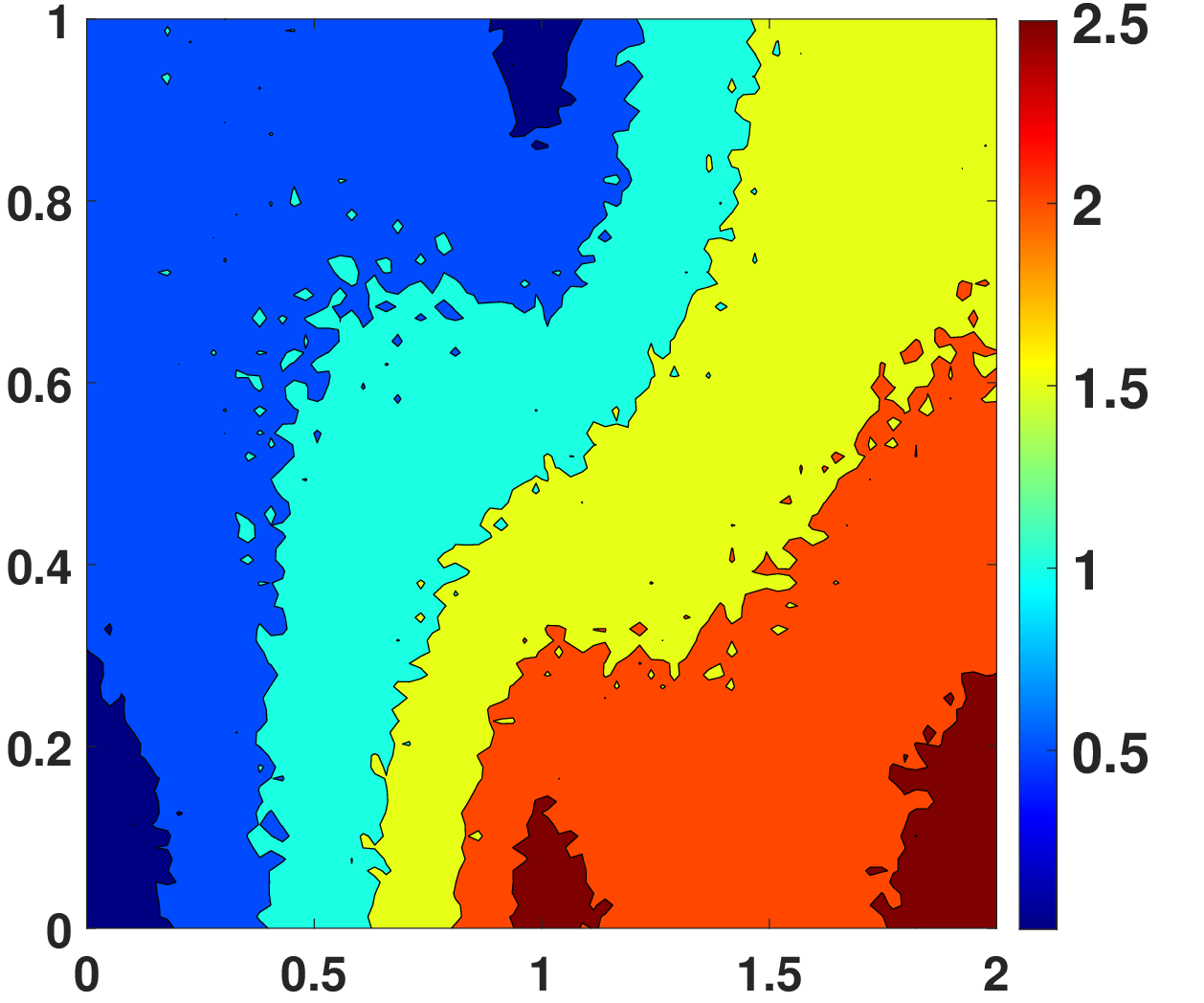}
\end{minipage}%
\begin{minipage}{0.23\textwidth}
\hspace{0.1cm}\includegraphics[scale=0.22]{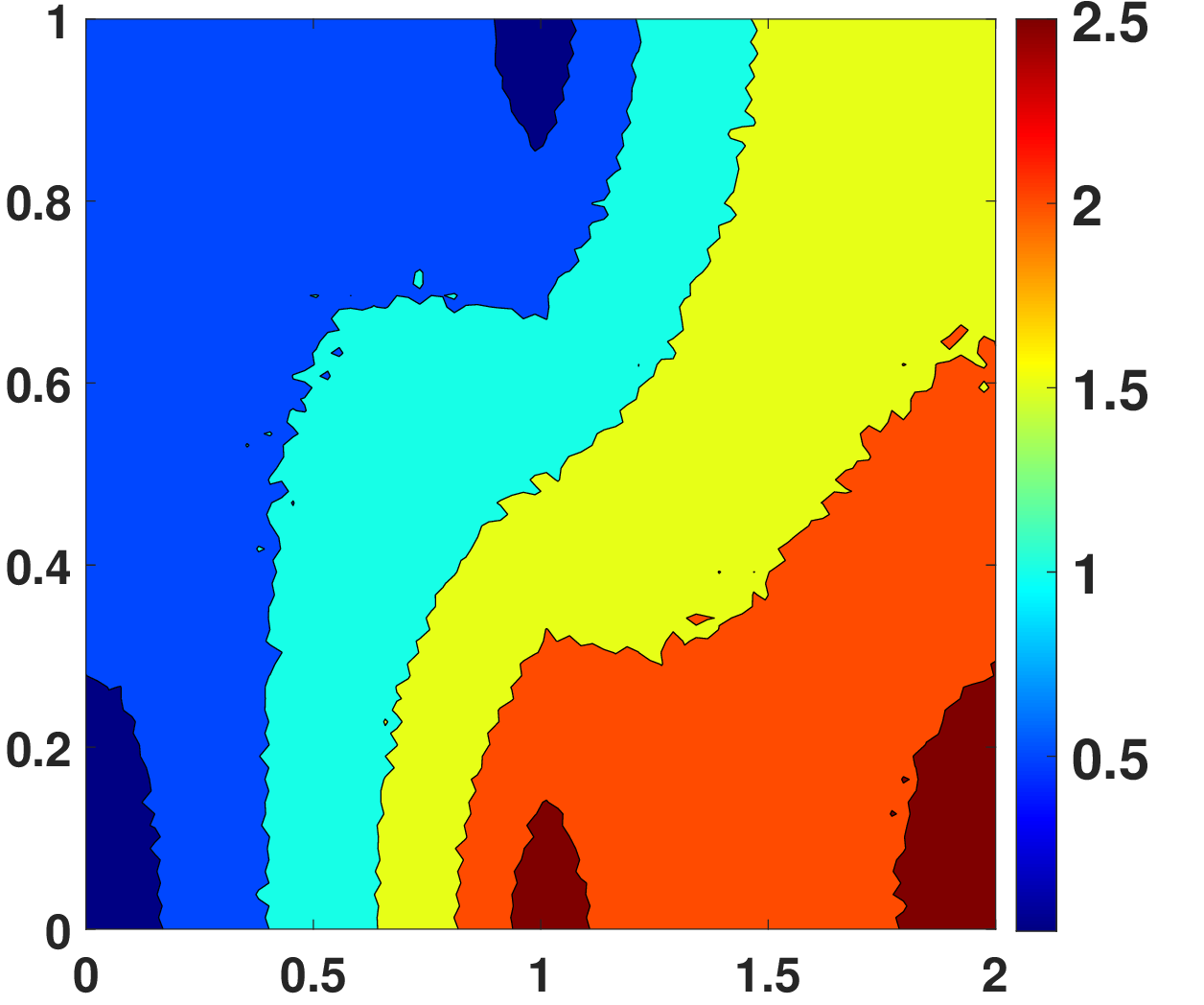}
\end{minipage}
\caption{[Test Case 1] Contour map illustrating the accuracy of the United Filter's pressure state estimation. (First) Reference pressure field at final time $T$. (Second) Estimated pressure field state with $100\%$ direct observation. (Third) With $75\%$ direct observation. (Fourth) With $50\%$ mixed observation.}
\label{Test1PresField_Contour}
\vspace{-0.2cm}
\end{figure}

\begin{figure}[h!]
\centering
\begin{minipage}{0.25\textwidth}
\includegraphics[scale=0.21]{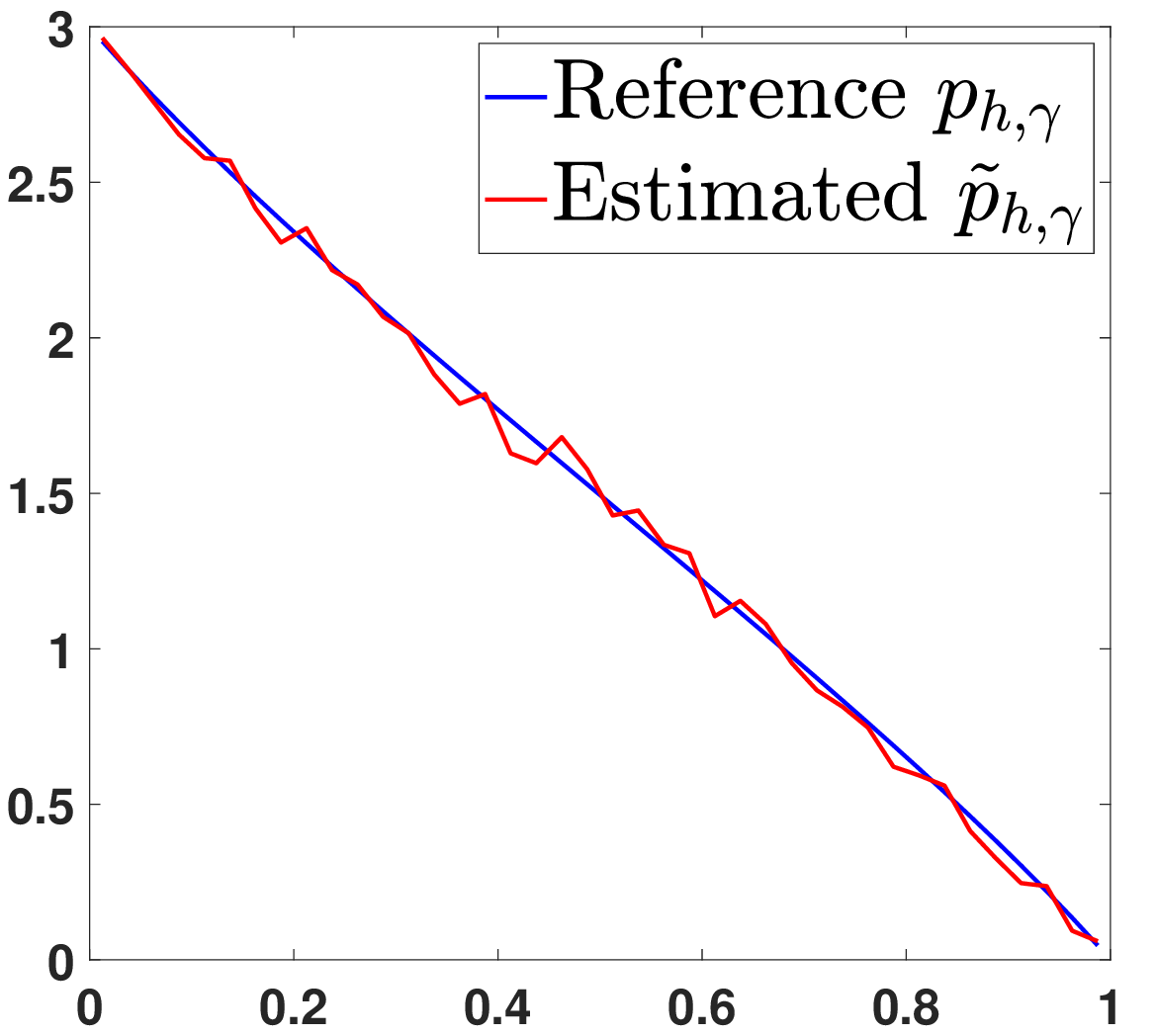}
\end{minipage}%
\begin{minipage}{0.25\textwidth}
\hspace{0.2cm}\includegraphics[scale=0.21]{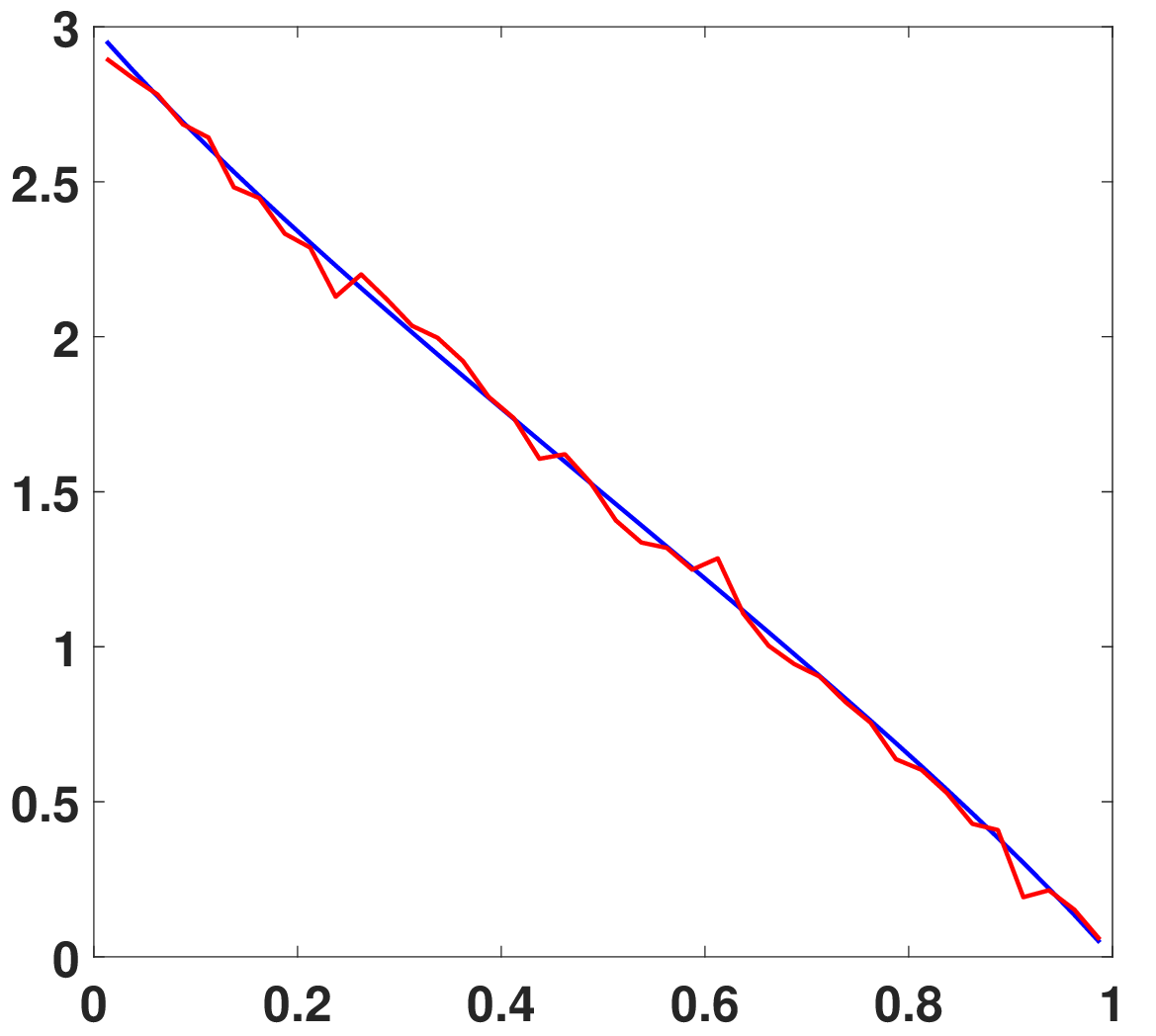}
\end{minipage}%
\begin{minipage}{0.25\textwidth}
\hspace{0.3cm}\includegraphics[scale=0.21]{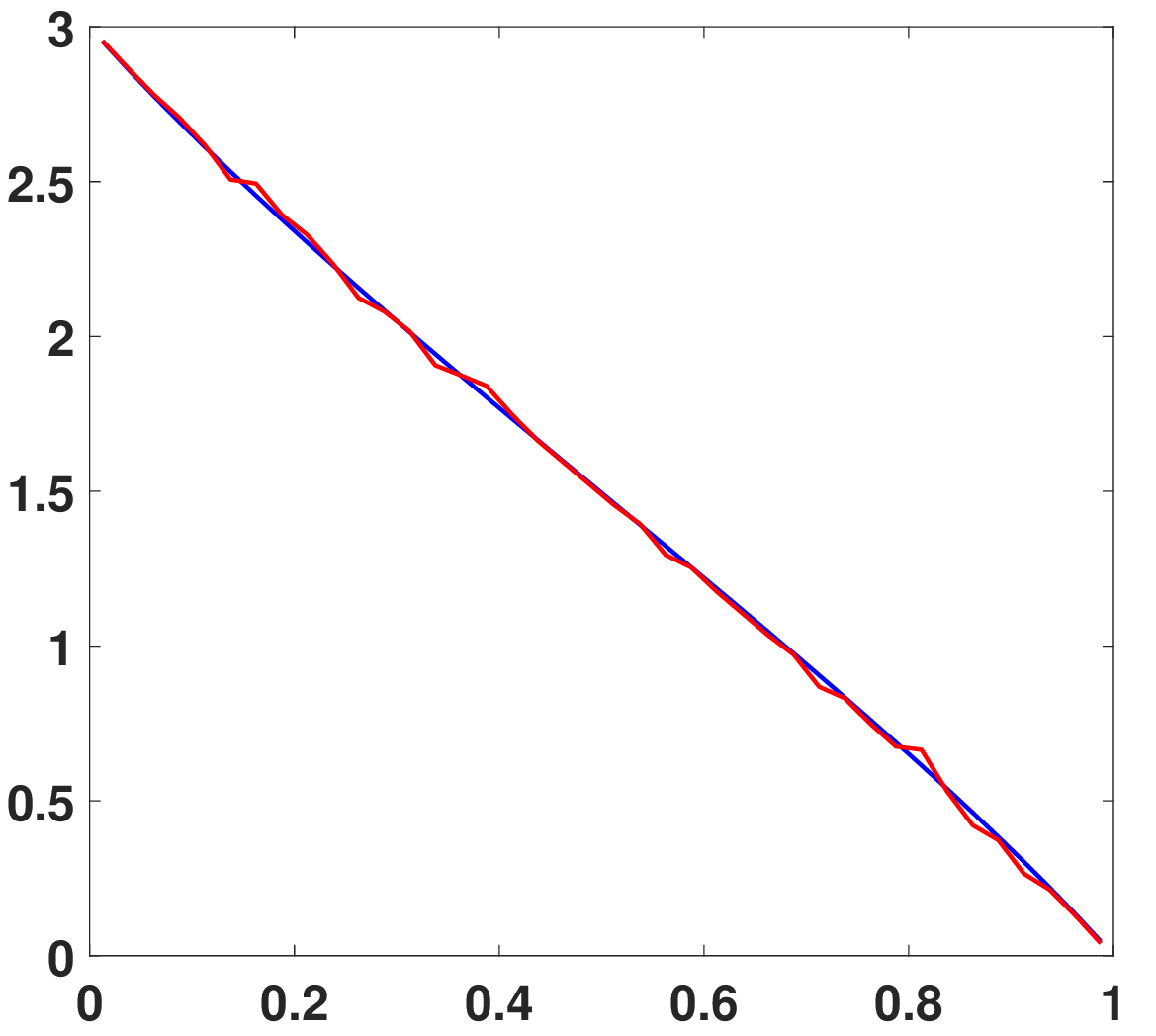}
\end{minipage}
\caption{[Test Case 1] Demonstration for 1D pressure estimation at fracture: (Left) With $100\%$ direct observation.
(Middle) With $75\%$ direct observation. (Right) With $50\%$ mixed observation.}
\label{Test1_PresFrt}
\end{figure}

\begin{figure}[h!]
\centering
\begin{minipage}{0.23\textwidth}
\includegraphics[scale=0.2]{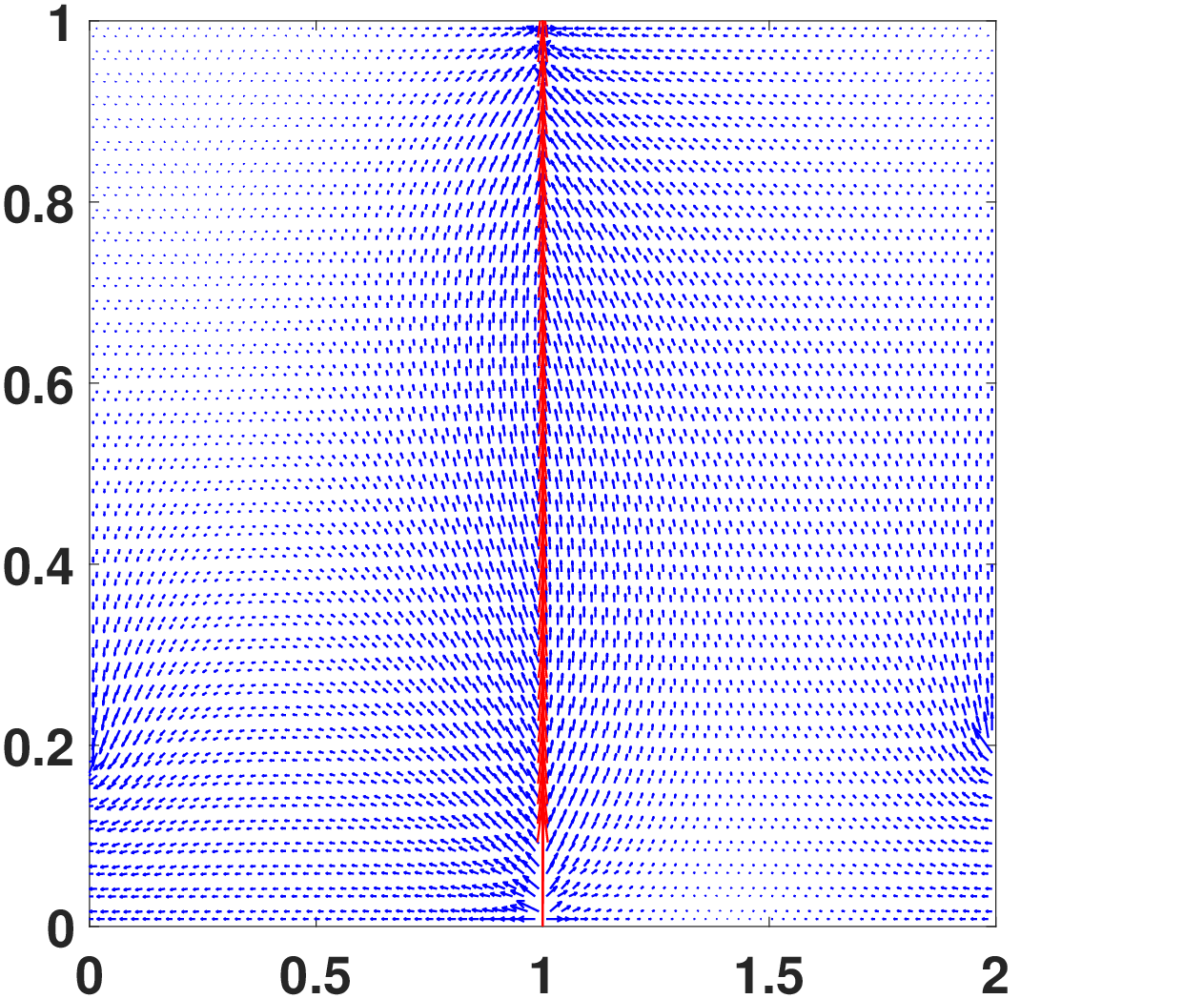}
\end{minipage}%
\begin{minipage}{0.23\textwidth}
\hspace{0.1cm}\includegraphics[scale=0.2]{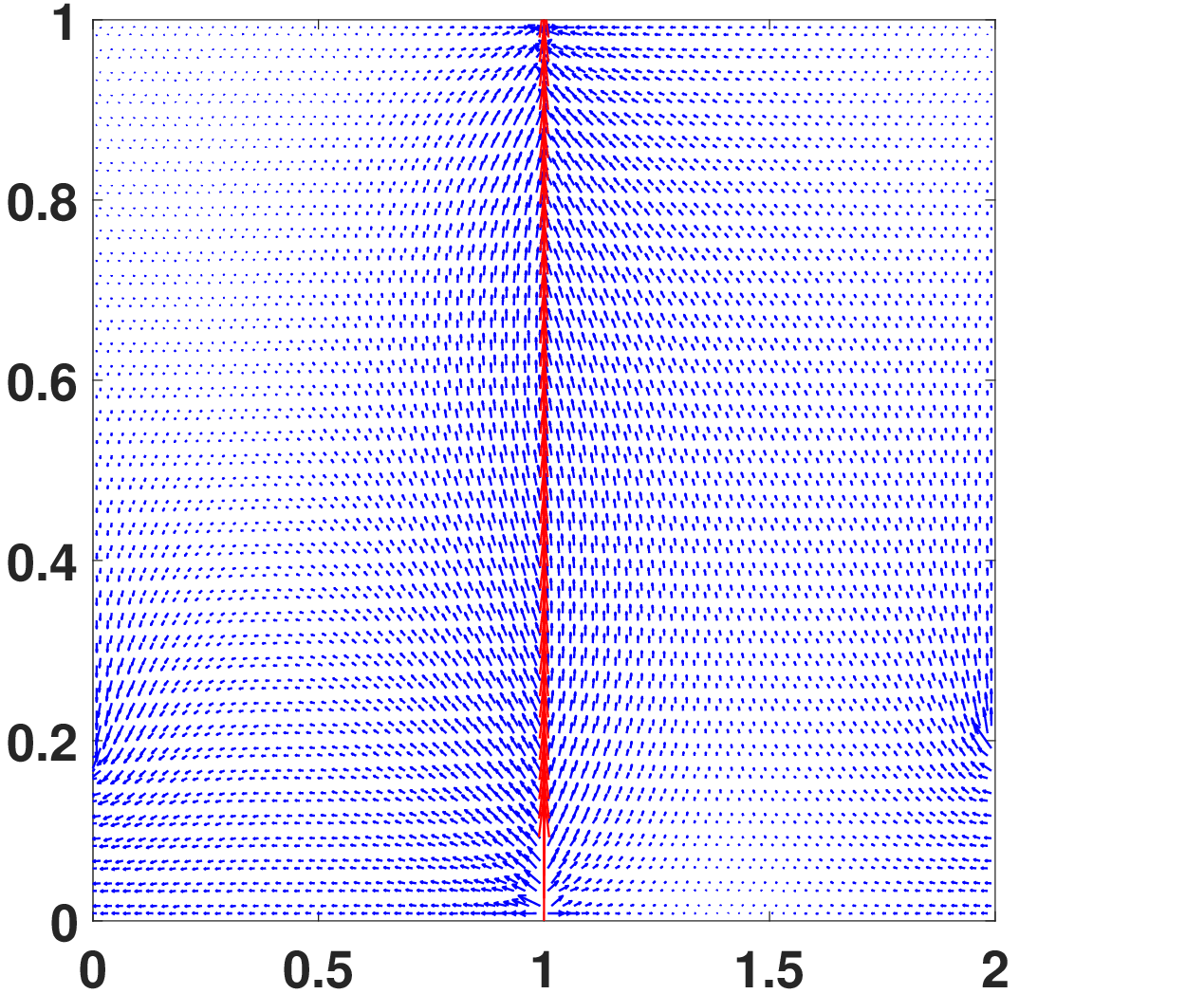}
\end{minipage}
\hspace{0.1cm}
\begin{minipage}{0.23\textwidth}
\includegraphics[scale=0.2]{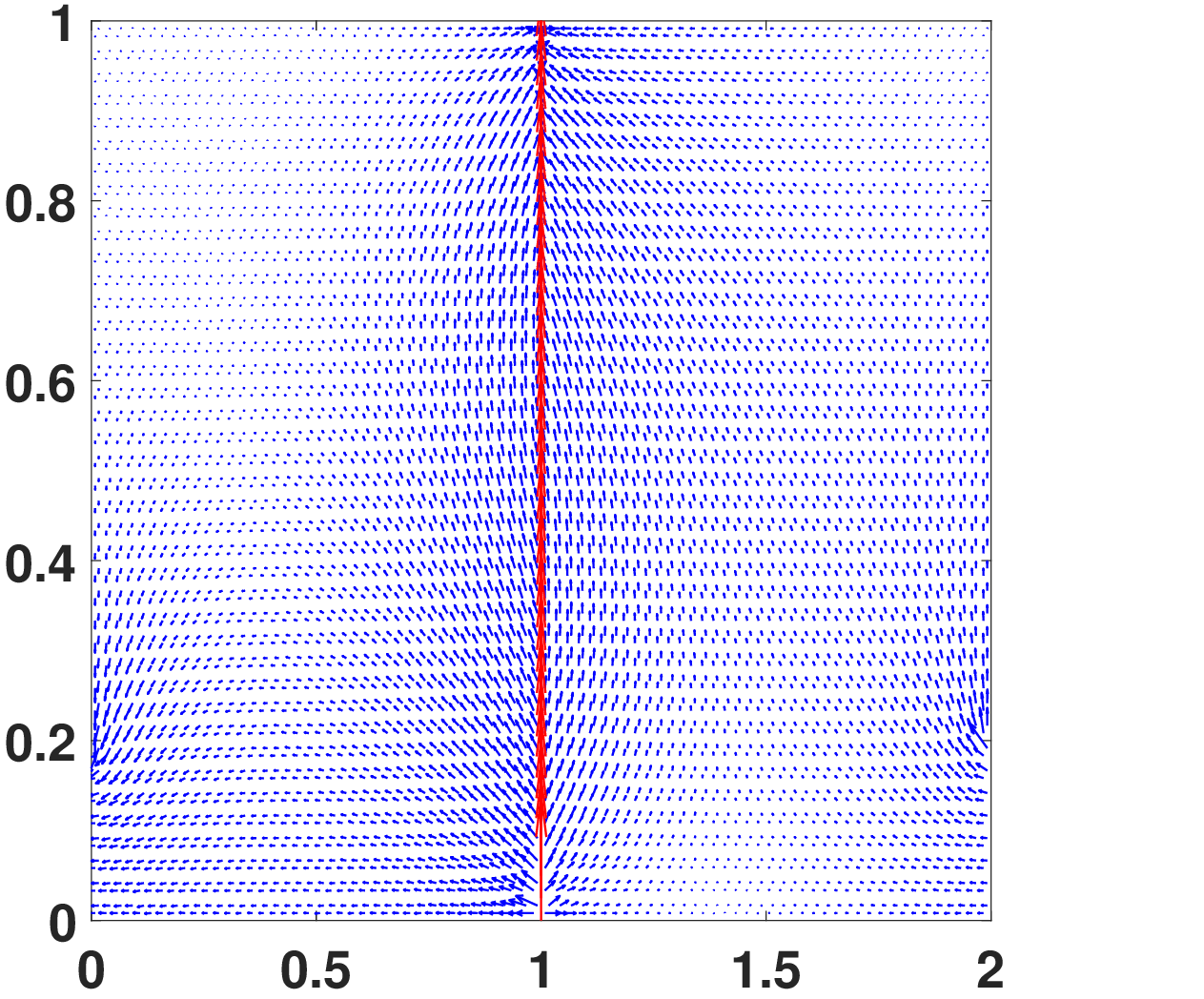}
\end{minipage}%
\begin{minipage}{0.23\textwidth}
\hspace{0.1cm}\includegraphics[scale=0.2]{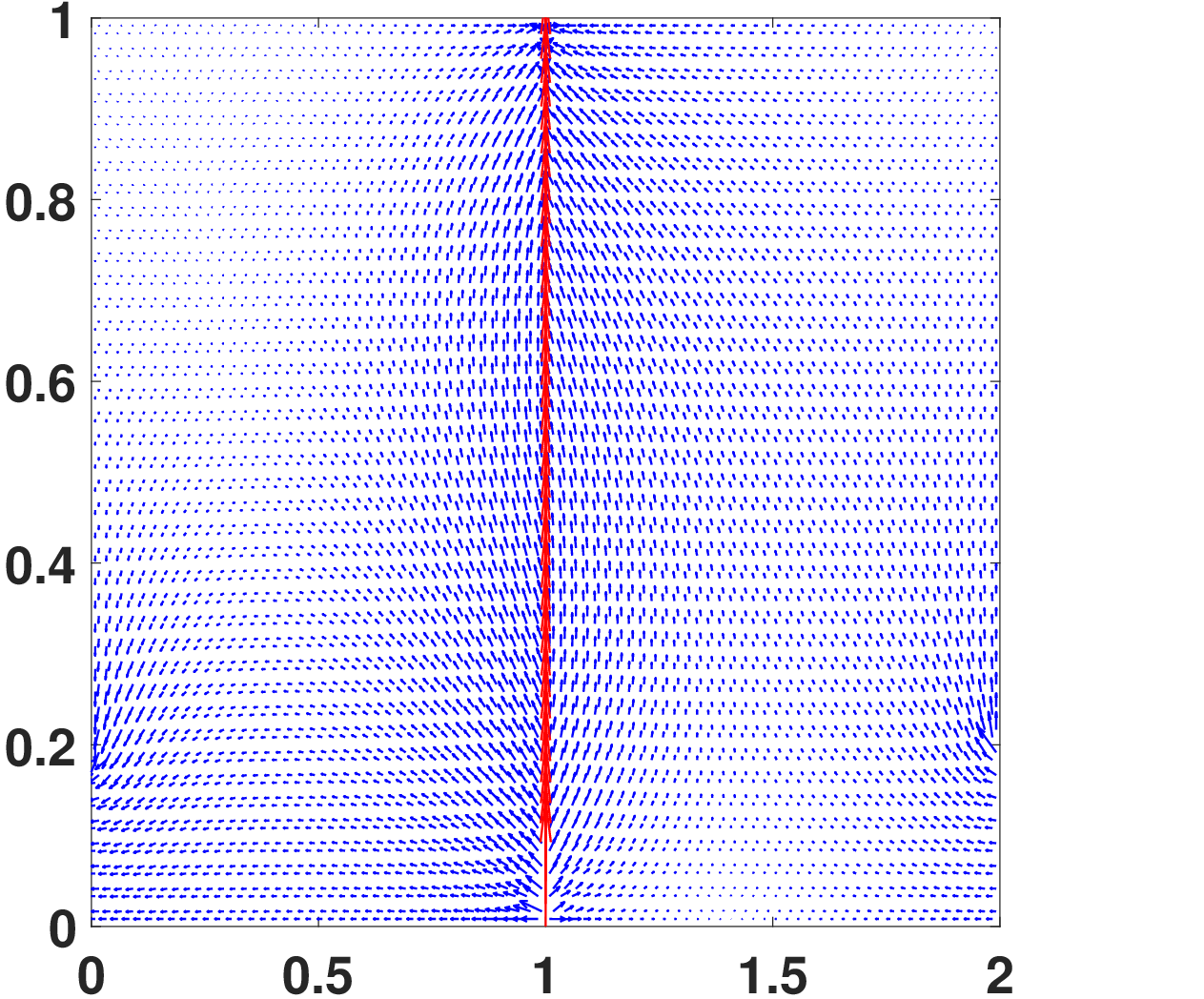}
\end{minipage}
\caption{[Test Case 1] Demonstration of the United Filter's accuracy in velocity state estimation. (First) Reference velocity field at final time $T$. (Second) Estimated velocity field state with $100\%$ direct observation.
(Third) With $75\%$ direct observation. (Fourth) With $50\%$ mixed observation.}
\label{Test1VeloField}
\vspace{-0.3cm}
\end{figure}
\begin{figure}[h!]
\centering
\begin{minipage}{0.4\textwidth}
\hspace{1cm}\includegraphics[scale=0.25]{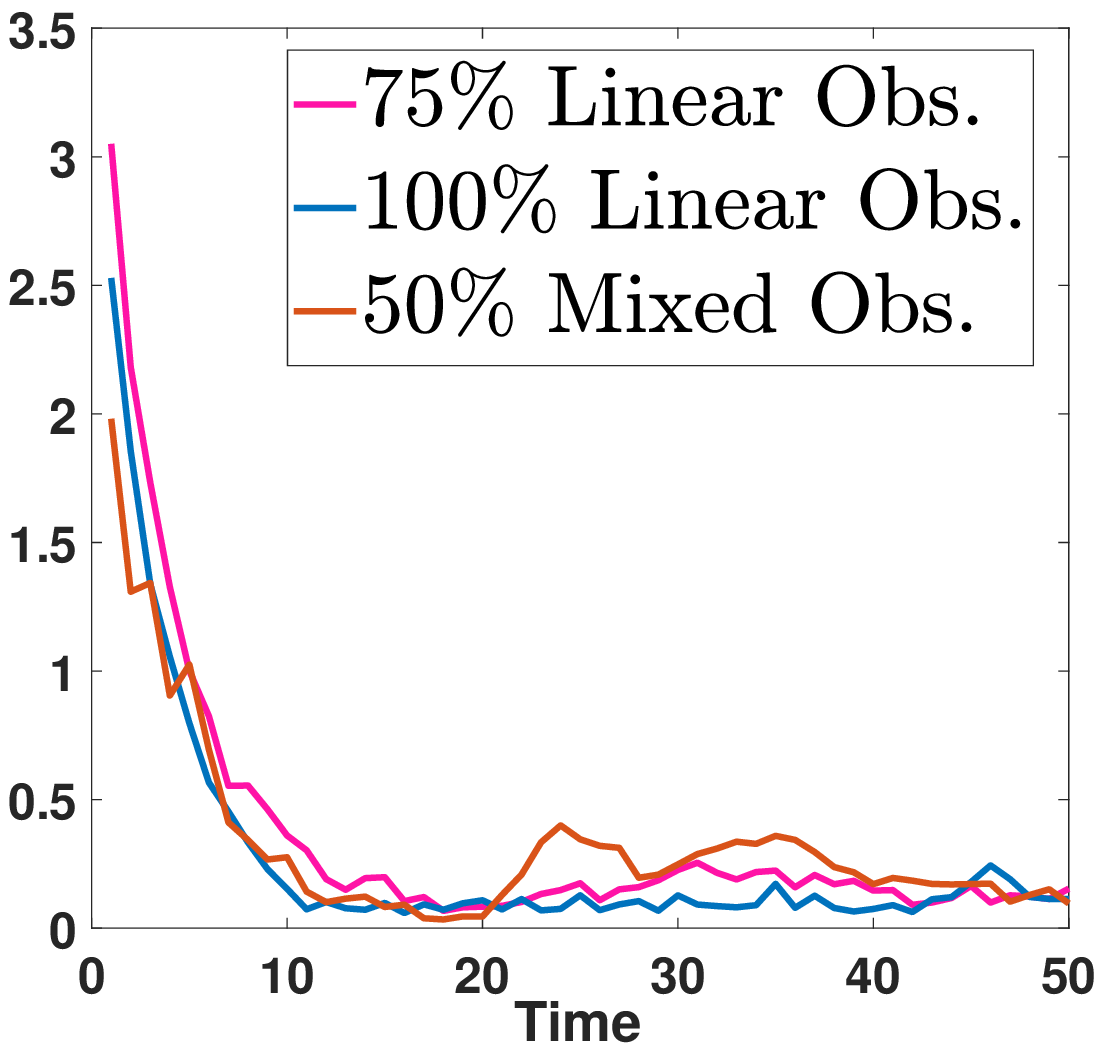}
\end{minipage}%
\caption{[Test Case 1] Root Mean Square Errors (RMSEs) for state estimation}
\label{Test1RMSE_Pure}
\vspace{-0.3cm}
\end{figure}
\begin{figure}[h!]
\centering
\begin{minipage}{0.33\textwidth}
\hspace{-0.1cm}\includegraphics[scale=0.21]{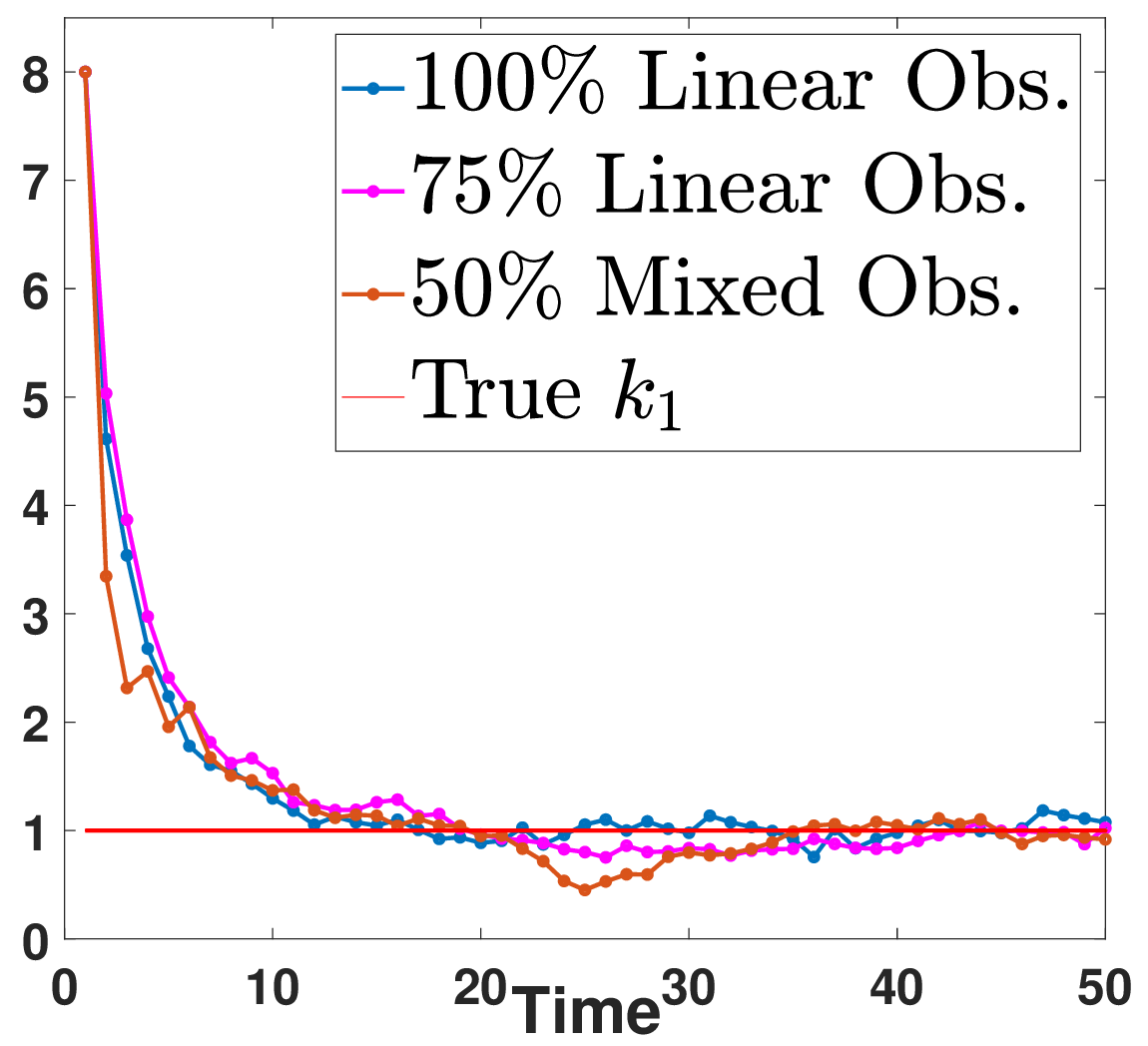}
\end{minipage}%
\begin{minipage}{0.33\textwidth}
\hspace{0.2cm}\includegraphics[scale=0.21]{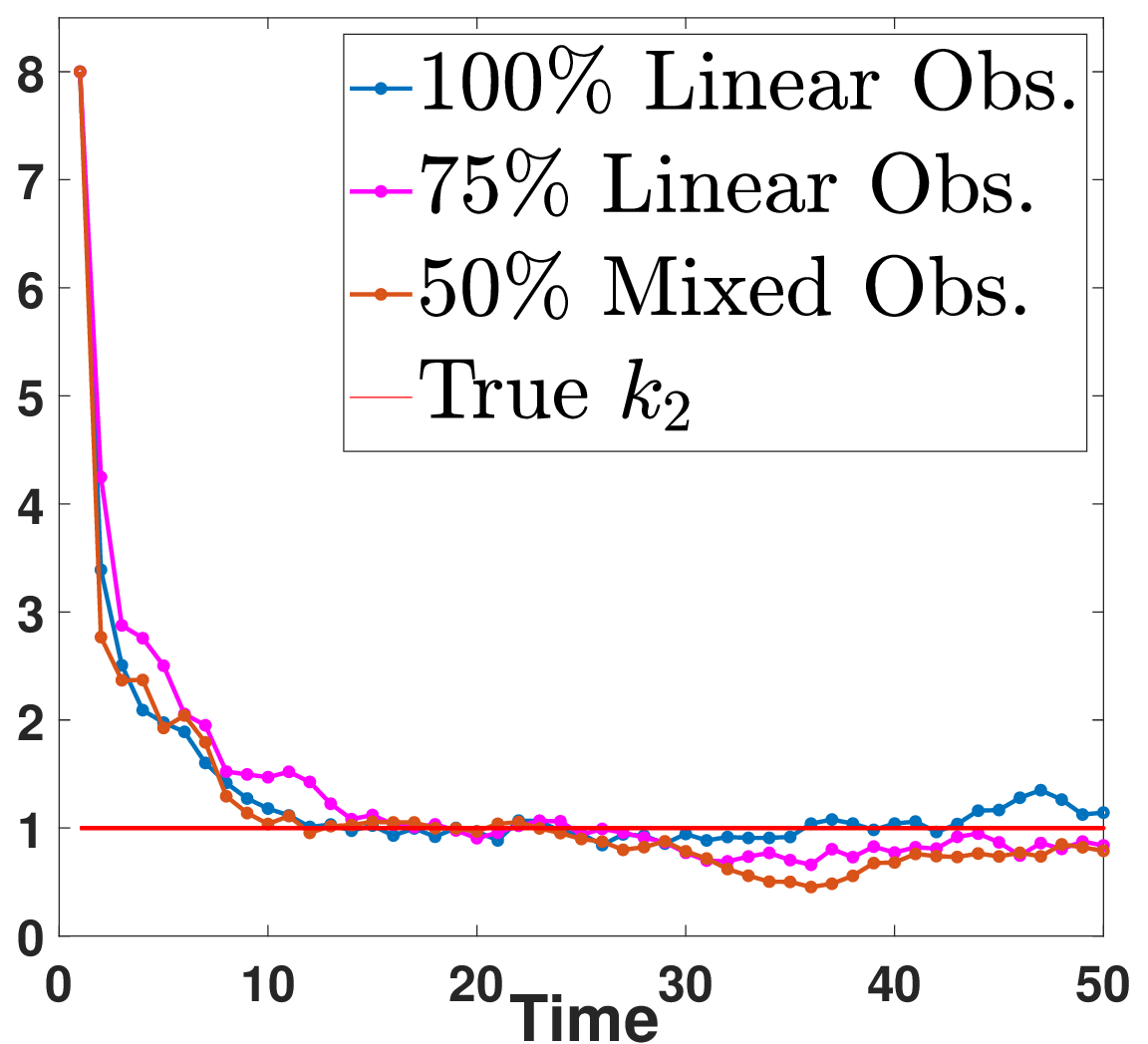}
\end{minipage}%
\begin{minipage}{0.33\textwidth}
\hspace{0.5cm}\includegraphics[scale=0.21]{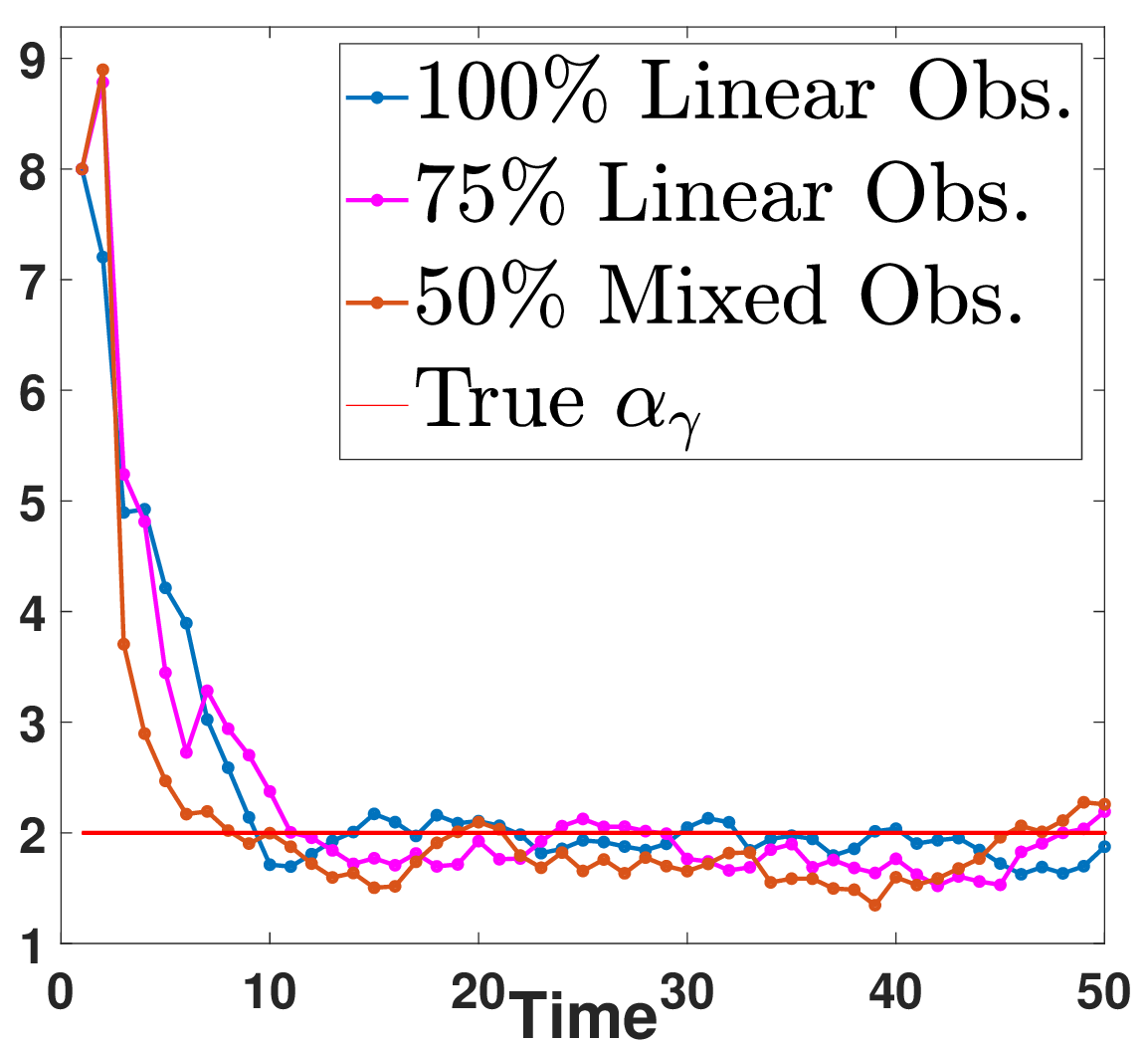}
\end{minipage}
\caption{[Test Case 1] The United Filter's performance of parameter estimation: (Left) Estimation for $k_1$. (Middle) Estimation for $k_2$. (Right)  Estimation for $\alpha_{\gamma}$.}
\label{ParaEst_Pure}
\vspace{-0.5cm}
\end{figure}

We first focus on the pressure part of the state estimation and present in Figure~\ref{Test1PresField_Heat} the reference 2D pressure field as well as the estimated pressure state obtained by the United Filter. To better observe the layers of the heat maps, the corresponding contour maps are shown in Figure~\ref{Test1PresField_Contour}. Even in the case that only 50$\%$ linear/nonlinear mixed observations are collected, the United Filter method still provides accurate estimates for the pressure field. The estimates for the pressure on the fracture are displayed in Figure~\ref{Test1_PresFrt}, which demonstrates that all cases yield good estimates of the exact state. Note that the graphs for both the reference pressure fields and the estimated pressure states are presented at the terminal time $T$.  

Next, we plot the exact and estimated velocity fields in Figure~\ref{Test1VeloField}. Nearly identical estimation performance are achieved across all observation scenarios. 
To better illustrate the accuracy of the United Filter in state estimation, we present the Root Mean Square Errors (RMSEs) of the estimated state (combining errors over pressure and velocity) in Figure~\ref{Test1RMSE_Pure} to provide a clear quantitative assessment of the United Filter's performance. These figures indicate that the RMSEs decrease overtime and eventually converge to nearly the same value.

Regarding the performance of the United Filter's parameter estimation, the estimates for all three parameters (with respect to data reception time) are shown in Figure~\ref{ParaEst_Pure}. We can see from this figure that the United Filter can provide accurate results for parameter estimation. 
Regardless of the amount of observational data, the particles approximating $k_1$ and $k_2$ (the first and second figures) quickly converge to the true values and remain stable. For the parameter $\alpha_{\gamma}$, the particles still converge relatively fast, but they exhibit mild oscillations as they approach the true parameter (see the last plot in Figure~\ref{ParaEst_Pure}). We will see in the next numerical test that the parameters on the fracture are also harder to approximate, but the United Filter algorithm still provides precise and consistent results.

\vspace{0.1cm}
To highlight the superior performance of the United Filter in joint state-parameter estimation, we compare it with the state-of-the-art joint state-parameter estimation method, the augmented Ensemble Kalman Filter (AugEnKF), for solving the same pure diffusion equation with a single fracture problem.

\begin{figure}[h!]
\centering
\begin{minipage}{0.35\textwidth}
\includegraphics[scale=0.22]{./figures_All/figures/figures_update/HeatMap_TruePres_Update_v2.eps}
\end{minipage}%
\begin{minipage}{0.35\textwidth}
\hspace{0.5cm}\includegraphics[scale=0.22]{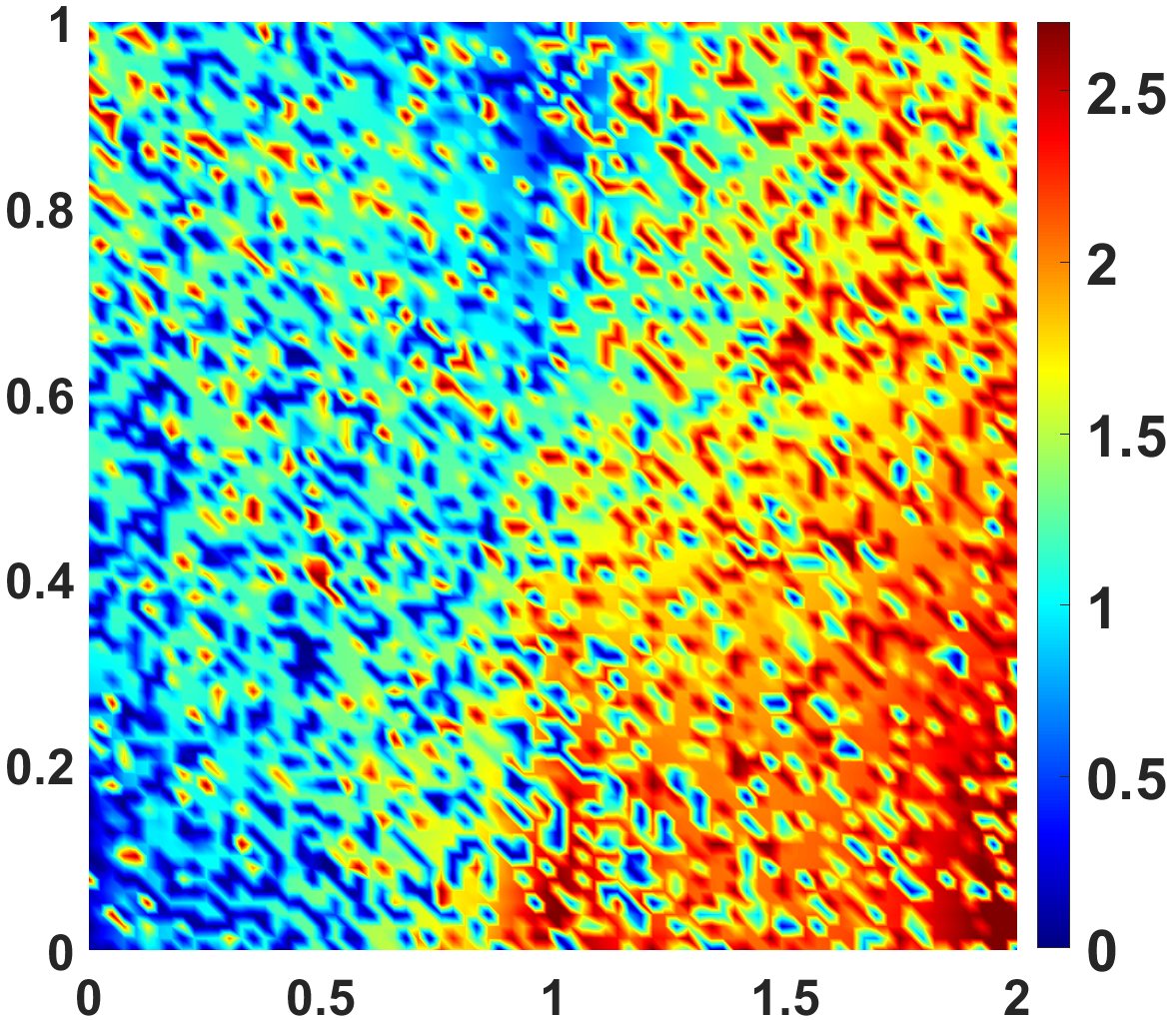}
\end{minipage}
\caption{[Test Case 1] Heat maps for the pressure field estimated by the AugEnKF: (Left) Reference pressure field. (Right) Estimation by the AugEnKF under 50\% of mixed observation. This is the comparison result by the AugEnKF corresponding to the United Filter's results in Figure \ref{Test1PresField_Heat}. }
\label{Test1HeatPresField_EnKF}
\vspace{-0.4cm}
\end{figure}
For the AugEnKF, the ensemble size is chosen as $100$, and we assume that 50$\%$ of the data is observed, half of which is received directly, while the other half is obtained through the arctangent function. For the Ensemble Kalman Filter method, we can either use the arctangent operator directly or its linearized version to compute the Kalman gain at each filtering step (see, e.g., \cite{roth2017}). In our work, we choose the latter approach as we aim to have a linear setting for the Kalman Filter. We reuse the reference solution computed in the setting of the United Filter method to construct the synthetic observational data. 
\begin{figure}[h!]
\vspace{-0.2cm}
\centering
\begin{minipage}{0.35\textwidth}
\includegraphics[scale=0.22]{./figures_All/figures/figures_update/ContourMap_TruePres_Update_v2.eps}
\end{minipage}%
\begin{minipage}{0.35\textwidth}
\hspace{0.5cm}\includegraphics[scale=0.22]{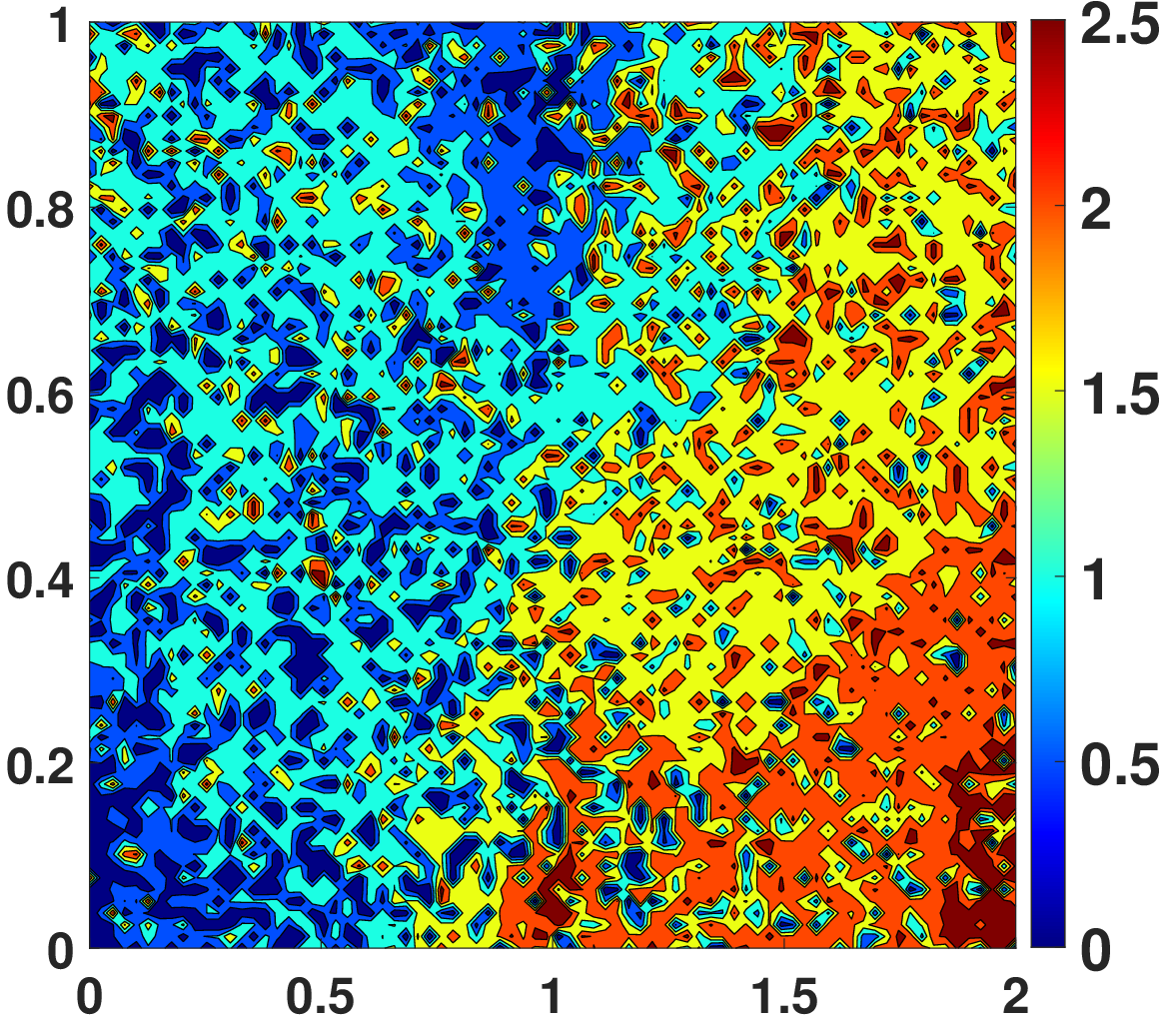}
\end{minipage}
\caption{[Test Case 1] Contour map estimated by the AugEnKF: (Left) Reference pressure field. (Right) Estimation by Augmented EnKF under 50\% of mixed observation. This is the comparison result by the AugEnKF corresponding to the United Filter's results in Figure \ref{Test1PresField_Contour}.}
\label{Test1ContourPresField_EnKF}
\vspace{-0.4cm}
\end{figure}
\begin{figure}[h!]
\centering
\includegraphics[scale=0.25]{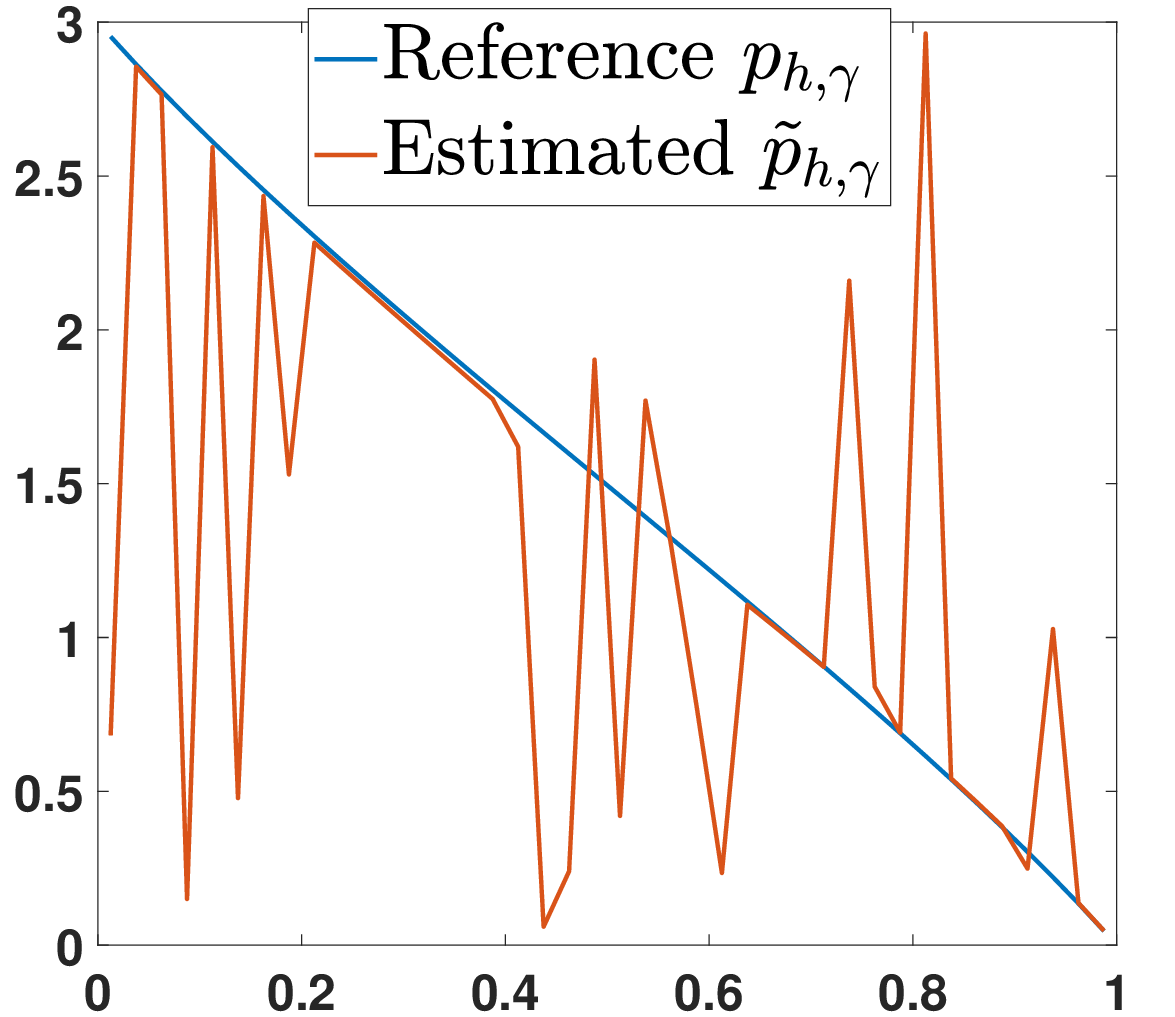}
\caption{[Test Case 1] 1D Concentration on the fracture estimated by the AugEnKF. This is the comparison result by the AugEnKF corresponding to the United Filter's results in Figure \ref{Test1_PresFrt}.}
\label{Test1_1DPres_EnKF}
\vspace{-0.4cm}
\end{figure}
\begin{figure}[h!]
\centering
\begin{minipage}{0.35\textwidth}
\includegraphics[scale=0.25]{./figures_All/figures/figures_update/TrueVelo_Test1.eps}
\end{minipage}%
\begin{minipage}{0.35\textwidth}
\hspace{0.5cm}\includegraphics[scale=0.25]{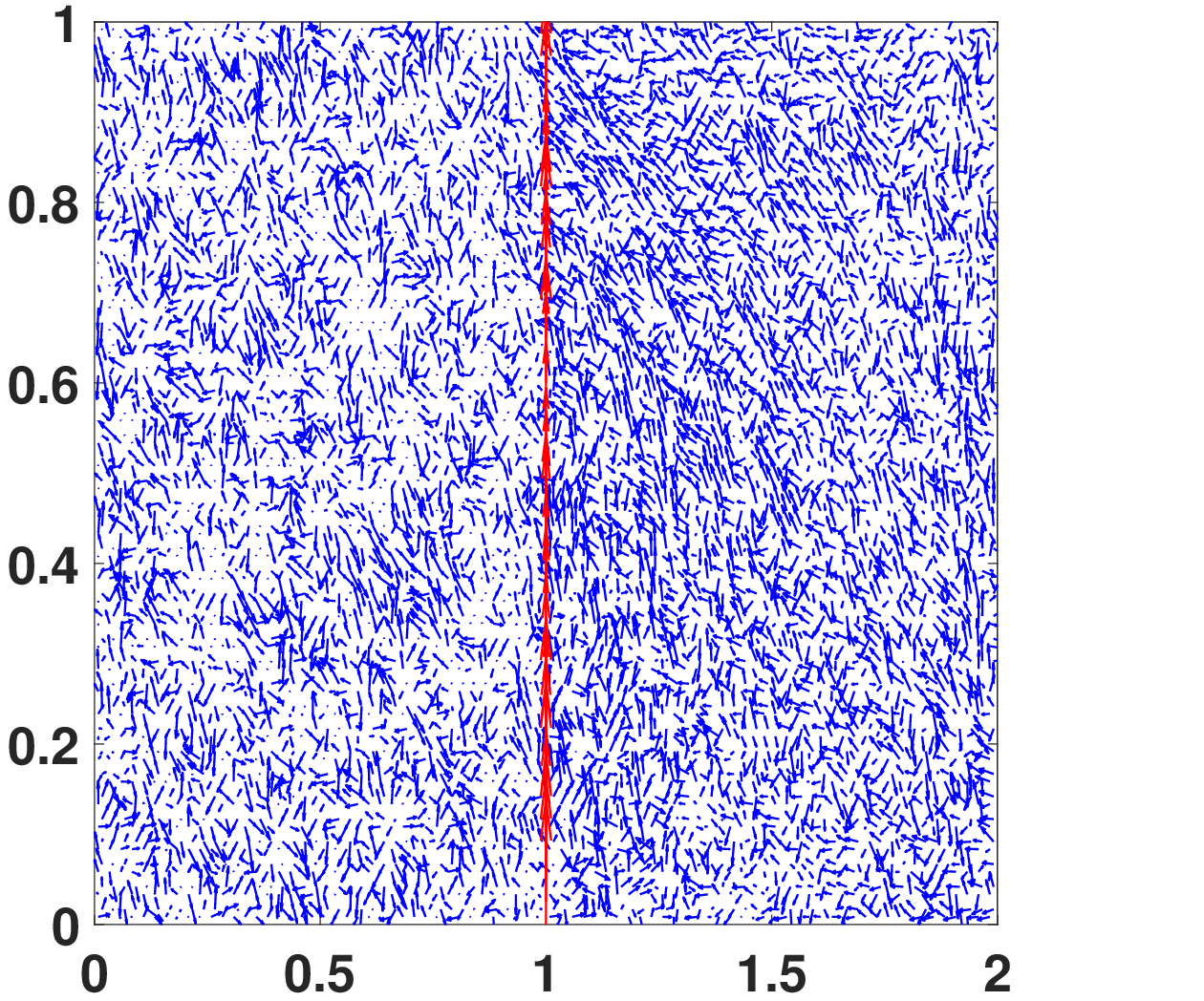}
\end{minipage}
\caption{[Test Case 1] Velocity field estimated by the AugEnKF: (Left) Reference field. (Right) Estimation by the AugEnKF. This is the comparison result by the AugEnKF corresponding to the United Filter's results in Figure \ref{Test1VeloField}. }
\label{Test1VeloField_EnKF}
\vspace{-0.3cm}
\end{figure}

We first present the true pressure field and the AugEnKF estimated pressure field in Figure~\ref{Test1HeatPresField_EnKF}; their corresponding contour maps are presented in Figure~\ref{Test1ContourPresField_EnKF}; and the velocity fields are plotted in Figure~\ref{Test1VeloField_EnKF}, while the 1D pressure on the fracture is displayed in Figure~\ref{Test1_1DPres_EnKF}. From Figures~\ref{Test1HeatPresField_EnKF}, \ref{Test1ContourPresField_EnKF} and \ref{Test1_1DPres_EnKF}, we can see that AugEnKF can only recovered a portion of the pressure from the sparse observation, and the approximate 2D and 1D pressure are very noisy. Similarly, from Figure~\ref{Test1VeloField_EnKF}, it is clear that the estimated velocity field is very chaotic and inaccurate compared to the fields obtained from the United Filter in Figure~\ref{Test1VeloField}. The poor performance of the AugEnKF also stems from the inherent difficulties that the EnKF typically encounters in handling nonlinear problems.

\begin{figure}[h!]
\centering
\begin{minipage}{0.3\textwidth}
\includegraphics[scale=0.21]{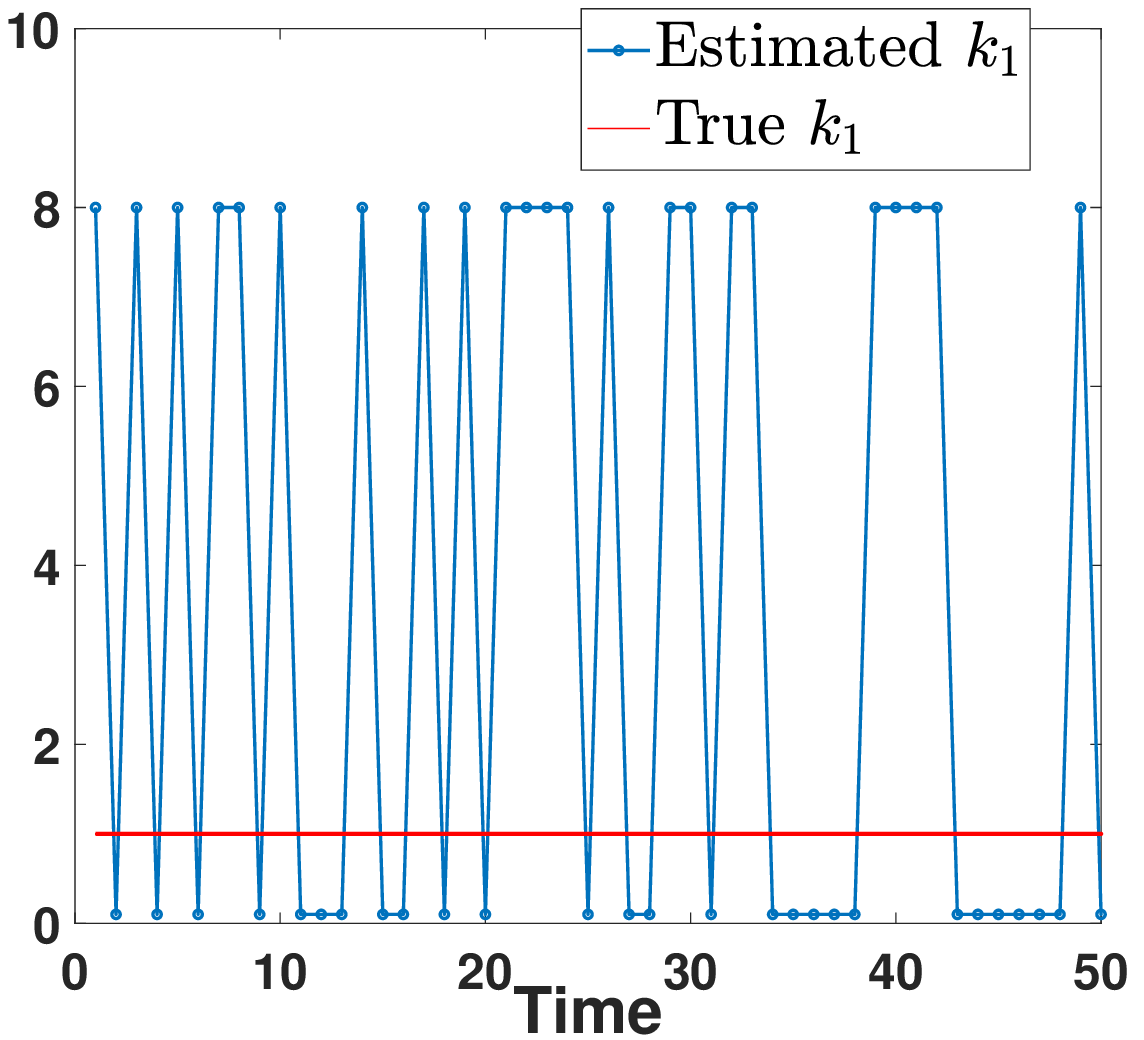}
\end{minipage}%
\begin{minipage}{0.3\textwidth}
\hspace{0.2cm}\includegraphics[scale=0.21]{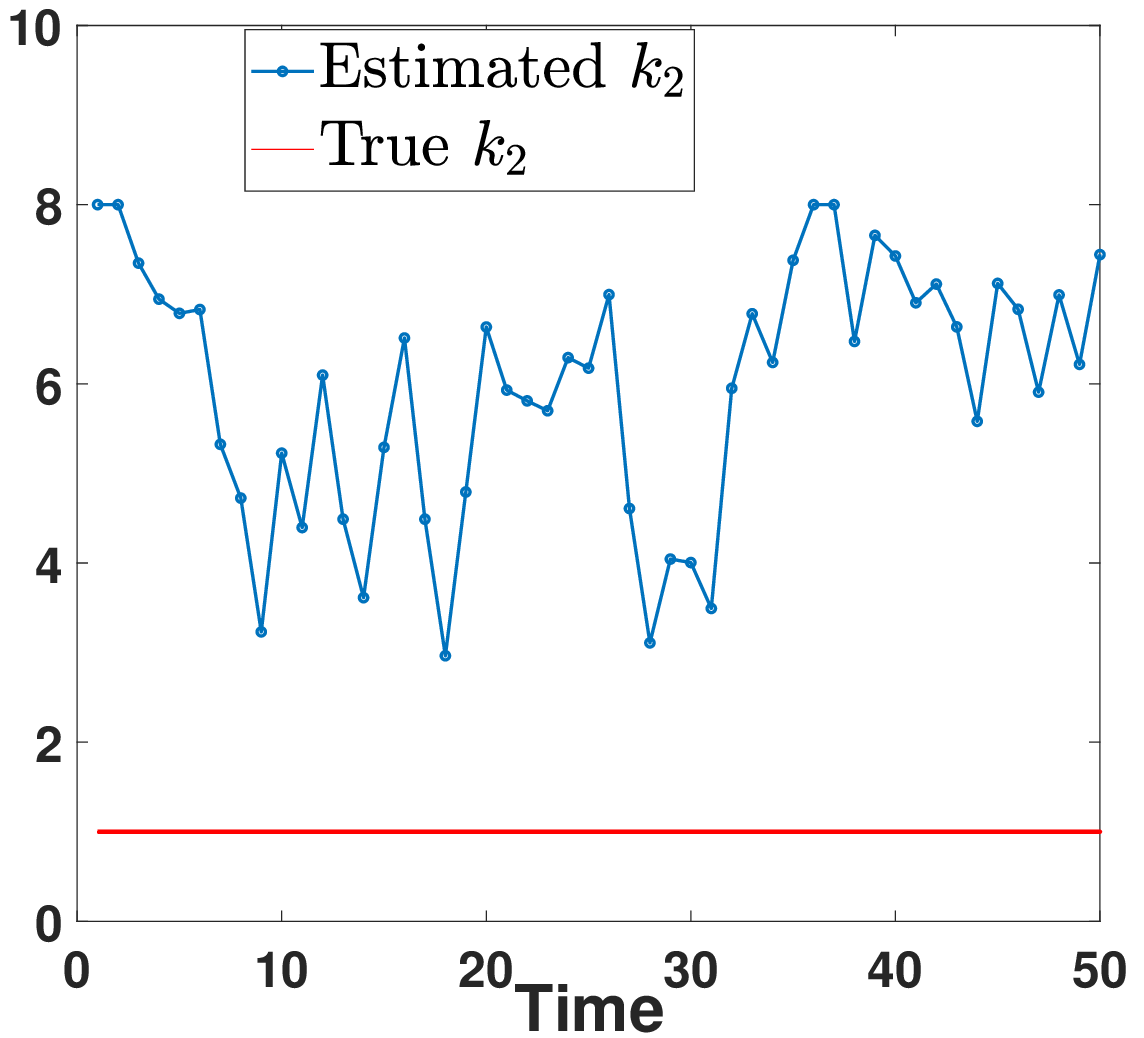}
\end{minipage} %
\begin{minipage}{0.3\textwidth}
\hspace{0.3cm}\includegraphics[scale=0.21]{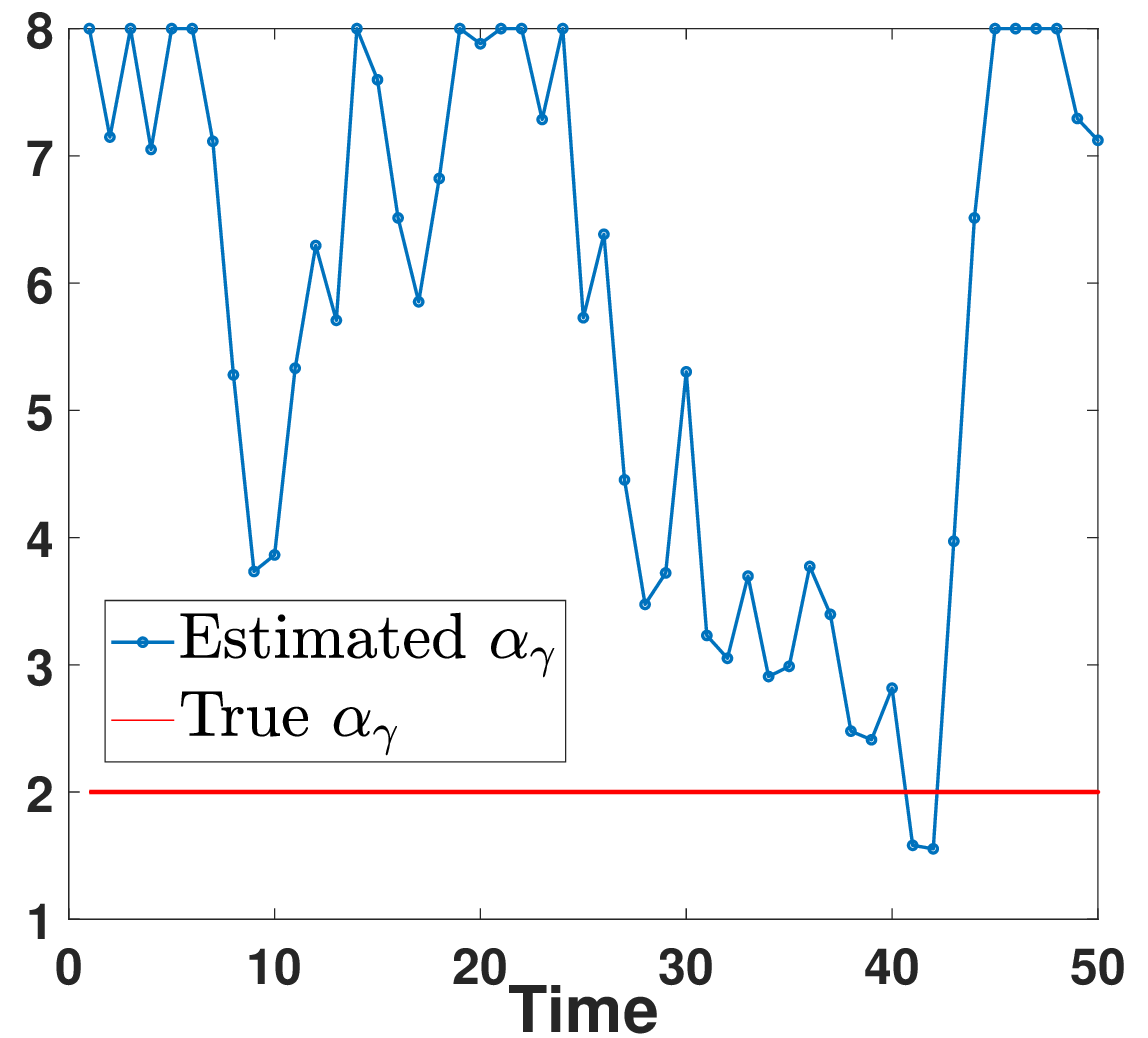}
\end{minipage}
\caption{[Test Case 1] Estimated Parameters obtained by the AugEnKF: (Left) Parameter $k_1$. (Middle) Parameter $k_2$. (Right) Parameter  $\alpha_{\gamma}$. This is the comparison result by the AugEnKF corresponding to the United Filter's results in Figure \ref{ParaEst_Pure}.}
\label{Test1_ParaEst_Pure_EnKF}
\vspace{-0.4cm}
\end{figure}

Then, in Figure~\ref{Test1_ParaEst_Pure_EnKF} we present the AugEnKF estimated parameters. It is apparent that the Kalman filter samples do not converge to the true values and remain unstable throughout the filtering process. Moreover, we also observed a strange behavior where the samples only take two values $0.1$ and $8$, rather than varying continuously, as shown in the first plot of Figure~\ref{Test1_ParaEst_Pure_EnKF}. This occurs because we enforce the particles to remain within the range [0, 8] to prevent divergence. Overall, the above results indicate that the United Filter method significantly outperforms AugEnKF in the setting of Test Case 1.  

\subsection{Test Case 2:  Pure diffusion equation with anisotropic and heterogeneous fracture permeability} \label{subsec4_2}
\begin{figure}[h!]
\vspace{-0.4cm}
\centering
\includegraphics[scale=0.3]{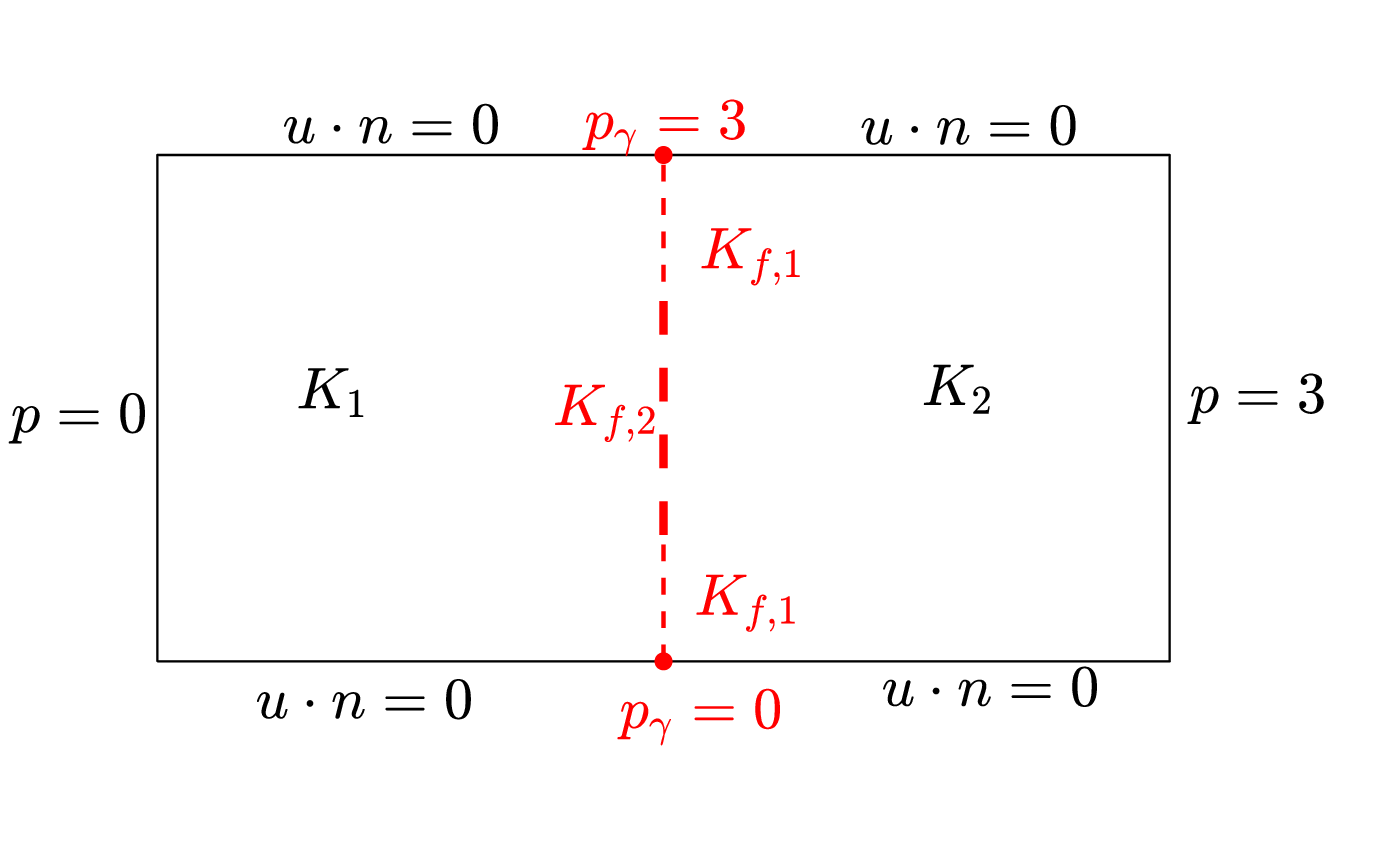}
\vspace{-0.8cm}
\caption{[Test Case 2] Setting and boundary conditions.}
\label{Test2BCs}
\end{figure}

For the second test case, we adapt an experiment introduced in \cite{Frih2012} where the permeability of the fracture is heterogeneous and anisotropic. More precisely, in the upper and lower quarters of the fracture, the fracture permeability is
$\begin{array}{l}
\pmb{K}_{f} = \pmb{K}_{f, 1} = \begin{bmatrix}
1 / k_f & 0 \vspace{0.1cm} \\
0 & k_f
\end{bmatrix},
\end{array}
$
while in the middle part, we have
$\begin{array}{l}
\pmb{K}_{f} = \pmb{K}_{f, 2} = \begin{bmatrix}
k_f & 0 \vspace{0.1cm} \\
0 & 1 / k_f
\end{bmatrix}.
\end{array}$
Assume that $K_i = k_i\pmb{I}, \; i=1, 2$, the true values of the parameters in this test case are $k_1 = k_2=1, k_f = 2000$. The boundary conditions for Test Case 2 are shown in Figure~\ref{Test2BCs}. For this case, the fluid can hardly flow along the fracture in the central zone but can easily cross it. However, in the two other parts of the fracture, the fluid tends to flow along the fracture, but cannot cross it \cite{Frih2012}. The reduced model \eqref{reduced_subdomain}-\eqref{reduced_fracture} is not applicable to Test Case 2 since it is only valid when the fracture has higher permeability than the rock matrix. Instead, we use a more general reduced model developed in \cite{Jaffre2005} as the forward solver for this test case. 

Assume that $\pmb{K}_f$ composes of a tangential part $\pmb{K}_{f, \tau}$  and a normal part $\pmb{K}_{f, \nu}$, the reduced fracture model of \eqref{original_problem} in this test case reads as follows \cite{Jaffre2005}:
\begin{equation}
\label{reduced_subdomain_Test2}
\begin{array}{rcll}
\partial_t{p_i}+\text{div }\bu_i&=&q_{i} &\text{ in } \Omega_i\times (0, T), \;i=1, 2, \\
\bu_i&=&-\bK_i\nabla{p_i} &\text{ in } \Omega_i\times (0, T), \; i=1, 2, \\
p_i&=&0 &\text{ on } \left(\partial\Omega_i \cap \partial\Omega\right) \times (0, T), \; i=1, 2, \\
\partial_t{p_{\gamma}}+\text{div}_{\tau }\bu_{\gamma}&=& q_{\gamma} +\sum\limits^{2}_{i=1}\left( \bu_i \cdot \bn_i\right)_{\vert \gamma} & \text{ in } \gamma \times (0, T), \\
\bu_{\gamma} &=&{-\bK_{\gamma}\nabla_{\tau}p_{\gamma}} & \text{ in } \gamma \times (0, T), \\
\kappa_{\gamma}(p_i-p_{\gamma})&=& \xi \bu_i\cdot\bn_i - \bar{\xi}\bu_{3-i}\cdot\bn_{3-i} &\text{ on } \gamma \times (0, T), i=1, 2, \\
p_{\gamma}&=&0 &\text{ on } \partial\gamma \times (0, T), \\
p_i(\cdot, 0)&=&p_{0, i} &\text{ in } \Omega_i, \; i=1, 2,\\
p_{\gamma}(\cdot, 0)&=&p_{0, \gamma} & \text{ in } \gamma,
\end{array} \vspace{-0.2cm}
\end{equation}
where $\pmb{K}_{\gamma} = \pmb{K}_{f, \tau}\delta$ and $\kappa_{\gamma} = \dfrac{2\pmb{K}_{f, \nu}}{\delta}$, $\xi$ is a parameter with $\xi > \dfrac{1}{2}$ and $\bar{\xi} = 1-\xi$. We fix $\xi = 1$ in \eqref{reduced_subdomain_Test2} for Test Case 2. Unlike Test Case 1, the solution for this test case is irregular, and the pressure on the fracture is discontinuous (see, e.g., \cite{Frih2012}). 

Due to the irregularity in the 2D and 1D pressure on the subdomains on the fracture, Test Case 2 is indeed more challenging than Test Case 1. Therefore, we increase the number of Direct Filter particles to $M = 50$, and the number of iteration for the United Filter to $R=4$. Moreover, due to the heterogeneity of $\pmb{K}_f$, we have $\pmb{K}_f = \pmb{K}_{f, 2}$ in the central part of the fracture, meaning $K_{f, \nu} =2000$ and $K_{f, \tau} = 1/2000$ , which implies that $K_{\gamma} = 5\times 10^{-7}$ and $\kappa_{\gamma} = 4\times 10^6$. On the other hand, in the upper and lower quarters of the fracture, the values are reversed: $K_{f, \tau} =2000$ and $\pmb{K}_{f, \nu} = 1/2000$, which implies that $K_{\gamma}=2$ and $\kappa_{\gamma} = 1$. Thus, to simplify the parameter estimation process, in addition to estimating the permeability $k_i, i=1, 2$ on the subdomains as in Test Case 1, we fix $\delta = 0.001$ and only aim to estimate $k_{f}$, which characterizes the permeability on the fracture, rather than the transmissivity of the fracture. To obtain the synthetic observational data, we solve the fully discretization problem of the reduced model with the spatial mesh size $h = 1/40$ and with the timestep $\Delta{t}_{\text{ref}} = T/800$ where $T=1$. The United Filter algorithm is performed on the same mesh size and with filtering timesteps given by $\Delta{t}_{\text{Filter}} = T/Nt_{\text{Filter}}$. Our aim in Test case 2 is to investigate the performance of the United Filter algorithm under varying noise levels. More specifically, we introduce three levels of disturbance in the state model \eqref{StateModel}: $\omega^1 = 0.001\sqrt{\Delta{t}_{\text{Filter}}}\varepsilon^1$, $\omega^2 = 0.01\sqrt{\Delta{t}_{\text{Filter}}}\varepsilon^2$, and $\omega^3 = 0.1\sqrt{\Delta{t}_{\text{Filter}}}\varepsilon^3$ where $\varepsilon^i \sim \mathcal{N}(0, \pmb{I}_l)$ for $i=1, 2, 3$.
\begin{figure}[h!]
\vspace{-0.3cm}
\centering
\begin{minipage}{0.25\textwidth}
\includegraphics[scale=0.22]{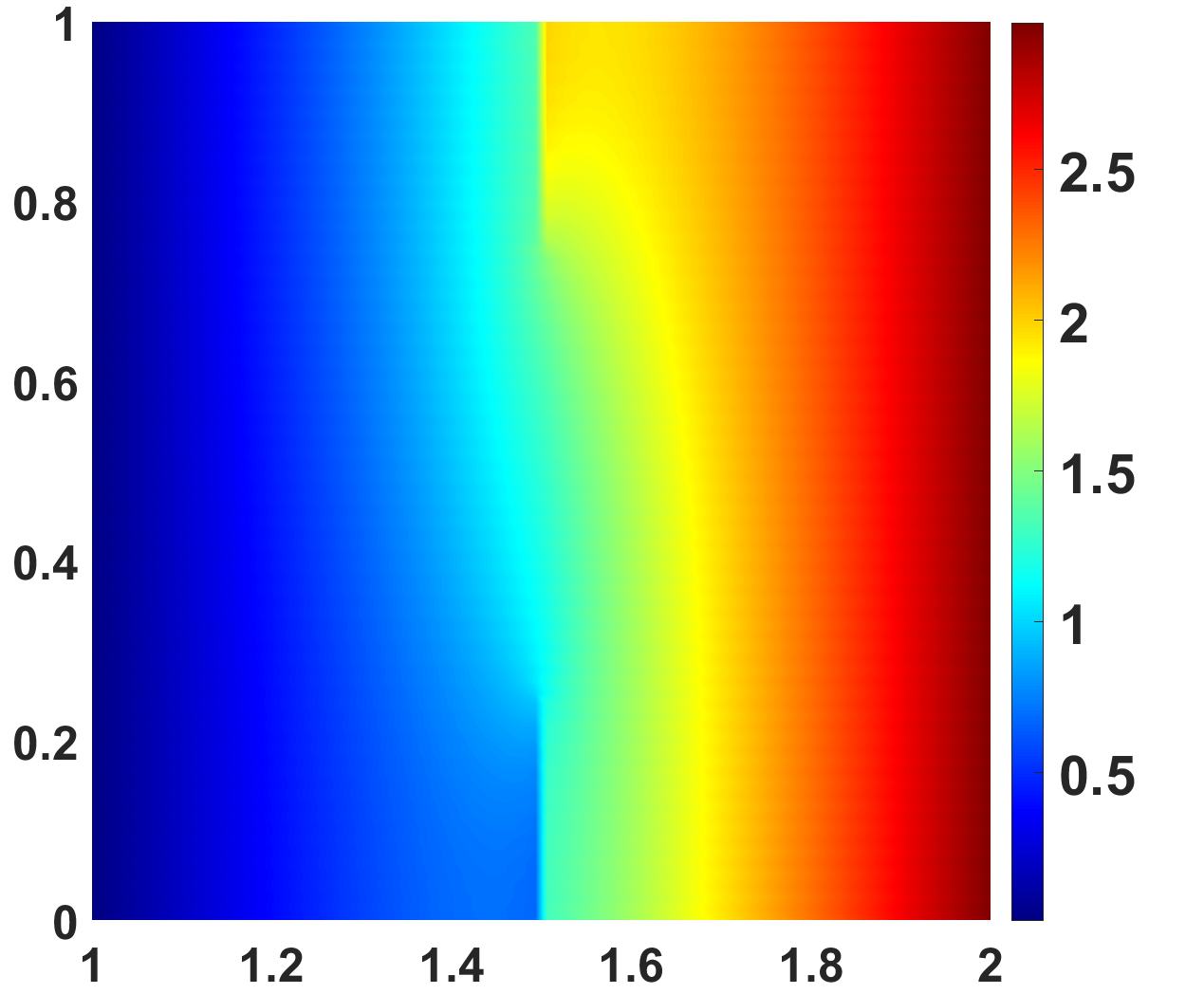}
\end{minipage}%
\begin{minipage}{0.24\textwidth}
\hspace{0.1cm}\includegraphics[scale=0.22]{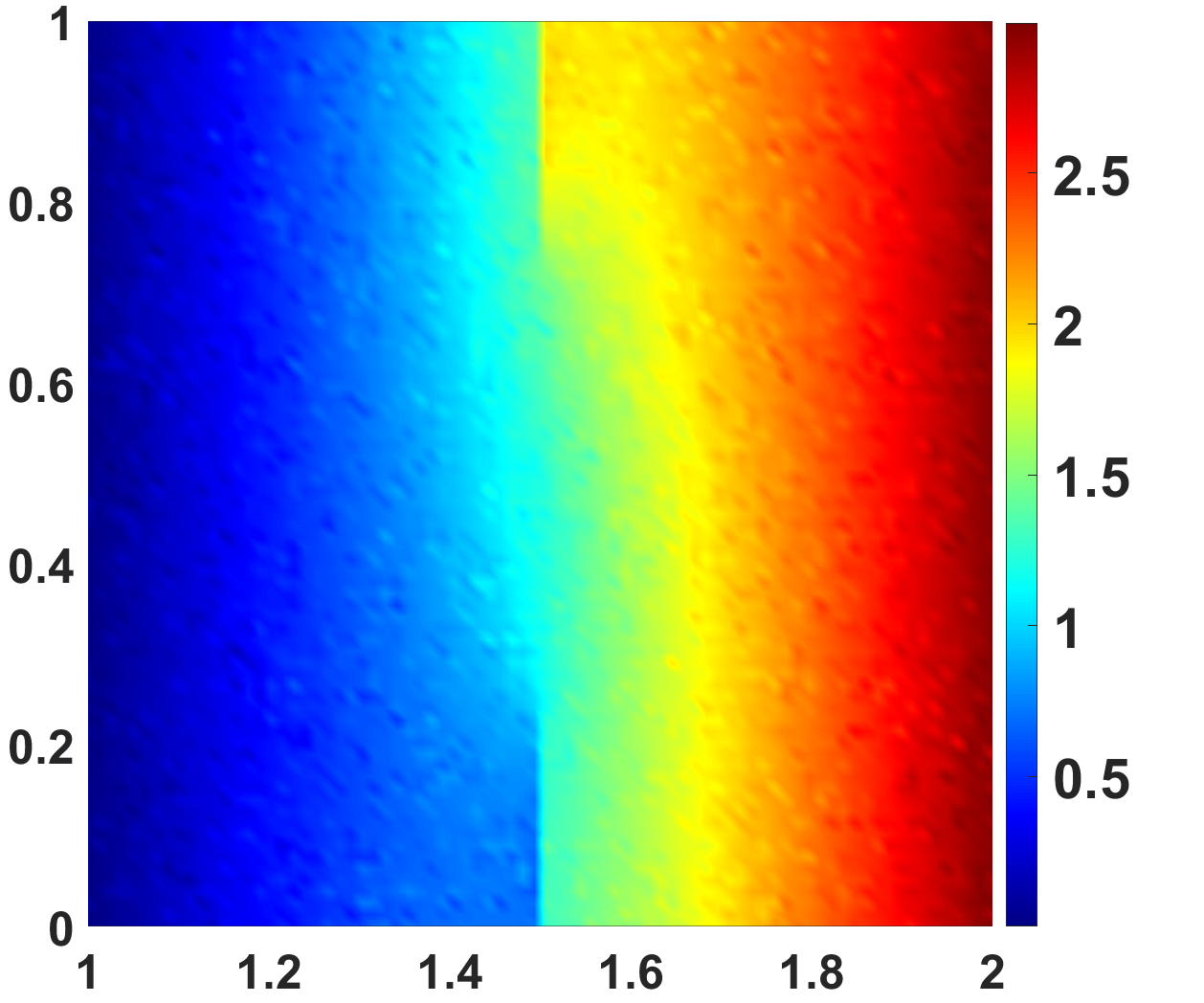}
\end{minipage} %
\begin{minipage}{0.226\textwidth}
\includegraphics[scale=0.22]{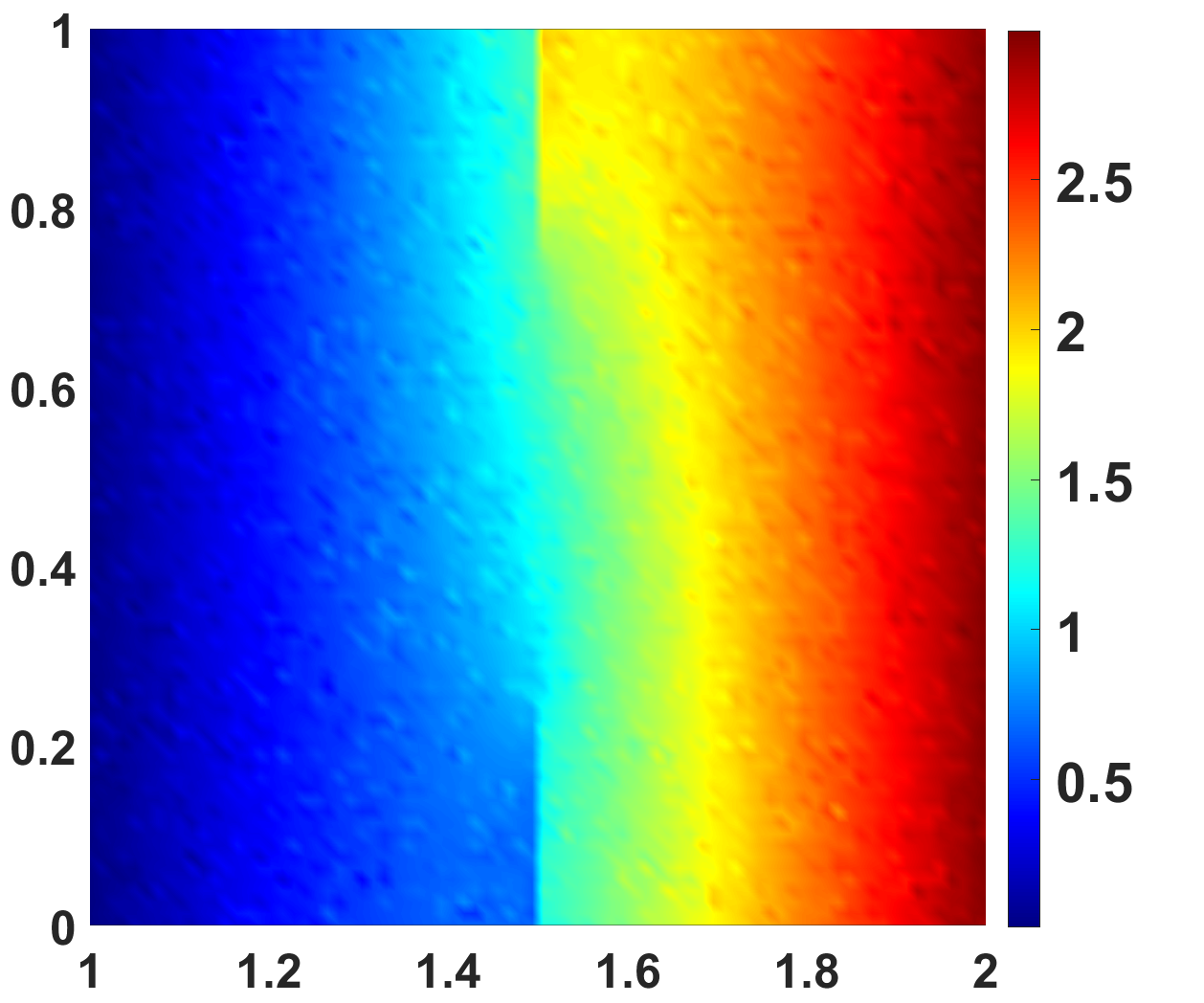}
\end{minipage}
 \hspace{0.1cm}
\begin{minipage}{0.23\textwidth}
\hspace{0.1cm}\includegraphics[scale=0.22]{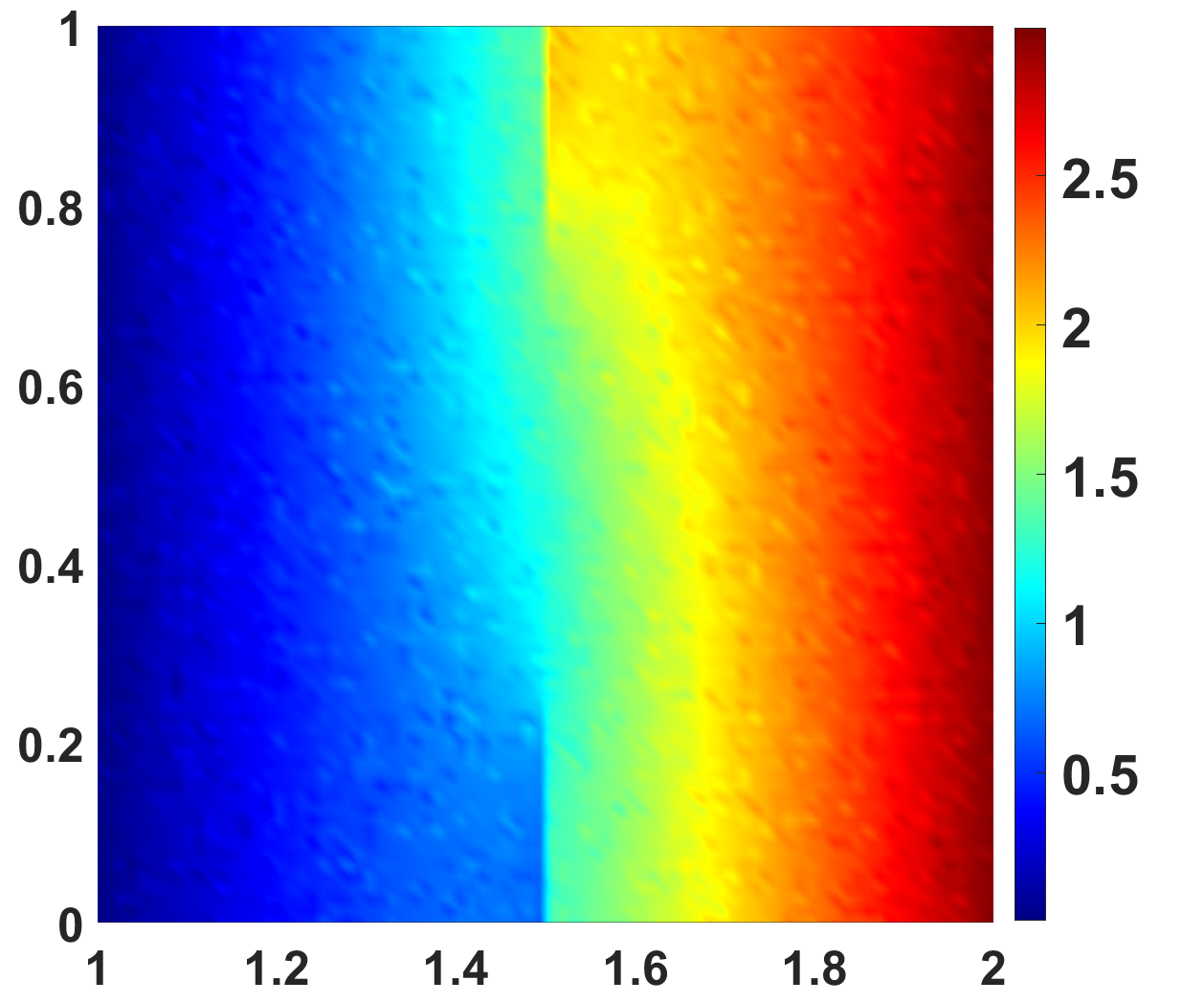}
\end{minipage}
\caption{[Test Case 2] Heat map illustrating the accuracy of the United Filter’s pressure state estimation. (First) Reference field. (Second) With $\omega_1$. (Third) With $\omega_2$. (Fourth) With $\omega_3$.}
\label{HeatPresField_WithNoise}
\vspace{-0.3cm}
\end{figure}
\begin{figure}[h!]
\vspace{-0.2cm}
\centering
\begin{minipage}{0.248\textwidth}
\includegraphics[scale=0.22]{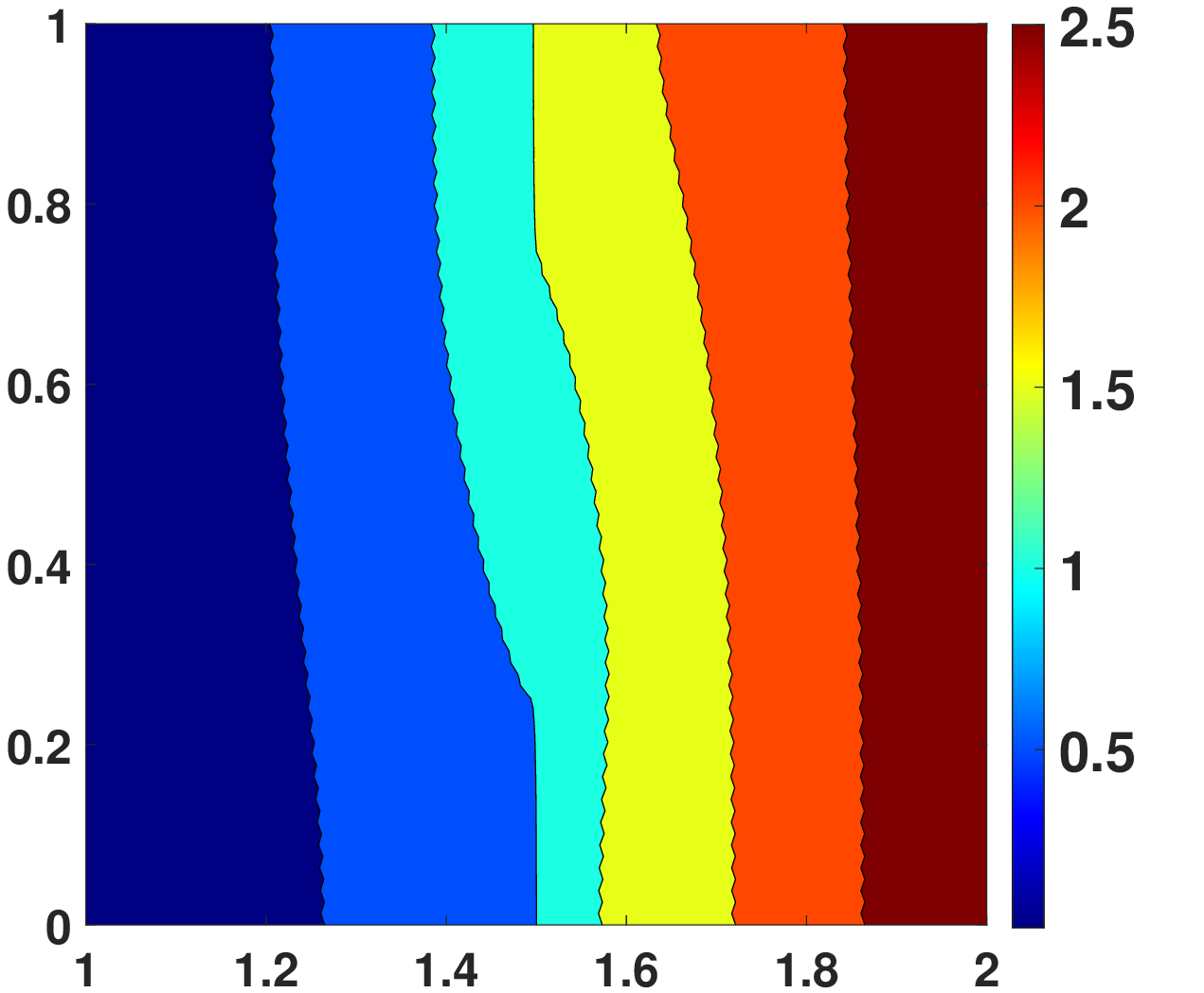}
\end{minipage}%
\begin{minipage}{0.226\textwidth}
\hspace{0.1cm}\includegraphics[scale=0.22]{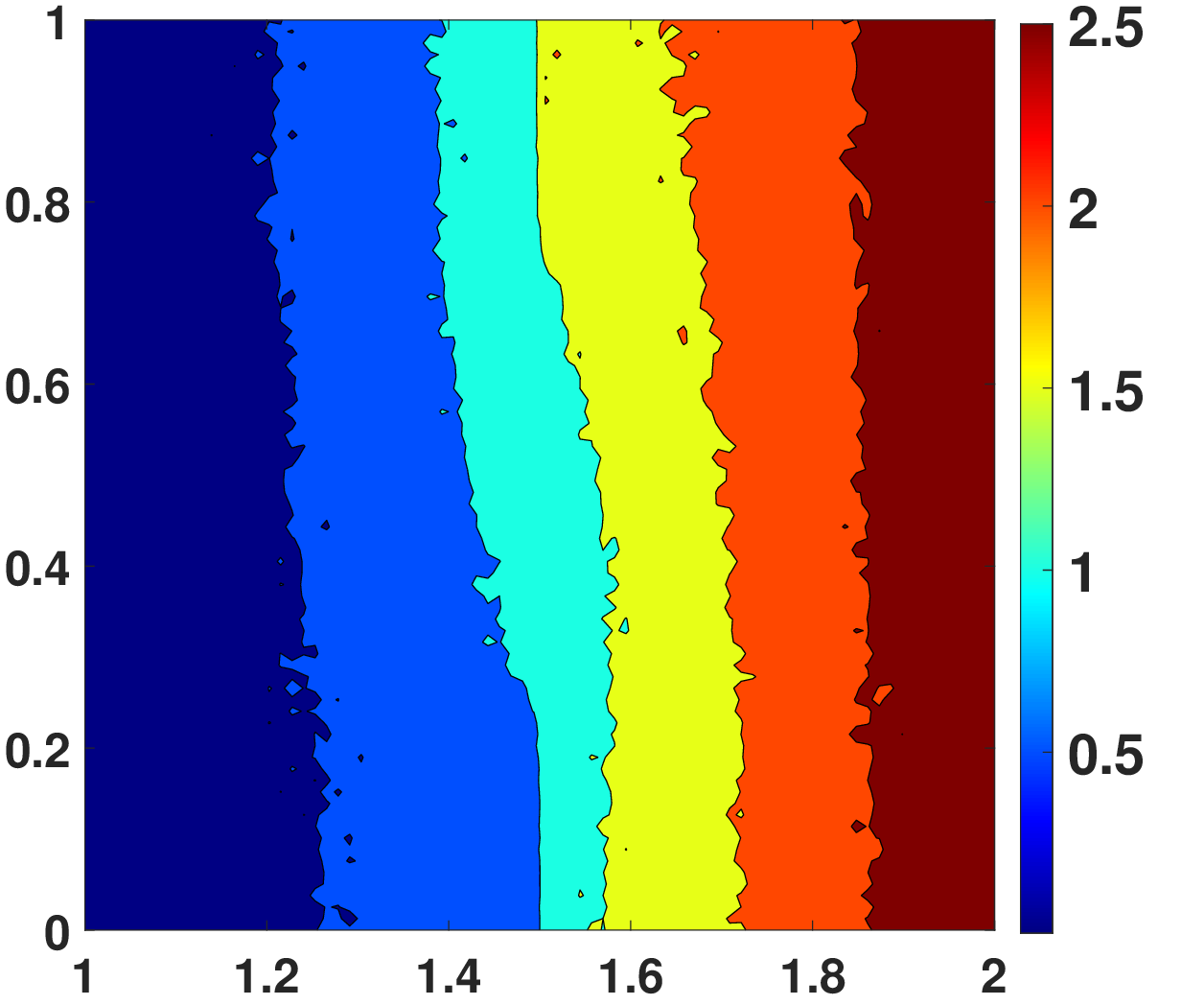}
\end{minipage}
\hspace{0.1cm}
\begin{minipage}{0.237\textwidth}
\hspace{0.1cm}\includegraphics[scale=0.22]{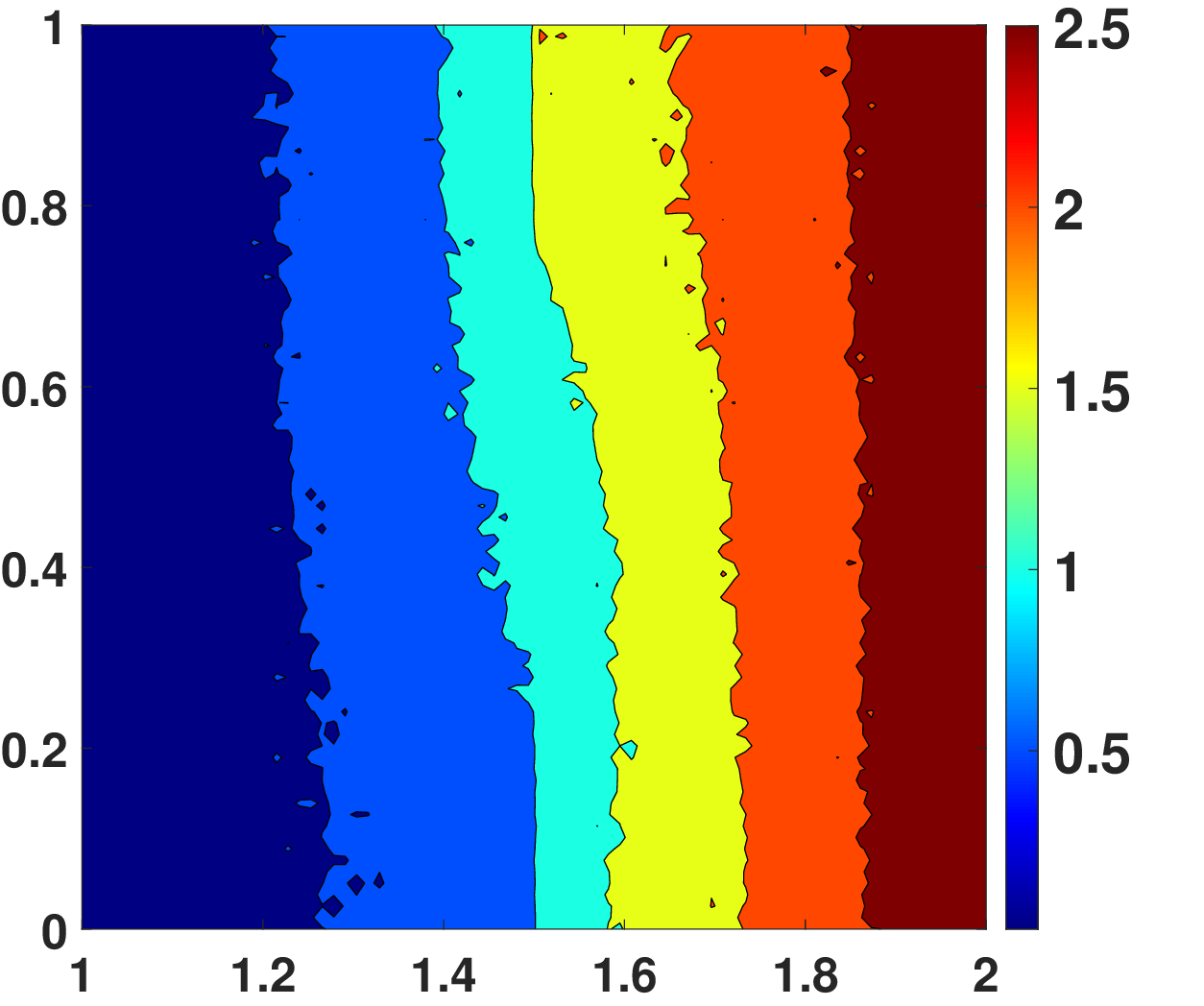}
\end{minipage} 
\begin{minipage}{0.23\textwidth}
\hspace{0.1cm}\includegraphics[scale=0.22]{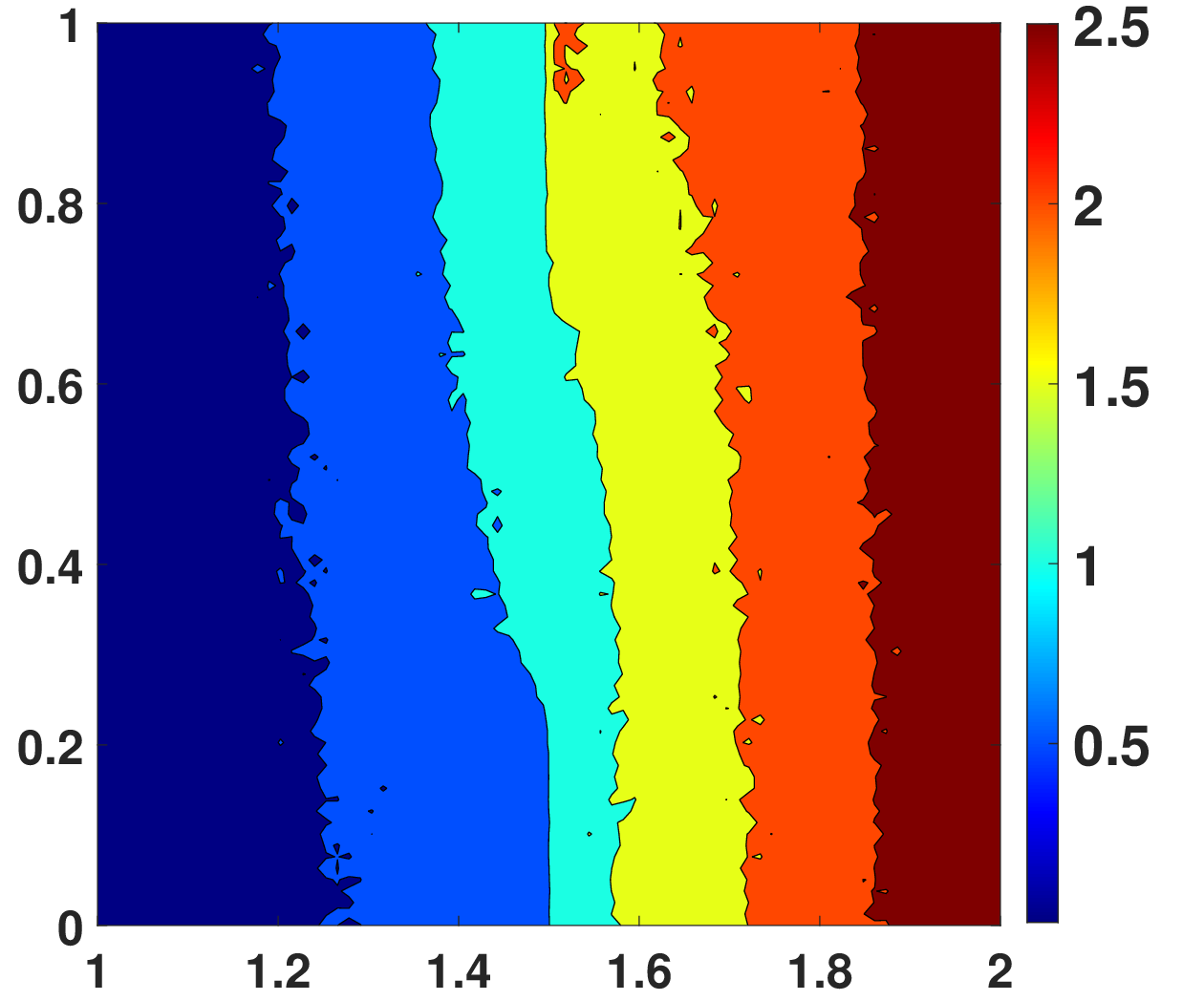}
\end{minipage} 
\caption{[Test Case 2] Contour map illustrating the accuracy of the United Filter’s pressure state estimation. (First) Reference field. (Second) With $\omega^1$. (Third) With $\omega^2$. (Fourth) With $\omega^3$.}
\label{ContourMap_UnitedF_WithNoise}
\vspace{-0.3cm}
\end{figure}
\begin{figure}[h!]
\centering
\begin{minipage}{0.25\textwidth}
\includegraphics[scale=0.21]{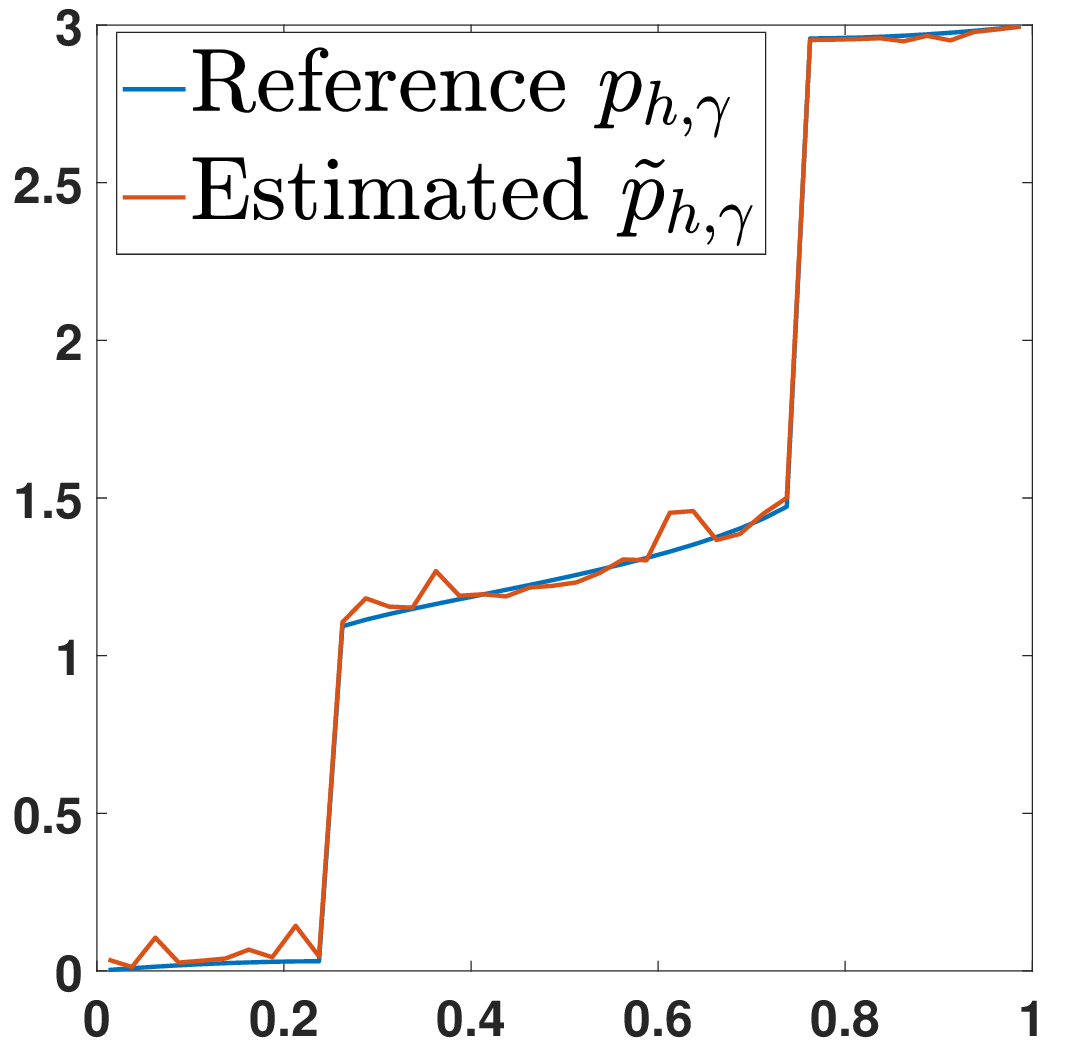}
\end{minipage}%
\begin{minipage}{0.25\textwidth}
\hspace{0.3cm}\includegraphics[scale=0.21]{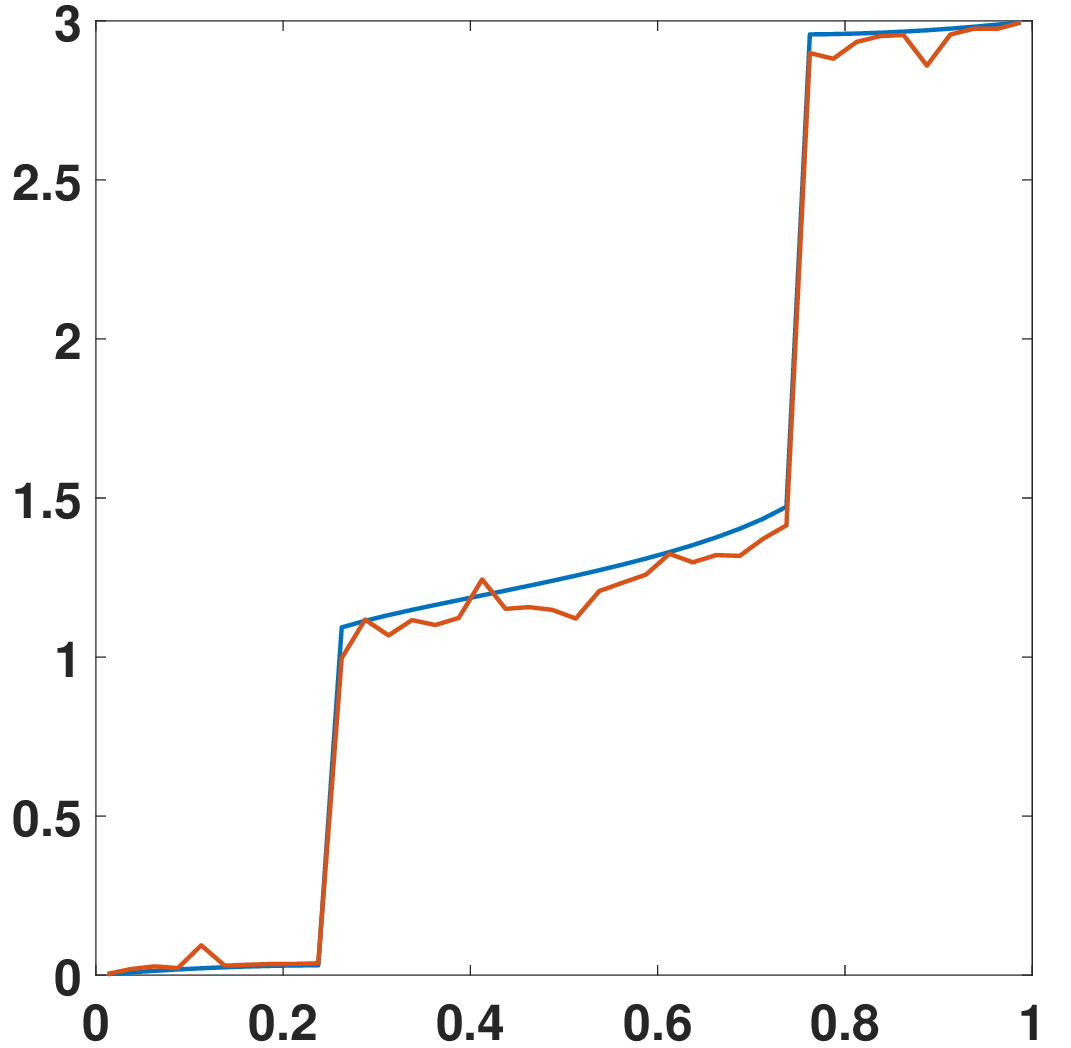}
\end{minipage} %
\begin{minipage}{0.25\textwidth}
\includegraphics[scale=0.21]{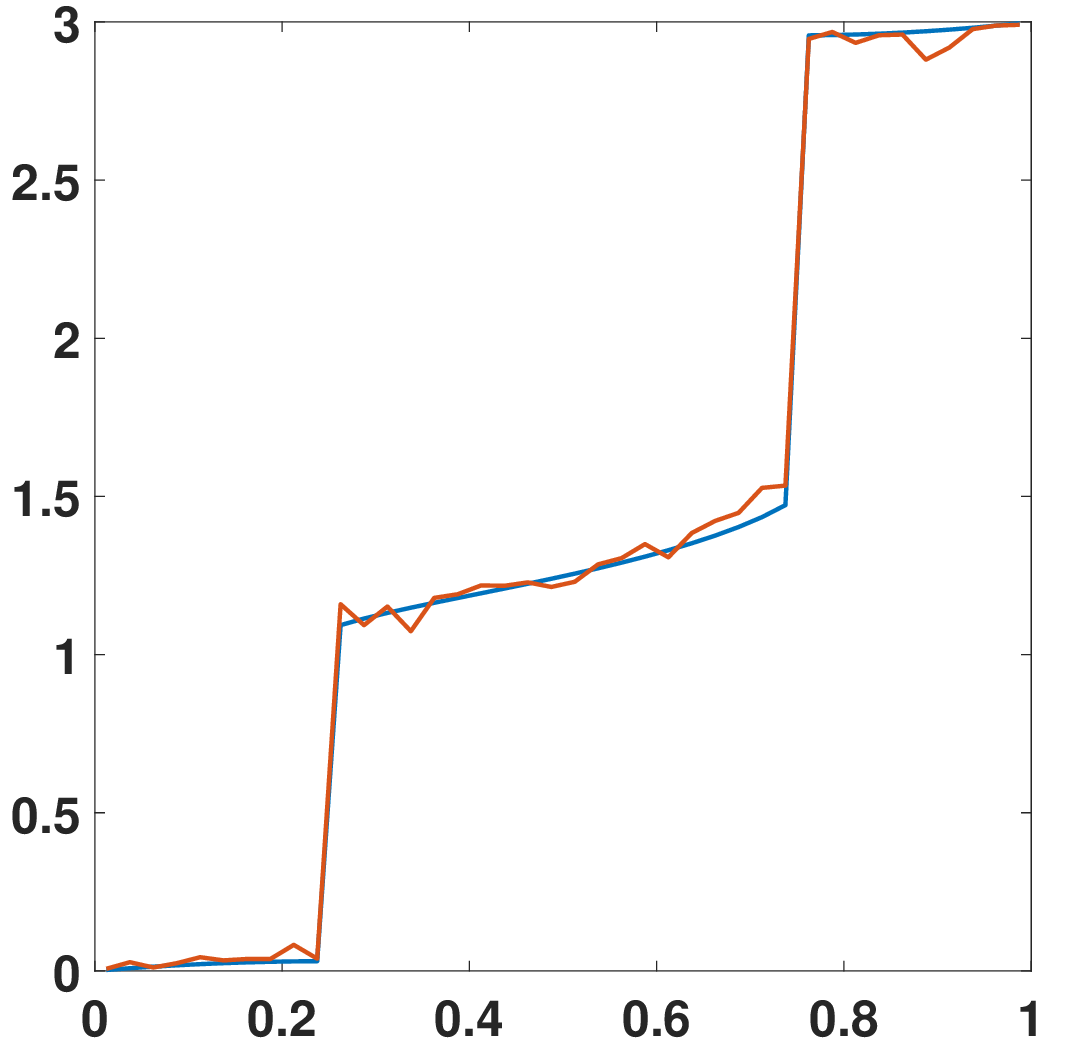}
\end{minipage}
\caption{[Test Case 2] 1D Pressure on the fracture estimated by the United Filter. [Left] With $\omega^1$. [Middle] With $\omega^2$. [Right] With $\omega^3$.}
\label{Test2_1DPresFrt_General_UnitedF_WithNoise}
\label{Test2_1DPresFrt_General_UnitedF_WithNoise}
\vspace{-0.3cm}
\end{figure}

\begin{figure}[h!]
\centering
\begin{minipage}{0.25\textwidth}
\includegraphics[scale=0.22]{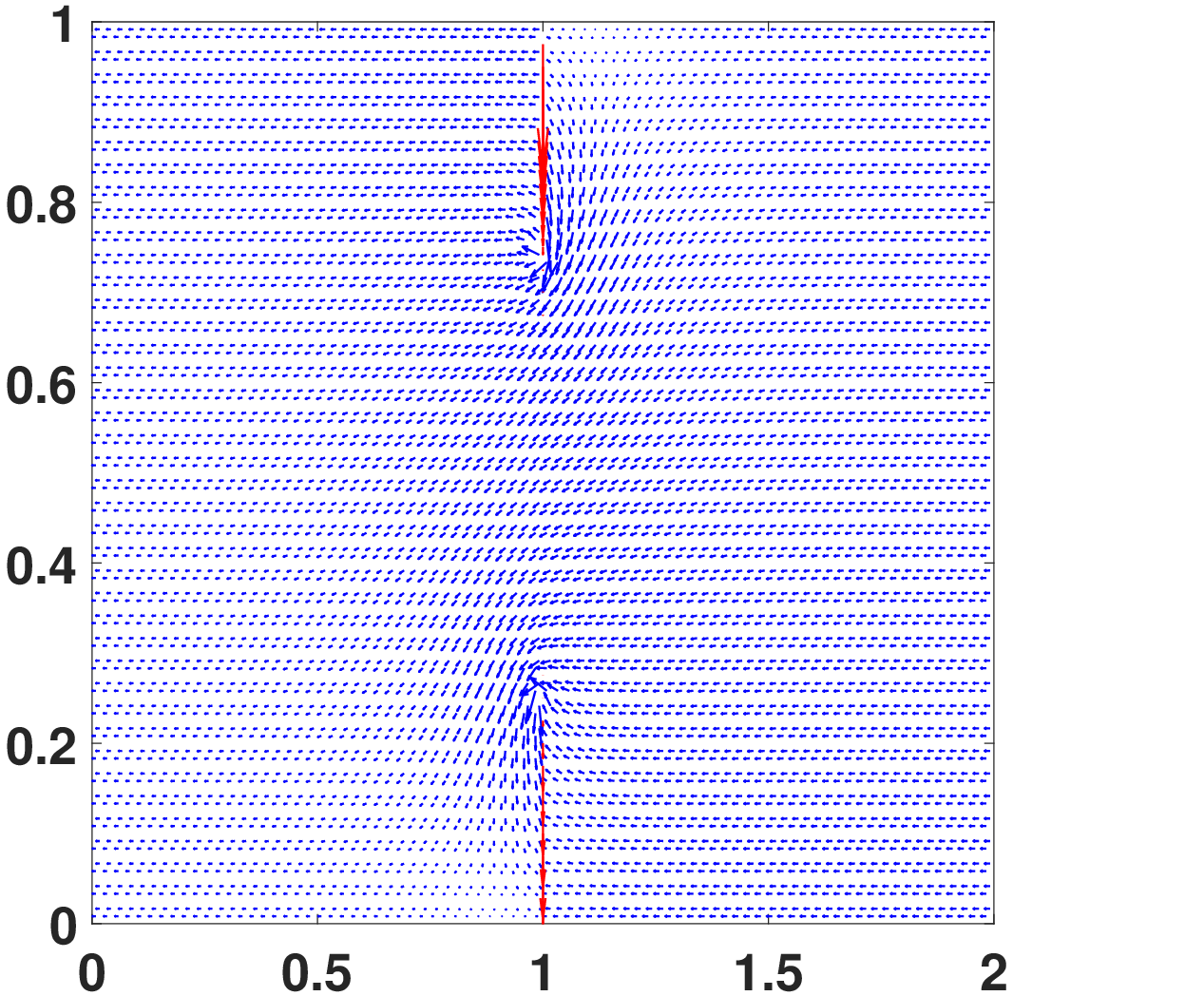}
\end{minipage}%
\begin{minipage}{0.23\textwidth}
\hspace{0.1cm}\includegraphics[scale=0.22]{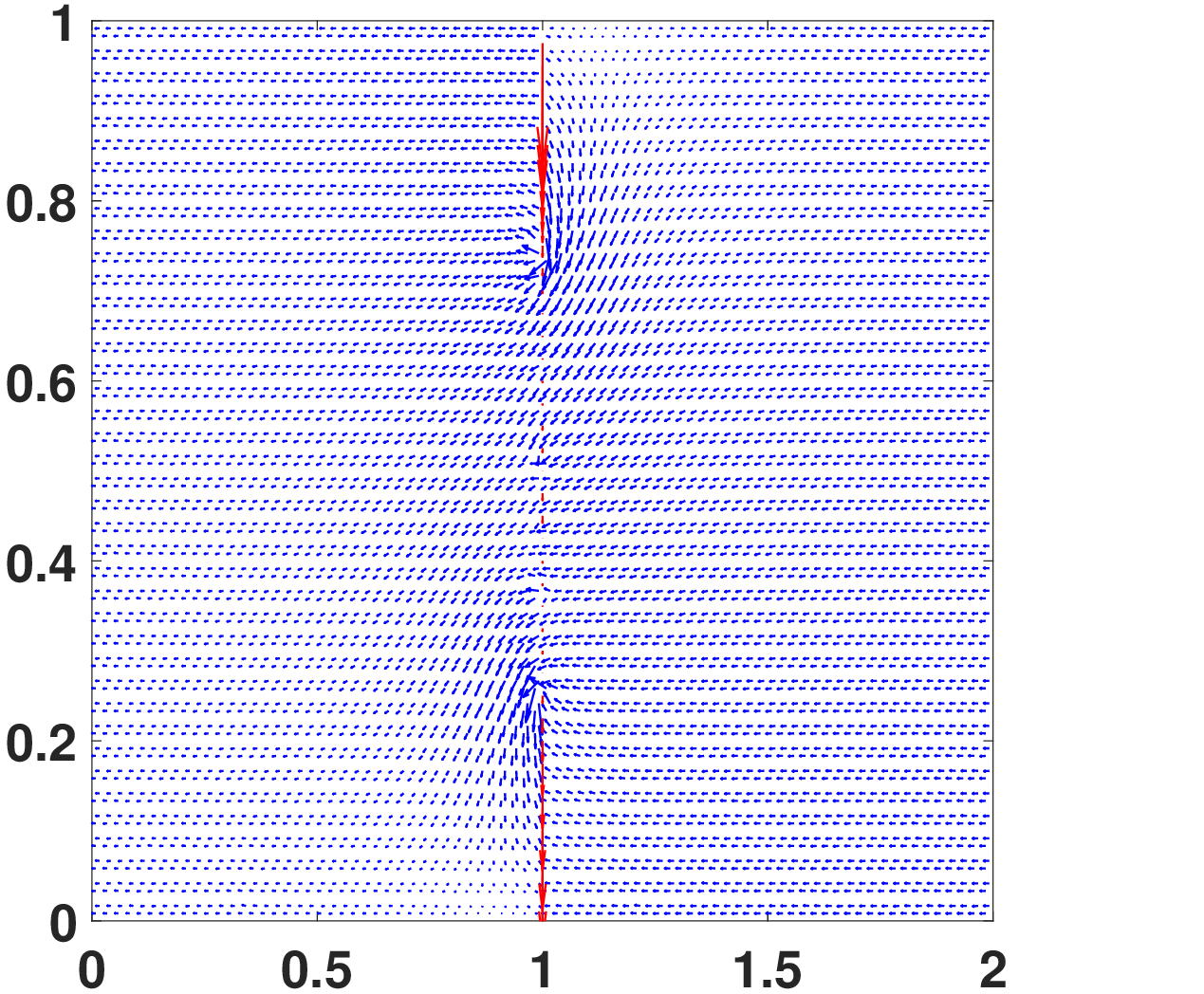}
\end{minipage} %
\hspace{0.1cm}
\begin{minipage}{0.25\textwidth}
\includegraphics[scale=0.22]{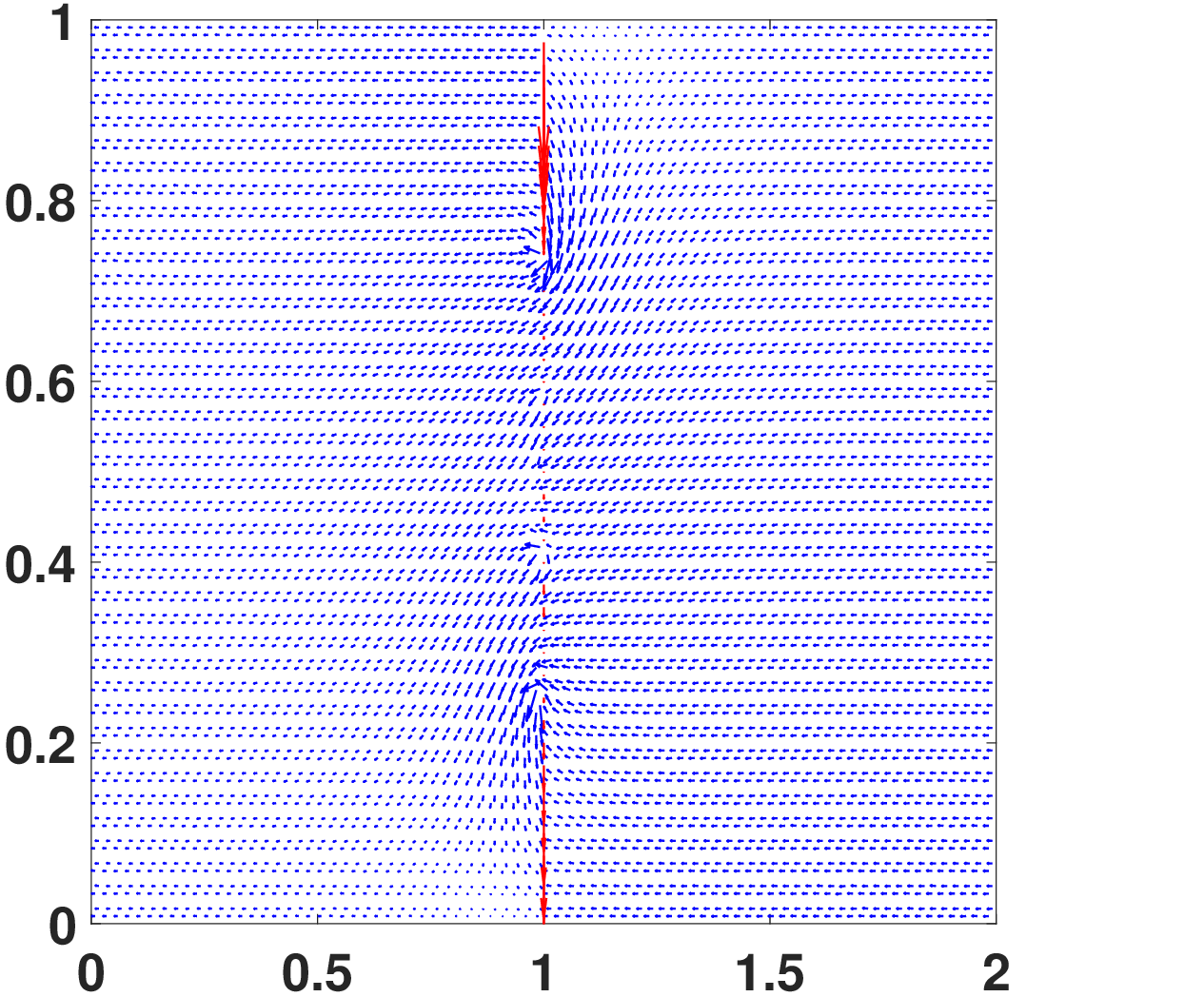}
\end{minipage}%
\begin{minipage}{0.23\textwidth}
\hspace{0.1cm}\includegraphics[scale=0.22]{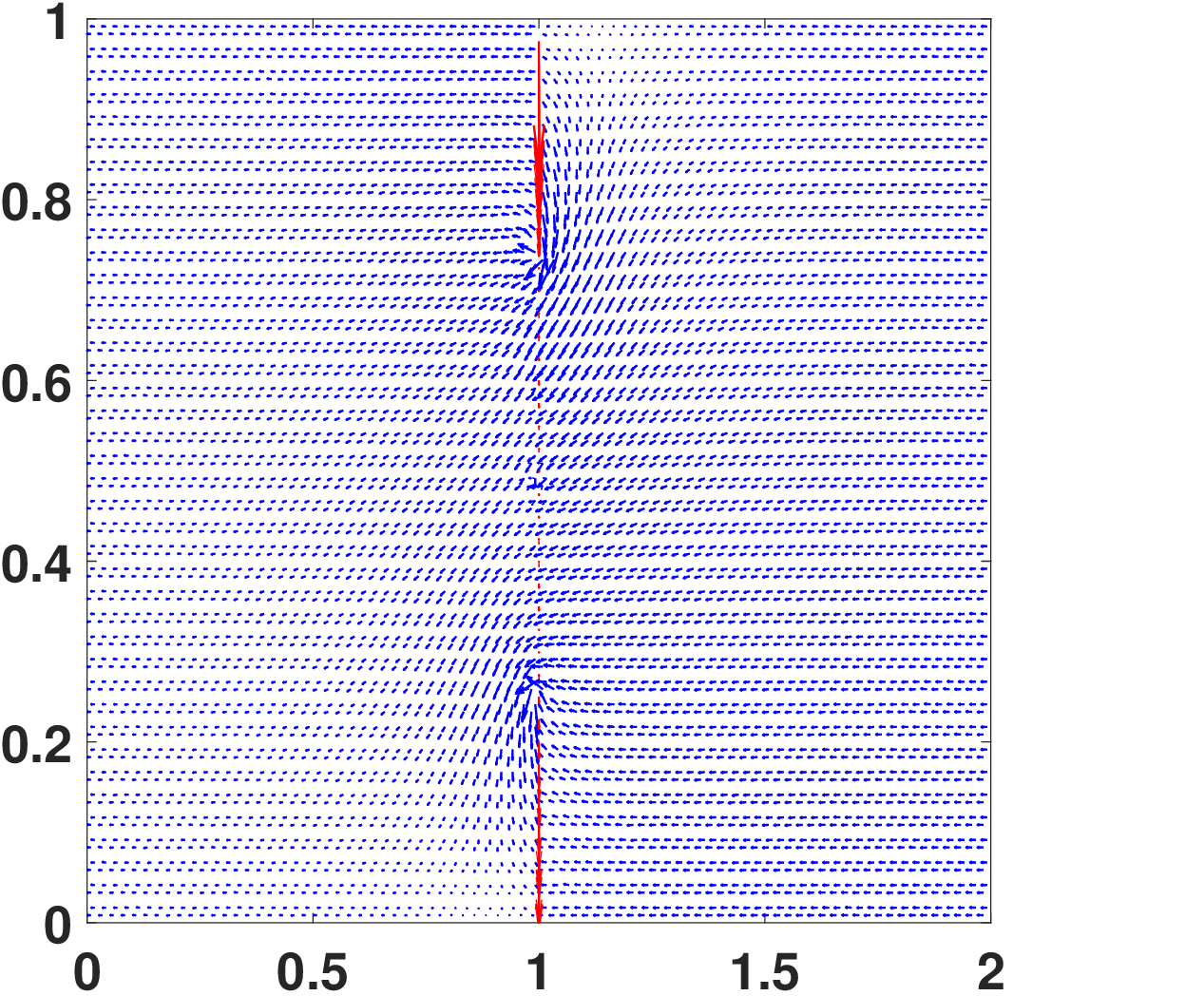}
\end{minipage} 
\caption{[Test Case 2] Demonstration of the United Filter’s accuracy in velocity state estimation. (First) Reference velocity field. (Second) With $\omega^1$. (Third) With $\omega^2$. (Fourth) With $\omega^3$.}
\label{VelField_General}
\vspace{-0.2cm}
\end{figure}
\begin{figure}[h!]
\centering
\begin{minipage}{0.4\textwidth}
\hspace{1cm}\includegraphics[scale=0.25]{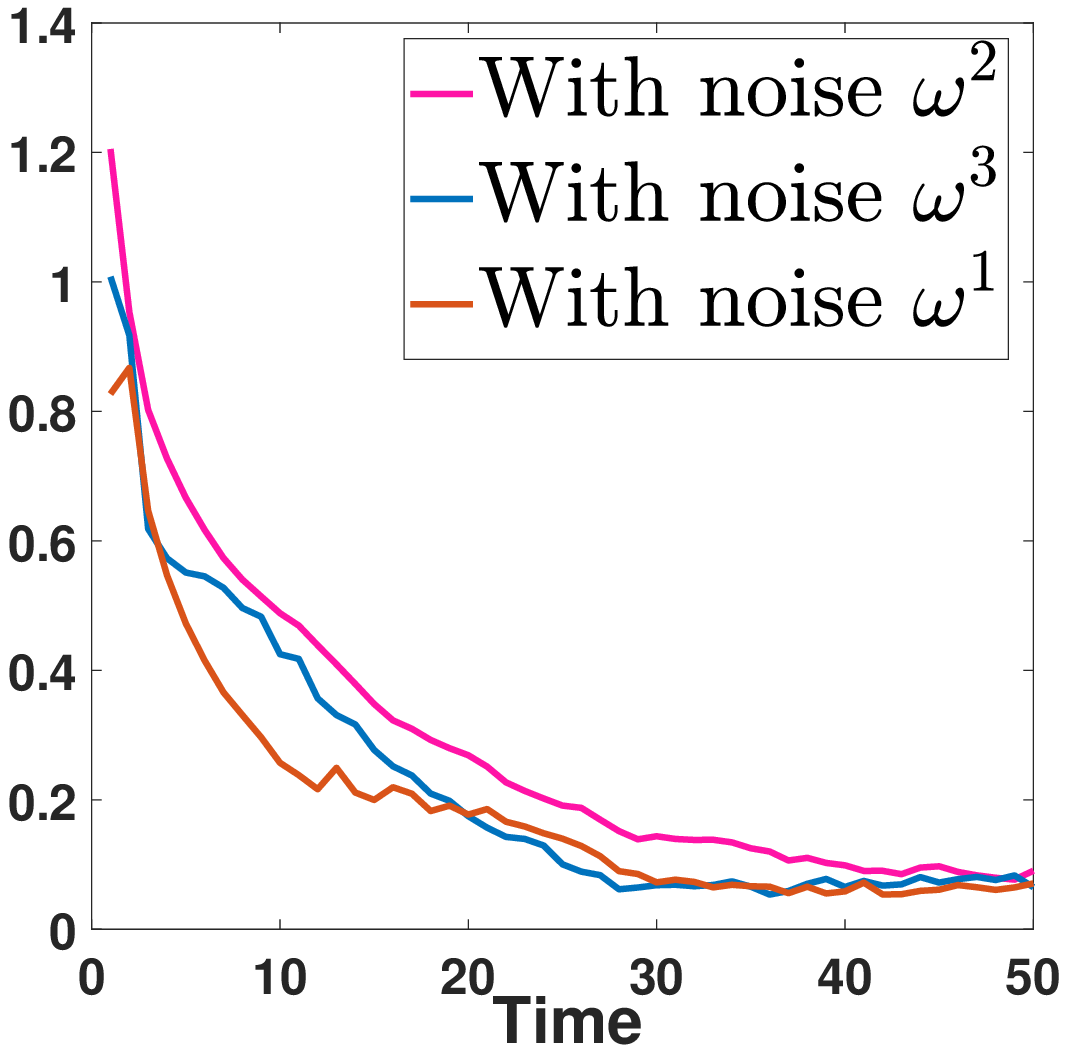}
\end{minipage}%
\caption{[Test Case 2] RMSEs of state estimation by the United Filter for all three types of noise. }
\label{RMSE_General_WithNoise}
\vspace{-0.4cm}
\end{figure}

We first focus on the approximate pressure in the state estimation. The 2D reference and approximate pressure fields, along with their corresponding contour maps, are presented in Figure~\ref{HeatPresField_WithNoise} and Figure~\ref{ContourMap_UnitedF_WithNoise}. Additionally, the 1D approximate pressure on the pressure is illustrated in Figure~\ref{Test2_1DPresFrt_General_UnitedF_WithNoise}. It is obvious that the filtering states yield accurate approximations for the pressure components, regardless of the types of noise. Moreover, as shown in Figure~\ref{ContourMap_UnitedF_WithNoise} and Figure~\ref{Test2_1DPresFrt_General_UnitedF_WithNoise}, the approximate states can effectively capture the discontinuous regions in the reference pressure with high precision. It is also the case for the velocity fields, which are displayed in Figure~\ref{VelField_General}. The magnitude of the approximate velocities closely match that of the reference one. Moreover, all four plots in Figure~\ref{VelField_General} exhibit similar behaviors around the middle zone as well as the lower and upper quarters of the fracture . Finally, we present in Figure~\ref{RMSE_General_WithNoise} the RMSEs for state estimation and observe that all three curves consistently decrease over the filtering steps and converge to nearly similar values.
\begin{figure}[h!]
\centering
\begin{minipage}{0.3\textwidth}
\includegraphics[scale=0.21]{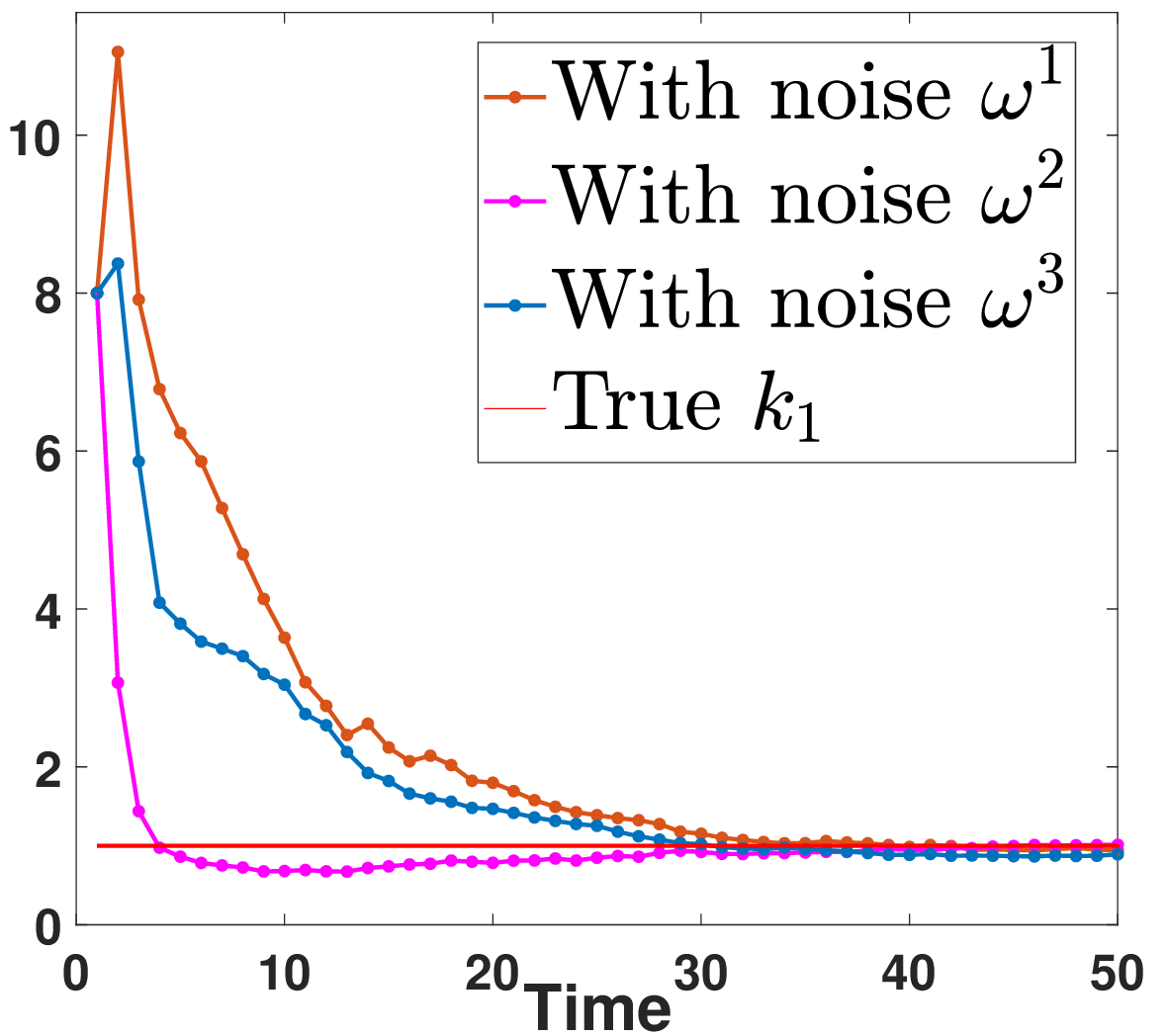}
\end{minipage}%
\begin{minipage}{0.3\textwidth}
\includegraphics[scale=0.21]{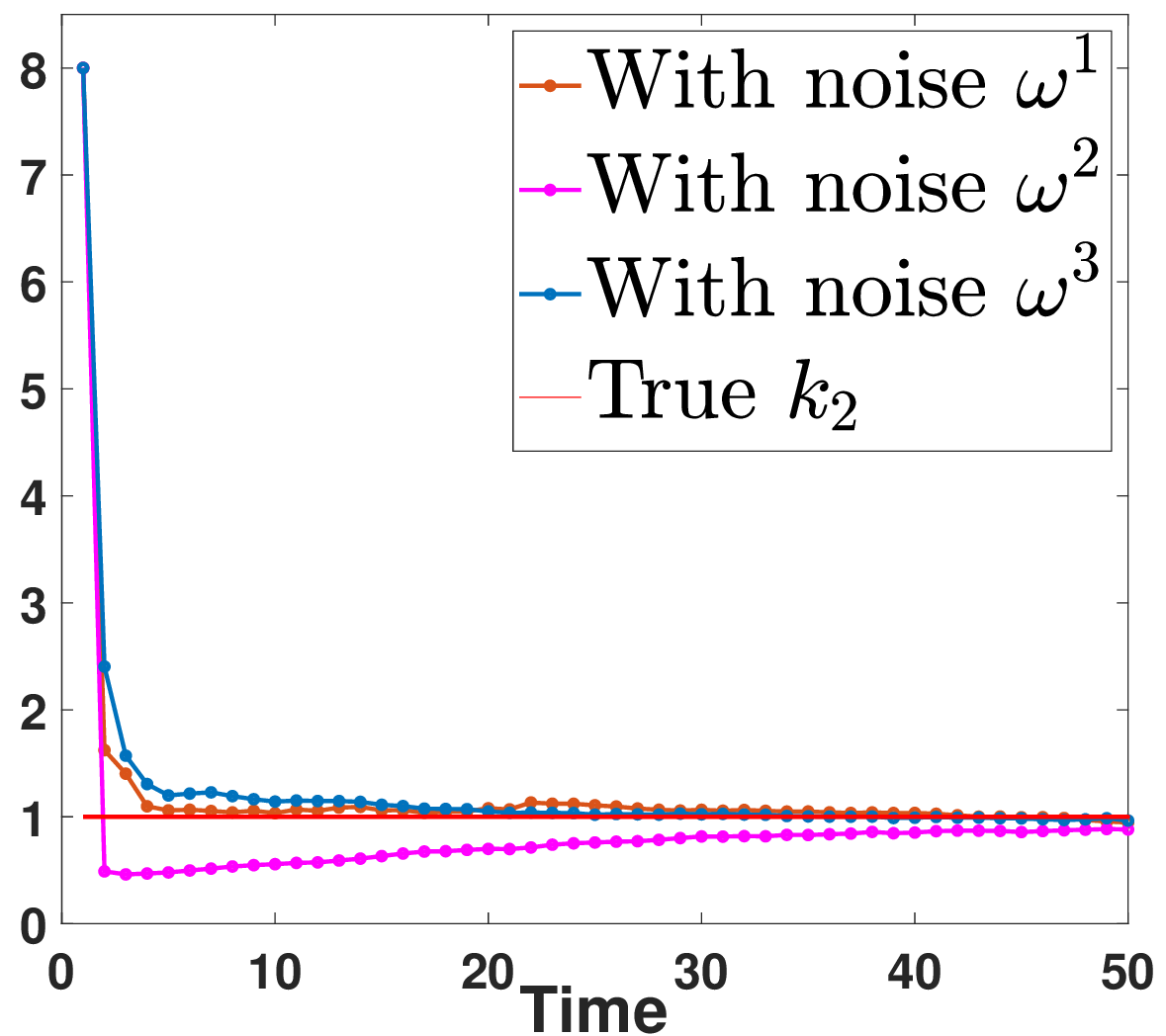}
\end{minipage} %
\begin{minipage}{0.3\textwidth}
\includegraphics[scale=0.21]{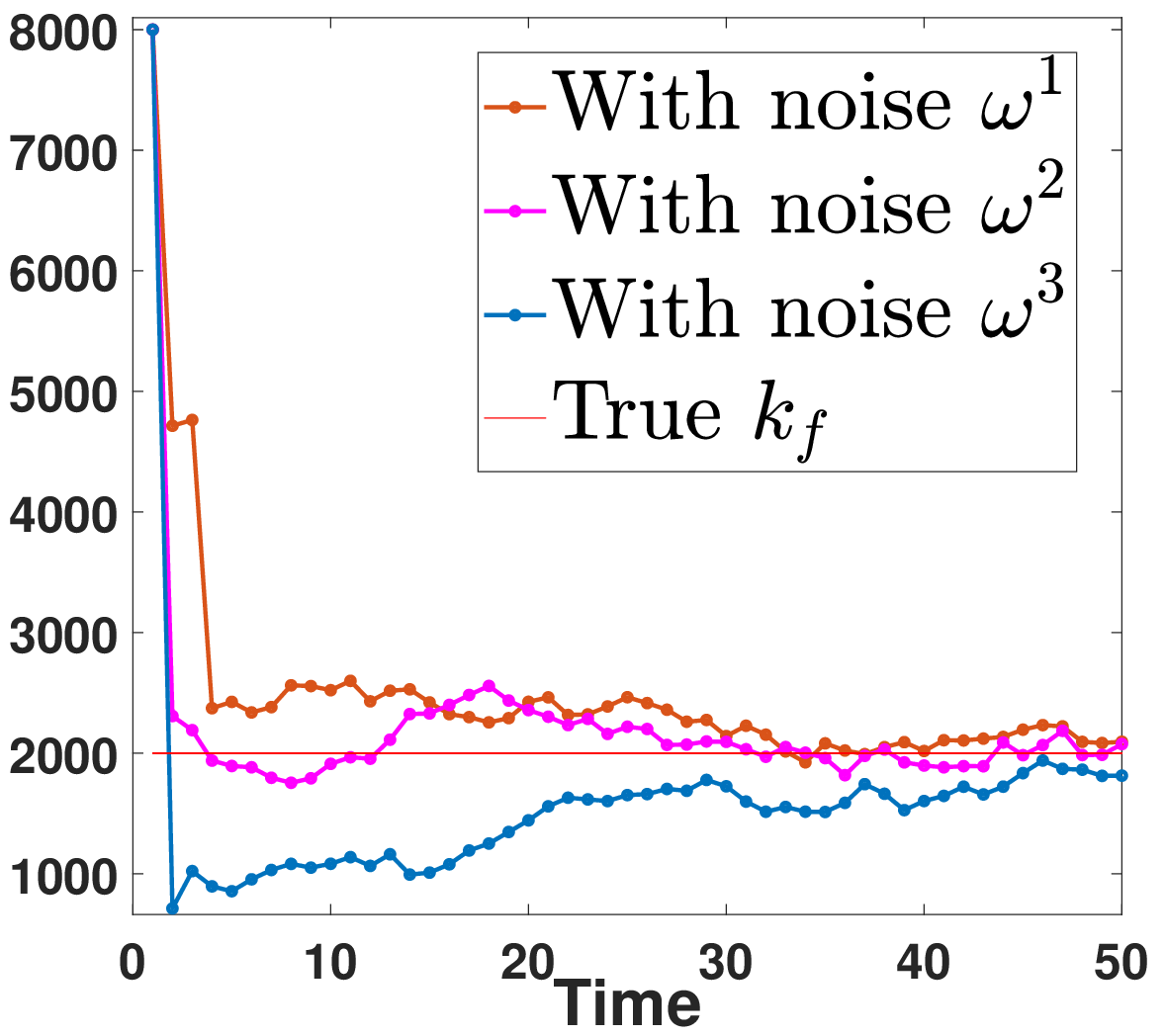}
\end{minipage}%
\caption{[Test Case 2] Parameter estimation by the United Filter. (Left) Estimation for $k_1$. (Middle) Estimation for $k_2$. (Right) Estimation for $k_{f}$.}
\label{ParaEst_UnitedF_General_WithNoise}
\vspace{-0.3cm}
\end{figure}

Finally, we present the parameter estimation results by the United Filter algorithm.  The estimates are plotted in Figure~\ref{ParaEst_UnitedF_General_WithNoise}. For all types of noise, the particles converge to nearly the same values of the true parameters. The parameters $k_1$ and $k_2$ show minimal sensitivity to the noise, as all three curves converge and stabilize relative fast. In contrast, the particles approximating $k_{f}$ show more frequent oscillations and require longer iterations to converge, particularly under the influence of $\omega^3$.
\begin{figure}[h!]
\vspace{-0.2cm}
\centering
\begin{minipage}{0.3\textwidth}
\includegraphics[scale=0.22]{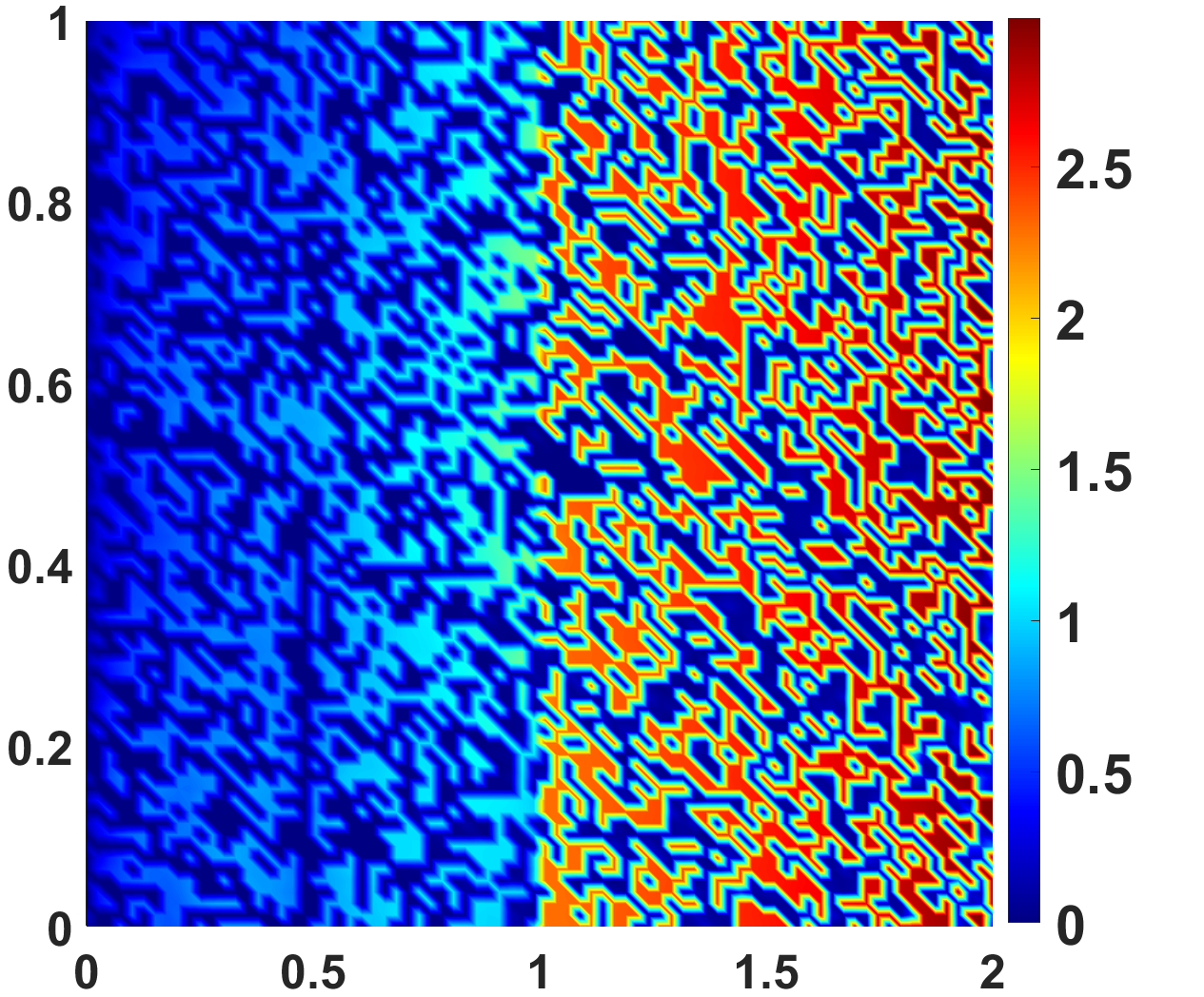}
\end{minipage} %
\hspace{0.1cm}
\begin{minipage}{0.3\textwidth}
\includegraphics[scale=0.22]{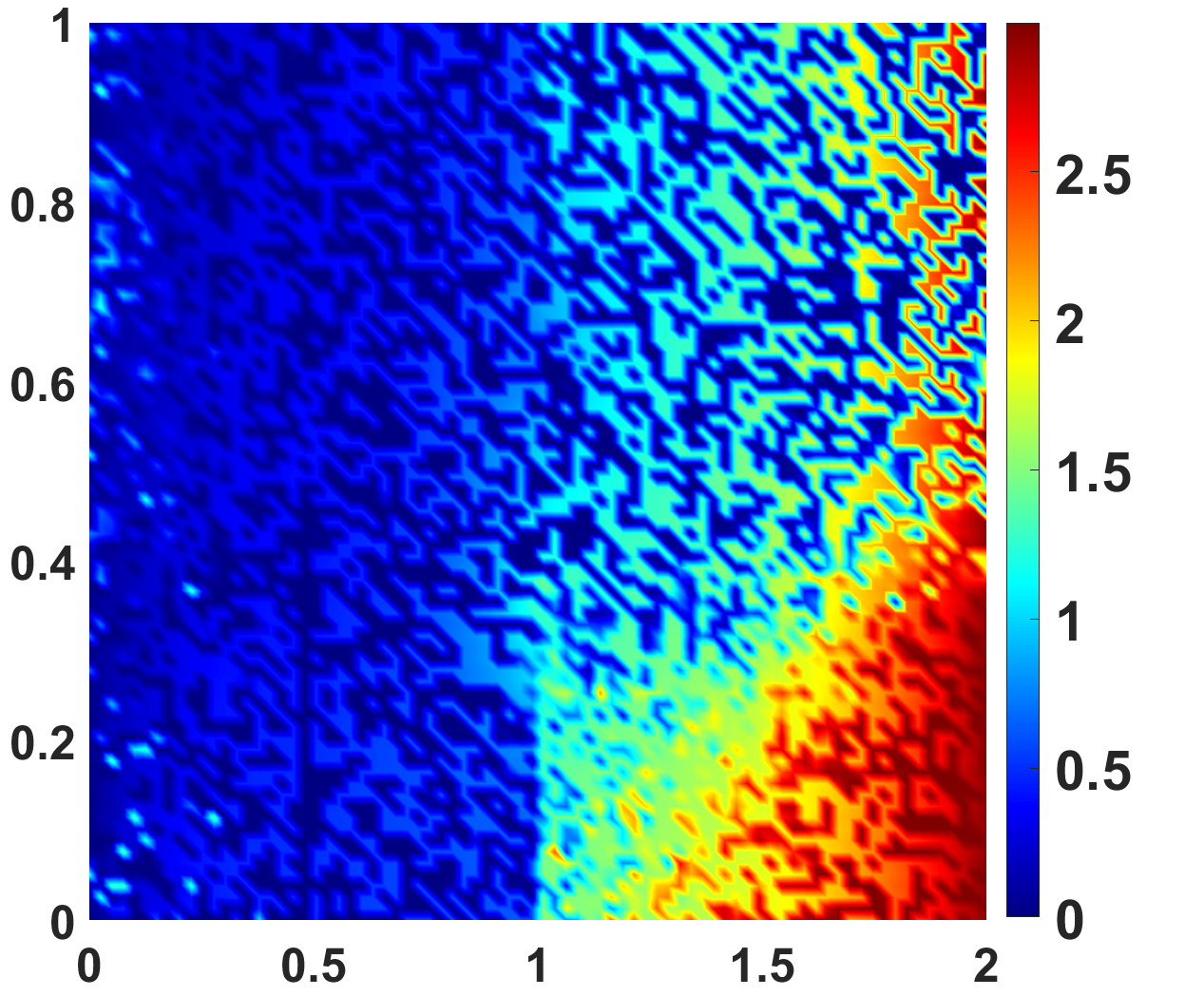}
\end{minipage}
\hspace{0.1cm}
\begin{minipage}{0.3\textwidth}
\includegraphics[scale=0.22]{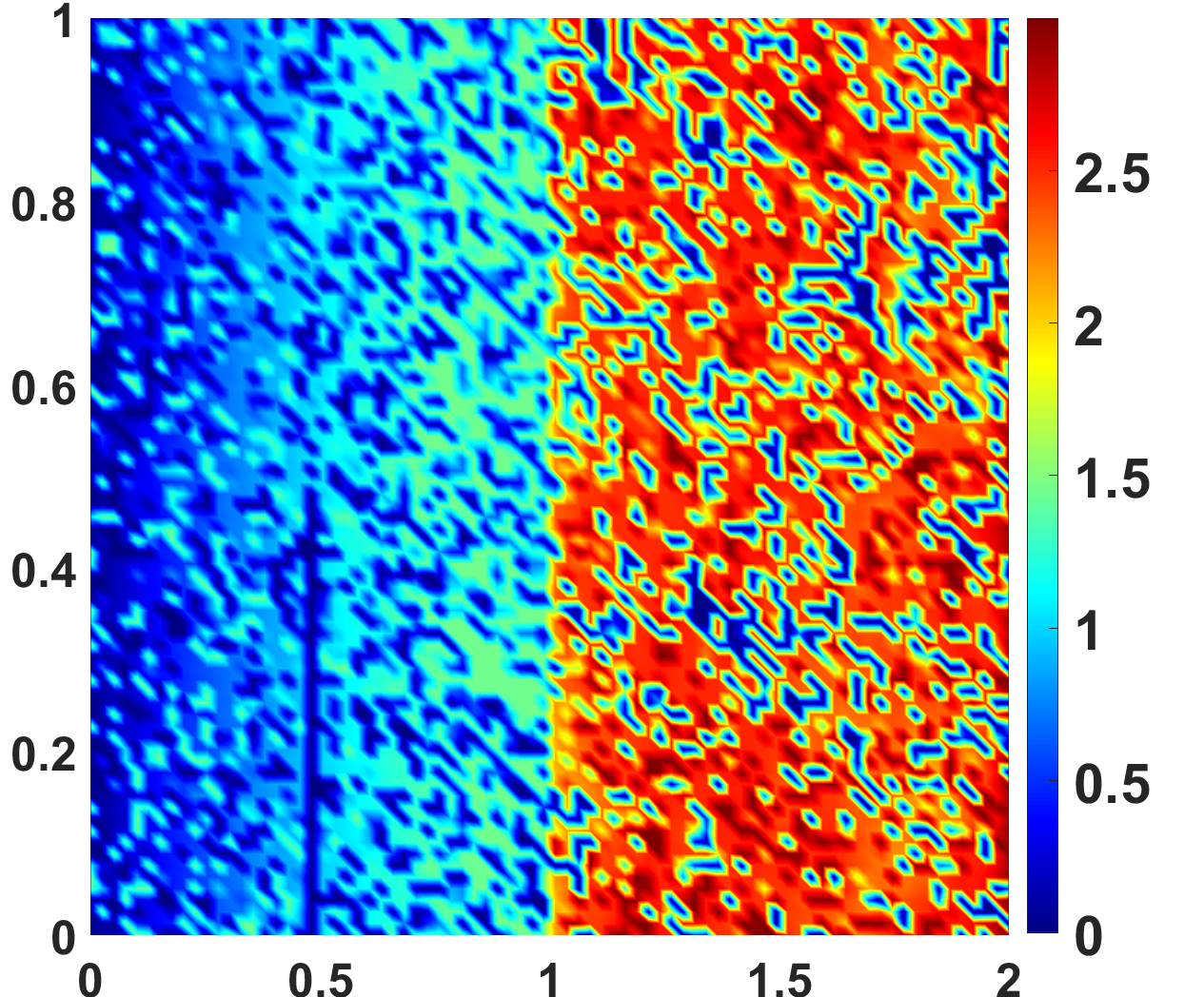}
\end{minipage}
\caption{[Test Case 2] Heat Maps for Pressure Fields estimated by AugEnKF. This is a comparison result to the United Filter estimates presented in Figure \ref{HeatPresField_WithNoise}. (Left) With $\omega^1$. (Middle) With $\omega^2$. (Right) With $\omega^3$.}
\label{Test2HeatPresField_WithNoise_EnKF}
\vspace{-0.2cm}
\end{figure}

In summary, the results for Test Case 2 demonstrate that the United Filter algorithm achieves high accuracy and remains mostly insensitive to the strength of the perturbation. In the last numerical test, we examine a similar scenario to study the impacts of the disturbance in a different setting. 
\begin{figure}[h!]
\vspace{-0.1cm}
\centering
\begin{minipage}{0.3\textwidth}
\includegraphics[scale=0.22]{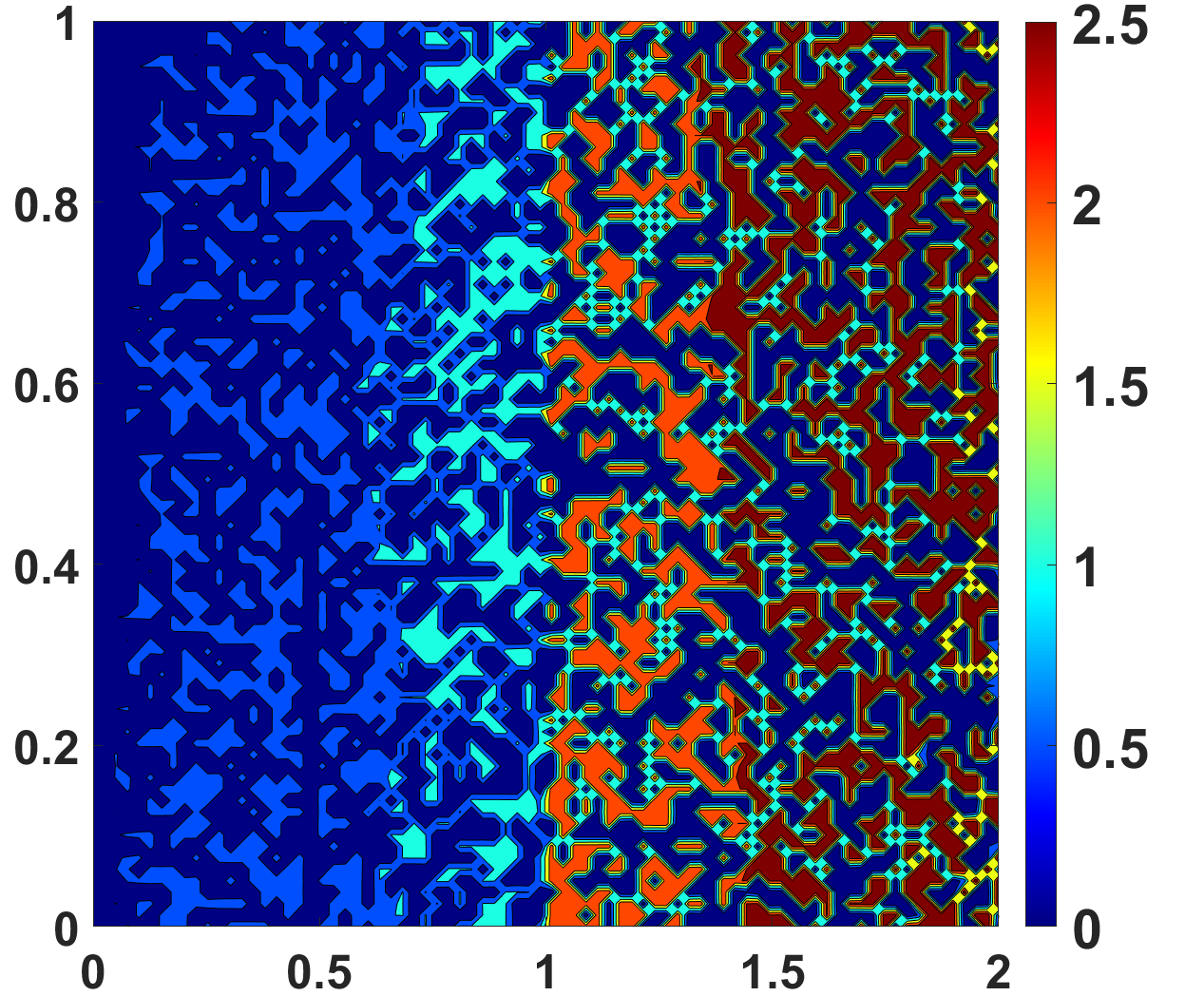}
\end{minipage}%
\hspace{0.1cm}
\begin{minipage}{0.3\textwidth}
\includegraphics[scale=0.22]{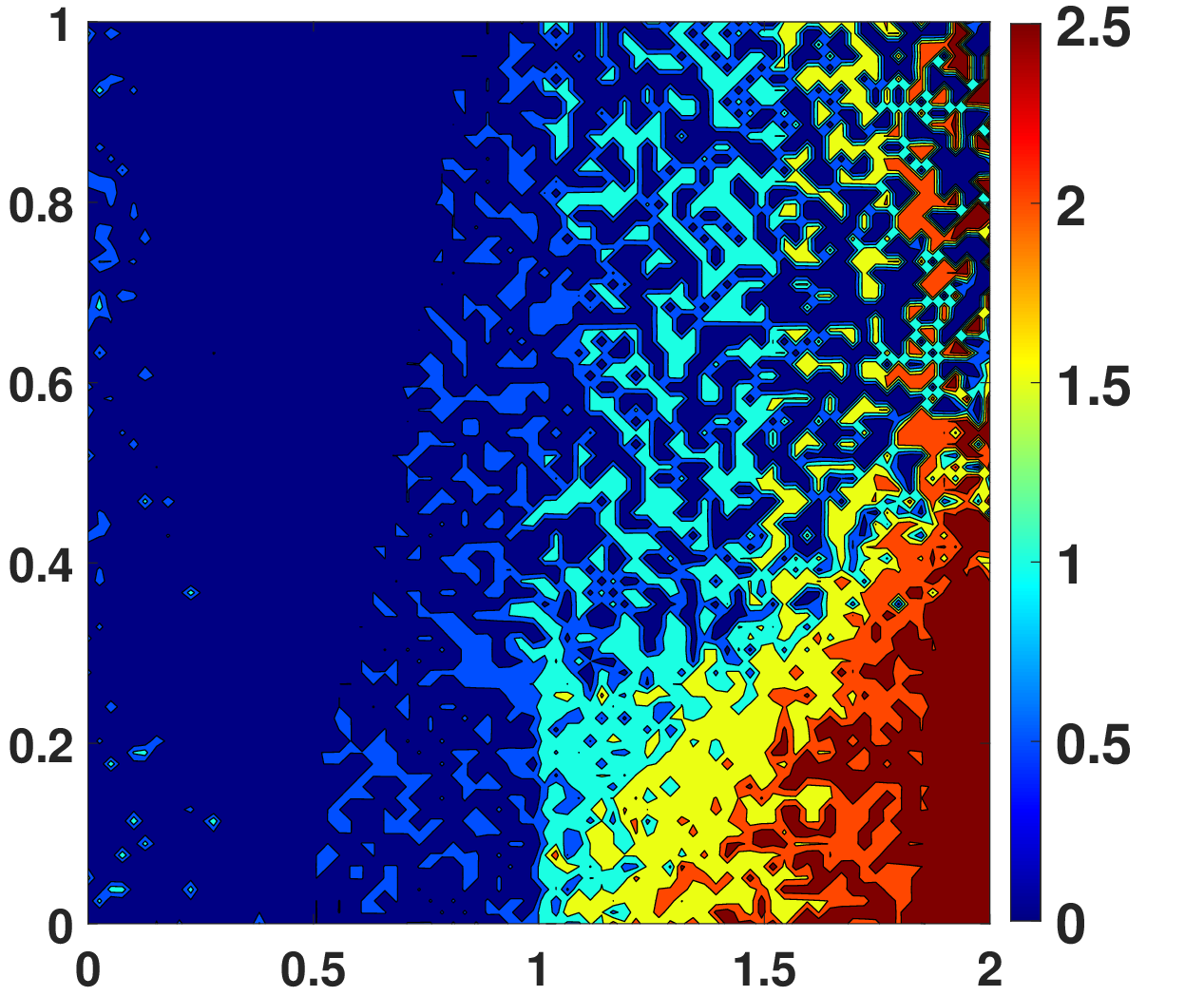}
\end{minipage} 
\hspace{0.1cm}
\begin{minipage}{0.3\textwidth}
\includegraphics[scale=0.22]{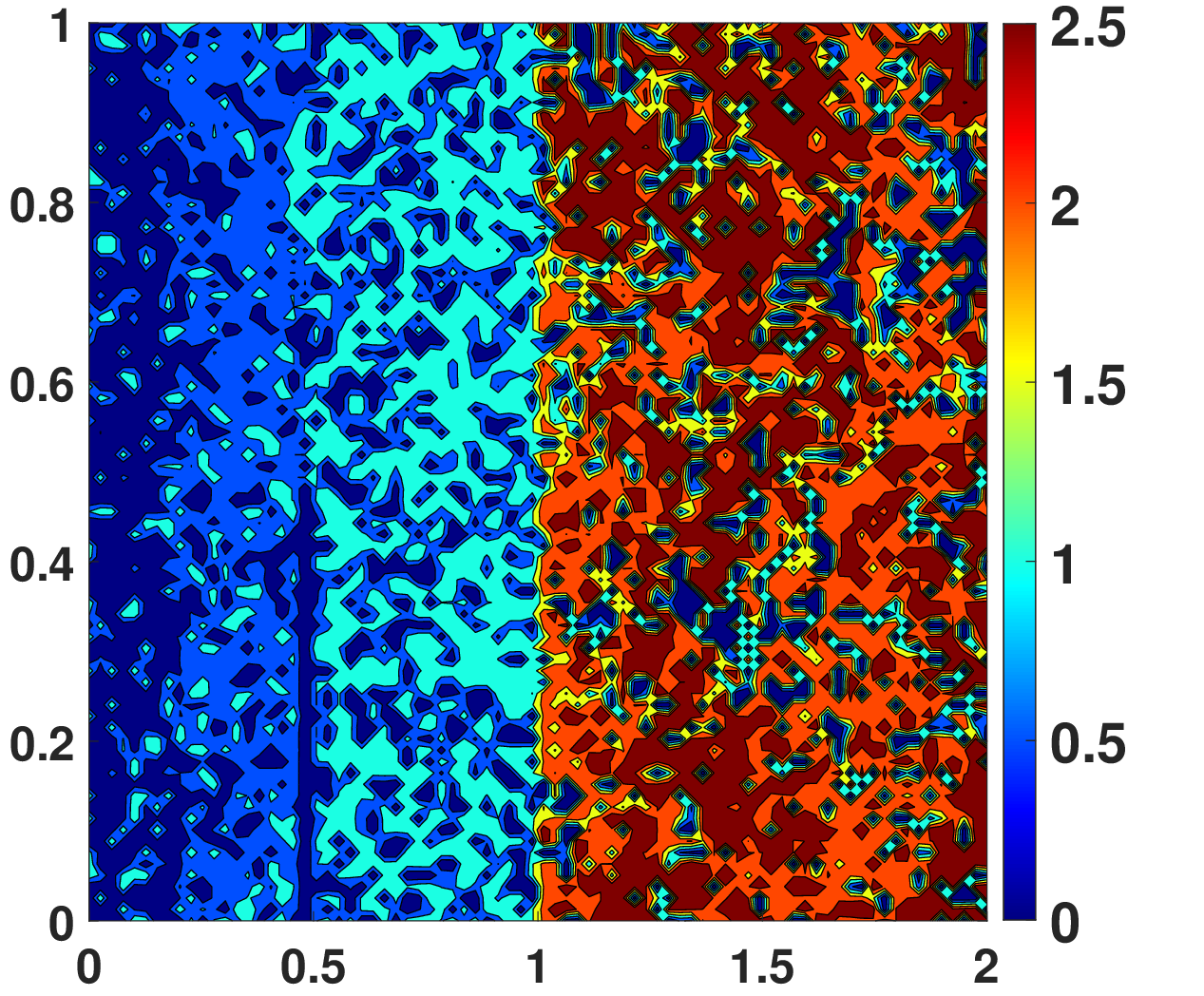}
\end{minipage} 
\caption{[Test Case 2] Pressure Contour Map estimated by AugEnKF. This is a comparison result to the United Filter estimates presented in Figure \ref{ContourMap_UnitedF_WithNoise}. (Left) With $\omega^1$. (Middle) With $\omega^2$. (Right) With $\omega^3$. }
\label{Test2Contour_WithNoise_EnKF}
\vspace{-0.3cm}
\end{figure}
\begin{figure}[h!]
\centering
\begin{minipage}{0.3\textwidth}
\includegraphics[scale=0.22]{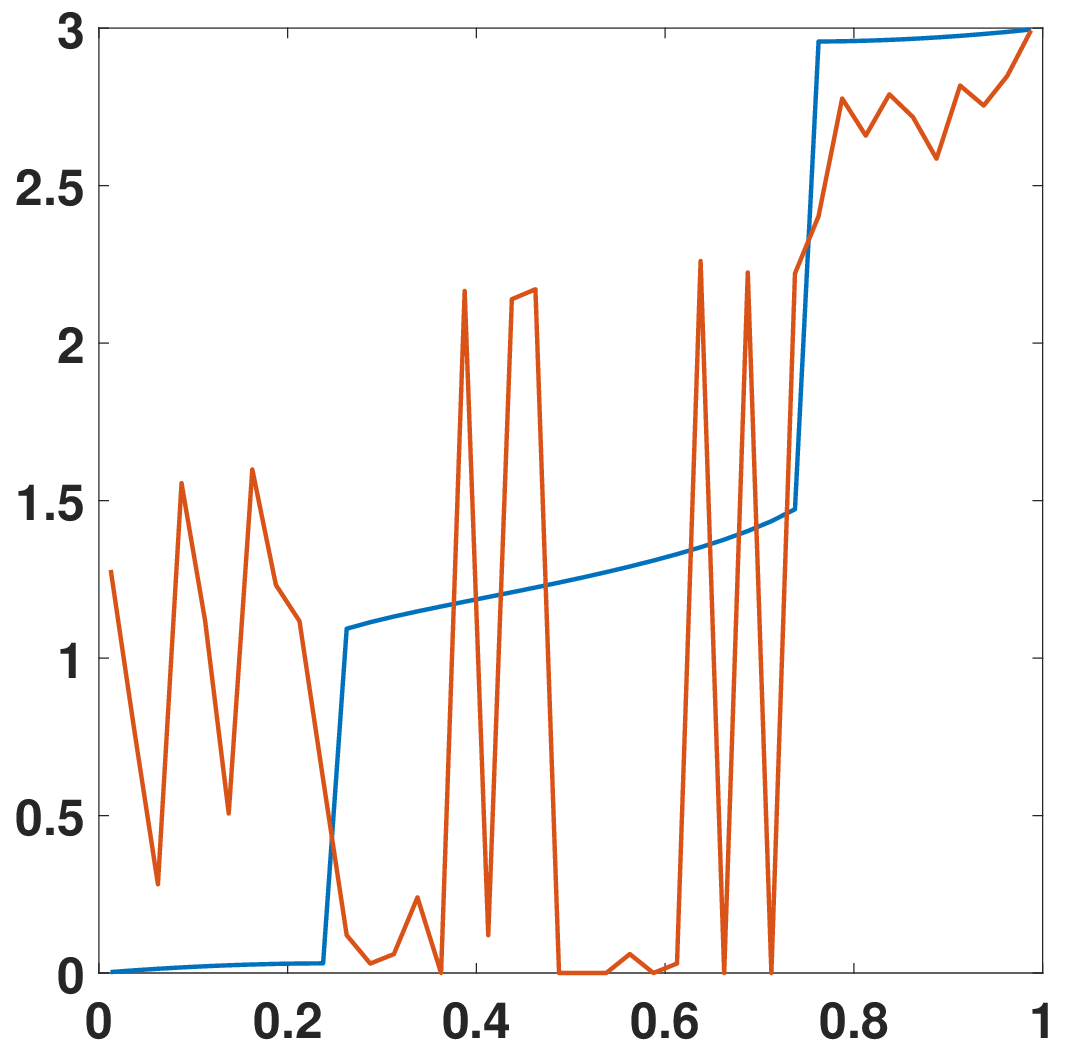}
\end{minipage}%
\hspace{0.1cm}
\begin{minipage}{0.3\textwidth}
\includegraphics[scale=0.22]{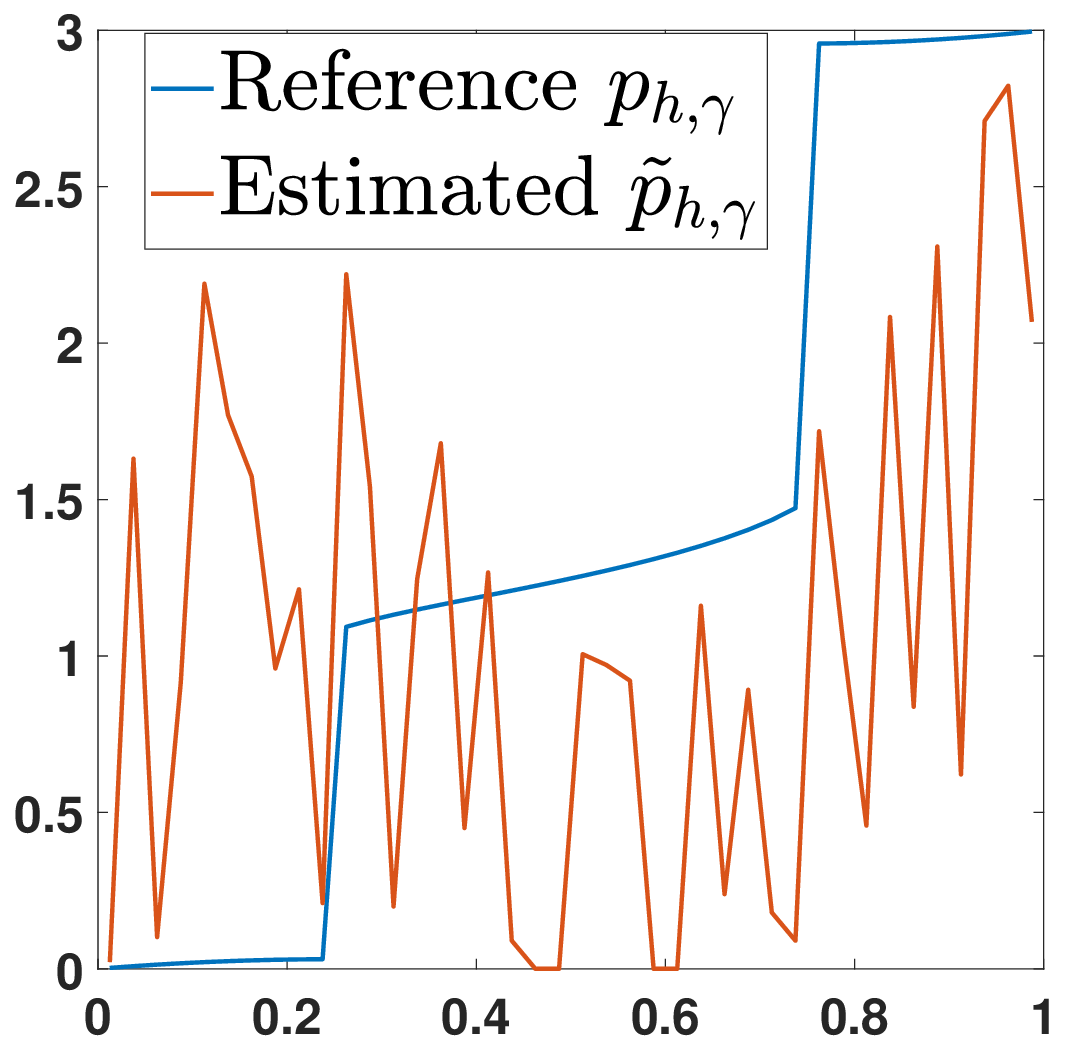}
\end{minipage} 
\hspace{0.1cm}
\begin{minipage}{0.3\textwidth}
\includegraphics[scale=0.22]{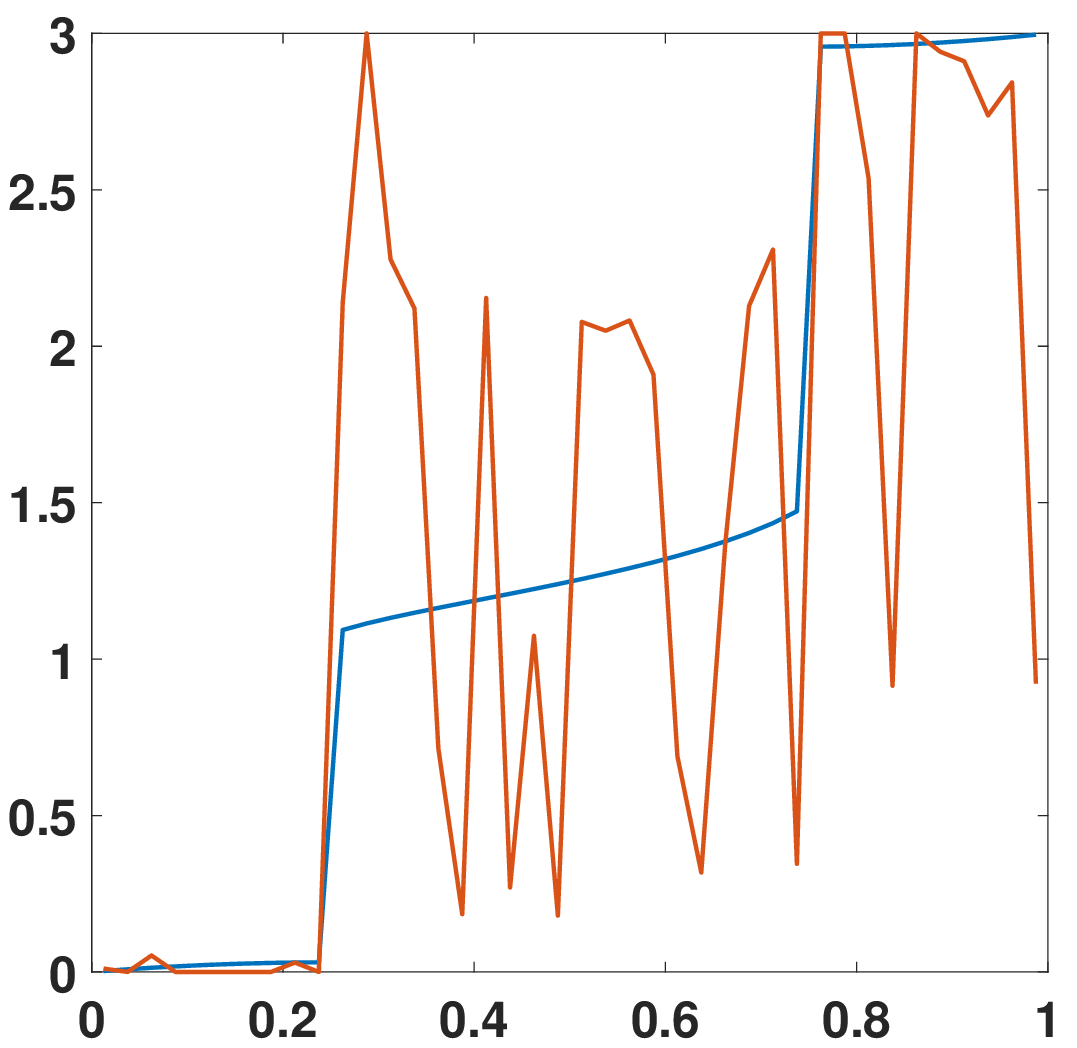}
\end{minipage}
\caption{[Test Case 2] 1D Pressure on the fracture estimated by the AugEnKF.  This is a comparison result to the United Filter's estimates presented in Figure \ref{Test2_1DPresFrt_General_UnitedF_WithNoise}. (Left) With $\omega^1$. (Middle) With $\omega^2$. (Right) With $\omega^3$.}
\label{Test2_1DPresFrt_General_EnKF}
\vspace{-0.3cm}
\end{figure}
\begin{figure}[h!]
\centering
\begin{minipage}{0.3\textwidth}
\includegraphics[scale=0.22]{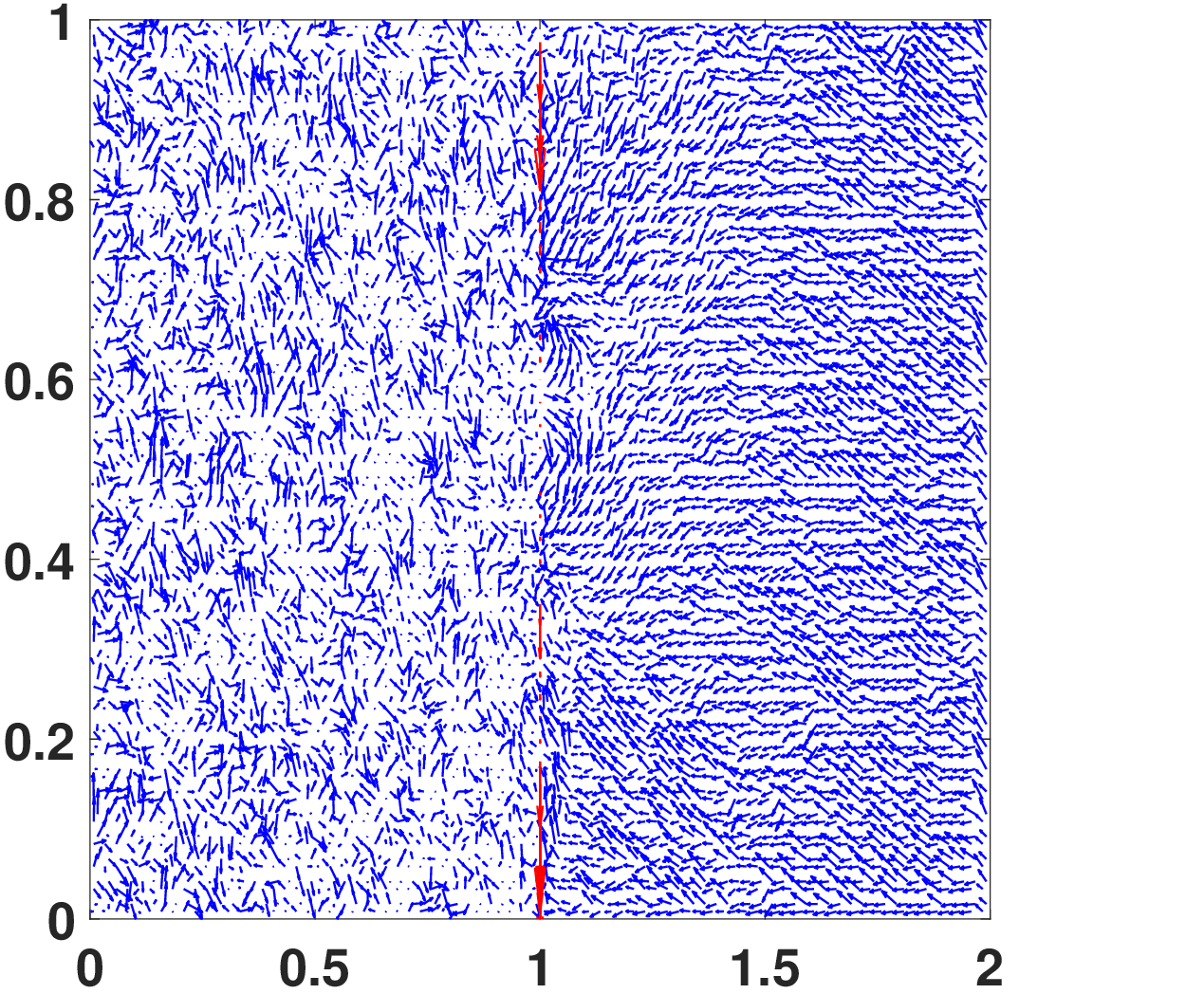}
\end{minipage} %
\hspace{0.1cm}
\begin{minipage}{0.3\textwidth}
\includegraphics[scale=0.22]{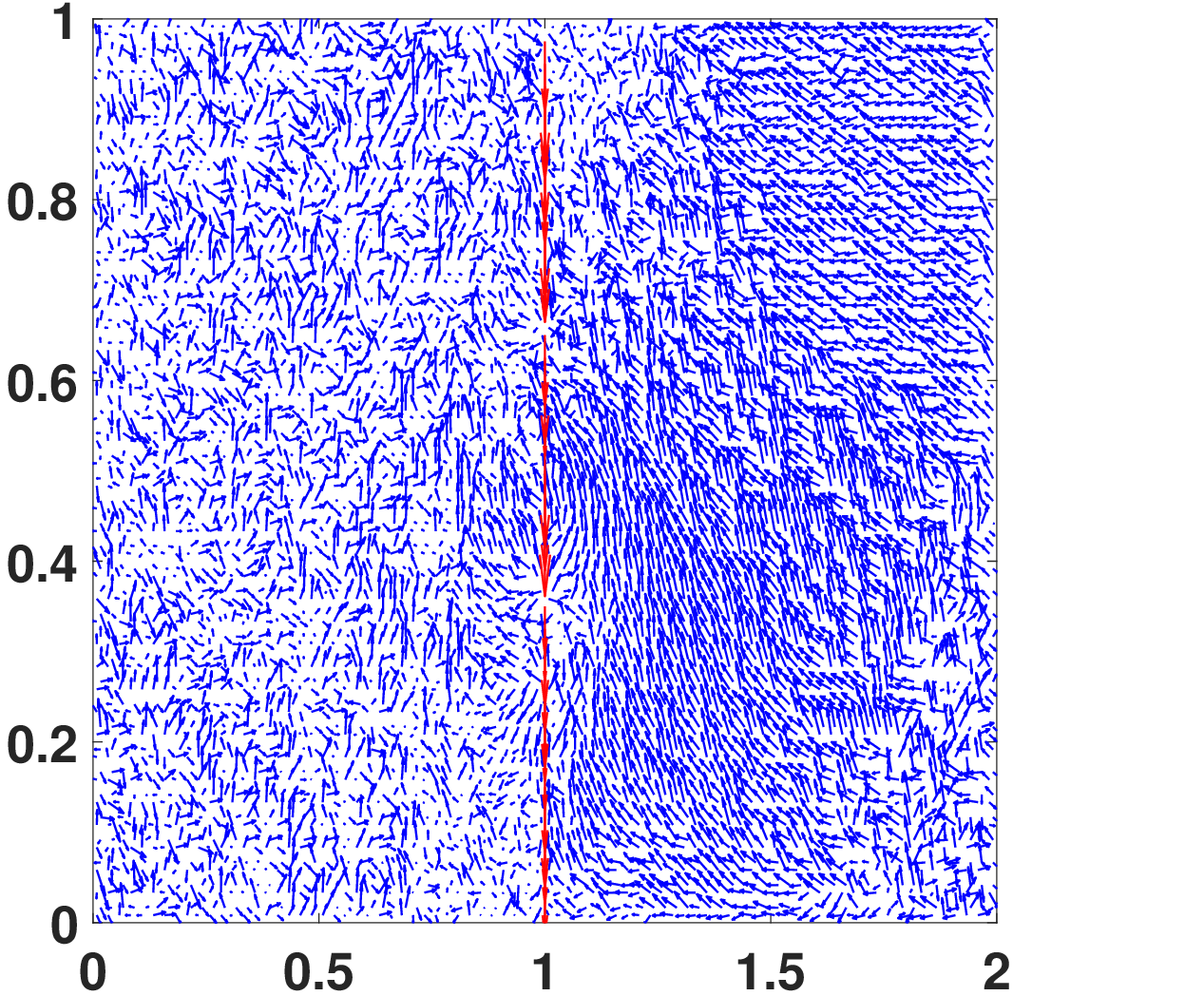}
\end{minipage}%
\hspace{0.1cm}
\begin{minipage}{0.3\textwidth}
\includegraphics[scale=0.22]{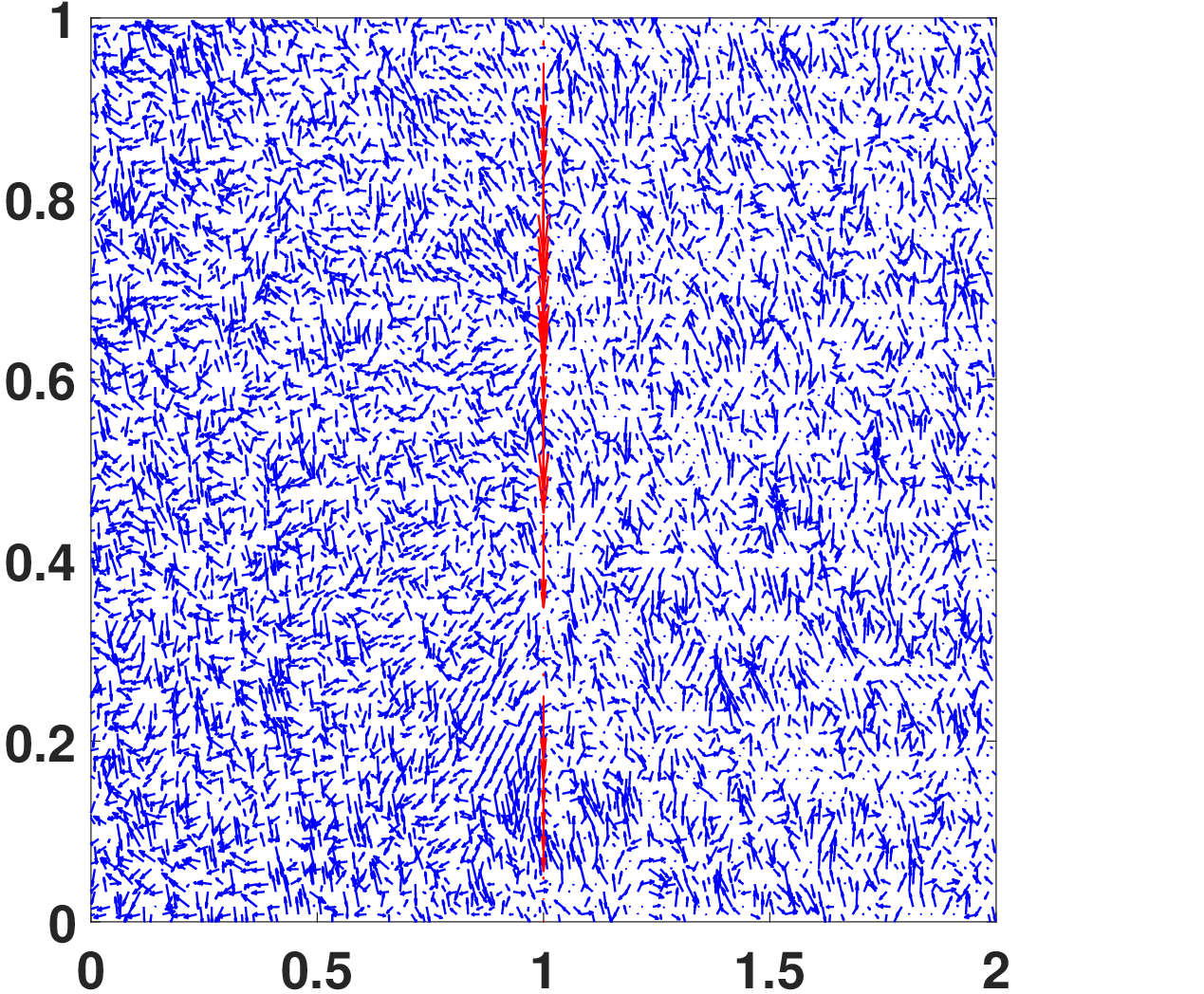}
\end{minipage} 
\caption{[Test Case 2] Velocity Fields estimated by the AugEnKF. This is a comparison result to the United Filter's estimates presented in Figure \ref{VelField_General}. (Left) With $\omega^1$. (Middle) With $\omega^2$. (Right) With $\omega^3$.}
\label{Test2VelField_General_EnKF}
\vspace{-0.3cm}
\end{figure}
\begin{figure}[h!]
\centering
\begin{minipage}{0.3\textwidth}
\includegraphics[scale=0.22]{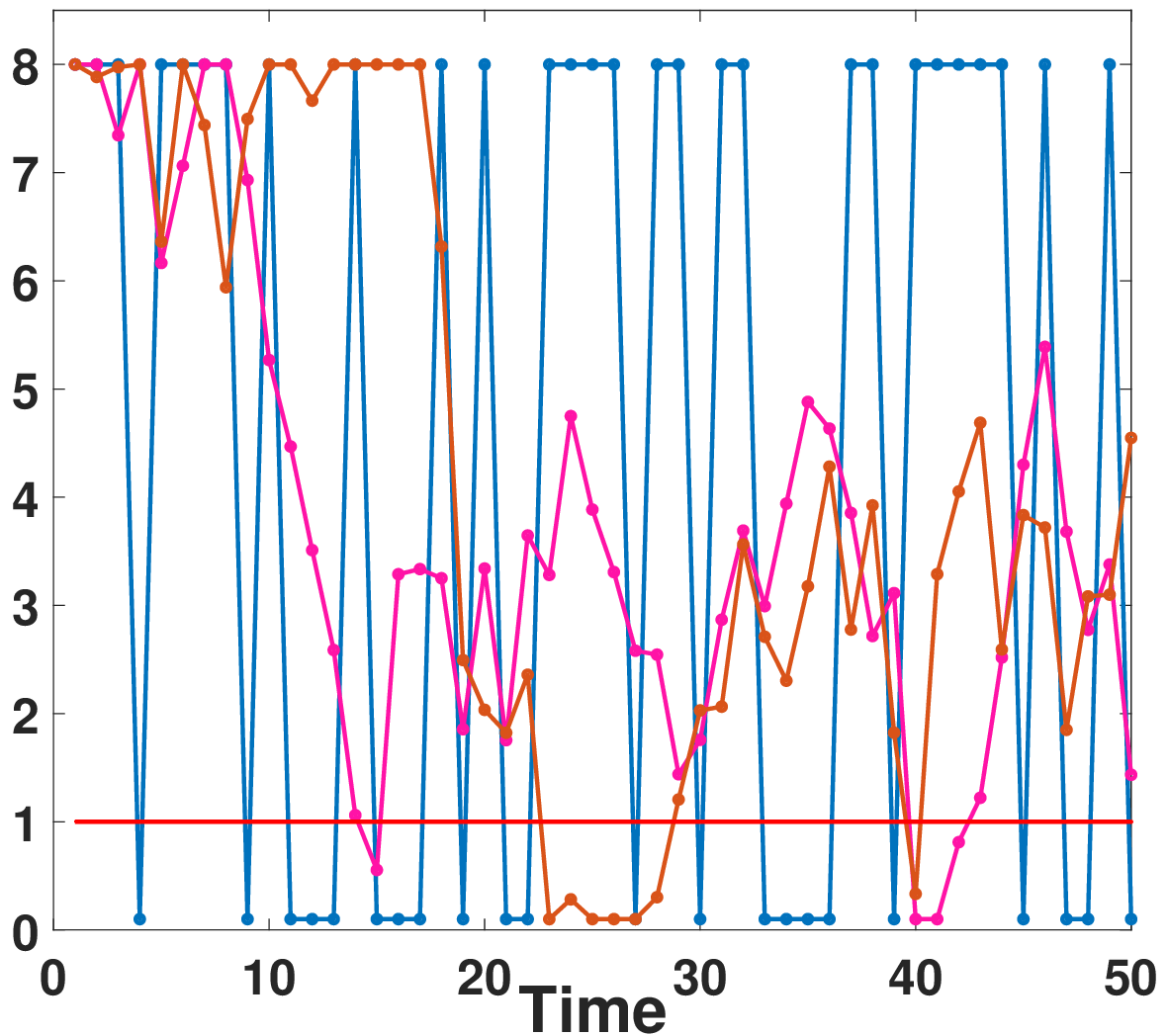}
\end{minipage}%
\hspace{0.1cm}
\begin{minipage}{0.3\textwidth}
\includegraphics[scale=0.22]{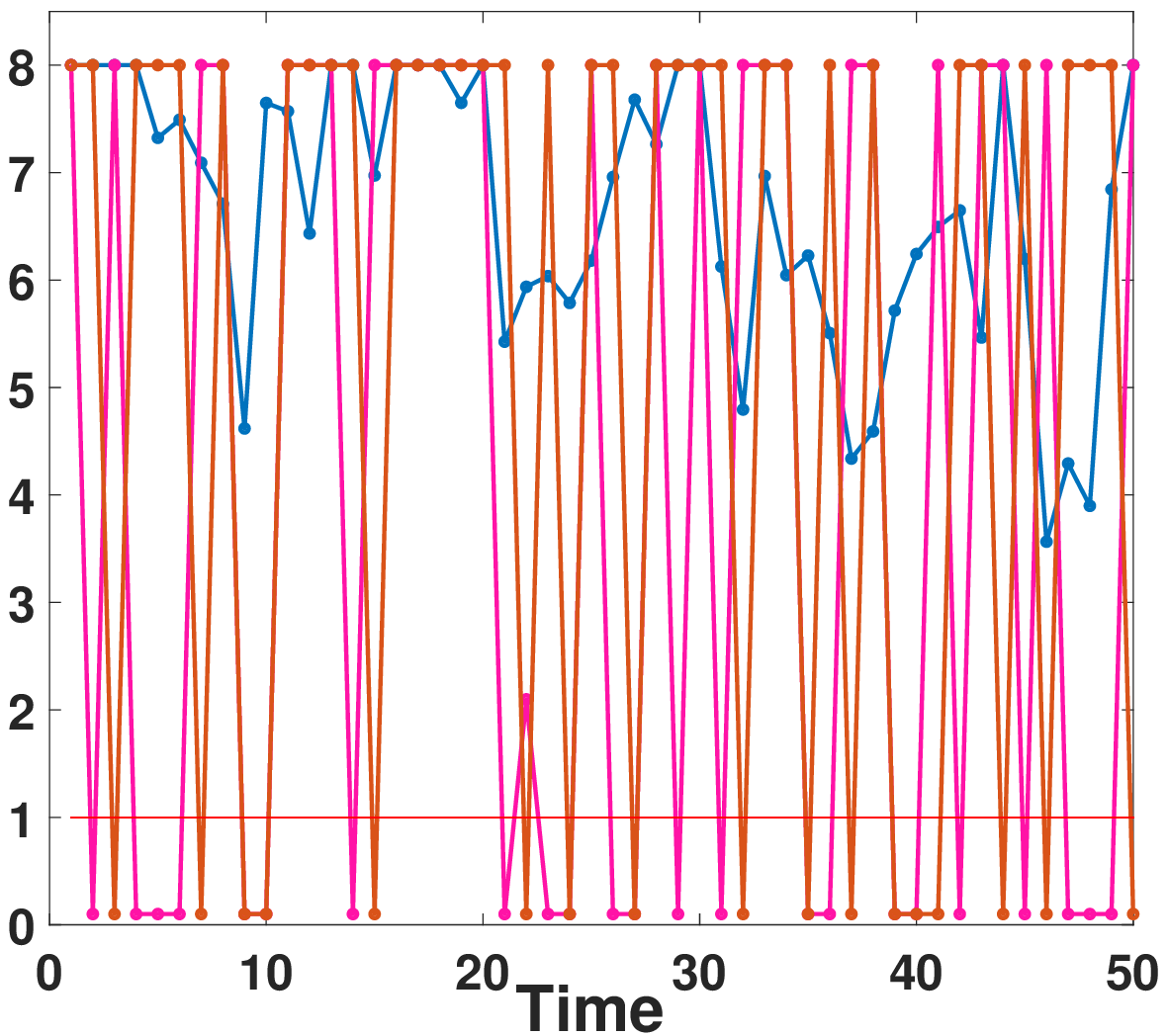}
\end{minipage} %
\hspace{0.1cm}
\begin{minipage}{0.33\textwidth}
\includegraphics[scale=0.22]{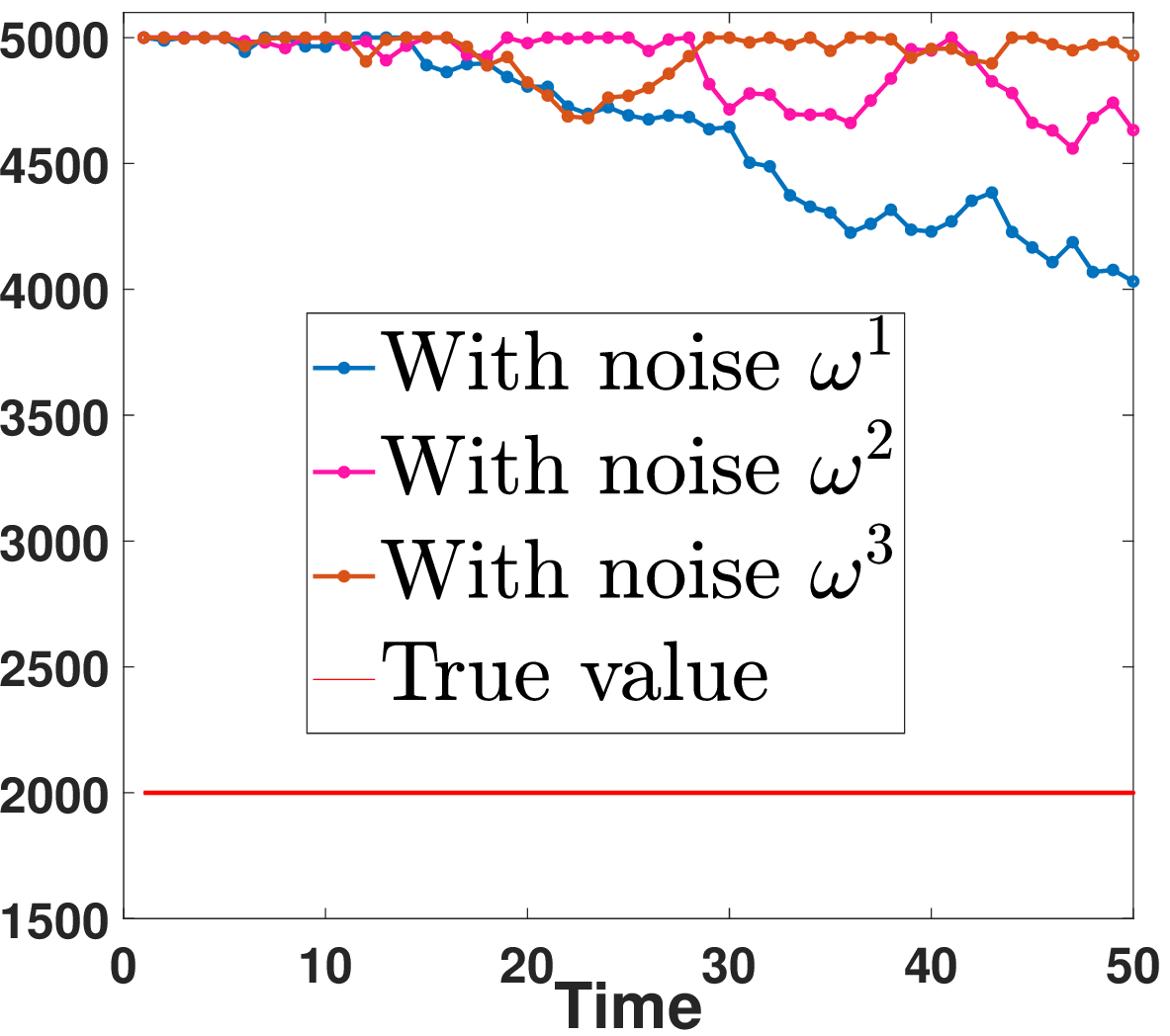}
\end{minipage}%
\caption{[Test Case 2] Parameters estimation by the AugEnK. This is a comparison result to the United Filter's estimates presented in Figure \ref{ParaEst_UnitedF_General_WithNoise}. (Left) Estimation for $k_1$. (Middle) Estimation for $k_2$. (Right) Estimation for $k_{f}$.}
\label{ParaEst_UnitedF_General_WithNoise_EnKF}
\vspace{-0.6cm}
\end{figure}



We now examine the performance of AugEnKF for Test Case 2. Similar to the setting for the United Filter, we consider three types of model uncertainty $\omega^1, \omega^2$ and $\omega^3$. Regarding the state estimation, the 2D and 1D pressure are shown in Figure~\ref{Test2HeatPresField_WithNoise_EnKF}, Figure~\ref{Test2Contour_WithNoise_EnKF} and Figure~\ref{Test2_1DPresFrt_General_EnKF}. By comparing with Figures~\ref{HeatPresField_WithNoise}, \ref{ContourMap_UnitedF_WithNoise}, and \ref{Test2_1DPresFrt_General_UnitedF_WithNoise}, we can see that the results from AugEnKF are inconsistent across different noise levels. Moreover, the estimated 2D pressures from AugEnKF fail to capture any region of the reference pressure field and are entirely incorrect. The estimated 1D pressures on fracture also fail to recover the discontinuities in the reference state. Similarly, from Figure~\ref{Test2VelField_General_EnKF}, it is clear that the approximate velocity fields from AugEnKF under all three noise levels deviate significantly from the reference velocity field. More precisely, the magnitude of all estimated 2D velocities is larger compared to that of the reference field. Moreover, since the middle zone of the fracture acts as a barrier in this case, the velocity along this section of the fracture should be negligible. However, the approximate velocity in this region is noticeably higher, which is physically inaccurate. 

Finally, we present in Figure~\ref{ParaEst_UnitedF_General_WithNoise_EnKF} the results of the parameter estimation by the AugEnKF algorithm. Similar to Test Case 1, AugEnKF failed to converge to the true parameters. Moreover, we have a similar behavior as in Figure~\ref{Test1_ParaEst_Pure_EnKF} where the particles only take two values due to the constraint enforced on the particles to avoid divergence. The above numerical results show that the United Filter outperformed the AugEnKF in Test Case 2.

\subsection{Test Case 3: Advection-diffusion equation}
Finally, we consider the following reduced model for the advection-diffusion equation (see, e.g., \cite{Alboin2002}) which consists of equations in the subdomains, \vspace{-0.2cm}
\begin{equation}
\label{reduced_adv_diff_subdom}
\left.\begin{array}{rcll}
\phi_i\partial_t{c_i}+\text{div }\pmb{\varphi}_{i}&=&q_{i} &\text{ in } \Omega_i\times (0, T), \vspace{0.05cm}\\
\pmb{\varphi}_i&=&\textbf{\textit{u}}_ic_i-\textbf{\textit{D}}_i\nabla{c_i} &\text{ in } \Omega_i\times (0, T), \vspace{0.05cm}\\
c_i&=&0 &\text{ on } \left(\partial\Omega_i \cap \partial\Omega\right) \times (0, T), \vspace{0.05cm}\\
c_i&=&c_{\gamma} &\text{ on } \gamma \times (0, T), \vspace{0.05cm}\\
c_i(\cdot, 0)&=&c_{0, i} &\text{ in } \Omega_i, 
\end{array}\right. 
\end{equation}
for $i=1, 2, $ coupled with the following equation in the one-dimensional fracture,
\vspace{-0.2cm}
\begin{equation}
\label{reduced_adv_diff_fracture}
\begin{array}{rcll}
\phi_{\gamma}\partial_t{c_{\gamma}}+\text{div}_{\tau}\pmb{\varphi}_{\gamma}&=& q_{\gamma} +\sum\limits^{2}_{i=1}\pmb{\varphi}_i\cdot\pmb{n}_{i\vert\gamma} & \text{ in } \gamma \times (0, T), \vspace{0.1cm} \\
\pmb{\varphi}_{\gamma} &=&\pmb{u}_{\gamma}c_{\gamma}-\textbf{\textit{D}}_{f, \tau}\delta\nabla_{\tau}c_{\gamma} & \text{ in } \gamma \times (0, T), \\
c_{\gamma}&=&0 &\text{ on } \partial\gamma \times (0, T), \\
c_{\gamma}(\cdot, 0)&=&c_{0, \gamma} & \text{ in } \gamma. \vspace{-0.1cm}
\end{array}
\end{equation}
where $c_i, i=1, 2, \gamma$ are the concentration on $\Omega_i, i=1, 2$ and $\gamma$, $q_i, i=1, 2, \gamma$ are source terms, $\phi_i, i=1, 2, \gamma$ are the porosity, $\textbf{\textit{D}}_i, i=1, 2$ and $\textbf{\textit{D}}_{f, \tau}$ are symmetric time-independent diffusion tensors, and $\textbf{\textit{u}}$ is the Darcy velocity given by solving the following coupled steady-state flow problem in the subdomains,
\begin{equation}
\label{Darcy_subdom}
\begin{array}{rll}
\text{div }\bu_i &= 0, &\text{in } \Omega_i \times (0, T), \vspace{0.1cm} \\
\bu_i &= -\bK_i\nabla{p}_i, &\text{in } \Omega_i \times (0, T), \vspace{0.1cm} \\
p_i &= g, &\text{on } \left(\partial{\Omega_i} \cap \partial{\Omega}\right) \times (0, T), \vspace{0.1cm} \\
p_i &= p_{\gamma}, &\text{on } \gamma \times (0, T), 
\end{array} i =1, 2,
\end{equation}
and on the fracture,
\begin{equation}
\label{Darcy_fracture}
\begin{array}{rll}
\text{div}_{\tau} \bu_{\gamma} &= \sum\limits^2_{i=1}\bu_i\cdot\bn_{i\vert{\gamma}}, &\text{in } \gamma \times (0, T), \vspace{0.1cm} \\
\bu_{\gamma} &= -\bK_{f, \tau}\delta\nabla_{\tau}p_{\gamma}, &\text{in } \gamma \times (0, T), \vspace{0.1cm} \\
p_{\gamma} &= g_{\gamma}, &\text{on } \partial{\gamma} \times (0, T).
\end{array}
\end{equation}

We assume further that $\bD_i = d_i\pmb{I}, \; i=1, 2$, $\bD_{f, \tau} = d_f\pmb{I}$ and denote $\alpha_{\gamma} := d_{f}\delta$. The problems \eqref{reduced_adv_diff_subdom}-\eqref{reduced_adv_diff_fracture} and \eqref{Darcy_subdom}-\eqref{Darcy_fracture} are solved numerically using the mixed hybrid finite element methods on the rectangular mesh (\cite{Brezzi1991}). This class of methods has been shown to be efficient for the case with strongly advection-dominated problems (\cite{Brunner2014}).

\begin{figure}[h!]
\vspace{-0.1cm}
\centering
\includegraphics[width=0.4\textwidth]{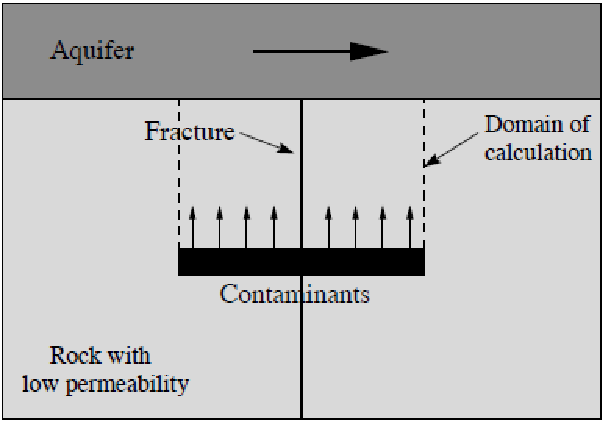}
\caption{{ [Test Case 3] A contaminant storage crossed by a fracture}
\label{Test_case1_image} }
\vspace{-0.2cm}
\end{figure}

We consider an adapted version of the test case used in \cite{Alboin2002} which describes a leaking  contaminant repository, located in a rock with low permeability (Figure~\ref{Test_case1_image}). The repository is crossed by a fracture and transported mostly upward. The rock is covered by an aquifer and the contaminant is assumed to be moved away instantly at the top boundary of the domain calculation so the boundary condition there is a vanishing concentration. We rescale the physical parameters from \cite{Alboin2002} and obtain the values shown in Table~\ref{physical_parametes_Testcase1}. Boundary conditions are as follows: for the velocity, we assume that there is no horizontal flow on the lateral sides of the domain while a pressure drop constant in time is given between the top and bottom boundaries.  At the top, the pressure is constant in space while at the bottom it is increasing slightly from the fracture toward the lateral sides.  For the concentration, it is given, constant, at the top and bottom boundaries, vanishing at the top. On the lateral sides we assume that there is no exchange with the outside.

\begin{table}[h!]
 \centering
\begin{tabular*}{\columnwidth}{@{\extracolsep{\fill}}lcc}
\hline
\textit{Parameters} &\textit{Subdomains} & \textit{Fracture} \\ \hline
Hydraulic conductivity  & 9.92$\times 10^{-6}$ & $3.15^{-5}$ \\ \hline
Molecular diffusion  & $3.15^{-4}$ & $9.92\times 10^{-3}$ \\ \hline
Porosity & $0.05$ & $0.1$ \\ \hline
Subdomains dimensions  & $1 \times 1$ & - \\ \hline
Fracture width  & - & $0.1$ \\ \hline
\end{tabular*}
\caption{[Test Case 3] Physical parameters for the experiment shown in Figure~\ref{Test_case1_image}\label{physical_parametes_Testcase1} }.
\vspace{-0.5cm}
\end{table}
\begin{figure}[h!]
\centering
\begin{minipage}{0.35\textwidth}
\includegraphics[scale=0.23]{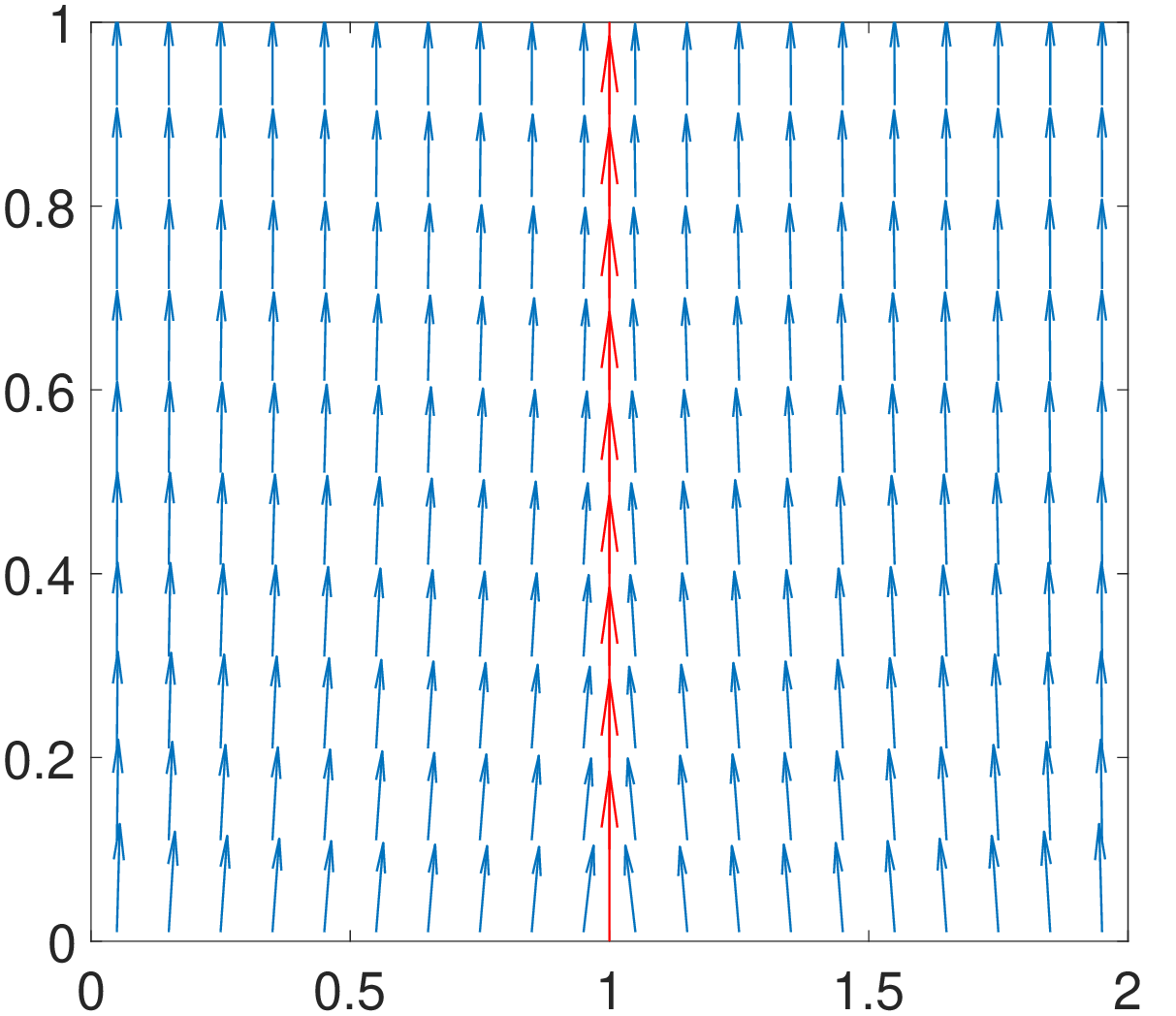}
\end{minipage}%
\begin{minipage}{0.35\textwidth}
\hspace{0.25cm}\includegraphics[scale=0.23]{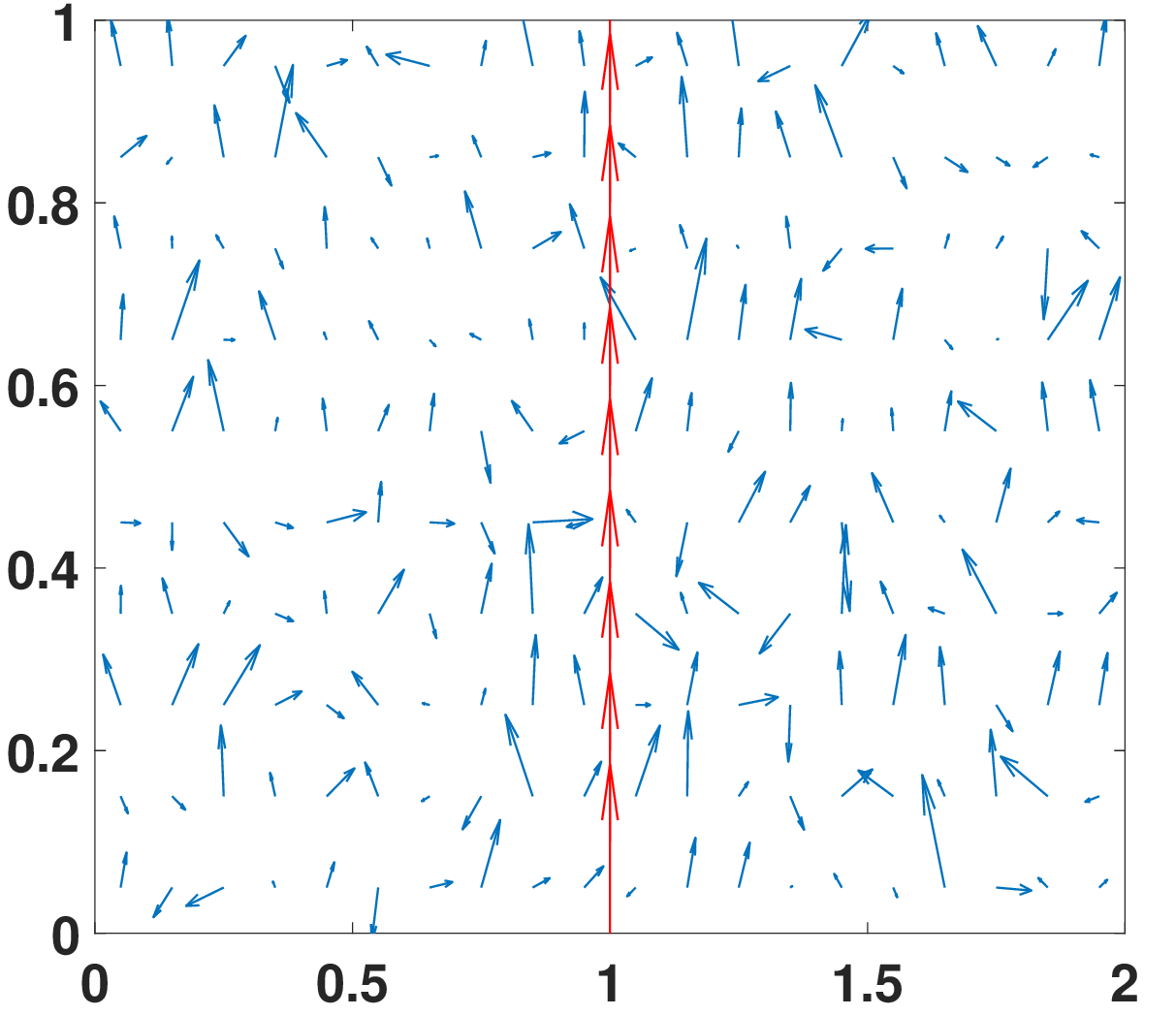}
\end{minipage}
\caption{[Test Case 3] (Left) Reference Darcy Velocity Field. (Right) Perturbed Velocity Field.}
\label{DarcyField}
\vspace{-0.3cm}
\end{figure} 
\begin{figure}[h!]
\centering
\begin{minipage}{0.35\textwidth}
\includegraphics[scale=0.24]{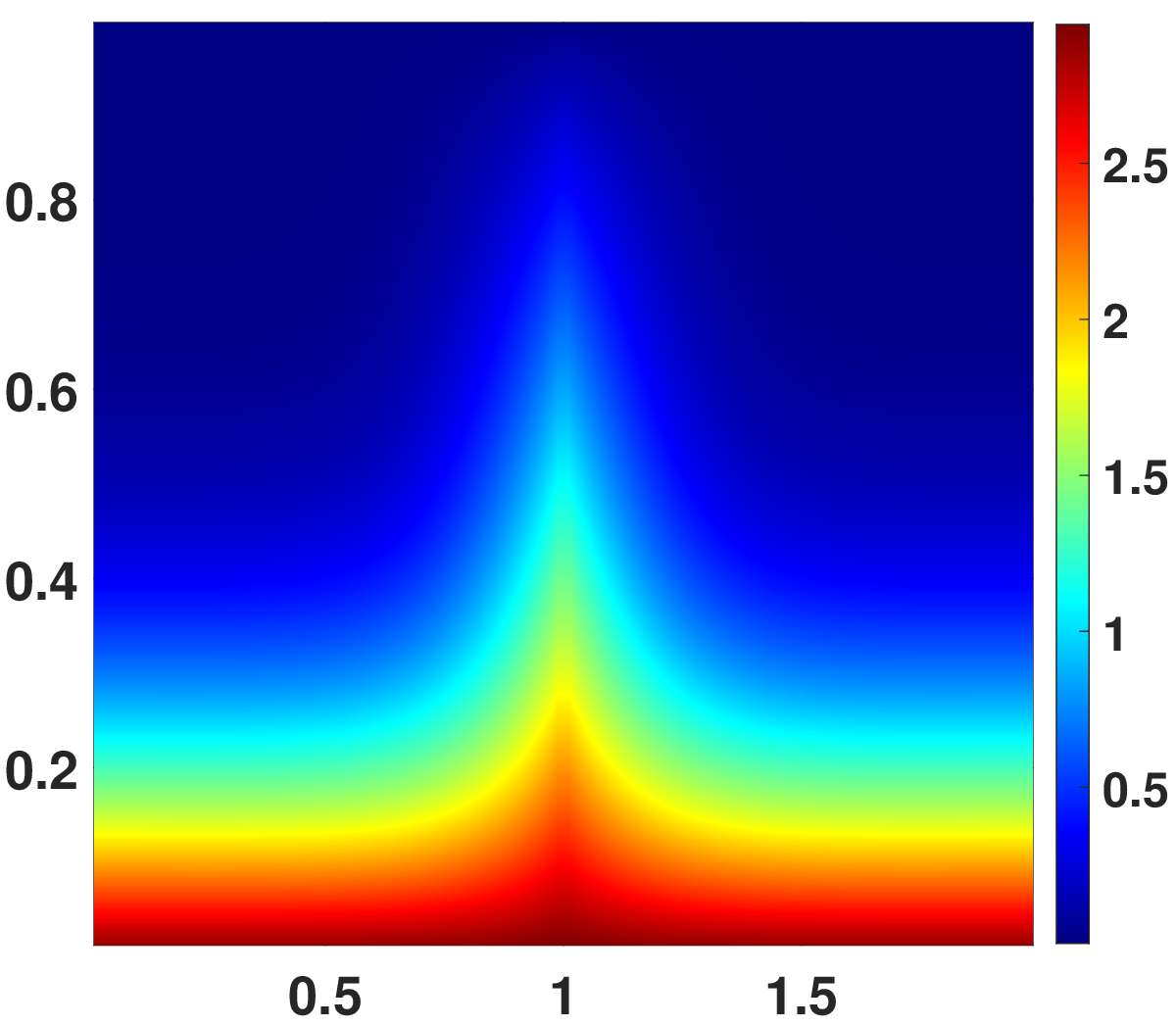}
\end{minipage}%
\begin{minipage}{0.35\textwidth}
\hspace{0.5cm}\includegraphics[scale=0.24]{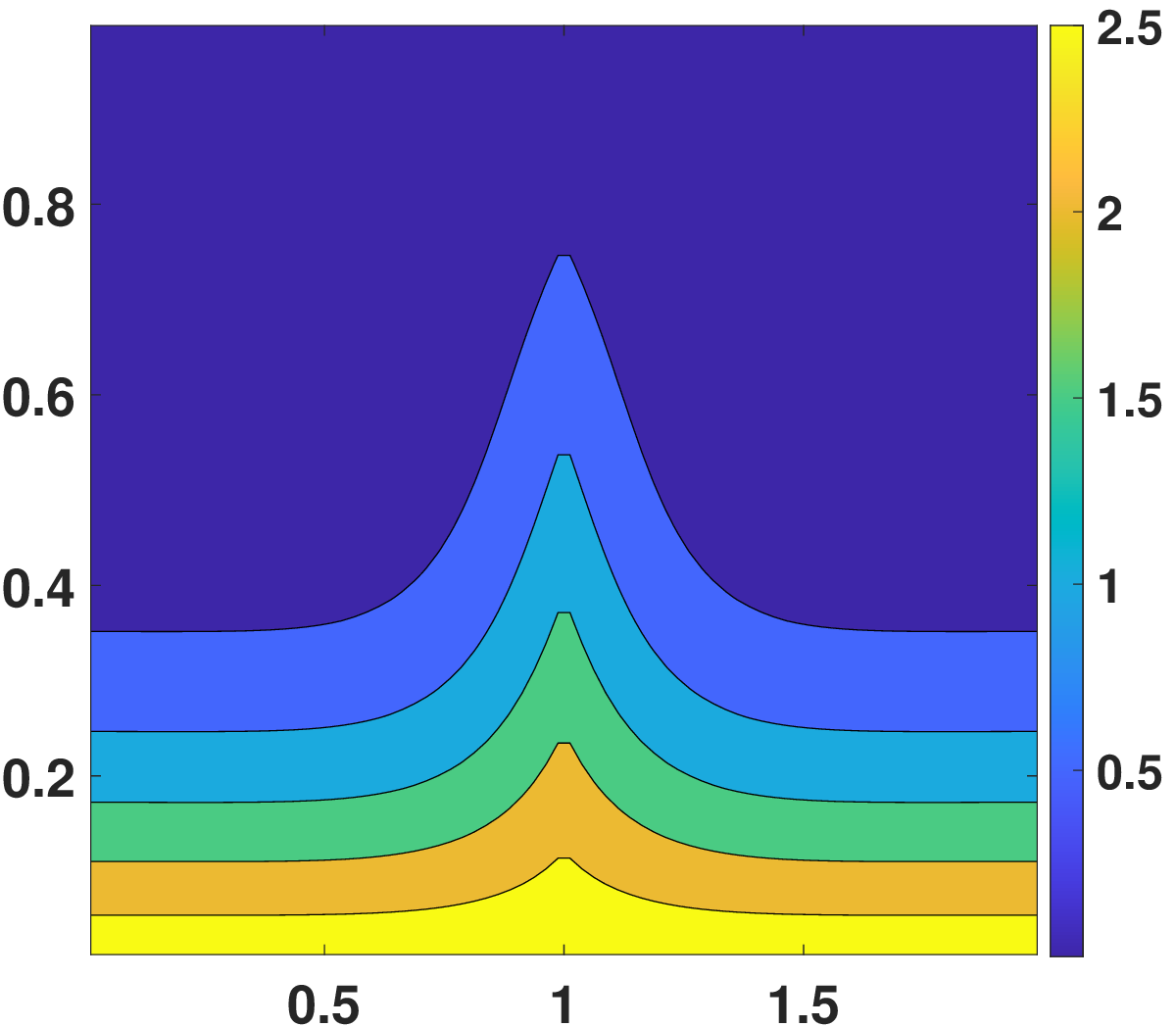}
\end{minipage}
\caption{\small [Test Case 3] Reference Concentration Field. (Left) Heat Map. (Right) Contour Map.}
\label{HeatMap_AdvDiff_Ref}
\vspace{-0.4cm}
\end{figure} 

To generate the observational data, we first solve the Darcy equations \eqref{Darcy_subdom}-\eqref{Darcy_fracture} to obtain the Darcy velocity fields $\bu = \left(\bu_1, \bu_2, \bu_f\right)$. Equations \eqref{reduced_adv_diff_subdom}-\eqref{reduced_adv_diff_fracture} are then solved with the computed velocity $\bu$ and the exact parameters in Table~\ref{physical_parametes_Testcase1} to obtain the reference solution. This is done on a spatial mesh size $h = 1/40$ and on a very fine temporal mesh $\Delta{t} = T/800$ with $T=5$. The data assimilation algorithm is performed over the same spatial mesh size and on a coarser temporal filtering steps  $\Delta{t}_{\text{Filter}} = T/ Nt_{\text{Filter}}$ where $Nt_{\text{Filter}}=50$, with the same assumption on the observational operator as in Test Case 2. Since concentration transport is driven by the Darcy velocity, we assume that there are some uncertainties in the process of solving the Darcy's flow system to increase the complexity of state estimation. Thus, we add a small perturbation $\xi_{\text{Darcy}} = 0.0001\epsilon_D$ to the reference Darcy velocities with $\epsilon_{\text{Darcy}} \sim N(0, \pmb{I})$  and input the disturbed velocity to the forward solver. The reference Darcy velocity field and one example of the perturbed velocity field are shown in Figure~\ref{DarcyField} on a very coarse mesh for the purpose of illustration. We can see that even with a small amount of noise, the perturbed Darcy field is significantly chaotic, which makes the task of state estimation challenging. 

Similar to Test Case 2, we aim to investigate the performance of the United Filter under various levels of perturbation. Thus, we consider two types of disturbed noise $\omega$ in \eqref{StateModel}: $\tilde{\omega}^1 = 0.001\sqrt{\Delta{t}_{\text{Filter}}}\tilde{\varepsilon}^1$ and $\tilde{\omega}^2 = 0.1\sqrt{\Delta{t}_{\text{Filter}}}\tilde{\varepsilon}^2$, where $\tilde{\varepsilon}^i \sim \mathcal{N}(0, \pmb{I}_l)$ for $i=1, 2$. 
Regarding the setting for the parameter estimation, as the exact values are small, we apply the assimilation algorithm to approximate $\rho_i = 1/d_i, i=1, 2$ and $\rho_{\gamma} = 1/\alpha_{\gamma}$. We choose the number of Direct Filter particles to be $M=40$, and the number of iterations in the United Filter to be $R=4$.
\begin{figure}[h!]
\vspace{-0.2cm}
\centering
\begin{minipage}{0.35\textwidth}
\includegraphics[scale=0.24]{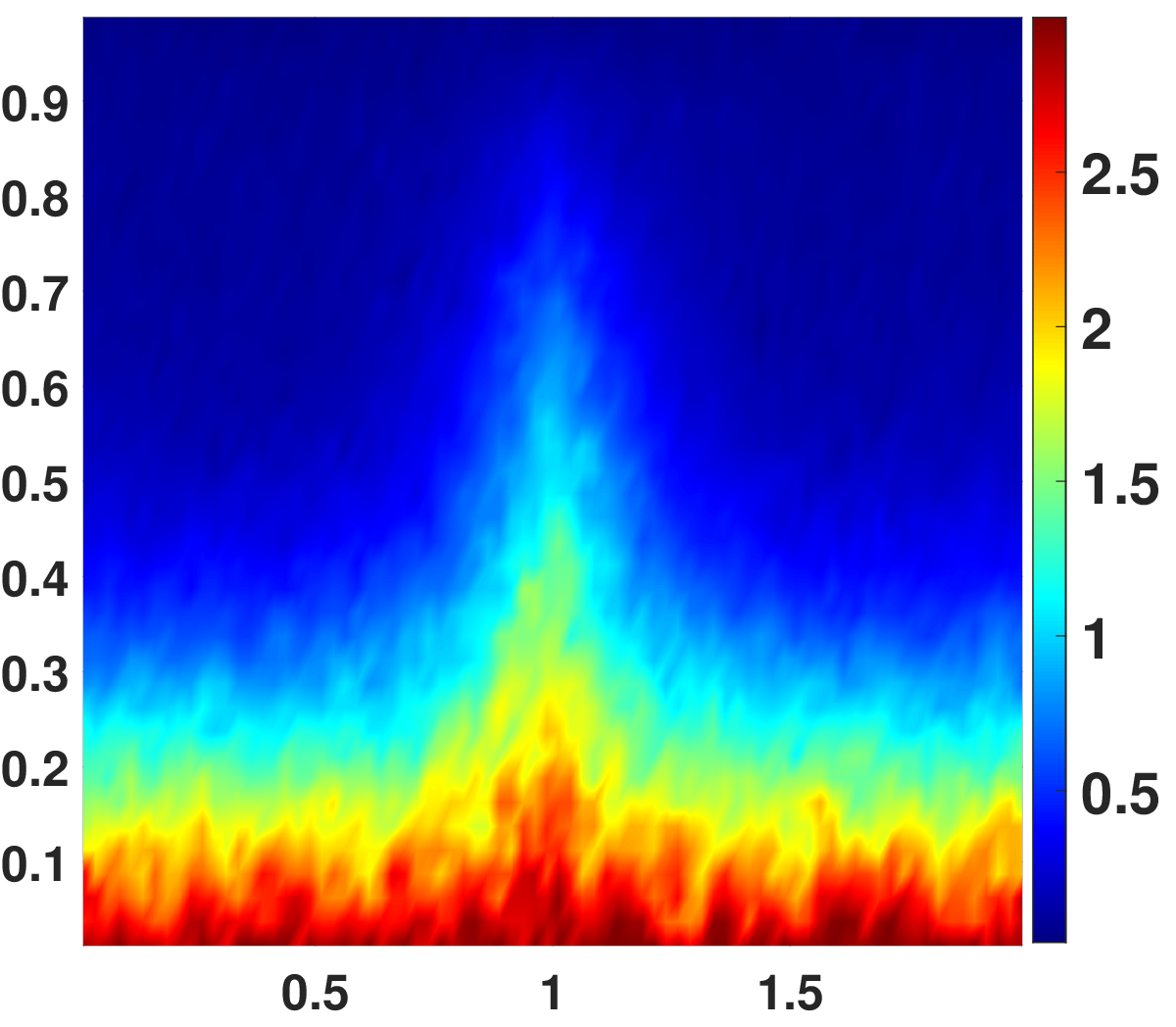}
\end{minipage}%
\begin{minipage}{0.35\textwidth}
\hspace{0.5cm}\includegraphics[scale=0.24]{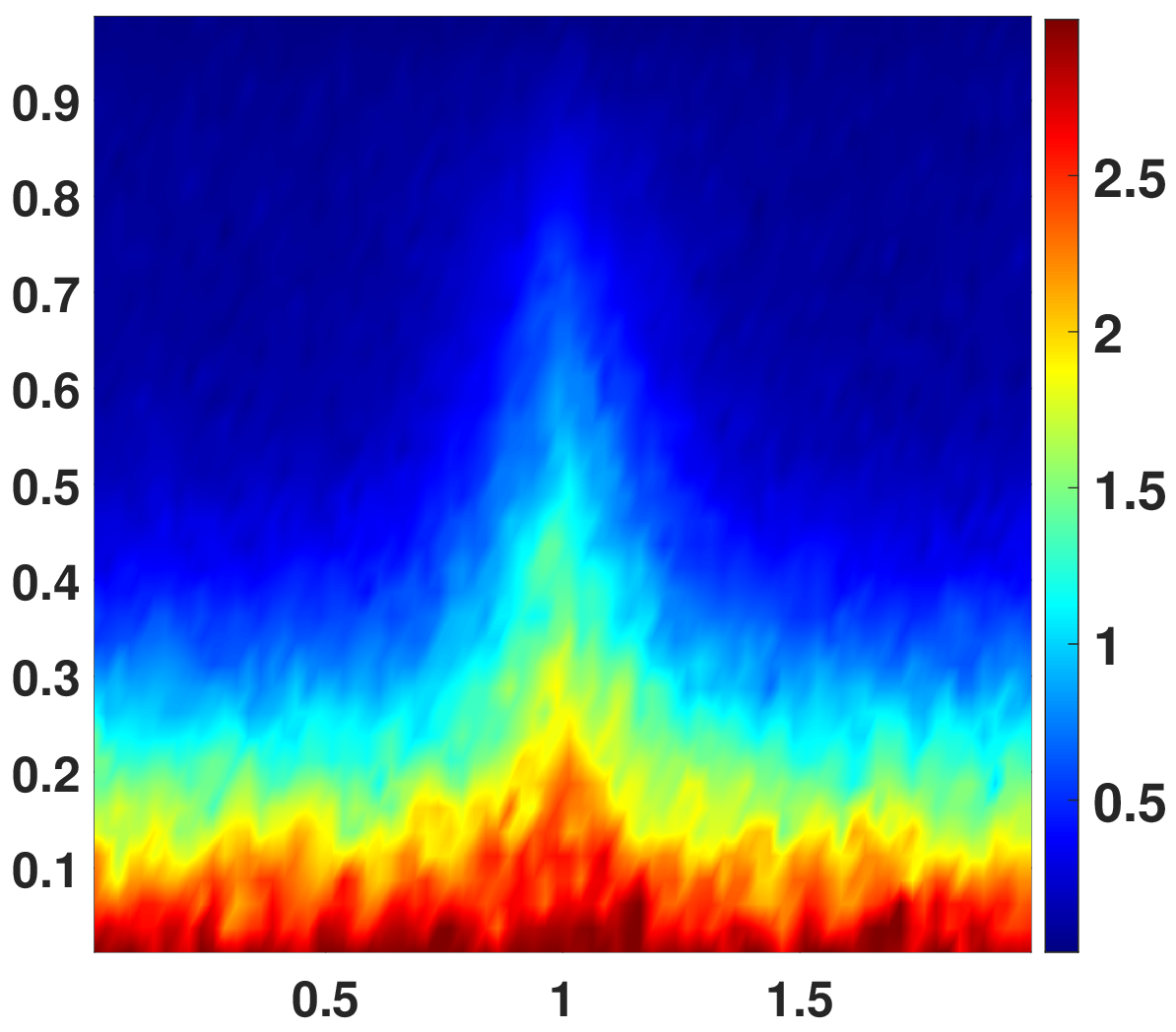}
\end{minipage}
\caption{[Test Case 3] Heat map of estimated concentration field by the United Filter: (Left) With noise $\tilde{\omega}^1$. (Right) With noise $\tilde{\omega}^2$.}
\label{HeatMap_AdvDiff_UnitedF}
\vspace{-0.4cm}
\end{figure}
\begin{figure}[h!]
\centering
\begin{minipage}{0.35\textwidth}
\includegraphics[scale=0.24]{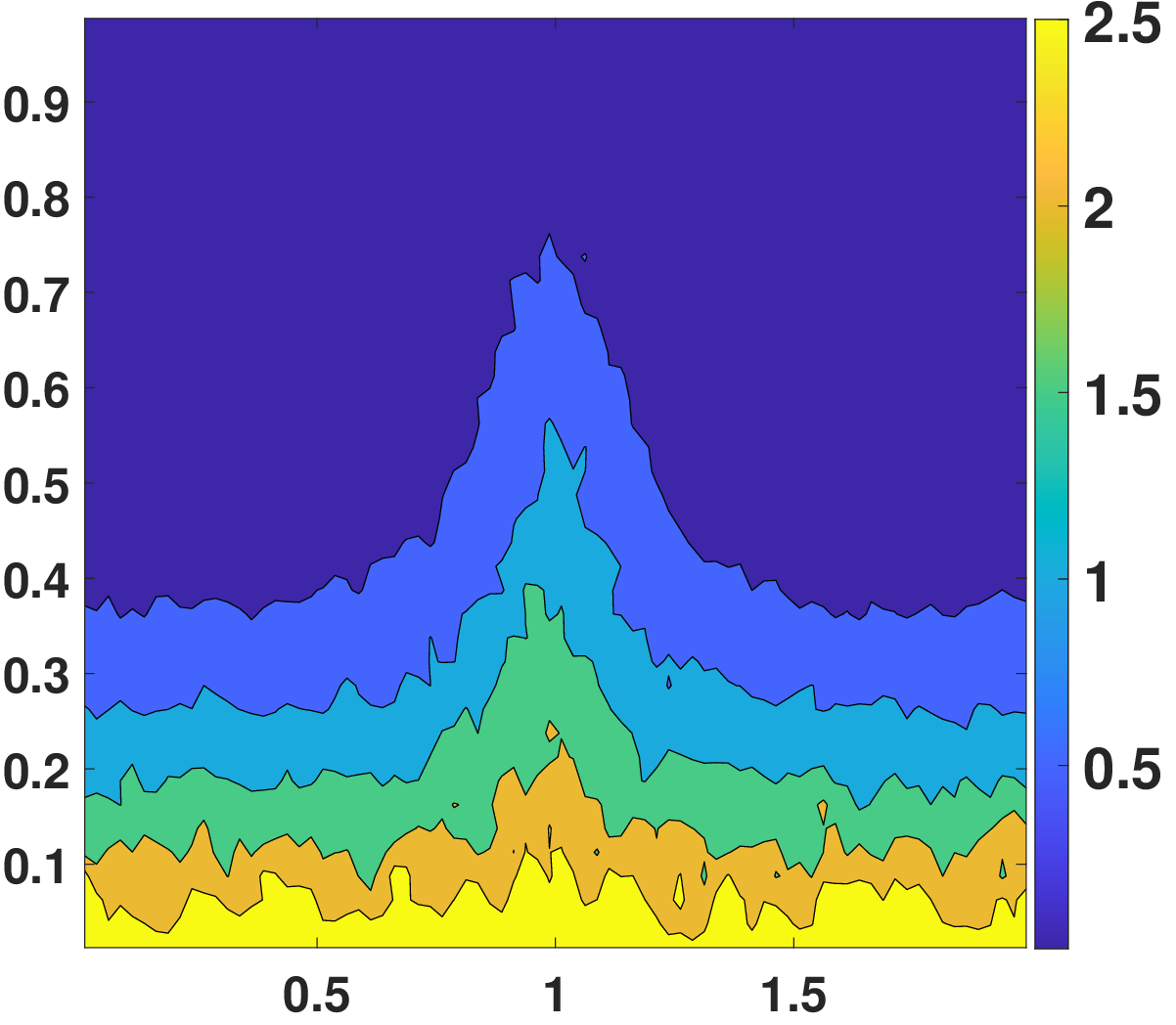}
\end{minipage}%
\begin{minipage}{0.35\textwidth}
\hspace{0.5cm}\includegraphics[scale=0.24]{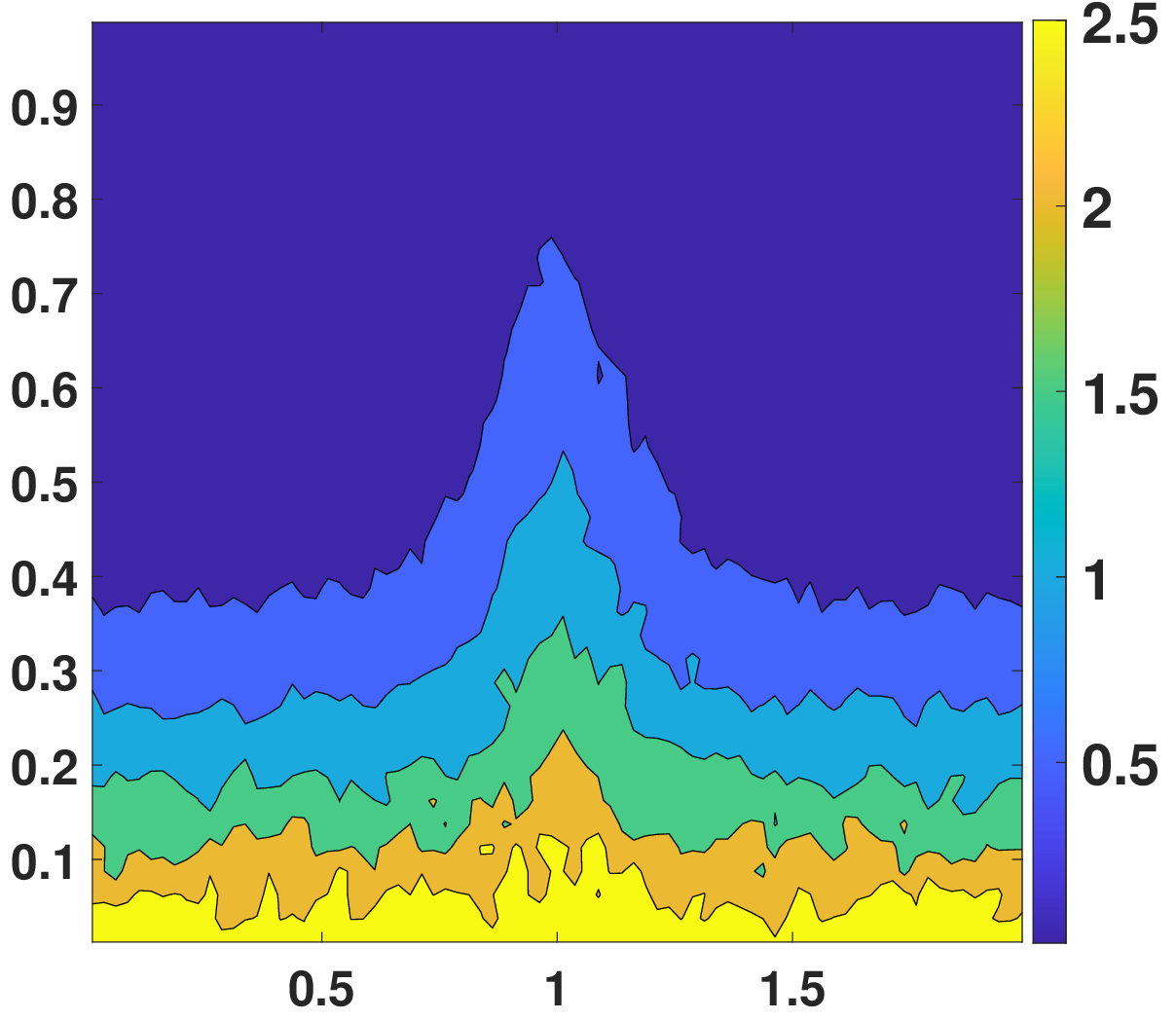}
\end{minipage}
\caption{[Test Case 3] Contour Map of estimated Concentration Field by the United Filter. (Left) With noise $\tilde{\omega}^1$. (Right) With noise $\tilde{\omega}^2$.}
\label{ContourMap_AdvDiff_UnitedF}
\vspace{-0.3cm}
\end{figure}
\begin{figure}[h!]
\centering
\begin{minipage}{0.35\textwidth}
\includegraphics[scale=0.24]{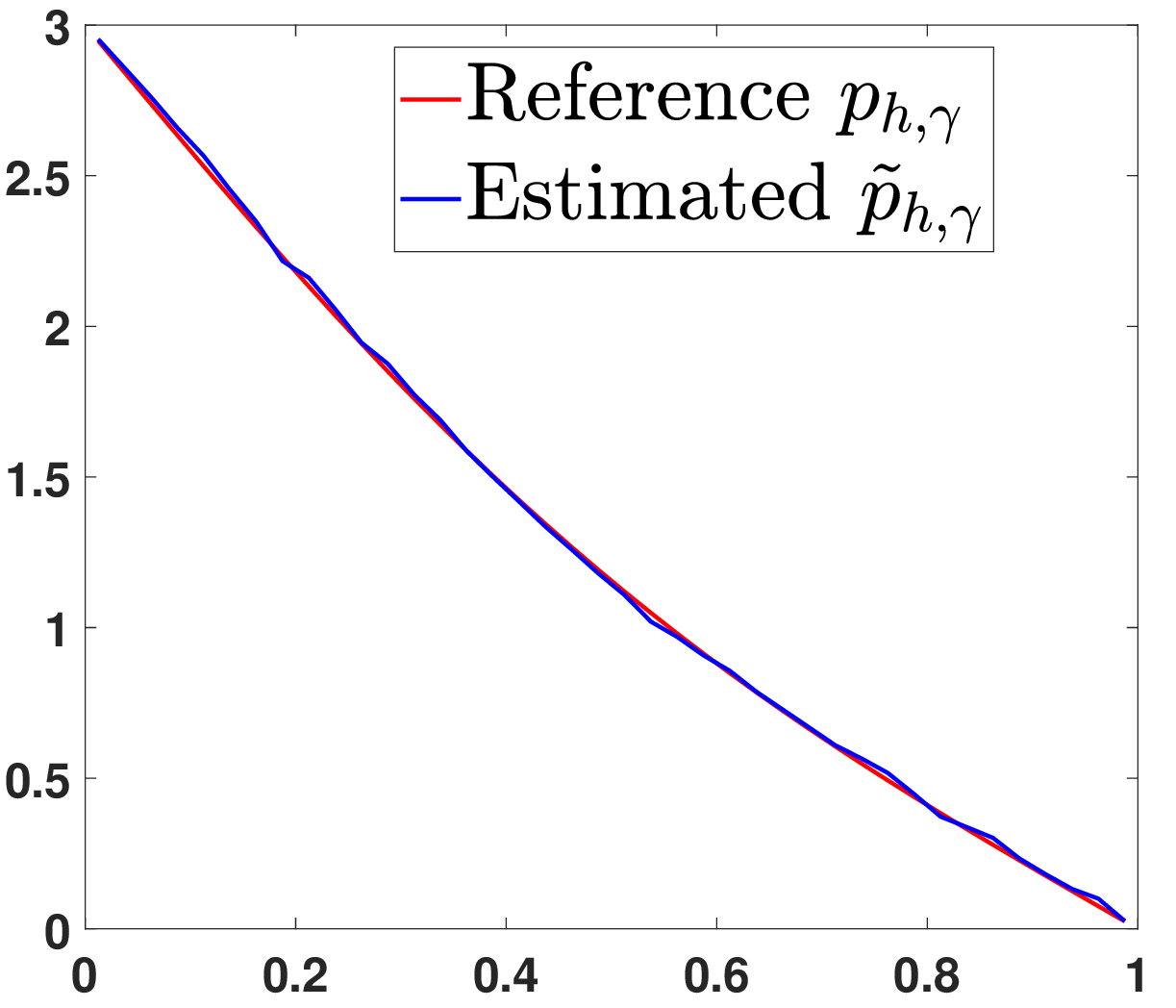}
\end{minipage}%
\begin{minipage}{0.35\textwidth}
\hspace{0.5cm}\includegraphics[scale=0.24]{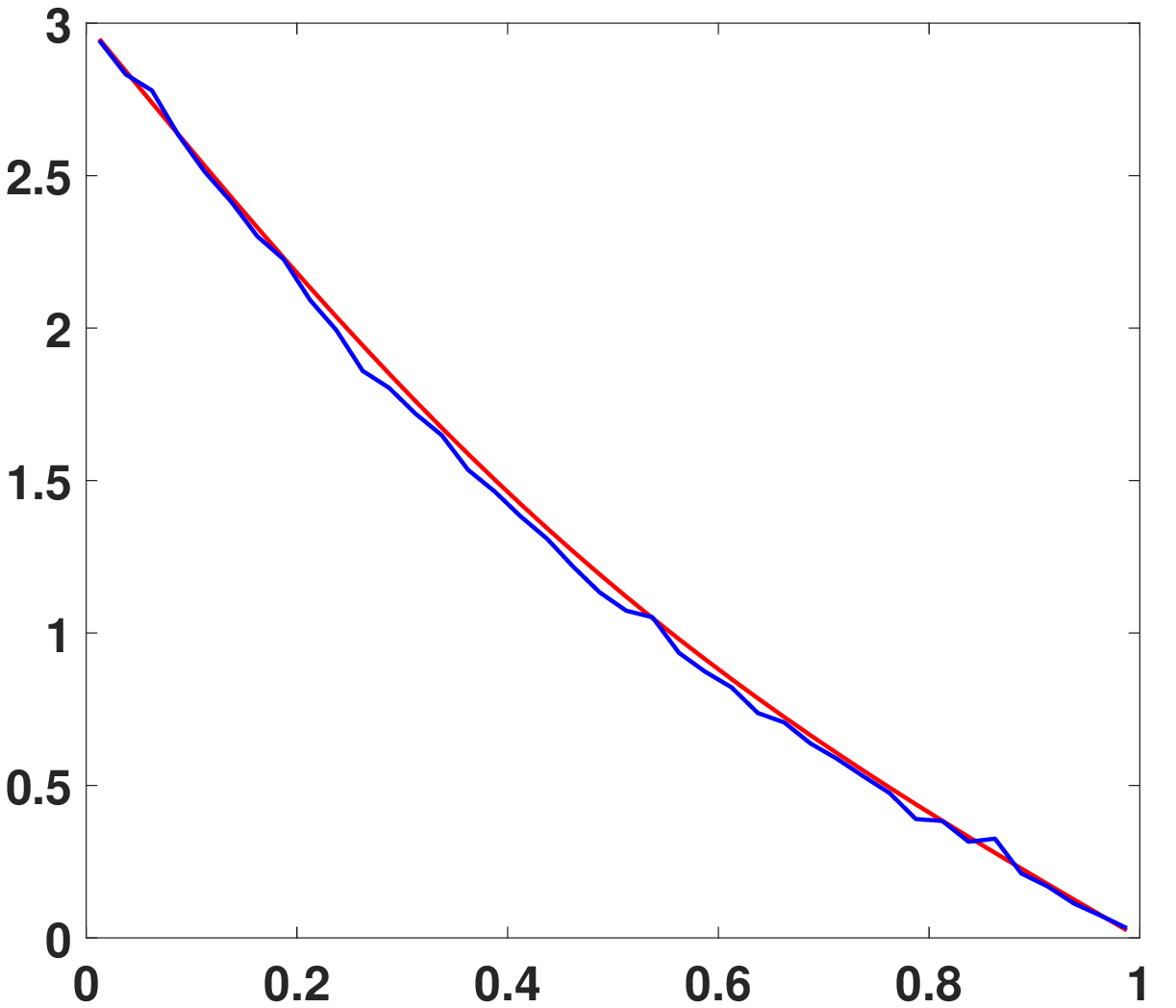}
\end{minipage}
\caption{[Test Case 3] Concentration on the fracture estimated by the United Filter. (Left) With noise $\tilde{\omega}^1$. (Right) With noise $\tilde{\omega}^2$.}
\label{1DConcentration_AdvDiff}
\vspace{-0.4cm}
\end{figure}

We begin by presenting the results of the state estimation process. Figure~\ref{HeatMap_AdvDiff_Ref} displays the true reference concentration fields, represented by the heat map (left) and the contour map (right). Subsequently, Figures~\ref{HeatMap_AdvDiff_UnitedF}-\ref{ContourMap_AdvDiff_UnitedF} illustrate the concentration fields estimated by the United Filter, shown as heat and contour maps, respectively, under varying levels of perturbation noise. The reference and the estimated 1D concentration on the fracture are shown in Figures~\ref{1DConcentration_AdvDiff}. These results demonstrate that, regardless of the noise level - whether due to model inaccuracies or numerical simulation errors - the United Filter consistently provides accurate estimates of the concentration fields.
\begin{figure}[h!]
\centering
\begin{minipage}{0.3\textwidth}
\hspace{0.1cm}\includegraphics[scale=0.24]{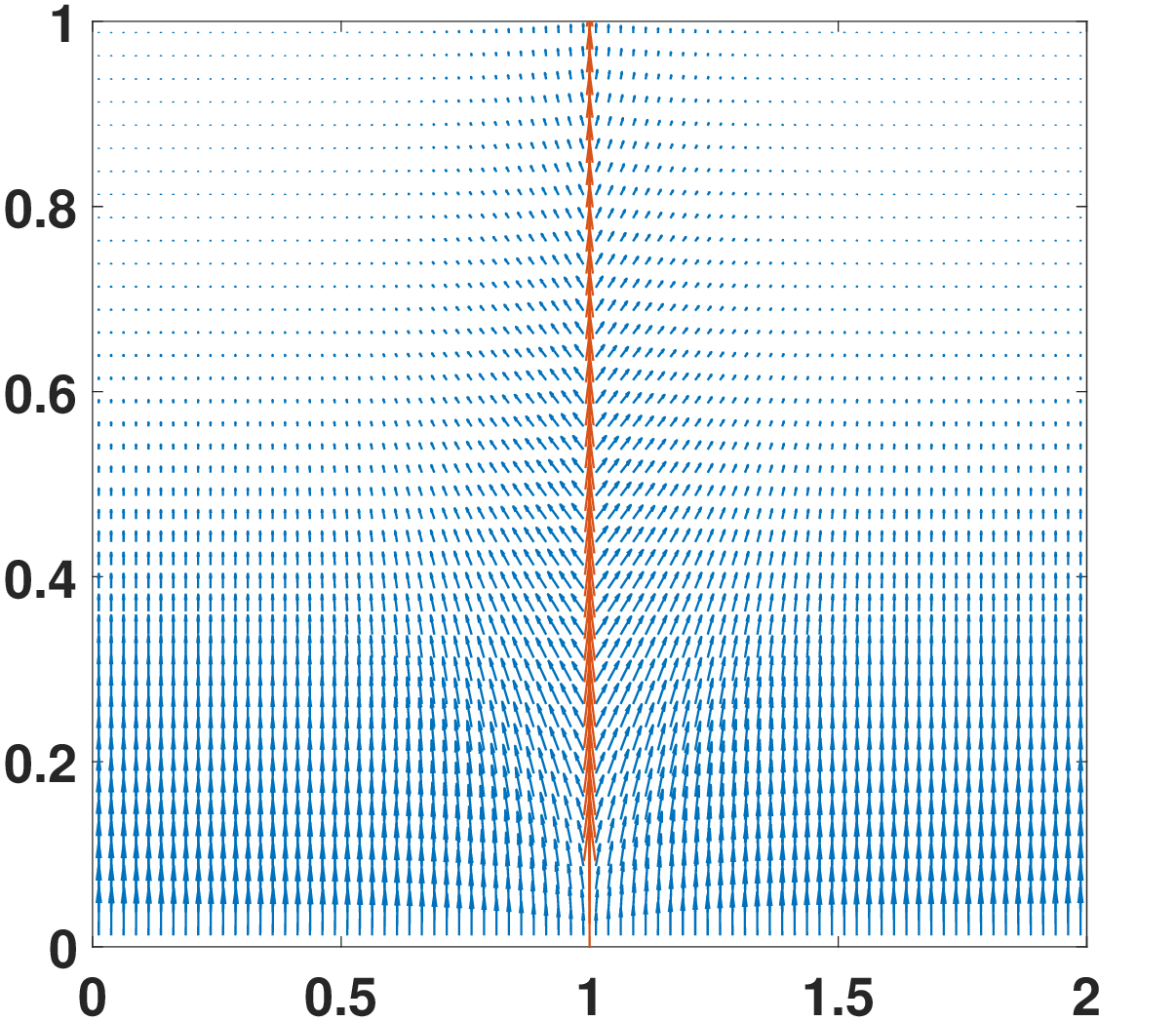}
\end{minipage}
\hspace{0.1cm}
\begin{minipage}{0.3\textwidth}
\includegraphics[scale=0.24]{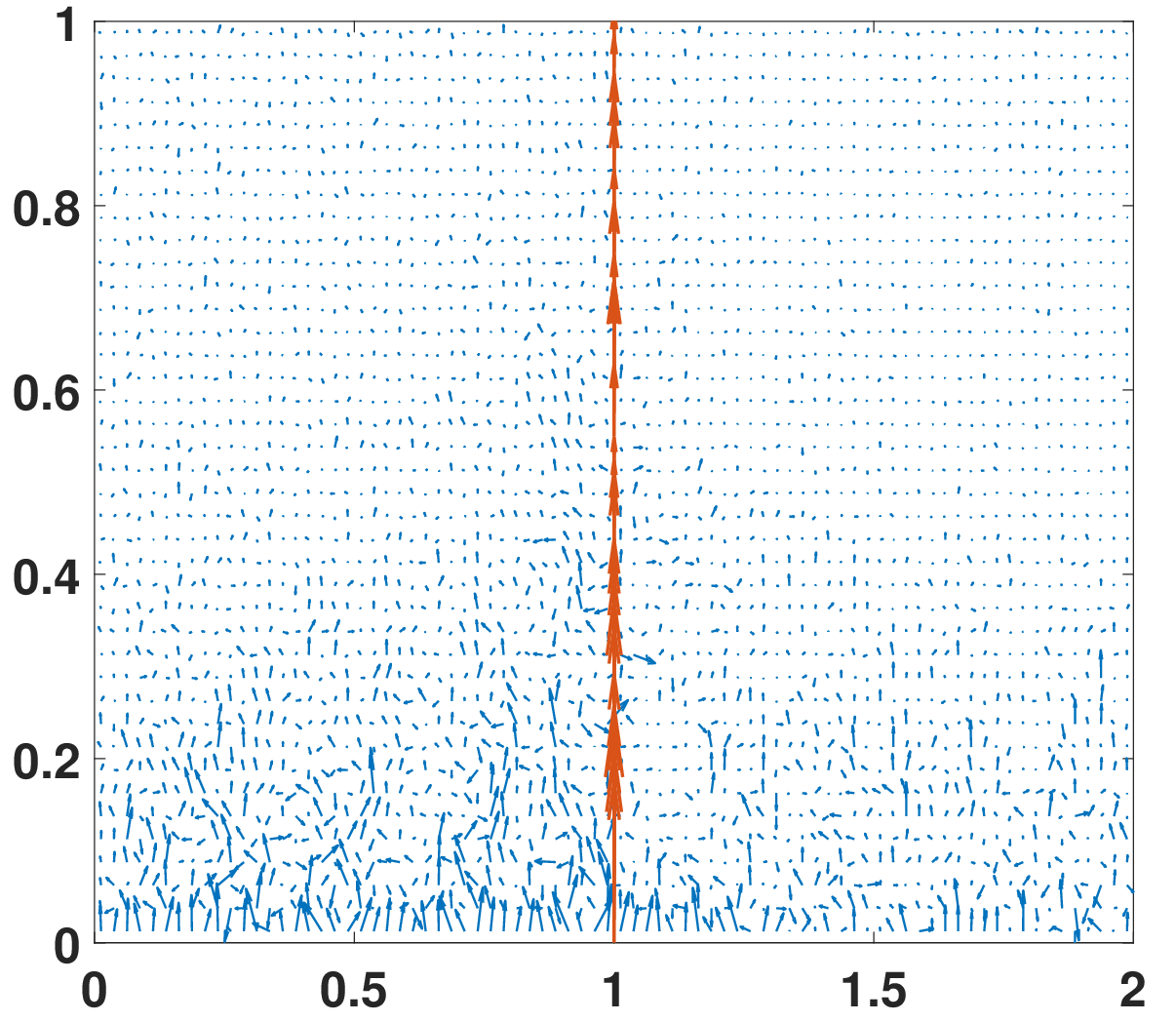}
\end{minipage}
\hspace{0.1cm}
\begin{minipage}{0.3\textwidth}
\includegraphics[scale=0.24]{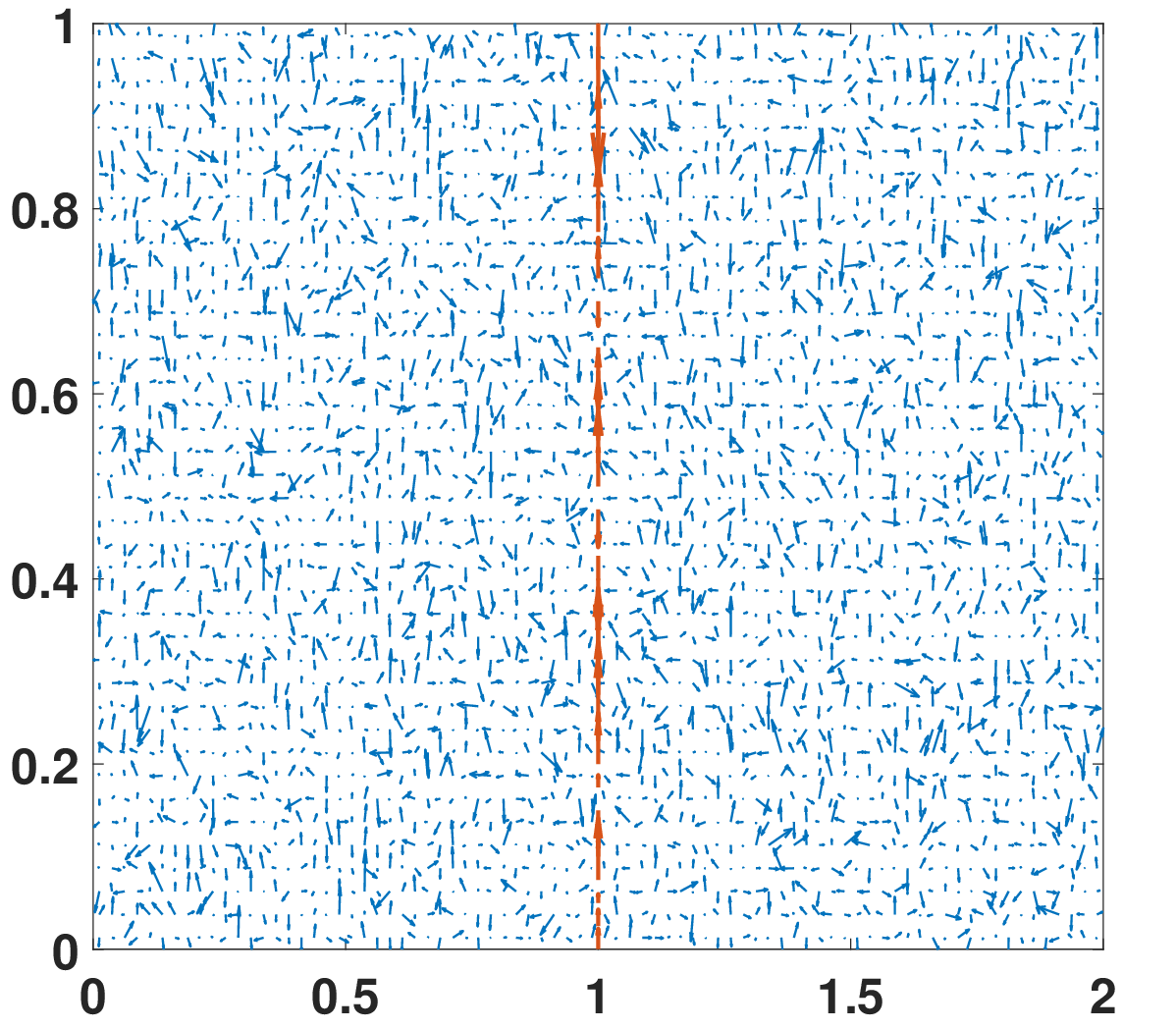}
\end{minipage}
\caption{[Test Case 3] Velocity Fields estimated by the United Filter. (Left) Reference Fields. (Middle) Estimation with noise $\tilde{\omega}^1$. (Right) Estimation with noise $\tilde{\omega}^2$.}
\label{Velo_AdvDiff_UnitedF}
\vspace{-0.4cm}
\end{figure} 
\begin{figure}[h!]
\centering
\begin{minipage}{0.4\textwidth}
\hspace{1cm}\includegraphics[scale=0.25]{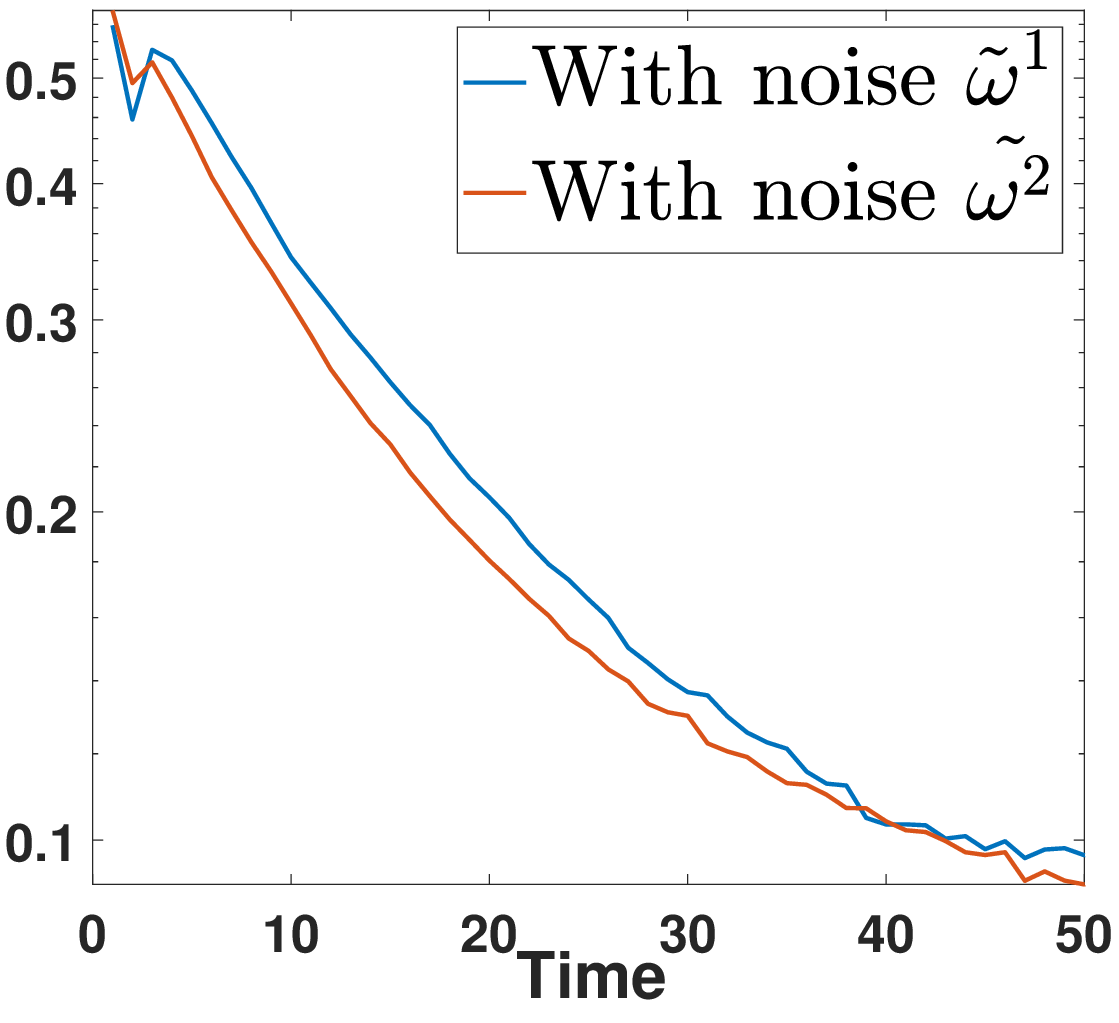}
\end{minipage}%
\caption{[Test Case 3] RMSEs for state estimation by the United Filter. }
\label{RMSEAdvDiff}
\vspace{-0.1cm}
\end{figure}

We next investigate the results for the estimation of the velocity fields under the two levels of perturbation noise. We display the reference and the estimated velocities in Figure~\ref{Velo_AdvDiff_UnitedF}. To better illustrate the accuracy of state estimation, we plot the RMSEs for state estimation in Figure~\ref{RMSEAdvDiff} using a logarithmic scale to enhance the visibility of velocity errors.
\begin{figure}[h!]
\vspace{-0.3cm}
\centering
\begin{minipage}{0.33\textwidth}
\includegraphics[scale=0.22]{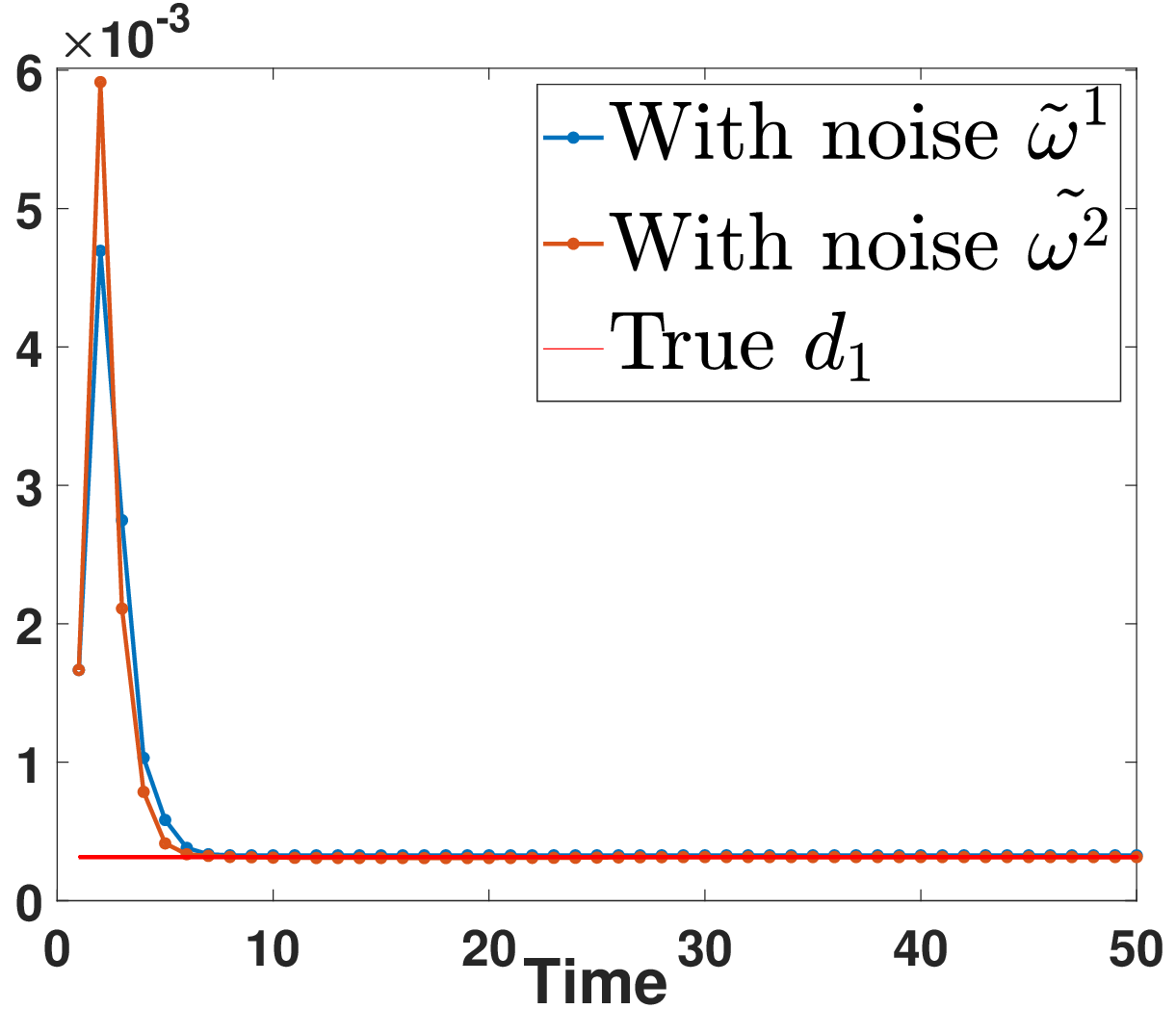}
\end{minipage}%
\begin{minipage}{0.33\textwidth}
\hspace{0.1cm}\includegraphics[scale=0.22]{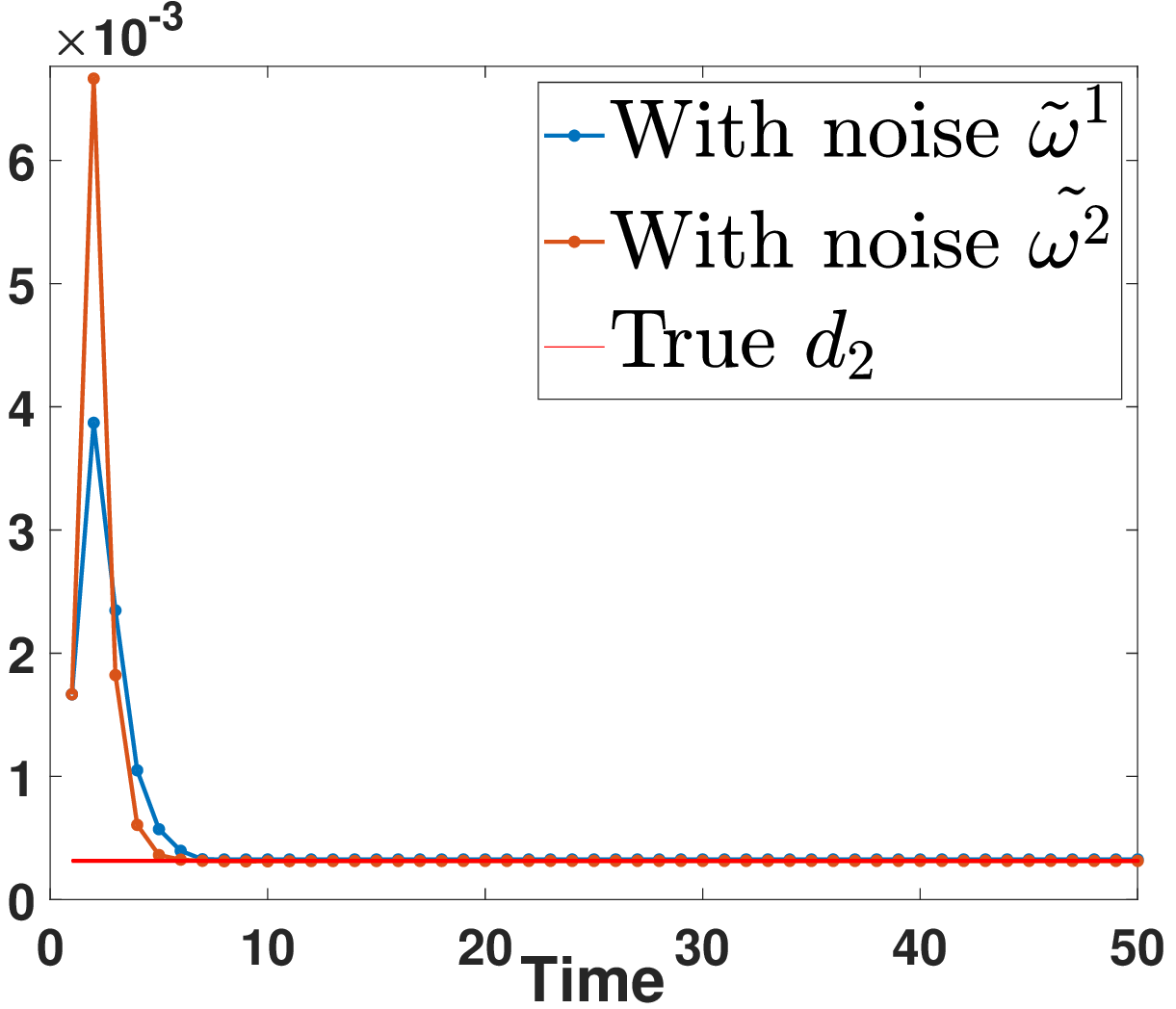}
\end{minipage}%
\begin{minipage}{0.33\textwidth}
\hspace{0.1cm}\includegraphics[scale=0.22]{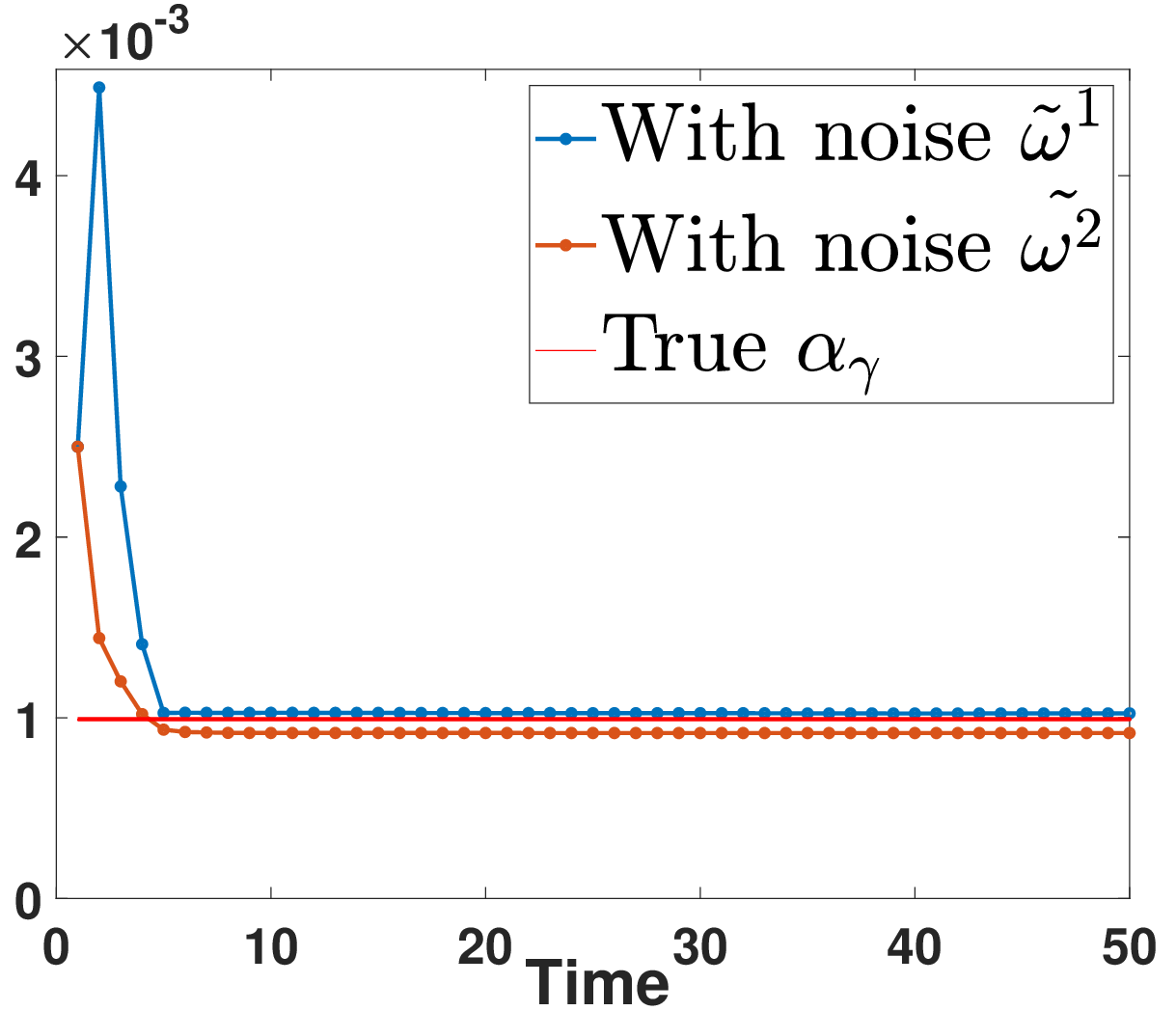}
\end{minipage}
\caption{[Test Case 3] Parameter estimation by the United Filter: (Left) Estimation of $d_1$
(Middle) Estimation of $d_2$. (Right) Estimation of $\alpha_{\gamma}$.}
\label{ParaEst_AdvDiff}
\vspace{-0.4cm}
\end{figure}

\begin{figure}[h!]
\centering
\begin{minipage}{0.35\textwidth}
\includegraphics[scale=0.23]{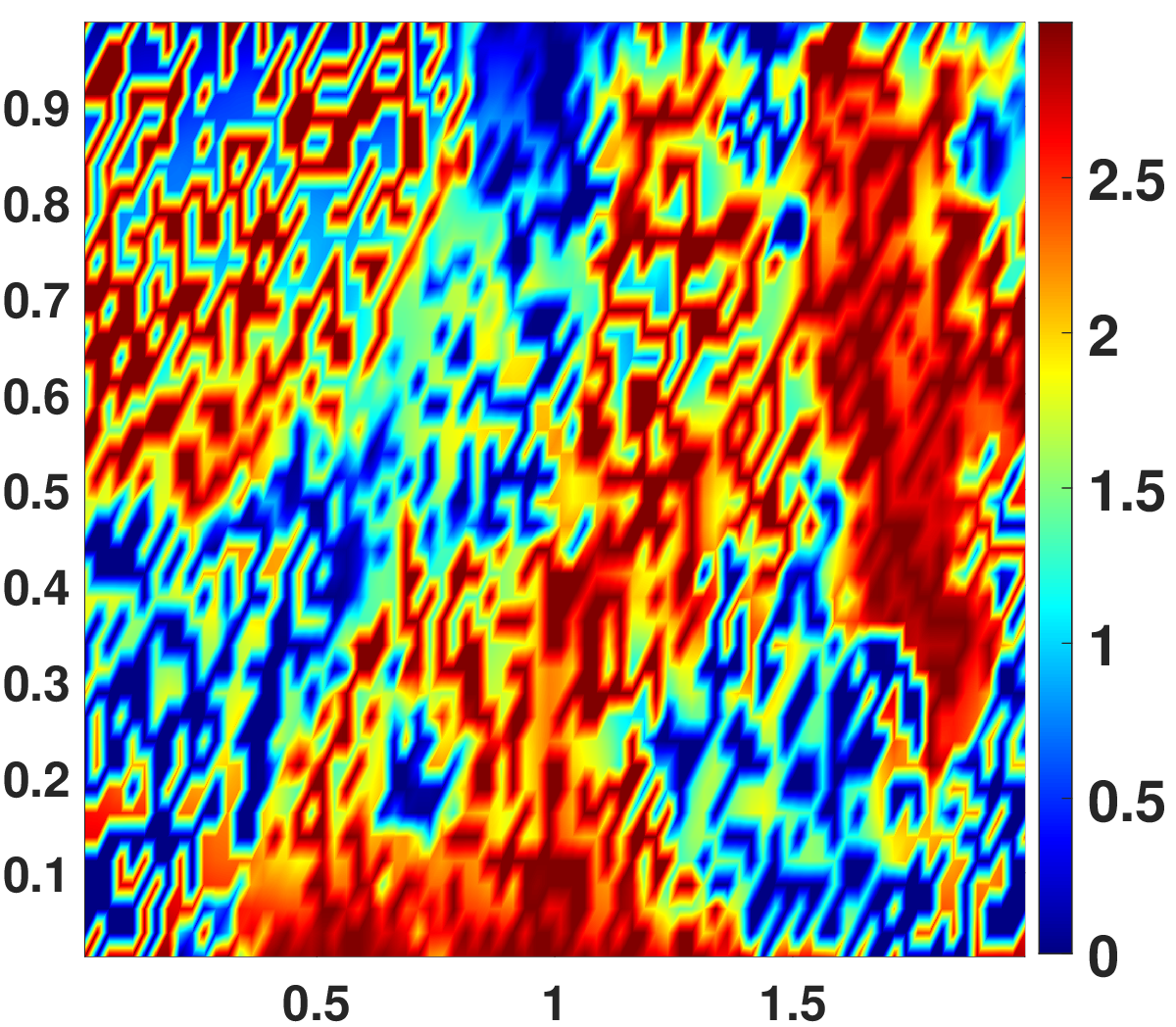}
\end{minipage}%
\begin{minipage}{0.35\textwidth}
\hspace{0.5cm}\includegraphics[scale=0.23]{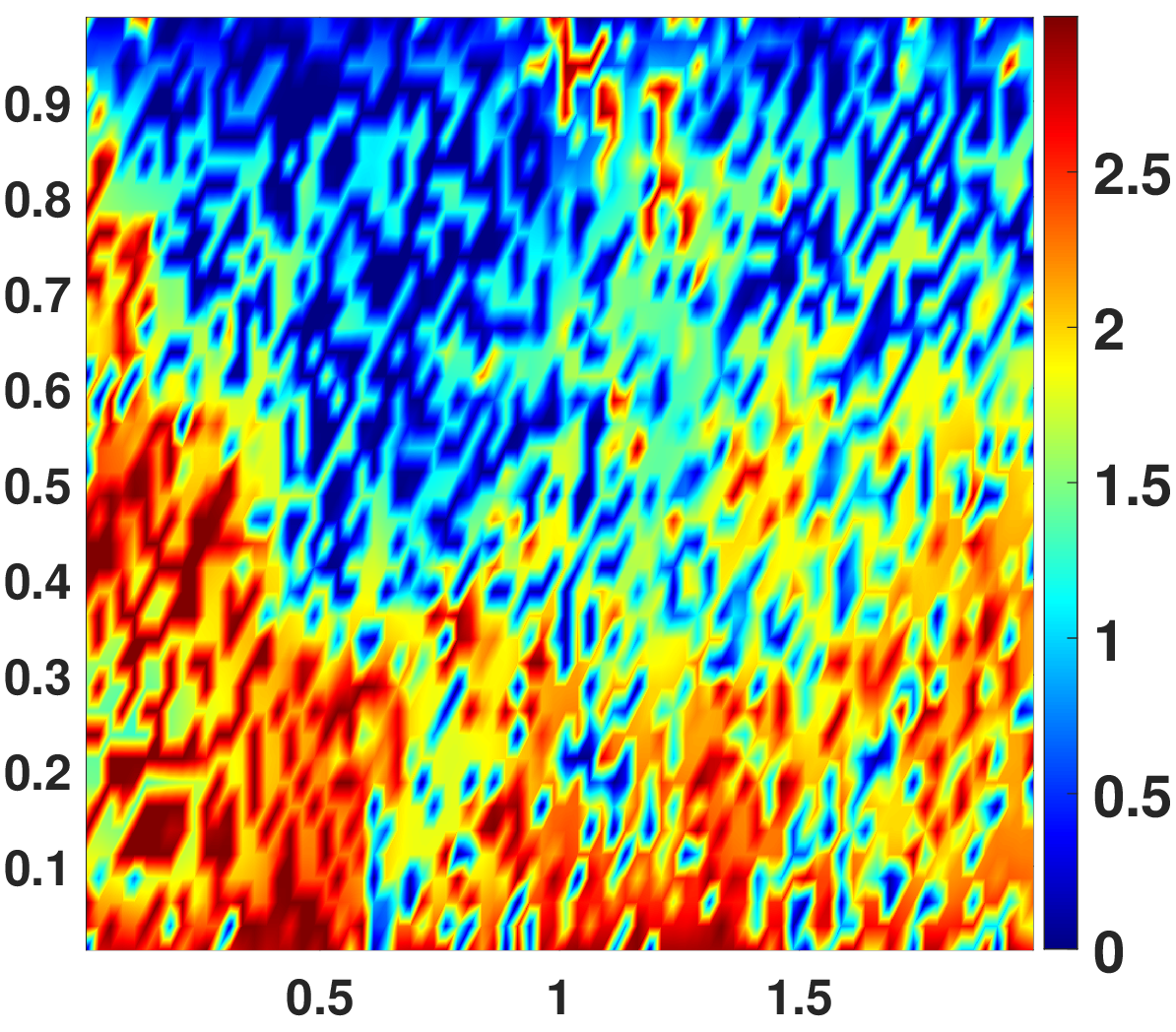}
\end{minipage}
\caption{[Test Case 3] Concentration heat map estimated by AugEnKF: (Left) With noise $\tilde{\omega}^1$; (Right) With noise $\tilde{\omega}^2$. This is a comparison result to the United Filter's estimates presented in Figure \ref{HeatMap_AdvDiff_UnitedF}. }
\label{HeatMap_AdvDiff_EnKF}
\vspace{-0.4cm}
\end{figure}

\begin{figure}[h!]
\centering
\begin{minipage}{0.35\textwidth}
\includegraphics[scale=0.23]{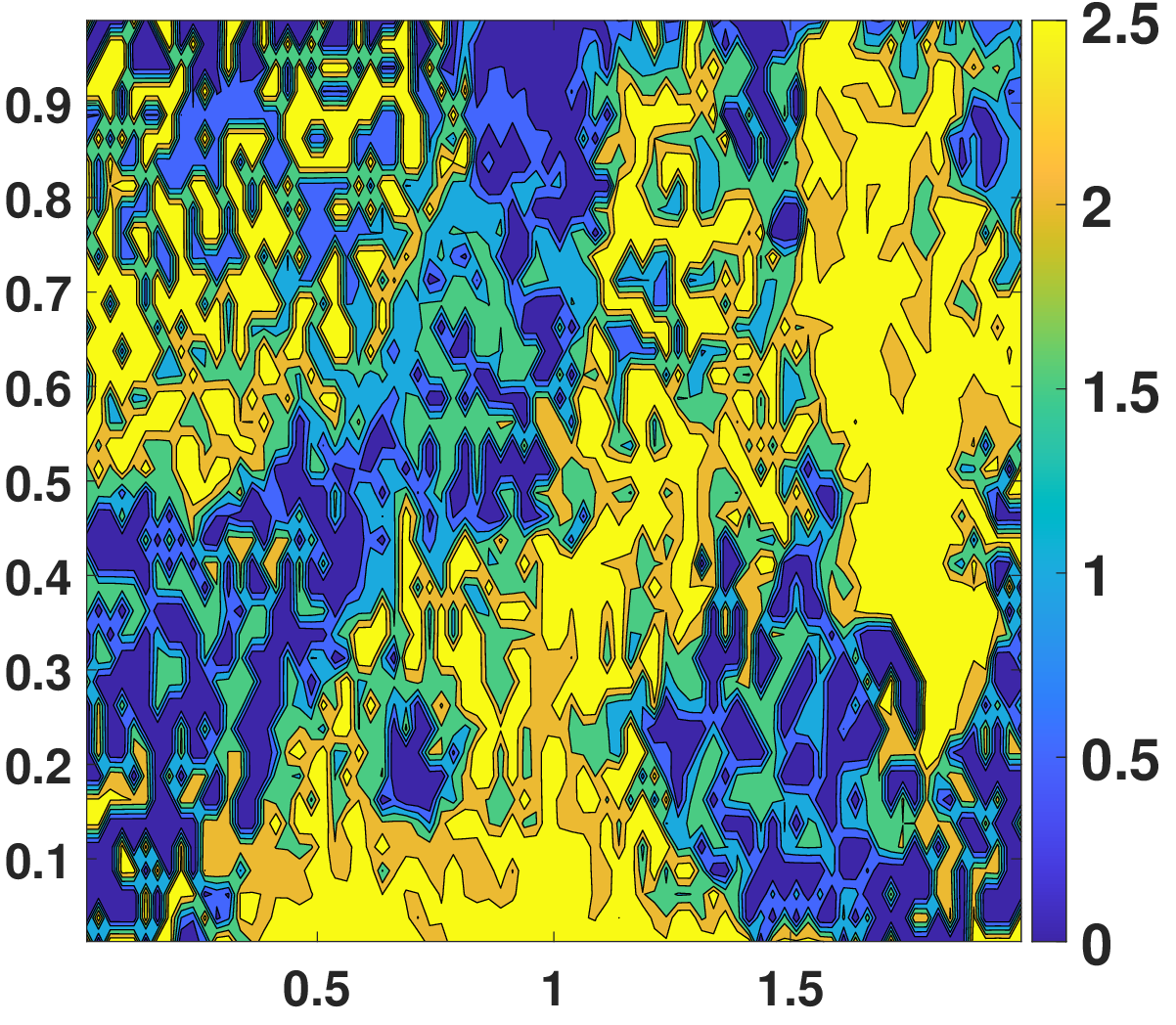}
\end{minipage}%
\begin{minipage}{0.35\textwidth}
\hspace{0.5cm}\includegraphics[scale=0.23]{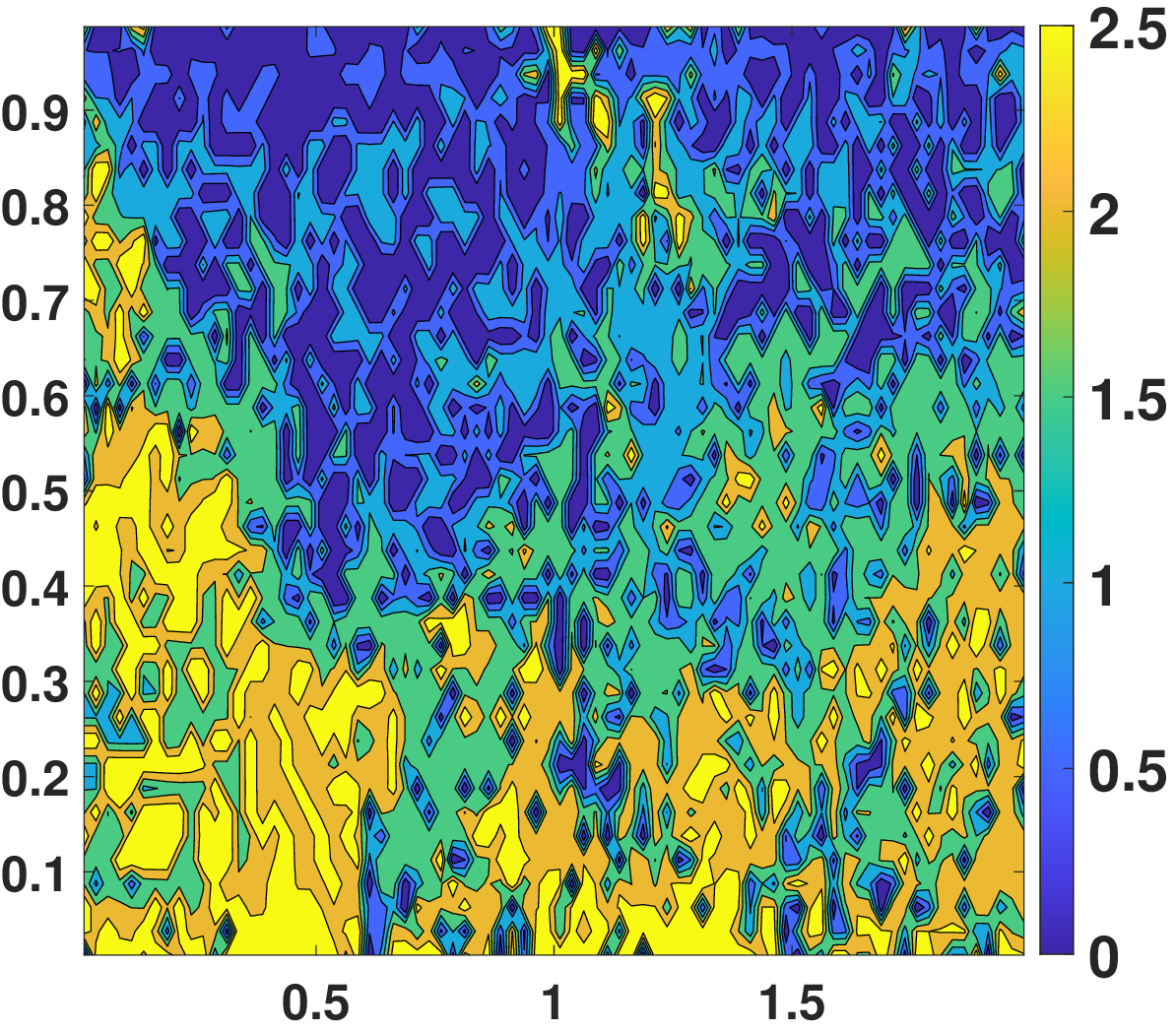}
\end{minipage}
\caption{[Test Case 3] Pressure contour map estimated by the AugEnKF. (Left) With noise $\tilde{\omega}^1$; (Right) With noise $\tilde{\omega}^2$. This is a comparison result to the United Filter's estimates presented in Figure \ref{ContourMap_AdvDiff_UnitedF}.}
\label{ContourMap_AdvDiff_EnKF}
\vspace{-0.4cm}
\end{figure}

\begin{figure}[h!]
\centering
\begin{minipage}{0.35\textwidth}
\includegraphics[scale=0.22]{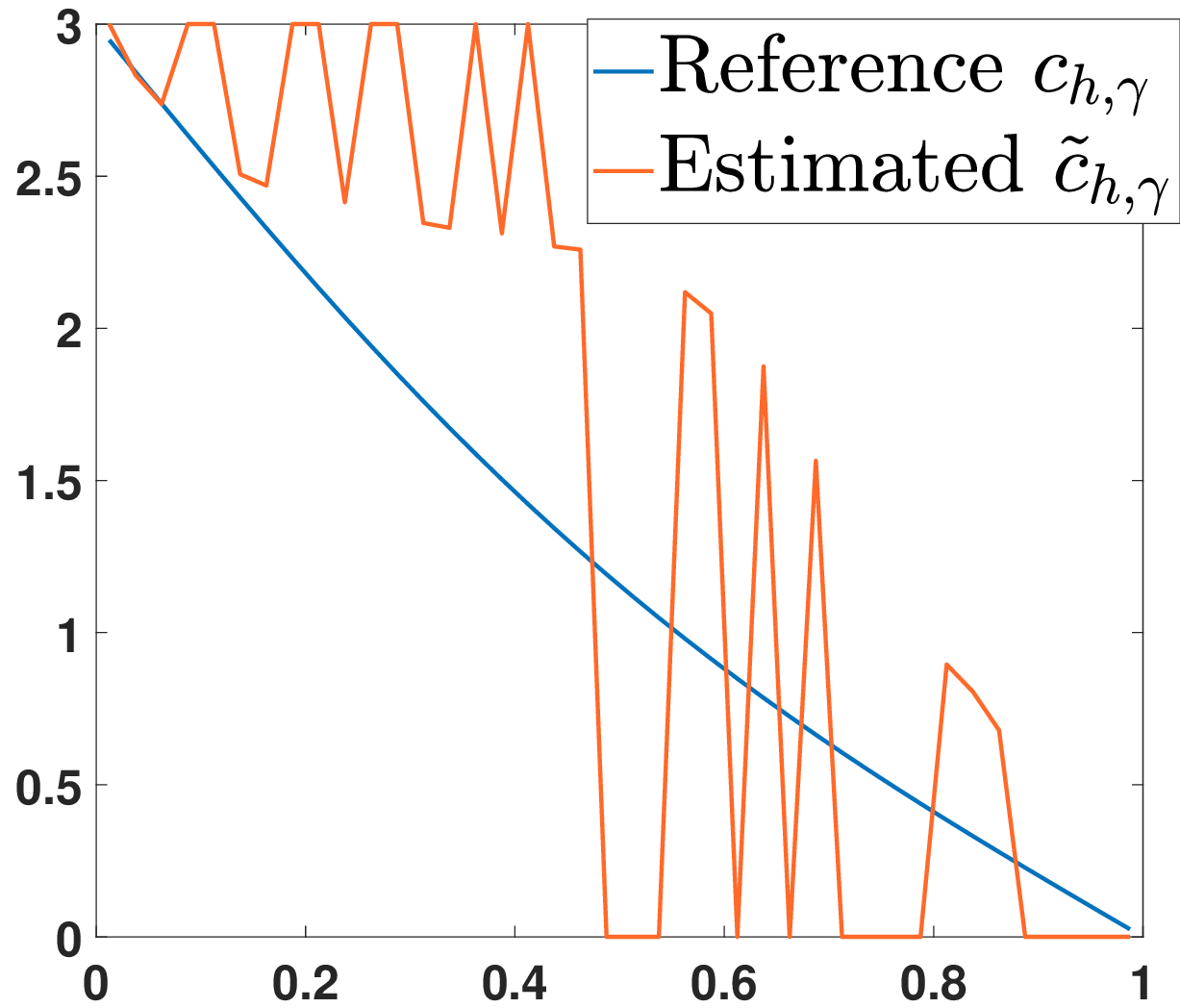}
\end{minipage}%
\begin{minipage}{0.35\textwidth}
\includegraphics[scale=0.22]{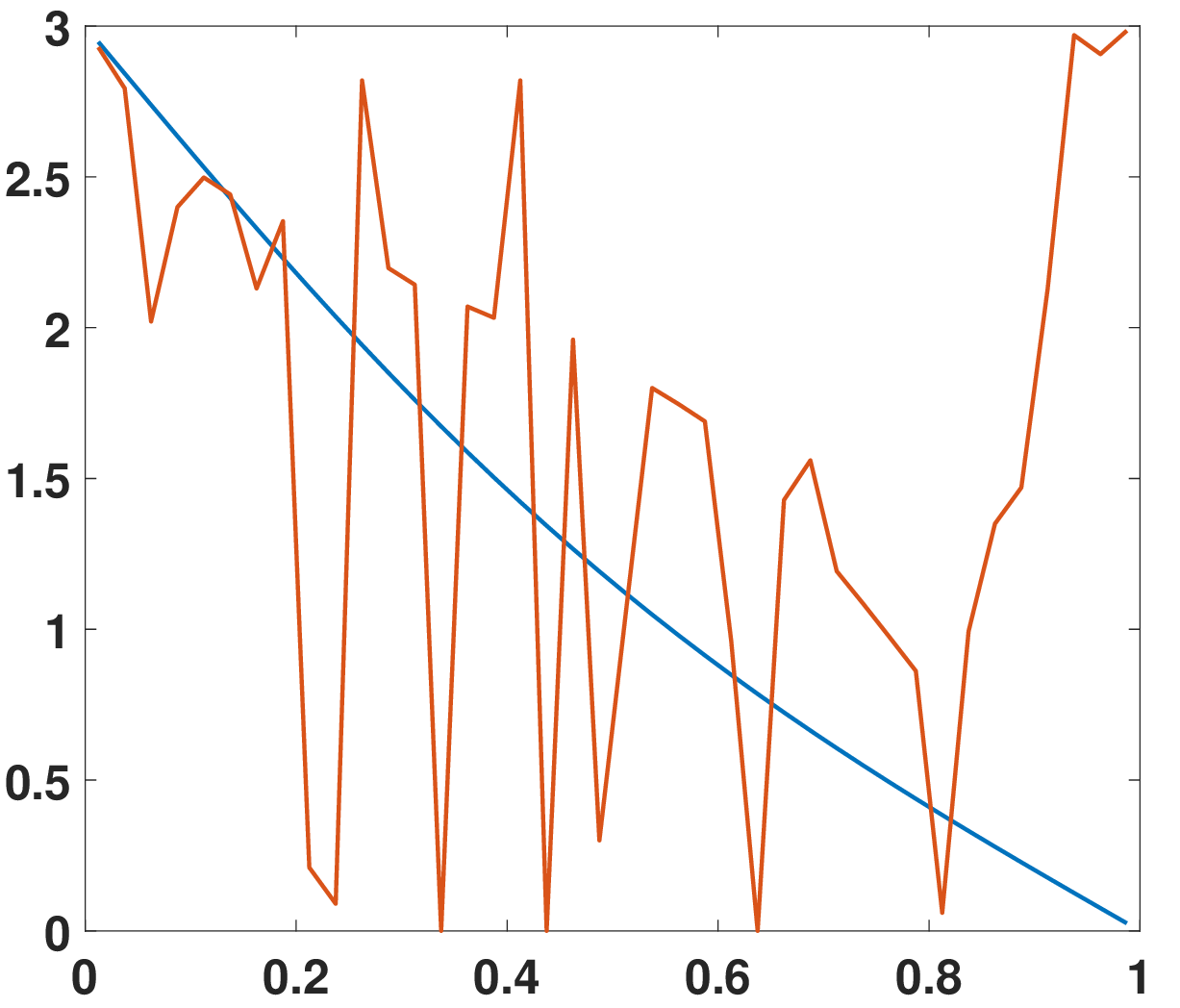}
\end{minipage}
\caption{[Test Case 3] Concentration estimation on the fracture by AugEnKF: (Left) With noise $\tilde{\omega^1}$; (Right) With noise $\tilde{\omega^2}$. This is a comparison result to the United Filter's estimates presented in Figure \ref{1DConcentration_AdvDiff}. }
\label{1DConcentration_AdvDiff_EnKF}\vspace{-0.4cm}
\end{figure}

\begin{figure}[h!]
\centering
\begin{minipage}{0.35\textwidth}
\includegraphics[scale=0.23]{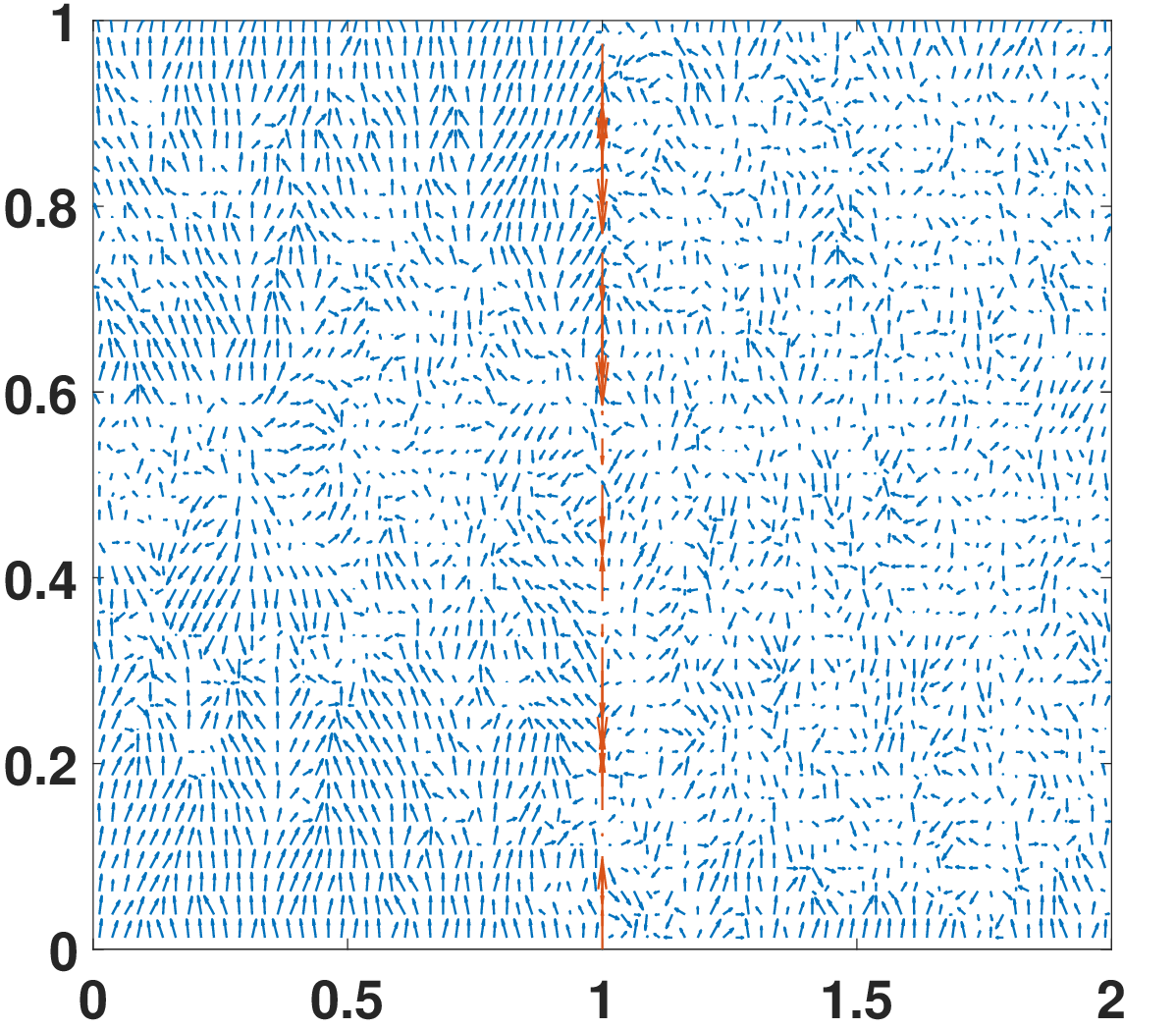}
\end{minipage}%
\begin{minipage}{0.35\textwidth}
\hspace{0.2cm}\includegraphics[scale=0.23]{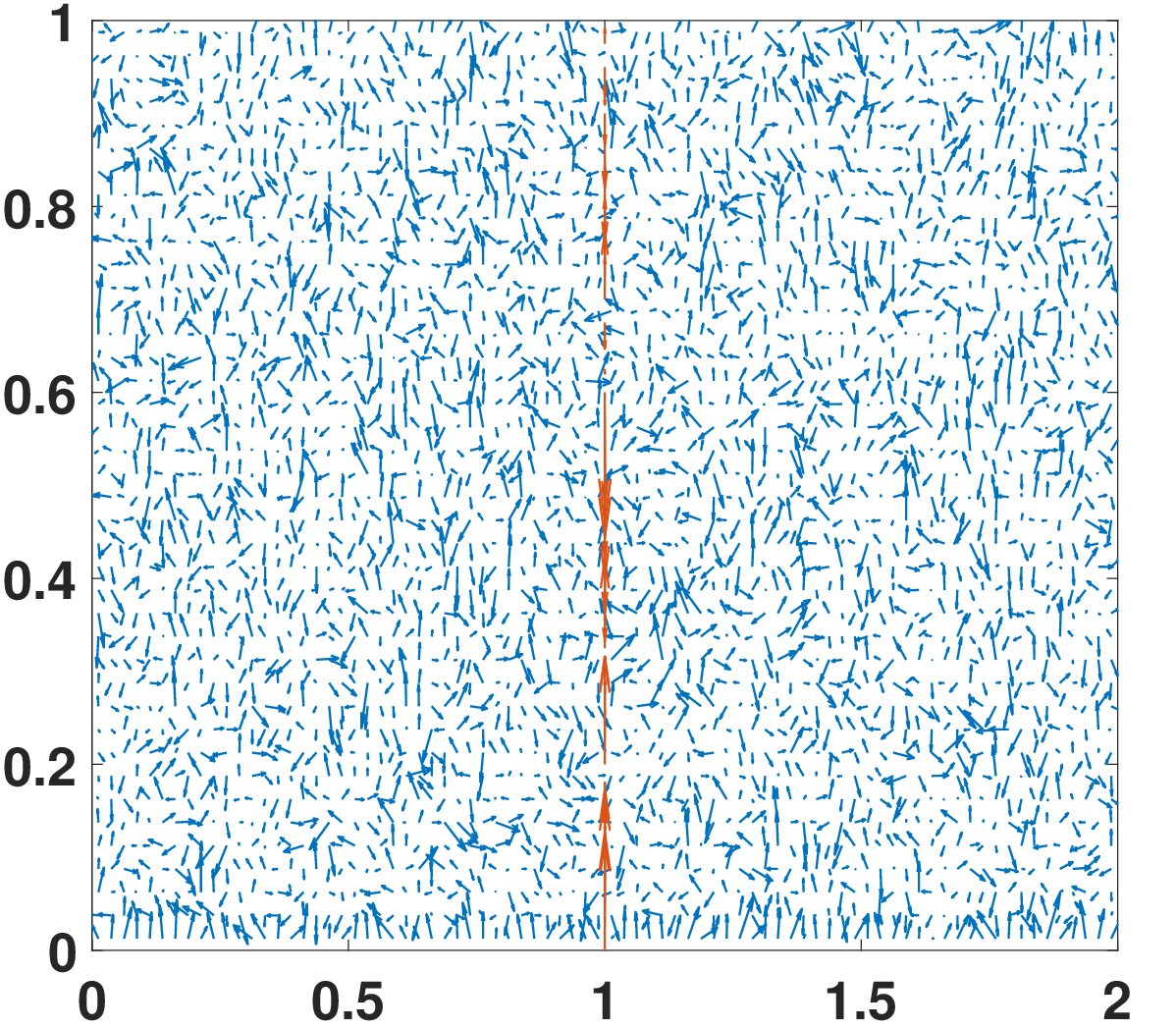}
\end{minipage}
\caption{[Test Case 3] Estimated velocity fields by the AugEnKF: (Left)  With noise $\tilde{\omega}^1$. (Right) With noise $\tilde{\omega}^2$. This is a comparison result to the United Filter's estimates presented in Figure \ref{Velo_AdvDiff_UnitedF}. }
\label{Velo_AdvDiff_EnKF}
\vspace{-0.4cm}
\end{figure}

\begin{figure}[h!]
\centering
\begin{minipage}{0.4\textwidth}
\includegraphics[scale=0.25]{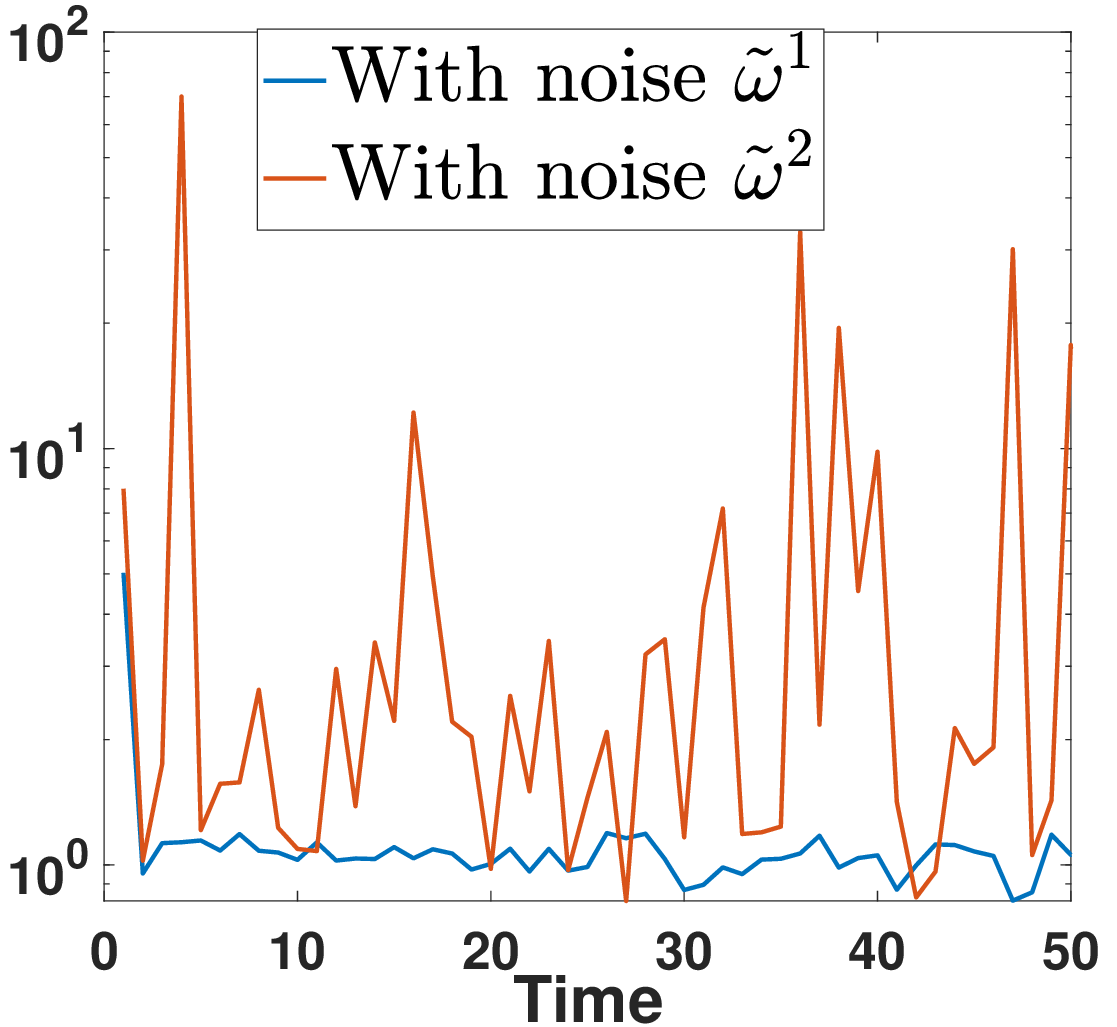}   
\end{minipage}%
\caption{[Test Case 3] RMSEs of the AugEnKF for state estimation with two types of noise. This is a comparison result to the United Filter's estimates presented in Figure \ref{RMSEAdvDiff}.  }
\label{RMSE_Velo_AdvDiff_EnKF}
\vspace{-0.4cm}
\end{figure}

We conclude our presentation of numerical results by showcasing the results for Test Case 3 obtained using the AugEnKF method. We also evaluate the performance of AugEnKF under two model uncertainties $\tilde{\omega}^1$ and $\tilde{\omega}^2$. The estimated concentration fields are shown in Figures~\ref{HeatMap_AdvDiff_EnKF}-\ref{ContourMap_AdvDiff_EnKF}, while the estimated concentration on the fracture is depicted in Figure~\ref{1DConcentration_AdvDiff_EnKF}. By comparing with the true reference concentration presented in Figure~\ref{HeatMap_AdvDiff_Ref}, we conclude that the AugEnKF failed to recover the concentration from sparse and noisy observations. The estimated velocity fields are presented in Figure~\ref{Velo_AdvDiff_EnKF}. The overall state estimation performance of the AugEnKF is presented in Figure~\ref{RMSE_Velo_AdvDiff_EnKF} in terms of  RMSEs. From this figure, we observe that the AugEnKF exhibits relatively low errors when the perturbation noise is at a lower level, i.e., $\tilde{\omega}^1$, although its performance remains notably worse than that of the United Filter. However, when the perturbation noise increases significantly, i.e., $\tilde{\omega}^2$, the AugEnKF essentially fails in the state estimation task.

Next, we consider the parameter estimation performance of the United Filter. The results for the estimated parameters are shown in Figure~\ref{ParaEst_AdvDiff}. From this figure, we observe that the United Filter produces very accurate estimation for the parameters - regardless of levels of perturbation noise.

\begin{figure}[h!]
\centering
\begin{minipage}{0.33\textwidth}
\includegraphics[scale=0.22]{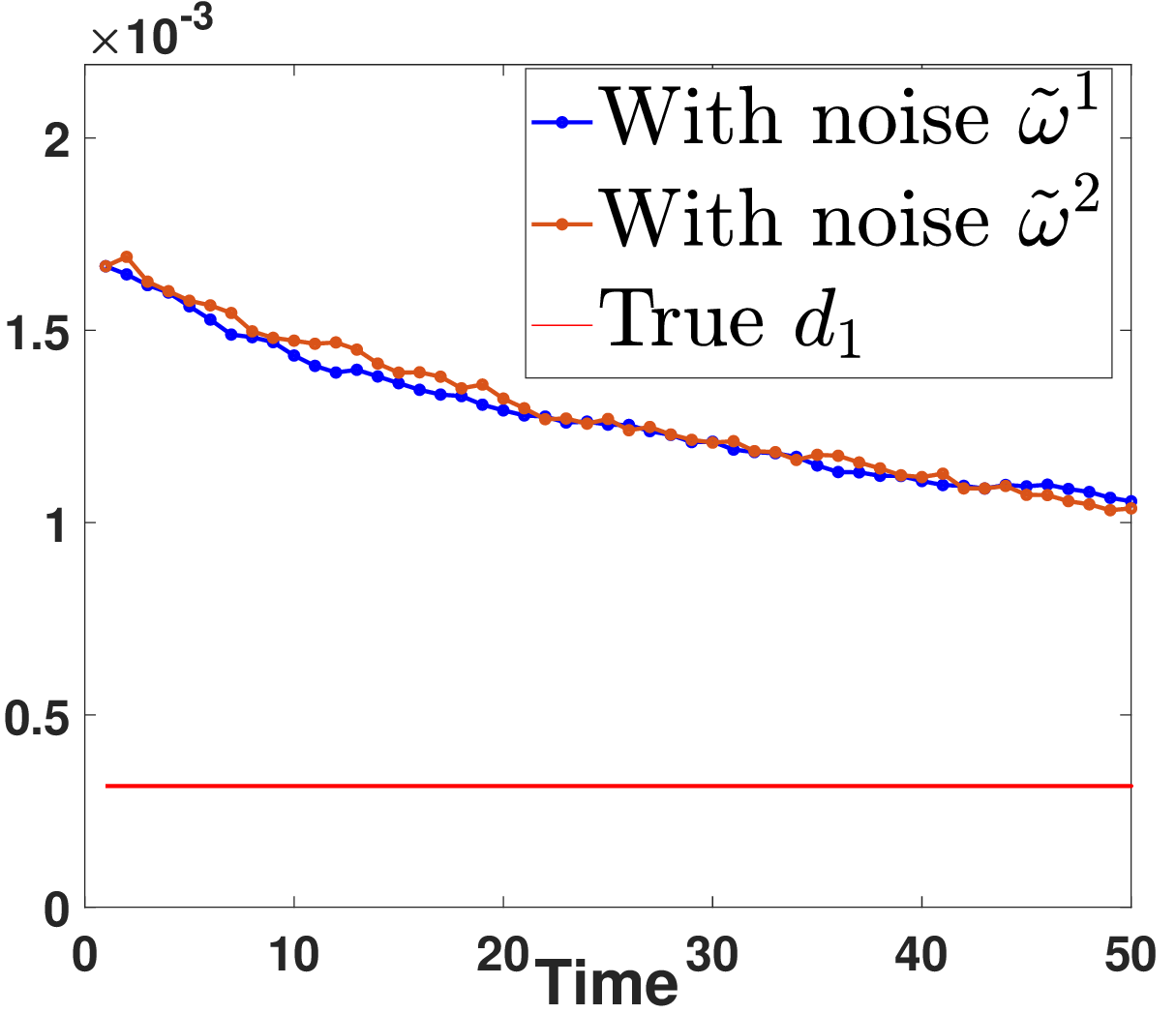}
\end{minipage}%
\begin{minipage}{0.33\textwidth}
\hspace{0.1cm}\includegraphics[scale=0.22]{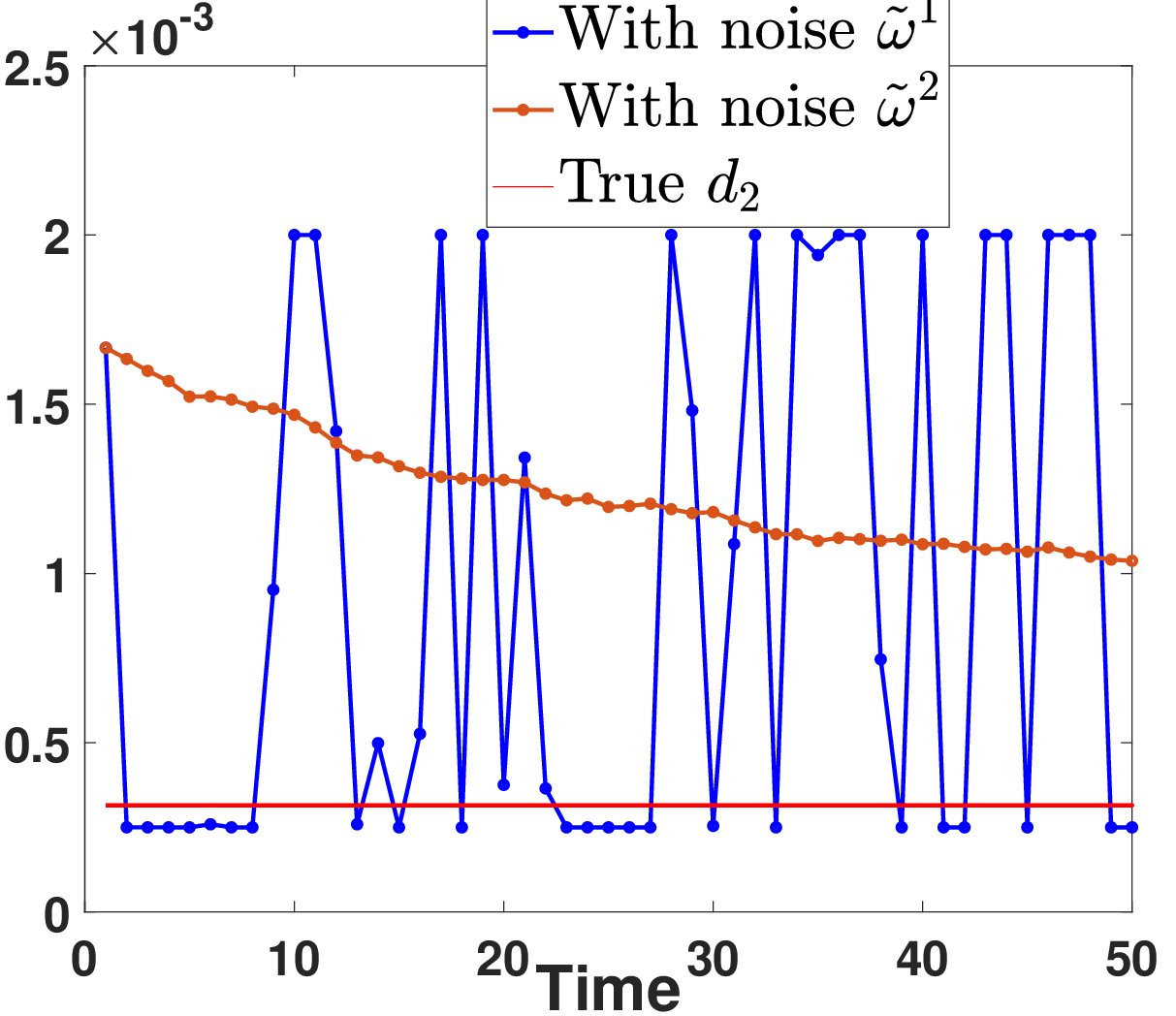}
\end{minipage}%
\begin{minipage}{0.33\textwidth}
\hspace{0.1cm}\includegraphics[scale=0.22]{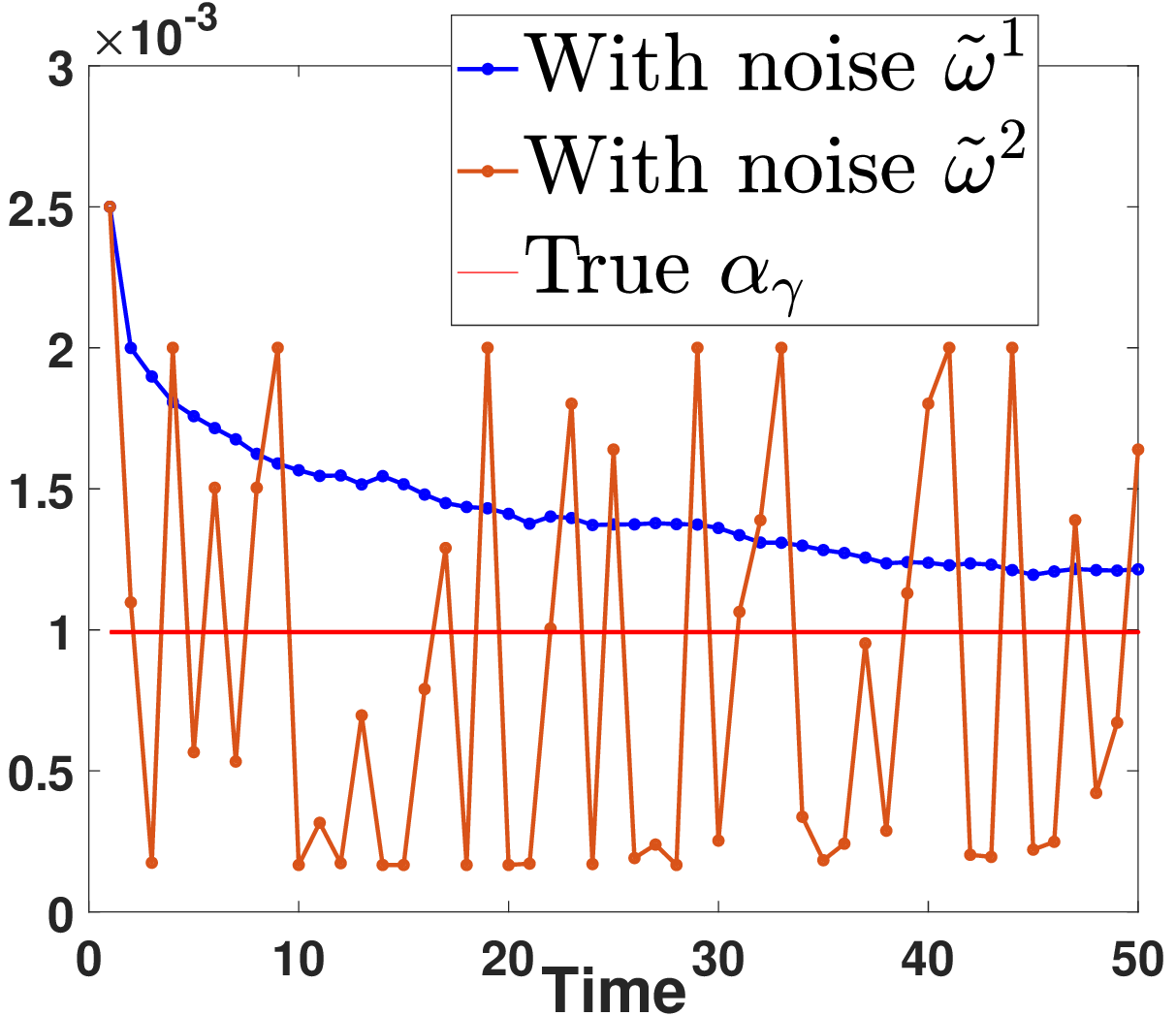}
\end{minipage}
\caption{[Test Case 3] Parameter estimation by the AugEnKF. (Left) Estimation for $d_1$ (Middle); Estimation for $d_2$; (Right) Estimation for $\alpha_{\gamma}$. This is a comparison result to the United Filter's estimates presented in Figure \ref{ParaEst_AdvDiff}.}
\label{ParaEst_AdvDiff_EnKF}
\vspace{-0.4cm}
\end{figure}

Finally, we present the results for the parameter estimation obtained by the AugEnKF in Figure~\ref{ParaEst_AdvDiff_EnKF}. 
Similar to the setting of the United Filter, we aim to approximate the reciprocals of the true parameters due to their very small values. We can observe from Figure~\ref{ParaEst_AdvDiff_EnKF} that the AugEnKF failed to capture the true parameters, regardless of the level of uncertainty, and exhibit significant instability.

As a conclusion for the numerical experiments, we verified that the United Filter algorithm outperforms the AugEnKF in solving the reduced fracture model for flow and transport problems. The United Filter method can give highly reliable results, even in challenging scenarios such as Test Case 2 and Test Case 3.

\section{Conclusion}
In this work, we investigated the data assimilation problem of jointly estimating the state of flow and transport in fractured porous media, along with unknown parameters representing the media's physical characteristics. The state - parameter estimation was performed in an online manner using noisy and sparse observations of fluid flow and contaminants in the rock matrix. To model the forward dynamics in the data assimilation problem, we employed a reduced fracture model and solved the estimation task using a novel United Filter method. This method iteratively approximates the state variables and updates the unknown parameters based on state predictions. State estimation was conducted using the Ensemble Score Filter (EnSF), which leverages a diffusion model framework, while parameter estimation was implemented using the Direct Filter method. To evaluate the performance, reliability, and stability of the United Filter, we presented three numerical experiments under different settings. Additionally, we applied the Augmented Ensemble Kalman Filter (AugEnKF) to perform similar tasks and demonstrated that our proposed algorithm outperforms AugEnKF.

For future study, we aim to tackle more challenging scenarios where observational data are extremely sparse. Additionally, we plan to explore cases where the physical properties of the state variables must be explicitly accounted for and preserved throughout the filtering process.

\section*{Acknowledgement}

This material is based upon work supported by the U.S. Department of Energy, Office of Science, Office of Advanced Scientific Computing Research, Applied Mathematics program under the contract ERKJ387 at the Oak Ridge National Laboratory, which is operated by UT-Battelle, LLC, for the U.S. Department of Energy under Contract DE-AC05-00OR22725. The corresponding author (FB) would like to acknowledge the support from U.S. National Science Foundation through project DMS-2142672 and the support from the U.S. Department of Energy, Office of Science, Office of Advanced Scientific Computing Research, Applied Mathematics program under Grants DE-SC0025412 and DE-SC0024703.

\bibliographystyle{plain} 
\bibliography{References}
\end{document}